\documentclass[a4paper,12pt]{amsart}
\usepackage{amsmath,amssymb,amsfonts,amsthm}
\usepackage{bm}
\usepackage[top=25mm,bottom=25mm,left=25mm,right=25mm]{geometry}
\usepackage[]{amsrefs}
\usepackage{ascmac}
\usepackage[all]{xy}
\usepackage{amsthm}
\usepackage{relsize}
\usepackage{tikz}
\usepackage{mathrsfs}
\usetikzlibrary{trees}
\usetikzlibrary{cd,arrows,matrix,backgrounds,shapes,positioning,calc,decorations.markings,decorations.pathmorphing,decorations.pathreplacing}
\usepackage{multirow}
\usetikzlibrary{intersections,arrows,calc,decorations.markings}
\newtheorem{theorem}{Theorem}[section]
\usepackage{mathtools,leftidx}

\newtheorem{proposition}[theorem]{Proposition}

\newtheorem{corollary}[theorem]{Corollary}
\newtheorem{lemma}[theorem]{Lemma}

\theoremstyle{definition}
\newtheorem{remark}[theorem]{Remark}
\newtheorem{example}[theorem]{Example}
\newtheorem{definition}[theorem]{Definition}

\newtheorem{question}[theorem]{Question}

\newtheorem*{theorem*}{Theorem}

\makeatletter
 
 \@addtoreset{equation}{section}
\makeatother
\def\bb{\mathbf{b}}
\def\cc{\mathbf{c}}

\def\ee{\mathbf{e}}

\def\gg{\mathbf{g}}

\def\vv{\mathbf{v}}

\def\xx{\mathbf{x}}
\def\yy{\mathbf{y}}

\def\TT{\mathbb{T}}

\def\PP{\mathbb{P}}
\def\ZZ{\mathbb{Z}}
\def\QQ{\mathbb{Q}}
\def\Acal{\mathcal{A}}

\def\Fcal{\mathcal{F}}
\def\Xcal{\mathcal{X}}

\def\QQsf{\mathbb{Q}_{\text{sf}}}
\def\trop{\mathrm{Trop}}

\allowdisplaybreaks[4]
\usepackage{xcolor}

\newcommand{\gb}{\mathbf{g}}
\newcommand{\xb}{\mathbf{x}}
\newcommand{\xhat}{\widehat{x}}
\newcommand{\fpq}{\frac{p}{q}}
\newcommand{\frs}{\frac{r}{s}}
\newcommand{\Ptilde}{\widetilde{\mathcal{P}}}
\newcommand{\calP}{\mathcal{P}}
\newcommand{\prin}{\mathrm{prin}}
\usepackage{mathdots}
\usepackage{fancybox}
\newcommand\scalemath[2]{\scalebox{#1}{\mbox{\ensuremath{\displaystyle #2}}}}
\usepackage[T1]{fontenc}
\usepackage{lmodern}

\title[Cluster algebraic interpretation]{Cluster algebraic interpretation of generalized Markov numbers and their matrixizations}

\author{Esther Banaian}
\address[Esther Banaian]{Department of Mathematics, University of California, Riverside, Riverside, CA, 92501, United States}
\email{estherb@ucr.edu}

\author{Yasuaki Gyoda}
\address[Yasuaki Gyoda]{Institute for Advanced Research, Nagoya University Furo-cho, Chikusa-ku, Nagoya-shi, 464-8601, Japan}
\email{ygyoda@math.nagoya-u.ac.jp}


\begin{document}
\begin{abstract}
Markov numbers, i.e. positive integers appearing in solutions to $x^2 + y^2 + z^2 = 3xyz$, can be viewed as specializations of cluster variables. The second author and Matsushita gave a generalization of the Markov equation, $x^2 + y^2 + z^2 + k_1yz + k_2xz + k_3xy = (3+k_1+k_2+k_3)xyz$, whose solutions can be viewed as specializations of cluster variables in generalized cluster algebras. We give two families of matrices in $SL(2,\mathbb{Z}[x_1^\pm,x_2^\pm,x_3^\pm])$ associated to these cluster structures. These matrix formulas relate to previous matrices appearing in the context of Markov numbers, including Cohn matrices and generalized Cohn matrices given by the second author, Maruyama, and Sato, as well as matrices appearing in the context of cluster algebras, including matrix formulas given by Kanatarc{\i} O\u{g}uz and Y{\i}ld{\i}r{\i}m. We provide a classification of the two families of matrices and exhibit an explicit family of each. The latter is done by realizing cluster variables in generalized Markov cluster algebras as weight-generating functions of order ideals in certain fence posets which are related to Christoffel words. An interesting observation is that these functions resemble Caldero-Chapoton functions for string modules, and a byproduct of our proofs is a new skein-like formula for such functions.
\end{abstract}
\maketitle
\tableofcontents

\section{Introduction}

Markov numbers are numbers which appear in positive integral solutions to the Markov equation \[
x^2 + y^2 + z^2 = 3xyz.
\]
These were first introduced in the context of Diophantine approximation \cite{markoff1879formes}, but Markov numbers have since become an area of considerable interest in many mathematical areas, including hyperbolic geometry, algebraic geometry, number theory, and combinatorics. Part of this interest stems from Frobenius' unicity conjecture, which posits that Markov triples, i.e. $(a,b,c) \in \mathbb{Z}_{>0}^3$ satisfying the Markov equation, are uniquely determined by their largest element. The unicity conjecture is still open today.

Another modern motivation to study Markov numbers is their connection to cluster algebras. It is easy to verify that if $(a,b,c)$ is a Markov triple, so is $(a,b,\frac{a^2 + b^2}{c})$, and similarly for $a$ and $b$. This process is an instance of Vieta jumping, and every Markov triple is reachable from $(1,1,1)$ using Vieta jumps \cite{hurwitz1963aufgabe}. This process is reminiscent of cluster mutation, and in fact Markov numbers are specializations of cluster variables in the \emph{Markov cluster algebra} \cite{beineke2011cluster,propp2005combinatorics}. This perspective has been useful in making progress towards the unicity conjecture \cite{lee2020ordering,rabideau2020continued}. The Markov cluster algebra is also interesting in its own right. This cluster algebra can be seen as arising from the once-punctured torus. While in general cluster algebras from surfaces enjoy many nice properties, those arising from closed surfaces can be more complicated. For instance, the Markov cluster algebra is a classic example of a cluster algebra which is not equal to its upper cluster algebra \cite{BFZ} and it is non-Noetherian \cite{muller2013locally}.

The second author and Matsushita introduced a family of generalizations of the Markov equation in \cite{GyodaMatsushita}. For any choice of a triple of nonnegative integers, $k_1,k_2,k_3$, this equation is given by \begin{equation}\label{eq:GM}
x^2 + y^2 + z^2 + k_1yz + k_2xz +k_3 xy = (3 + k_1 + k_2 + k_3)xyz.
\end{equation}

Here, a Vieta jump which replaces the third component in a positive integral triple $(a,b,c)$ satisfying Equation \eqref{eq:GM} is of the form \[
(a,b,c) \mapsto \bigg(a,b, \frac{a^2 + k_3 ab + b^2}{c}\bigg),
\]
and all positive integral solutions are again reachable from $(1,1,1)$ using Vieta jumps (in any of the three positions) \cite{GyodaMatsushita}.
When $k_i$ is nonzero, this replacement now is of the form of \emph{generalized cluster mutation}, in the sense of Chekhov and Shapiro \cite{chekhov2014teichmuller}. Solutions to \eqref{eq:GM} for all values of $k_i$ are specializations of cluster variables in \emph{generalized cluster algebras}, as was explained in \cite{GyodaMatsushita}. An initial motivation for the definition of generalized cluster algebras was the Teichm{\"u}ller space of an orbifold. Triples satisfying Equation \eqref{eq:GM} in the case $k_1 = k_2 = k_3 = 1$ can be seen to be specializations of cluster variables in a generalized cluster algebra from a once-punctured sphere with three orbifold points of order three. The first author and Sen independently studied this special case of the Gyoda-Matsushita equation in \cite{banaian2024generalization}.

One line of research on Markov numbers, distinct from the cluster algebraic approach, involves studying elements of the modular group $SL(2,\mathbb Z)$ that contain Markov numbers as $(1,2)$-components and analyzing their associated modular forms. A major turning point in this direction was the introduction of \emph{Cohn matrices} by Cohn \cites{Cohn1,Cohn2} building on work by Frobenius \cite{frobenius1968uber}. The definition of Cohn matrices, which sit in \emph{Cohn triples}, is inspired by the
resemblance between the \emph{Fricke identity} \cite{fricke1896uber} in the modular group
\[\mathrm{tr}(A)^2+\mathrm{tr}(B)^2+\mathrm{tr}(AB)^2=\mathrm{tr}(A)\mathrm{tr}(B)\mathrm{tr}(AB)+\mathrm{tr}(ABA^{-1}B^{-1})+2\]
and the Markov equation. When substituting a Cohn triple with traces $x, y, z$ into the Fricke identity, the equation takes the form
\begin{align}\label{eq:second-Markov}
x^2 + y^2 + z^2 = xyz.
\end{align}
Solutions to Equation \eqref{eq:second-Markov} are known to be in one-to-one correspondence with the set of positive integer solutions to the Markov equation. For this reason, Equation \eqref{eq:second-Markov} is sometimes called the \emph{second Markov equation}.

Overall Cohn matrices provide a strong link between Markov numbers and hyperbolic geometry and have been used as a tool towards the unicity conjecture \cite{aigner2013Markov,zhang2006elementary}
Cohn matrices have also been studied in relation to \emph{Christoffel words} \cite{MR4265545}, and they now occupy a central place in modern research on Markov numbers.

More recently, Cohn matrices have been generalized by the second author, Maruyama, and Sato \cite{gyoda2024sl}. In the process of generalizing the Markov equation to the form \eqref{eq:GM}, they observed that a compatible generalization of the second Markov equation, \eqref{eq:second-Markov}, is
\begin{align}\label{eq:second-Markov-gen}
x^2+y^2+z^2+(2k_1+k_2k_3)x+(2k_2+k_1k_3)y+(2k_3+k_1k_2)z+k_1^2+k_2^2+k_3^2+2k_1k_2k_3=xyz.
\end{align}
The interaction between Equations \eqref{eq:GM} and \eqref{eq:second-Markov-gen} served as inspiration for the notion of a generalized Cohn triple. We remark that \cite{gyoda2024sl} worked in the case  $k_1 = k_2 = k_3$; here we allow the values $k_i$ to be distinct.

The authors of \cite{gyoda2024sl} also discovered that \eqref{eq:second-Markov-gen}  coincides with a special case of an identity in $(X,Y,Z) \in SL(2,\mathbb{C})^3$ previously established by Luo \cite{luo1998geodesic} or Nakanishi and N\"a\"at\"anen \cite{nakanishi2001areas},
\begin{align}
x^2 + y^2 + z^2 + (ad+bc)x + (bd+ca)y + (cd+ab)z + a^2 + b^2 + c^2 + d^2 + abcd - 4 = xyz,
\end{align}
where $x=-\mathrm{tr}(YZ)$, $y=-\mathrm{tr}(ZX)$, $z=-\mathrm{tr}(XY)$, $a=-\mathrm{tr}(X)$, $b=-\mathrm{tr}(Y)$, $c=-\mathrm{tr}(Z)$, $d=-\mathrm{tr}(XYZ)$. Building on this, they constructed a new type of matrix and investigated its relationship with the generalized Cohn matrix. These matrices are referred to as \emph{Markov-monodromy matrices}.

In this paper, we introduce the ``cluster matrixization'' of generalized Markov numbers, which unifies the two aforementioned lines of research on clusterization and matrixization. Specifically, we define matrices called the \emph{cluster generalized Cohn matrix} and the \emph{cluster Markov-monodromy matrix}, and investigate their properties. These matrices are $2 \times 2$ matrices whose entries are Laurent polynomials of three variables, with the $(1,2)$-entry being a cluster variable of a generalized Markov cluster algebra. By substituting $x_1 = x_2 = x_3 = 1$ into the three variables of the Laurent polynomials, these matrices reduce to the generalized Cohn matrix and the Markov-monodromy matrix, respectively. In this sense, they can be regarded as \emph{clusterizations} of the generalized Cohn matrix and the Markov-monodromy matrix. Furthermore, these matrices retain many of the properties of the generalized Cohn matrix and the Markov-monodromy matrix studied by \cite{gyoda2024sl}, leading to the conclusion that the combinatorial structures inherent in the generalized Cohn matrix and the Markov-monodromy matrix originate from these cluster matrices. In the second half, we pay special attention to one specific tree for each family of matrices. We summarize the main results.

\begin{itemize}
    \item We classify all cluster generalized Cohn matrices and cluster Markov-monodromy matrices. In particular, we organize both families of matrices into binary trees and show that each matrix belongs to a unique tree (Theorems \ref{thm:all-cohn-triple}, \ref{thm:all-cohn-triple2}, \ref{thm:all-Markov-triple}, and \ref{thm:all-markov-triple2}). This characterization can be seen as an analogue of the family of trees which encompass all (integer) Cohn triples  \cite{aigner2013Markov} and the generalized analogue in \cite{gyoda2024sl}.
    \item We exhibit an isomorphism $\Psi_g$ between cluster generalized Cohn trees and cluster Markov-monodromy trees (Theorem \ref{thm:BT-CT2}). This isomorphism is a parameter-generalized cluster lift of the correspondence $\Psi$ introduced by Gyoda, Maruyama, and Sato \cite[Proposition 5.4 and Theorem 5.6]{gyoda2024sl}: it allows distinct $k_1,k_2,k_3$ and an arbitrary Laurent-polynomial parameter $g$. Their map is recovered by specializing $x_1=x_2=x_3=1$, setting $k_1=k_2=k_3=k$, and taking $g=k$.
    \item We give a new combinatorial expansion formula for cluster variables in the generalized cluster algebras considered here (Theorem \ref{thm:CorrectnessOfPosetFormula}). This formula uses weighted order ideals of a family of fence posets.
    \item We construct one tree of cluster generalized Cohn matrices and one tree of cluster Markov-monodromy matrices, i.e. the \emph{combinatorial cluster generalized Cohn tree} and the \emph{combinatorial cluster Markov-monodromy tree}, using these fence posets (Theorem \ref{thm:Combinatorial_Cohn_Matrix_and_Triple} and Section \ref{sec:CombinatorialMM}).
\end{itemize}

The combinatorial expansion formula is as follows. It is well-known that Markov numbers (and solutions to the Gyoda-Matsushita equation in general) can be indexed by rational numbers. To each rational number $\fpq$, we associate a poset $\calP_\fpq$. As a fence poset, each $\calP_\fpq$ bears resemblance to the poset of join irreducibles in the lattice of perfect matchings of a \emph{snake graph}. Snake graphs are weighted graphs whose generating function of perfect matchings provides a combinatorial expansion formula for a cluster variable in a surface cluster algebra \cite{musiker2009cluster,musiker2011positivity, Wilson}.
Perfect matchings of snake graphs can be endowed with a partial order which turns out to form a \emph{distributive lattice} \cite{musiker2013bases,Wilson}. This allows these formulas to be recast in the language of order ideals of posets \cite{ezgieminecluster2024,pilaud2023posets}.
Similarly, Theorem \ref{thm:CorrectnessOfPosetFormula} could easily be phrased in terms of snake graphs. We remark that the specific generalized cluster algebras considered here are not covered by previous work since they come from closed orbifolds.

Our specific matrices in each combinatorial tree have each entry given as a weight-generating function of subsets of order ideals in $\calP_\fpq$; these are inspired by the matrices in \cite{ouguz2025oriented,ezgieminecluster2024}.

Snake graph calculus, introduced in \cite{canakci2013snake}, is a combinatorial procedure on snake graphs which produces multiplication formulas for surface-type cluster algebras. With Kang and Kelley, the first author adapted and extended snake graph calculus to the setting of posets \cite{banaian2024skein}. These formulas are crucial in showing that our combinatorial families of matrices indeed are cluster generalized Cohn matrices and cluster Markov-monodromy matrices.   Since we are working with generalized cluster algebras, unlike \cite{banaian2024skein} which focused on the ordinary case, we must provide a new skein relation. Proposition \ref{prop:KissingReverseOverlap} is the statement needed for our purposes. Similar relations have been discussed in the context of Caldero-Chapoton (CC) functions (discussed in Remark \ref{rem:CC}), and our Laurent polynomial associated to a poset resembles the calculation of a CC function from a string module. In Section \ref{subsec:ExtendReverseKissingToXs}, and in particular in Theorem \ref{thm:DecomposeVariables}, we discuss how Proposition \ref{prop:KissingReverseOverlap} can be extended in this setting.

We conclude our introduction by summarizing the contents of the remainder of the article. Section \ref{sec:GenCA} includes all relevant background on generalized cluster algebras, and Section \ref{sec:GenMarkovCA} introduces the specific cluster algebras we study here. In Sections \ref{sec:ClusterCohn} and \ref{sec:ClusterMM}, we introduce and characterize cluster generalized Cohn matrices and cluster Markov-monodromy matrices, respectively. The latter characterization is proved in Section \ref{sec:Relations}, where we introduce two relations between the two sets of matrices.

In Section \ref{sec:Posets}, up to symmetry, we construct a poset for each generalized Markov number using a common indexing via rational numbers. In Section \ref{sec:MarkovPosets}, we prove that these posets provide an expansion formula for the cluster variables in the generalized cluster algebras considered. The latter section also provides a discussion of skein relations on posets. In Sections \ref{sec:CombinatorialCohnMatrices} and \ref{sec:CombinatorialMM}, we exhibit specific families of cluster generalized Cohn matrices and cluster Markov-monodromy matrices. Finally, in Section \ref{sec:RelationshipsToOthers}, we give a brief discussion on the relationships between our matrices and other matrix formulas related to Markov numbers and cluster algebras.

\section*{Acknowledgements}
We thank Gregg Musiker for helpful discussions concerning the relationship amongst these matrix formulas. We thank Tomoki Nakanishi, Ryota Akagi and Zhichao Chen for helpful comments and an anonymous reviewer for useful suggestions regarding the exposition. This work was supported by JSPS KAKENHI Grant Number JP25K17224.

\section{Generalized Cluster Algebras}\label{sec:GenCA}
In this section, we introduce generalized cluster algebras. We start with recalling definitions of seed mutations and generalized cluster patterns according to \cites{chekhov2014teichmuller,nakanishi2015structure,nakanichi2016companion}.

\subsection{Generalized Cluster Algebras}

A \emph{semifield} $\mathbb P$ is an abelian multiplicative group equipped with an addition $\oplus$ which is distributive over the multiplication. We particularly make use of the following two semifields.

Let $\mathbb Q_{\text{sf}}(u_1,\dots,u_{\ell})$ be the set of rational functions in $u_1,\dots,u_{\ell}$ which have subtraction-free expressions. This set can be given the structure of a semifield under the usual multiplication and addition. It is called the \emph{universal semifield} of $u_1,\dots,u_{\ell}$ (\cite[Definition 2.1]{fomin2007cluster}).

Let $\text{Trop}(u_1,\dots, u_\ell)$ be the abelian multiplicative group freely generated by the elements $u_1,\dots,u_\ell$. Given two elements, $\prod_{j=1}^\ell u_j^{a_j}$ and $\prod_{j=1}^\ell u_j^{b_j}$ in $\text{Trop}(u_1,\dots, u_\ell)$, define their sum as
\begin{align}
\prod_{j=1}^\ell u_j^{a_j} \oplus \prod_{j=1}^{\ell} u_j^{b_j}=\prod_{j=1}^{\ell} u_j^{\min(a_j,b_j)}.
\end{align}
The set $\text{Trop}(u_1,\dots, u_\ell)$  is a semifield under $\oplus$ and the usual multiplication, often called  the \emph{tropical semifield} of $u_1,\dots,u_\ell$ (\cite[{Definition 2.2}]{fomin2007cluster}).

For any semifield $\PP$ and $p_1, \dots, p_{\ell}\in\PP$, there exists a unique semifield homomorphism $\pi$ such that
\begin{align} \label{qsfuniv}
	\pi:\QQsf(y_1, \dots, y_{\ell}) &\longrightarrow \PP\\
	y_i &\longmapsto p_i. \nonumber
\end{align}
For $F(y_1,\dots,y_\ell ) \in \QQsf(y_1, \dots, y_{\ell})$, we denote the \emph{evaluation} of $F$ at $p_1, \dots, p_{\ell}$ by
\begin{align}
	F|_{\PP}(p_1, \dots, p_{\ell}):=\pi(F(y_1, \dots, y_\ell)).
\end{align}
We fix a positive integer $n$ and a semifield $\PP$. Let $\mathbb{ZP}$ be the group ring of $\mathbb{P}$ as a multiplicative group. Since $\mathbb{ZP}$ is a domain (\cite[{Section 5}]{fomin2002cluster}), its total quotient ring is a field $\mathbb{Q}(\mathbb P)$. Let $\mathcal{F}$ be the field of rational functions in $n$ indeterminates with coefficients in $\mathbb{Q}(\mathbb P)$.
Let $n\in \ZZ_{\geq1}$ and $\Fcal$ be a rational function field of $n$ indeterminates. A \emph{labeled seed} is a quadruplet $(\mathbf{x}, \mathbf{y},B,\mathbf{Z})$, where
\begin{itemize}\setlength{\leftskip}{-15pt}
\item $\mathbf{x}=(x_1, \dots, x_n)$ is an $n$-tuple of elements of $\mathcal F$ forming a free generating set of $\mathcal F$,
\item $\mathbf{y}=(y_1, \dots, y_n)$ is an $n$-tuple of elements of $\mathbb P$,
\item $B=(b_{ij})$ is an $n \times n$ integer matrix which is \emph{skew-symmetrizable}, that is, there exists a positive integer diagonal matrix $S$ such that $SB$ is skew-symmetric,
\item $\mathbf{Z}=(Z_1,\dots,Z_n)$ is an $n$-tuple of non-constant polynomials with the coefficient in $\ZZ_{\geq 0}\mathbb P$ \footnote{In \cites{nakanishi2015structure, nakanichi2016companion}, $z_{i,j}$ is an element in $\mathbb P$, but we can take $z_{i,j}$ as an element in $\mathbb Z_{\geq 0} \mathbb P$. See \cite{nakanishi2024addendum}.}
\[Z_i(u)=z_{i,0}+z_{i,1}u+\cdots+z_{i,d_i}u^{d_i}\]
satisfying $z_{i,0}=z_{i,d_i}=1$.
\end{itemize}

We call $S$ in the above definition a \emph{skew-symmetrizer} of $B$.
We say that $\xx$ is a \emph{(labeled) cluster} and $\yy$ is a \emph{coefficient tuple}, and we refer to $x_i$, $B$ and $Z_i$ as the \emph{cluster variable}, the \emph{exchange matrix (\textrm{or} $B$-matrix)} and the \emph{exchange polynomial}, respectively. Furthermore, we set $D=\text{diag} (d_1,\dots, d_n)$; in particular, $D$ is a positive integer diagonal matrix of rank $n$.

Given an integer $b$,  let $[b]_+=\max\{b,0\}$.
Let $(\mathbf{x},\mathbf{y}, B,\mathbf{Z})$ be a labeled seed, and let $k \in\{1,\dots, n\}$. The \emph{seed mutation $\mu_k$ in direction $k$} transforms $(\mathbf{x},\mathbf{y}, B,\mathbf{Z})$ into another labeled seed $\mu_k(\mathbf{x},\mathbf{y}, B,\mathbf{Z})=(\mathbf{x}',\mathbf{y}', B',\mathbf{Z}')$ defined as follows:
\begin{itemize}\setlength{\leftskip}{-15pt}
\item The cluster variables $\mathbf{x'}=(x'_1, \dots, x'_n)$ are given by
\begin{align}\label{eq:x-mutation}
x'_j=\begin{cases}\frac{1}{x_k}\left(\mathop{\prod}\limits_{i=1}^{n} x_i^{[-b_{ik}]_+}\right)^{d_k}\frac{Z_k\left(y_k\mathop {\prod}\limits_{i=1}^{n} x_i^{b_{ik}}\right)}{Z_k|_{\mathbb P}(y_k)} &\text{if $j=k$,} \vspace{4mm}\\
x_j &\text{otherwise.}
\end{cases}
\end{align}
\item The coefficient tuple $\yy'=(y_1',\dots,y'_n)$ are given by
\begin{align}\label{eq:y-mutation}
y'_j=\begin{cases}
y_{k}^{-1} &\text{if $j=k$,} \\
y_j y_k^{d_k[b_{kj}]_+}Z_k|_{\mathbb P}(y_k)^{-b_{kj}} &\text{otherwise.}
\end{cases}
\end{align}
\item The entries of $B'=(b'_{ij})$ are given by
\begin{align} \label{eq:matrix-mutation}
b'_{ij}=\begin{cases}-b_{ij} &\text{if $i=k$ or $j=k$,} \\
b_{ij}+d_k\left(\left[ b_{ik}\right] _{+}b_{kj}+b_{ik}\left[ -b_{kj}\right]_+\right) &\text{otherwise.}
\end{cases}
\end{align}
\item The exchange polynomials $\mathbf{Z'}=(Z'_1, \dots, Z'_n)$ are given by
\begin{align}\label{eq:Z-mutation}
Z'_j(u)=\begin{cases}u^{d_k}Z_k(u^{-1})  &\text{if $j=k$,}\\
Z_j(u) &\text{otherwise.}
\end{cases}
\end{align}
\end{itemize}

One can check that mutation is an involution on labeled seeds. Let $\mathbb{T}_n$ be the \emph{$n$-regular tree} whose edges are labeled by the numbers $1, \dots, n$ such that the $n$ edges emanating from each vertex have different labels. We write
$\begin{xy}(0,0)*+{t}="A",(10,0)*+{t'}="B",\ar@{-}^k"A";"B" \end{xy}$
to indicate that vertices $t,t'\in \mathbb{T}_n$ are joined by an edge labeled by $k$. We fix an arbitrary vertex $t_0\in \TT_n$, which is called the \emph{rooted vertex}.
A \emph{generalized cluster pattern} is an assignment of a labeled seed $\Sigma_t=(\xx_t,\yy_t,B_t,\mathbf Z_t)$ to every vertex $t\in \mathbb{T}_n$ such that the labeled seeds $\Sigma_t$ and $\Sigma_{t'}$ assigned to the endpoints of any edge
$\begin{xy}(0,0)*+{t}="A",(10,0)*+{t'}="B",\ar@{-}^k"A";"B" \end{xy}$
are obtained from each other by the seed mutation in direction $k$. When the initial seed is $\Sigma_{t_0}=(\xx,\yy, B,\mathbf Z)$, we denote by $CP_{(\xx,\yy,B,\mathbf Z)}\colon t\mapsto \Sigma_t$ this assignment.
The degree $n$ of the regular tree $\TT_n$ is called the \emph{rank} of a generalized cluster pattern $CP_{(\xx,\yy,B,\mathbf Z)}$.

In this paper, we call a \emph{non-labeled cluster} a permutation of $(x_{1;t},\dots,x_{n;t})$.

\begin{definition}
A \emph{generalized cluster algebra} $\Acal(B,\mathbf Z)$ associated with a cluster pattern $CP_{(\xx,\yy,B,\mathbf Z)}$ is the $\ZZ\PP$-subalgebra of $\Fcal$ generated by $\Xcal=\{x_{i;t}\}_{1\leq i\leq n, t\in \TT_n}$.
\end{definition}

\begin{remark}
A generalized cluster algebra in which $D$ is an identity matrix is called a \emph{(ordinary) cluster algebra}. In this case, the symbol $\mathbf Z$ is often omitted from the notation of cluster patterns or cluster algebras because all polynomials $Z_i(u)$ are identically $1+u$, and these are invariant under mutation.

Since this article largely concerns generalized cluster algebras,  we may at times drop the ``generalized'' adjective.  We will instead stress when we are dealing with an ordinary cluster algebra.
\end{remark}

\subsection{$c$-vectors, $g$-vectors, $F$-polynomials}

Next, we define two important families of vectors associated to a cluster algebra, following \cites{fomin2007cluster}. We first need to introduce principal coefficients.
\begin{definition}
We say that a generalized cluster pattern $v\mapsto \Sigma_v$ or a generalized cluster algebra $\Acal$ of rank $n$ has \emph{principal coefficients} at the rooted vertex $t_0$ if $\mathbb{P}=\text{Trop}(y_1,\dots,y_n)$ and $\mathbf{y}_{t_0}=(y_1,\dots,y_n)$. In this case, we denote $\Acal=\Acal^{\prin}(B,\mathbf Z)$.
\end{definition}
First, we define the $c$-vectors. For $\bb=(b_1,\dots,b_n)^{\top}$, we use the notation $[\bb]_+=([b_1]_+,\dots,[b_n]_+)^{\top}$, where $\top$ stands for transpose.
\begin{definition}
Let $\Acal^{\prin}(B,\mathbf Z)$ be a generalized cluster algebra with principal coefficients at $t_0$. We define the \emph{$c$-vector} $\cc_{j;t}$ as the degree of $y_i$ in $y_{j;t}$, that is, if $y_{j;t}=y_1^{c_{1j;t}}\cdots y_n^{c_{nj;t}}$, then
\begin{align}
\cc_{j;t}^{B;t_0}=\cc_{j;t}=\begin{bmatrix}c_{1j;t}\\ \vdots \\ c_{nj;t} \end{bmatrix}.
\end{align}
We define the \emph{$C$-matrix} $C_t^{B;t_0}$ as
\begin{align}
C_t^{B;t_0}:=(\cc_{1;t},\dots,\cc_{n;t}).
\end{align}
\end{definition}
The $c$-vectors are the same as those defined by the following recursion. For all $j \in\{1,\dots,n\}$, define
\begin{align*}
\cc_{j;t_0}=\ee_j\quad \text{(canonical basis)},
\end{align*}
and for any \begin{xy}(0,1)*+{t}="A",(10,1)*+{t'}="B",\ar@{-}^k"A";"B" \end{xy},  define
\begin{align*}
\cc_{j;t'} =
\begin{cases}
-\cc_{j;t} & \text{if $j=k$;} \\[.05in]
\cc_{j;t} + [d_kb_{kj;t}]_+ \ \cc_{k;t} +d_kb_{kj;t} [-\cc_{k;t}]_+
 & \text{if $j\neq k$}.
 \end{cases}
\end{align*}
See \cite{nakanichi2016companion} for the equivalence of these definitions. Since the recursion formula only depends on exchange matrices,
we can regard $c$-vectors as vectors associated with vertices of $\TT_n$. In this way, we remark that it is possible to define $c$-vectors for a cluster algebra that does not have principal coefficients.

Nakanishi shows that the $c$-vectors in a generalized cluster algebra are equivalent to those in a certain ordinary cluster algebra.

\begin{proposition}[{\cite[Proposition 3.16]{nakanishi2015structure}}]\label{prop:NakanishiCvectors}
The $c$-vectors of the generalized cluster pattern $CP_{(\xx,\yy,B,\mathbf{Z})}$ coincide with the $c$-vectors of the cluster pattern associated to the ordinary cluster algebra $CP_{(\xx,\yy,DB)}$.
\end{proposition}

Every cluster algebra has a dual set of vectors, called $g$-vectors.  These vectors come from a homogenous multigrading we can place on a cluster algebra with We can regard cluster variables in cluster algebras principal coefficients.

\begin{theorem}[{\cite[Proposition 3.9]{nakanishi2015structure}}]
Let $\Acal^{\prin}(B,\mathbf Z)$ be a generalized cluster algebra with principal coefficients at $t_0$. Each cluster variable $x_{i;t}$ is a homogeneous Laurent polynomial in $x_1,\dots,x_n,y_1,\dots,y_n$ by the following $\ZZ^n$-grading:
\begin{align}\label{grading}
\deg x_{i}=\ee_i,\quad \deg y_{i}=-\mathbf{b}_i,
\end{align}
where $\ee_i$ is the $i$th canonical basis vector of $\ZZ^n$ and $\mathbf{b}_i$ is the $i$th column vector of $B$.
\end{theorem}
We denote by $(g_{1j;t},\dots,g_{nj;t})^\top$ the $\ZZ^n$-grading of $x_{j;t}$.
\begin{definition}
Let $\Acal^{\prin}(B,\mathbf Z)$ be a generalized cluster algebra with principal coefficients at $t_0$. We define the \emph{$g$-vector} $\gg_{j;t}$ as the degree of the homogeneous Laurent polynomial $x_{j;t}$ by the $\ZZ^n$-grading \eqref{grading}, that is,
\begin{align}
\gg_{j;t}^{B;t_0}=\gg_{j;t}=\begin{bmatrix}g_{1j;t}\\ \vdots \\ g_{nj;t} \end{bmatrix}.
\end{align}
We define the \emph{$G$-matrix} $G_t^{B;t_0}$ as
\begin{align}
G_t^{B;t_0}:=(\gg_{1;t},\dots,\gg_{n;t}).
\end{align}
\end{definition}
The $g$-vectors are the same as those defined by the following recursion: For any $j\in\{1,\dots,n\}$,
\begin{align*}
\gg_{j;t_0}=\ee_j\quad \text{(canonical basis)},
\end{align*}
and for any \begin{xy}(0,1)*+{t}="A",(10,1)*+{t'}="B",\ar@{-}^k"A";"B" \end{xy},
\begin{align}\label{g-recursion}
\gg_{j;{t'}}&=\begin{cases}
\gg_{j;t} \ \ & \text{if } j\neq k;\\
-\gg_{k;t}+\mathop\sum\limits_{i=1}^{n}[d_kb_{ik;t}]_+\gg_{i;t}-\mathop\sum\limits_{i=1}^{n}[d_kc_{ik;t}]_+\mathbf{b}_k &\text{if } j=k
\end{cases}
\end{align}
See \cite{nakanichi2016companion} for the equivalence of these definitions.  As with $c$-vectors,  through this recursive definition we can define $g$-vectors even for cluster algebras with coefficient systems other than principal coefficients.

Nakanishi gives a parallel result to Proposition \ref{prop:NakanishiCvectors} for $g$-vectors.

\begin{proposition}[{\cite[Proposition 3.17]{nakanishi2015structure}}]\label{prop:NakanishiGvectors}
The $g$-vectors of the generalized cluster pattern $CP_{(\xx,\yy,B,\mathbf{Z})}$ coincide with the $g$-vectors of the cluster pattern associated to the ordinary cluster algebra $CP_{(\xx,\yy,BD)}$.
\end{proposition}

We will see that set of $g$-vectors along with the following family of polynomials uniquely determine the set of cluster variables.

\begin{definition}\label{def:FPoly}
Let $\mathcal{A}^{\prin}(B,\mathbf Z)$ be a generalized cluster algebra of rank $n$ with principal coefficients at $t_0$. We define the \emph{$F$-polynomial} $F_{j;t}$ to be the result of specializing all initial variables $x_{1;t_0},\ldots,x_{n;t_0}$ to 1 in $x_{j;t}$.

The $F$-polynomials can be alternatively defined as follows.  For all $j \in \{1,\ldots,n\}$, \[
F_{j;t_0} = 1
\]
and for any $ \begin{xy}(0,1)*+{t}="A",(10,1)*+{t'}="B",\ar@{-}^k"A";"B" \end{xy}$,\[
F_{j;t'} = \begin{cases} F_{j;t}^{-1}\bigg(\prod\limits_{i=1}^n y_i^{[c_{jk;t}]_+}F_{i;t}^{[b_{jk;t}]_+}\bigg)^{d_j} \sum\limits_{s=0}^{d_j} z_{j,s}\big(\prod\limits_{i=1}^n  y_i^{c_{jk;t}}F_{i;t}^{b_{jk;t}}\big)^s & \text{if $j = k$;}\\
F_{j;t} &\text{if $j \neq k$.}\end{cases}
\]

\end{definition}

The fact that the second characterization of an $F$-polynomial is equivalent to the first comes from Proposition 3.12 in \cite{nakanishi2015structure}. The convenience of defining $F$-polynomials comes from the separation formula,  shown for ordinary cluster algebras in \cite{fomin2007cluster} and for generalized cluster algebras in \cite{nakanishi2015structure}.

\begin{theorem}[{\cite[Theorem 3.23]{nakanishi2015structure}}]\label{thm:SeparationFormulaNakanishi}
Let $\mathcal{A}(B,\mathbf Z)$ be a generalized cluster algebra of rank $n$. Let $\xx$ denote the initial cluster of  $\mathcal{A}(B,\mathbf Z)$.  Define $\widehat{y_i} = y_i \prod_{j=1}^n x_j^{b_{ji}}$. Then, for any seed $t$ and $i \in \{1,\ldots,n\}$,  we have \[
x_{i;t} = \xx^{\gg_{i;t}} \frac{F_{i;t}\vert_{\mathcal{F}}(\widehat{y}_1,\ldots,\widehat{y}_n)}{F_{i;t}\vert_{\mathbb{P}}(y_1,\ldots,y_n)}
\]
\end{theorem}

By using Theorem \ref{thm:SeparationFormulaNakanishi}, we have the following proposition.

\begin{proposition}\label{thm:ximpliesg}
For any generalized cluster algebra, $x_{i;t}=x_{j;t'}$ implies $\gg_{i;t}=\gg_{j;t'}$.
\end{proposition}

\begin{proof}
We follow the implication $x_{i;t}=x_{j;t'}\Rightarrow
\gg_{i;t}=\gg_{j;t'}$ in \cite[II. Theorem 7.2]{nakanishi2023scattering}. Regard $t$ as the initial vertex, and use a superscript $t$ for the corresponding $F$-polynomials and $g$-vectors. The generalized separation formula gives
\[
x_{i;t}=x_{j;t'}
=\xx_t^{\gg^{\,t}_{j;t'}}
\frac{F^{\,t}_{j;t'}\vert_{\mathcal F}
(\widehat y_{1;t},\ldots,\widehat y_{n;t})}
{F^{\,t}_{j;t'}\vert_{\mathbb P}
(y_{1;t},\ldots,y_{n;t})}.
\]
By \cite[Proposition 3.19]{nakanishi2015structure} and sign-coherence of $c$-vectors \cite[Theorem 8.13]{amanda2025positivity}, the relative $F$-polynomial $F^{\,t}_{j;t'}$ has constant term $1$. Laurent positivity \cite[Theorem 7.9]{amanda2025positivity} prevents cancellation among the Laurent monomials in the separation formula. Since the left-hand side is the single initial cluster variable $x_{i;t}$, the same argument as in \cite[II. Theorem 7.2]{nakanishi2023scattering} yields
\[
F^{\,t}_{j;t'}=1,
\qquad
\gg^{\,t}_{j;t'}=\ee_i=\gg^{\,t}_{i;t}.
\]
It remains to return to the original initial vertex $t_0$. By Proposition \ref{prop:NakanishiGvectors}, the $g$-vectors under consideration coincide with those of the corresponding ordinary companion cluster pattern. Therefore, the change-of-initial-vertex bijection \cite[II. Proposition 2.25]{nakanishi2023scattering} carries the preceding equality to
\[
\gg^{\,t_0}_{j;t'}=\gg^{\,t_0}_{i;t}.
\]
This is the desired equality.
\end{proof}

\begin{remark}
Proposition \ref{thm:ximpliesg} allows us to define the $g$-vector $\gg_{i;t}$ in a cluster algebra with non-principal coefficients, which was initially determined by the vertex $t$ of the cluster pattern and the position $i$ at that vertex, as an invariant associated with the cluster variable $x_{i;t}$. In other words, the \emph{$g$-vector associated with a cluster variable $x$} is well-defined.
\end{remark}

\begin{remark}
For the reciprocal generalized cluster algebras considered in this paper, the converse of Proposition \ref{thm:ximpliesg} also holds. In the principal-coefficient case, Cheung, Kelley, and Musiker show that every generalized cluster monomial is a theta function \cite[Theorem 4.8]{cheung2023cluster}, and identify the theta-function parameter with the usual $g$-vector \cite[Corollary 4.13]{cheung2023cluster}. Hence two cluster variables with the same $g$-vector are the same theta function. The statement for an arbitrary coefficient system of the same generalized cluster pattern follows by coefficient specialization, equivalently by the separation formula. See also \cite[Theorem 8.12]{mou2024scattering} for a parallel scattering-diagram construction in a slightly different coefficient framework.
\end{remark}

\section{Generalized Markov Cluster Algebras and Fraction Labelings}\label{sec:GenMarkovCA}
\subsection{Setting}\label{subsec:Setting}
In this paper, we will consider the following setting: we fix $k_1,k_2,k_3\in \mathbb Z_{\geq{0}}$, and assume that $k_1\leq k_2\leq k_3$ without loss of generality. For each choice of $k_1,k_2,k_3$, we will define a \emph{generalized Markov cluster algebra}, denoted $\mathcal{A}(k_1,k_2,k_3)$. These will be the generalized cluster algebras  $\mathcal A(B,\mathbf Z)$ with initial seeds specified as below. For all choices of $k_i$, we  set $\PP=\{1\}$ (therefore, we have $\yy=(1,1,1)$). Now, we specify the other pieces of data.

\begin{itemize}\setlength{\leftskip}{-15pt}
\item [(I)] If $k_1=k_2=k_3=0$, we set
\[\xx=(x_1,x_2,x_3),\ B=\begin{bmatrix}
    0&-2&2\\ 2&0&-2\\-2&2&0
\end{bmatrix},\ \begin{cases}Z_1(u)=1+u\\Z_2(u)=1+u\\Z_3(u)=1+u\end{cases},\]
 \item [(II)] If $k_1=k_2=0$ and $k_3\neq 0$, we set
\[\xx=(x_1,x_2,x_3),\ B=\begin{bmatrix}
    0&-2&1\\ 2&0&-1\\-2&2&0
\end{bmatrix}.\ \begin{cases}Z_1(u)=1+u\\Z_2(u)=1+u\\Z_3(u)=1+k_3u+u^2\end{cases}.\]
 \item [(III)] If $k_1=0$ and $k_2,k_3\neq 0$, we set
\[\xx=(x_1,x_2,x_3),\ B=\begin{bmatrix}
    0&-1&1\\ 2&0&-1\\-2&1&0
\end{bmatrix},\ \begin{cases}Z_1(u)=1+u\\Z_2(u)=1+k_2u+u^2\\Z_3(u)=1+k_3u+u^2\end{cases}.\]
 \item [(IV)] If $k_1,k_2,k_3\neq 0$, we set
\[\xx=(x_1,x_2,x_3),\ B=\begin{bmatrix}
    0&-1&1\\ 1&0&-1\\-1&1&0
\end{bmatrix},\ \begin{cases}Z_1(u)=1+k_1u+u^2\\Z_2(u)=1+k_2u+u^2\\Z_3(u)=1+k_3u+u^2\end{cases}.\]
\end{itemize}

Note that if all $k_i =0$, we recover the (ordinary) Markov cluster algebra.
These generalized cluster algebras have a good symmetry. For any $i\in\{1,2,3\}$, we have

\[\mu_i(\pm B)=\mp B,\quad \mu_i(\mathbf Z)=\mathbf Z.\]

In the ordinary case, such matrices are called \emph{sign-equivalent exchange matrices} \cite{chen2025cluster}. This sign-equivalence means that the mutation rule of clusters does not depend on the current seed, and it is always as follows: for $(x_{1;t},x_{2;t},x_{3;t})$,
\begin{align*}
\mu_1(x_{1;t},x_{2;t},x_{3;t})&=\left(\frac{x_{2;t}^2+k_1x_{2;t}x_{3;t}+x_{3;t}^2}{x_{1;t}},x_{2;t},x_{3;t}\right),\\
\mu_2(x_{1;t},x_{2;t},x_{3;t})&=\left(x_{1;t},\frac{x_{1;t}^2+k_2x_{1;t}x_{3;t}+x_{3;t}^2}{x_{2;t}},x_{3;t}\right),\\
\mu_3(x_{1;t},x_{2;t},x_{3;t})&=\left(x_{1;t},x_{2;t},\frac{x_{1;t}^2+k_3x_{1;t}x_{2;t}+x_{2;t}^2}{x_{3;t}}\right).
\end{align*}

\begin{remark}
It is possible to recover the structures (I), (II), and (III) from (IV) by setting various $k_i$ to 0. However, in Sections \ref{sec:Posets} and \ref{sec:MarkovPosets} we will enhance these cluster algebras to have principal coefficients. When we include non-trivial coefficients, there is a difference between, say, mutation structure (I) and mutation structure (IV) with all $k_i = 0$.
\end{remark}

We set
\[\mathscr M(k_1,k_2,k_3) :=\frac{x_1^2+x_2^2+x_3^2+k_1x_2x_3+k_2x_1x_3+k_3x_1x_2}{x_1x_2x_3},\]
and we call it the \emph{$(k_1,k_2,k_3)$-generalized Markov invariant}. When $k_1,k_2,$ and $k_3$ are understood, or are arbitrary values, $\mathscr M(k_1,k_2,k_3)$ is simply denoted as $\mathscr M$.  The following can easily be verified for two clusters one mutation apart, which implies it is true for all clusters.

\begin{proposition}
The $(k_1,k_2,k_3)$-generalized Markov invariant $\mathscr M $ is an invariant on mutations of $\mathcal A(k_1,k_2,k_3)$.
\end{proposition}

It will be important for us to keep track of which position in a cluster each cluster variable appears in.

\begin{theorem}\label{position-variable}
For any cluster variable $x$ in $\mathcal A(k_1,k_2,k_3)$, the position on which $x$ appears in clusters is uniquely determined.
\end{theorem}

For example, by Theorem \ref{position-variable}, if a cluster variable $x$ appears in the first coordinate of a cluster $\mathbf x_t$, then $x$ does not appears the second or third coordinates of any other cluster. We will prove this theorem in the next section. This theorem guarantees that the following definition is well-defined.

\begin{definition}
For a cluster variable $x$, let $i\in \{1,2,3\}$ be the position of $x$ in a cluster in $\mathcal{A}(k_1,k_2,k_3)$, that is, $x$ is the $i$th component of some cluster. Then, we set $k_x:=k_i$ and $k_x$ is called the \emph{parity} of $x$.
\end{definition}
By Theorem \ref{position-variable}, $k_x$ is well-defined. By definition, if a cluster variable $a$ satisfies $k_a=k_i$, $a\in \mathbf x_t$ and $\xx_t=\{a,b,c\}$, then we have $k_\frac{b^2+k_ibc+c^2}{a}=k_i$.

We denote by $\mathbb T(k_1,k_2,k_3)$ the tree given from the correspondence $t\mapsto \xx_t\in\mathcal A(k_1,k_2,k_3)$ for every $t\in\mathbb T_3$.

\begin{example}
We draw the first few vertices of $\mathbb T(0,1,2)$ below.
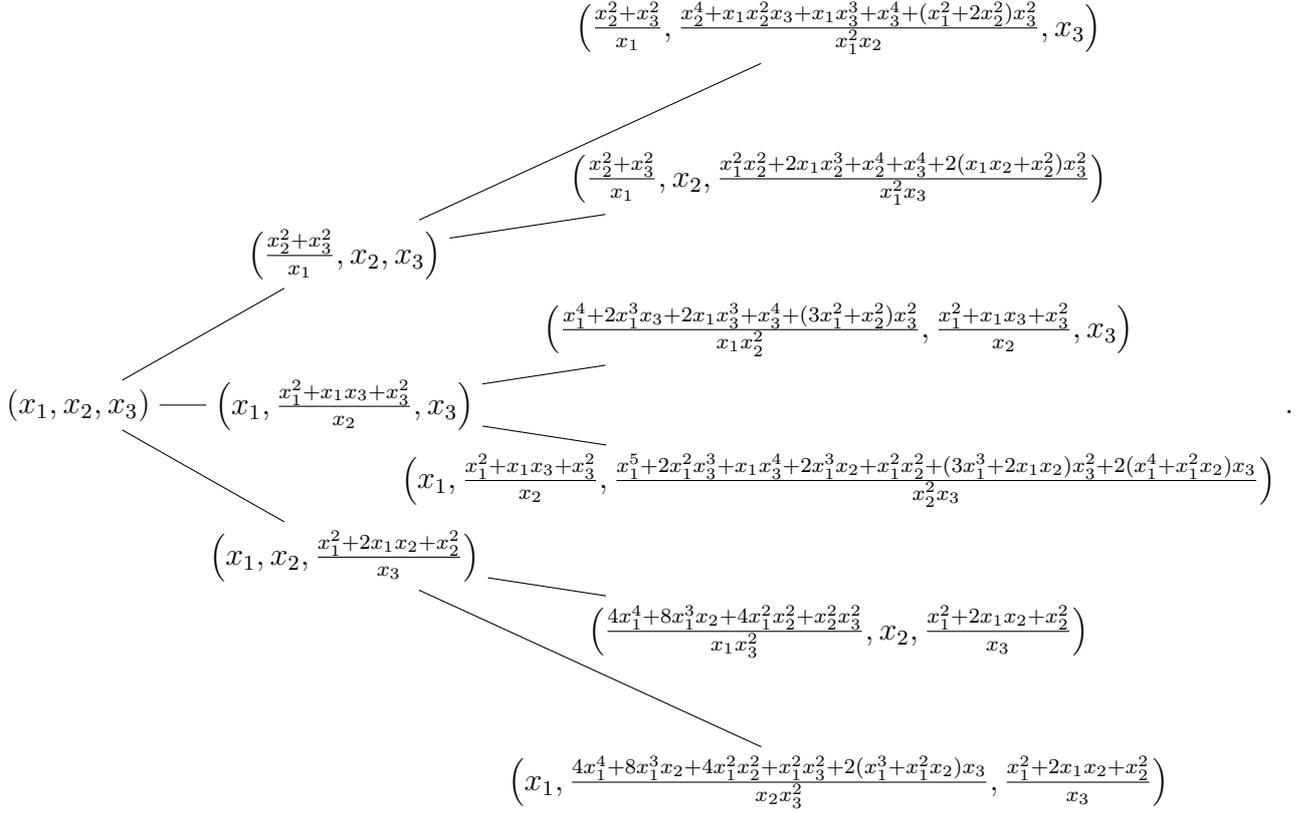
\begin{align*}\label{tree}
\begin{xy}(0,0)*+{(x_1,x_2,x_3)}="0",(35,20)*+{\left(\frac{x_2^2+x_3^2}{x_1},x_2,x_3\right)}="1",(35,0)*+{\left(x_1,\frac{x_1^2+x_1x_3+x_3^2}{x_2},x_3\right)}="1'",(35,-20)*+{\left(x_1,x_2,\frac{x_1^2+2x_1x_2+x_2^2}{x_3}\right)}="1''",(100,50)*+{\left(\frac{x_2^2+x_3^2}{x_1},\frac{x_2^4 + x_1x_2^2x_3 + x_1x_3^3 + x_3^4 + (x_1^2 + 2x_2^2)x_3^2}{x_1^2x_2},x_3\right)}="20",(100,30)*+{\left(\frac{x_2^2+x_3^2}{x_1},x_2,\frac{x_1^2x_2^2 + 2x_1x_2^3 + x_2^4 + x_3^4 + 2(x_1x_2 + x_2^2)x_3^2}{x_1^2x_3}\right)}="21",(100,10)*+{\left(\frac{x_1^4 + 2x_1^3x_3 + 2x_1x_3^3 + x_3^4 + (3x_1^2 + x_2^2)x_3^2}{x_1x_2^2},\frac{x_1^2+x_1x_3+x_3^2}{x_2},x_3\right)}="22",(100,-10)*+{\left(x_1,\frac{x_1^2+x_1x_3+x_3^2}{x_2},\frac{x_1^5 + 2x_1^2x_3^3 + x_1x_3^4 + 2x_1^3x_2 + x_1^2x_2^2 + (3x_1^3 + 2x_1x_2)x_3^2 + 2(x_1^4 + x_1^2x_2)x_3}{x_2^2x_3}\right)}="23",(100,-30)*+{\left(\frac{4x_1^4 + 8x_1^3x_2 + 4x_1^2x_2^2 + x_2^2x_3^2}{x_1x_3^2},x_2,\frac{x_1^2+2x_1x_2+x_2^2}{x_3}\right)}="24",(100,-50)*+{\left(x_1,\frac{4x_1^4 + 8x_1^3x_2 + 4x_1^2x_2^2 + x_1^2x_3^2 + 2(x_1^3 + x_1^2x_2)x_3}{x_2x_3^2},\frac{x_1^2+2x_1x_2+x_2^2}{x_3}\right)}="25", \ar@{-}"0";"1"\ar@{-}"0";"1'"\ar@{-}"0";"1''"\ar@{-}"1";"20"\ar@{-}"1";"21"\ar@{-}"1'";"22"\ar@{-}"1'";"23"\ar@{-}"1''";"24"\ar@{-}"1''";"25"
\end{xy}.
\end{align*}
\end{example}

\begin{definition}
The full subtree of $\mathbb{T}(k_1,k_2,k_3)$ whose root is $\left(\frac{x_2^2+k_1x_2x_3+x_3^2}{x_1},x_2,x_3\right)$ (resp. $\left(x_1,\frac{x_1^2+k_2x_1x_3+x_3^2}{x_2},x_3\right)$,$\left(x_1,x_2,\frac{x_1^2+k_3x_1x_2+x_2^2}{x_3}\right)$) is called the \emph{$k_1$(resp. $k_2$ and $k_3$)-branch of $\mathbb{T}(k_1,k_2,k_3)$}. We denote by $k_1\mathbb{T}(k_1,k_2,k_3)$, $k_2\mathbb{T}(k_1,k_2,k_3)$, $k_3\mathbb{T}(k_1,k_2,k_3)$ these subtrees, respectively.
\end{definition}

\subsection{Cluster Generalized Markov Trees}
To prove Theorem \ref{position-variable}, we will introduce the \emph{$(k_1,k_2,k_3)$-cluster generalized Markov tree} (abbreviated as \emph{$(k_1,k_2,k_3)$-CGM tree}).
Let $\mathfrak S_3$ be the symmetry group of rank $3$. We consider the left action of $\mathfrak S_3$ on $\{1,2,3\}$. We set a permutation $\sigma\in \mathfrak S_3$, and we set
\[\widetilde{\sigma(\xx,\mathbf k)}:=\left((x_{\sigma(1)},k_{\sigma(1)}),\left(\frac{x_{\sigma(1)}^2+k_{\sigma(2)}x_{\sigma(1)}x_{\sigma(3)}+x_{\sigma(3)}^2}{x_{\sigma(2)}},k_{\sigma(2)}\right),(x_{\sigma(3)},k_{\sigma(3)})\right).\]

We define the \emph{$(k_1,k_2,k_3)$-cluster generalized Markov tree} $\mathrm{CM}\mathbb T(k_1,k_2,k_3,\sigma)$ for $\sigma\in \mathfrak{S}_3$ as follows.
\begin{itemize}\setlength{\leftskip}{-15pt}
\item [(1)] The root vertex is $\widetilde{\sigma(\xx,\mathbf k)}$.
\item [(2)] Every vertex $((a,\alpha),(b,\beta),(c,\gamma))$ has the following two children.
\[\begin{xy}(0,0)*+{((a,\alpha),(b,\beta),(c,\gamma))}="1",(-40,-15)*+{\left((a,\alpha),\left(\frac{a^2+\gamma ab+b^2}{c},\gamma\right),(b,\beta)\right)}="2",(40,-15)*+{\left((b,\beta),\left(\frac{b^2+\alpha bc+c^2}{a},\alpha\right),(c,\gamma)\right).}="3", \ar@{-}"1";"2"\ar@{-}"1";"3"
\end{xy}\]
\end{itemize}

\begin{example}\label{ex:CMT(0,1,2,id)}
The first 3 vertices of $\mathrm{CM}\mathbb T(0,1,2,\textrm{id})$ are as follows:
\begin{align*}\label{tree}
\begin{xy}(35,0)*+{\left(x_1,\frac{x_1^2+x_1x_3+x_3^2}{x_2},x_3\right)}="1'",(100,10)*+{\left(\frac{x_1^2+x_1x_3+x_3^2}{x_2},\frac{x_1^4 + 2x_1^3x_3 + 2x_1x_3^3 + x_3^4 + (3x_1^2 + x_2^2)x_3^2}{x_1x_2^2},x_3\right)}="22",(100,-10)*+{\left(x_1,\frac{x_1^5 + 2x_1^2x_3^3 + x_1x_3^4 + 2x_1^3x_2 + x_1^2x_2^2 + (3x_1^3 + 2x_1x_2)x_3^2 + 2(x_1^4 + x_1^2x_2)x_3}{x_2^2x_3},\frac{x_1^2+x_1x_3+x_3^2}{x_2}\right)}="23", \ar@{-}"1'";"22"\ar@{-}"1'";"23"
\end{xy}.
\end{align*}
\end{example}

For $\tau\in \mathfrak S_3$, we denote by $(a,b,c)^\tau$ the triple with $(a,b,c)$ rearranged so that $a$ is the $\tau(1)$th component, $b$ is the $\tau(2)$th component, and $c$ is the $\tau(3)$th component. Given a tree  $\mathbb{T}$, let $V(\mathbb{T})$ denote its set of vertices. We define the map $\pi:V(\mathrm{CM}\mathbb T(k_1,k_2,k_3,\sigma))\to \mathbb Z[x_1^{\pm1},x_2^{\pm1},x_3^{\pm1}]^3$ as follows.
\begin{itemize}\setlength{\leftskip}{-15pt}
    \item We set $\pi(\widetilde{\sigma(\xx,\mathbf k)}) = \left(x_{\sigma(1)},\frac{x_{\sigma(1)}^2+k_{\sigma(2)}x_{\sigma(1)}x_{\sigma(3)}+x_{\sigma(3)}^2}{x_{\sigma(2)}},x_{\sigma(3)}\right)^{\sigma}$.
    \item If $\pi((a,\alpha),(b,\beta),(c,\gamma))=(a,b,c)^\tau$, then
    \begin{align*}
        \pi\left((a,\alpha),\left(\frac{a^2+\gamma ab+b^2}{c},\gamma\right),(b,\beta)\right)= \left(a,\frac{a^2+\gamma ab+b^2}{c},b\right)^{\tau\circ(2\ 3)}, \text{ and} \\
         \pi \left((b,\beta),\left(\frac{b^2+\alpha bc+c^2}{a},\alpha\right),(c,\gamma)\right)= \left(b,\frac{b^2+\alpha bc+c^2}{a},c\right)^{\tau\circ(1\ 2)}.
    \end{align*}
\end{itemize}

In this paper, the composition of permutations is assumed to act from right to left.

\begin{proposition}\label{prop:isomorphism-CMT-kT}
The map $\pi$ induces the following graph isomorphisms.
\begin{align*}
\mathrm{CM}\mathbb T(k_1,k_2,k_3,(1\ 2)) &\simeq k_{1}\TT(k_1,k_2,k_3)\\
\mathrm{CM}\mathbb T(k_1,k_2,k_3,(1\ 3\ 2)) &\simeq k_{1}\TT(k_1,k_2,k_3)\\
\mathrm{CM}\mathbb T(k_1,k_2,k_3,\mathrm{id}) &\simeq k_{2}\TT(k_1,k_2,k_3)\\
\mathrm{CM}\mathbb T(k_1,k_2,k_3,(1\ 3)) &\simeq k_{2}\TT(k_1,k_2,k_3)\\
\mathrm{CM}\mathbb T(k_1,k_2,k_3,(2\ 3)) &\simeq k_{3}\TT(k_1,k_2,k_3)\\
\mathrm{CM}\mathbb T(k_1,k_2,k_3,(1\ 2\ 3)) &\simeq k_{3}\TT(k_1,k_2,k_3)\\
\end{align*}
\end{proposition}

\begin{proof}
We will only prove $\mathrm{CM}\mathbb T(k_1,k_2,k_3,(1\ 2)) \simeq k_{1}\TT(k_1,k_2,k_3)$. We can directly check that $\pi$ gives the correspondence between the roots of $\mathrm{CM}\mathbb T(k_1,k_2,k_3,(1\ 2))$ and $ k_{1}\TT(k_1,k_2,k_3)$. We fix the vertex $((a,\alpha),(b,\beta),(c,\gamma))$ in $\mathrm{CM}\mathbb T(k_1,k_2,k_3,(1\ 2))$, and we assume $\pi$ maps it to $(a,b,c)$ (that is, $\tau$ in the definition of $\pi$ is $\textrm{id}$). Then, by the definition of $\pi$, we have
\begin{align*}
        \pi\left((a,\alpha),\left(\frac{a^2+\gamma ab+b^2}{c},\gamma\right),(b,\beta)\right)&=\left(a,b,\frac{a^2+\gamma ab+b^2}{c}\right),\\
         \pi\left((b,\beta),\left(\frac{b^2+\gamma bc+c^2}{a},\gamma\right),(c,\gamma)\right)&= \left(\frac{b^2+\alpha bc+c^2}{a},b,c\right).
    \end{align*}
The parent of $((a,k_h),(b,k_i),(c,k_j))$ in $\mathrm{CM}\mathbb T(k_1,k_2,k_3,(1\ 2))$ is
\[\left((a,\alpha),(c,\gamma),\left(\frac{a^2+\beta ac+c^2}{b},\beta\right)\right)\  \text{or}\ \left(\left(\frac{a^2+\beta ac+c^2}{c},\beta\right),(a,\alpha),(c,\gamma)\right) \]
By definition of $\pi$, in both cases, the image of $\pi$ of this vertex is  $\left(a,\frac{a^2+\beta ac+c^2}{b},c\right)$.
By checking the image of $\pi$ in the same way for cases where $\tau$ is not the identity permutation,  we can inductively show that the image of $\pi$ matches the generation rule of $k_1\TT(k_1,k_2,k_3)$. One can also show that $\alpha,\beta,\gamma$ appearing in the first (resp. second, third) entry of the image of $\pi$ are $k_1$ (resp. $k_2,k_3$) by inducting on a vertex's distance from the root.
\end{proof}
\begin{remark}\label{rem:mirror}
The isomorphism $\pi$ does not specify left-right for the child, since the tree $k_{i}\TT(k_1,k_2,k_3)$ is a left-right indistinguishable tree with respect to the child. Indeed, there exist distinguished pairs $\sigma,\sigma'$ such that $\sigma(2)=\sigma'(2)$, and the trees $\mathrm{CM}\mathbb T(k_1,k_2,k_3,\sigma)$ and $\mathrm{CM}\mathbb T(k_1,k_2,k_3,\sigma')$ are mirror images of each other.
\end{remark}

For cluster variable $a$, we will define $|a|$ as
$a| _{x_1=x_2=x_3=1}$.  Moreover, for a non-labeled cluster $(a,b,c)$, if $|b|>\max\{|a|,|c|\}$, we have
\[\left|\frac{a^2+k_cab+b^2}{c}\right|> \max\{|a|,|b|\}\quad \text{and}\quad  \left|\frac{b^2+k_abc+c^2}{a}\right|> \max\{|b|,|c|\}\] (see \cite[Proposition 2.4]{gyoda2024sl} or \cite[Proposition 4]{GyodaMatsushita}). By using this fact, Proposition \ref{prop:isomorphism-CMT-kT}, and Remark \ref{rem:mirror}, we have the following corollary.

\begin{corollary}\label{cor:unique-non-labeled-cluster}
Let $(a,b,c)$ be a non-labeled cluster in $\mathcal{A}(k_1,k_2,k_3)$, and its parities are $(k_i,k_j,k_h)$. We assume that $|b|>\max\{|a|,|c|\}.$ Then, there exist a unique $\sigma\in \mathfrak{S}_3$ and a unique vertex $v$ in $\mathrm{CM}\mathbb T(k_1,k_2,k_3,\sigma)$ such that $v=((a,k_i),(b,k_h),(c,k_j))$.
\end{corollary}

\subsection{Inverse Cluster Generalized Markov Trees}

We will introduce the \emph{inverse $(k_1,k_2,k_3)$-cluster generalized Markov tree} (abbreviated as \emph{inverse $(k_1,k_2,k_3)$-CGM tree}). We set a permutation $\sigma\in \mathfrak S_3$, and we set \[\sigma(\xx,\mathbf k):=\left((x_{\sigma(1)},k_{\sigma(1)}),\left(x_{\sigma(2)},k_{\sigma(2)}\right),(x_{\sigma(3)},k_{\sigma(3)})\right).\]

We consider a binary tree, the \emph{$(k_1,k_2,k_3)$-cluster generalized Markov tree} $\mathrm{CM}\mathbb T^\dag(k_1,k_2,k_3,\sigma)$ for $\sigma\in \mathfrak{S}_3$.
\begin{itemize}\setlength{\leftskip}{-15pt}
\item [(1)] The root vertex is $\sigma(\xx,\mathbf k)$.
\item [(2)] Every vertex $((a,\alpha),(b,\beta),(c,\gamma))$ has the following two children.
\[\begin{xy}(0,0)*+{((a,\alpha),(b,\beta),(c,\gamma))}="1",(-40,-15)*+{\left((a,\alpha),(c,\gamma),\left(\frac{a^2+\beta ac+c^2}{b},\beta\right)\right)}="2",(40,-15)*+{\left(\left(\frac{a^2+\beta ac+c^2}{b},\beta\right),\left(a,\alpha\right),(c,\gamma)\right).}="3", \ar@{-}"1";"2"\ar@{-}"1";"3"
\end{xy}\]
\end{itemize}
\begin{example}
The first 7 vertices of $\mathrm{CM}\mathbb T^\dag(0,1,2,\textrm{id})$ are as follows (we omit parities of cluster variables):
\begin{align*}\label{tree}
\begin{xy}(0,0)*+{(x_1,x_2,x_3)}="0",(35,20)*+{\left(\frac{x_1^2+x_1x_3+x_3^2}{x_2},x_1,x_3\right)}="1",(35,-20)*+{\left(x_1,x_3, \frac{x_1^2+x_1x_3+x_3^2}{x_2}\right)}="1'",(100,30)*+{\left(\frac{x_1^4 + 2x_1^3x_3 + 2x_1x_3^3 + x_3^4 + (3x_1^2 + x_2^2)x_3^2}{x_1x_2^2},\frac{x_1^2+x_1x_3+x_3^2}{x_2},x_3\right)}="20",(100,10)*+{\left(\frac{x_1^2+x_1x_3+x_3^2}{x_2},x_3,\frac{x_1^4 + 2x_1^3x_3 + 2x_1x_3^3 + x_3^4 + (3x_1^2 + x_2^2)x_3^2}{x_1x_2^2}\right)}="21",(100,-10)*+{\left(\frac{x_1^5 + 2x_1^2x_3^3 + x_1x_3^4 + 2x_1^3x_2 + x_1^2x_2^2 + (3x_1^3 + 2x_1x_2)x_3^2 + 2(x_1^4 + x_1^2x_2)x_3}{x_2^2x_3},x_1,\frac{x_1^2+x_1x_3+x_3^2}{x_2}\right)}="22",(100,-30)*+{\left(x_1,\frac{x_1^2+x_1x_3+x_3^2}{x_2},\frac{x_1^5 + 2x_1^2x_3^3 + x_1x_3^4 + 2x_1^3x_2 + x_1^2x_2^2 + (3x_1^3 + 2x_1x_2)x_3^2 + 2(x_1^4 + x_1^2x_2)x_3}{x_2^2x_3}\right)}="23",\ar@{-}"0";"1"\ar@{-}"0";"1'"\ar@{-}"1";"20"\ar@{-}"1";"21"\ar@{-}"1'";"22"\ar@{-}"1'";"23"
\end{xy}.
\end{align*}
\end{example}

In the same way as \cite[Proposition 3.5]{gyoda2024sl}, we have the relation between $\mathrm{CM}\mathbb T(k_1,k_2,k_3,\sigma)$ and $\mathrm{CM}\mathbb T^\dag(k_1,k_2,k_3,\sigma)$. Before describing the relation, we will introduce the \emph{canonical graph isomorphism} between two trees.

\begin{definition}
Let $\mathbb T$ and $\mathbb T'$ be full planar binary trees. Suppose $f \colon \mathbb T \to \mathbb T'$ is a graph isomorphism such that, for all $v \in V(\mathbb T)$, if $\ell$ and $r$ are the left and right children of $f$ respectively, then $f(\ell)$ and $f(r)$ are the left and right children of $f(v)$. In this case, we say $f$ is the \emph{canonical graph isomorphism}.
\end{definition}

\begin{proposition}\label{pr:rho-mor}
The correspondence \[\mu\colon((a,\alpha),(b,\beta),(c,\gamma))\mapsto\left((a,\alpha),\left(\frac{a^2+\beta ac+c^2}{b},\beta\right),(c,\gamma)\right)\] induces the canonical graph isomorphism from $\mathrm{CM}\mathbb T(k_1,k_2,k_3,\sigma)$ to $\mathrm{CM}\mathbb T^\dag(k_1,k_2,k_3,\sigma)$.
\end{proposition}

\begin{remark}
Since the correspondence $\mu$ is a convolution, it also induces the isomorphism from $\mathrm{CM}\mathbb T^{\dag}(k_1,k_2,k_3,\sigma)$ to $\mathrm{CM}\mathbb T(k_1,k_2,k_3,\sigma)$.
\end{remark}

By using Proposition \ref{pr:rho-mor}, we have the following corollary:

\begin{corollary}\label{cor:unique-non-labeled-cluster-2}
Let $(a,b,c)$ be a non-labeled cluster in $\mathcal{A}(k_1,k_2,k_3)$ and their parities are $(k_h,k_i,k_j)$. We assume that $|b|\leq\max\{|a|,|c|\}.$ Then, there exist a unique $\sigma\in \mathfrak{S}_3$ and a unique vertex $v$ in $\mathrm{CM}\mathbb T^\dag(k_1,k_2,k_3,\sigma)$ such that $v=((a,k_h),(b,k_i),(c,k_j))$.
\end{corollary}

\begin{proof}
Since $|b|\leq\max\{|a|,|c|\}$, we have $\left|\frac{a^2+k_iac+c^2}{b}\right|> \max\{|a|,|c|\}$ (see \cite[Proposition 2.4]{gyoda2024sl}). Therefore, by Corollary \ref{cor:unique-non-labeled-cluster}, there exist a unique $\sigma \in \mathfrak S_3$ and a unique vertex $v$ in $\mathrm{CM}\mathbb T(k_1,k_2,k_3,\sigma)$ such that $v=\mu((a,k_h),(b,k_i),(c,k_j))$. Since
\[((a,k_h),(b,k_i),(c,k_j))=\mu\circ\mu((a,k_h),(b,k_i),(c,k_j)),\]
we conclude the statement by Proposition \ref{pr:rho-mor}.
\end{proof}
Combining Corollaries \ref{cor:unique-non-labeled-cluster} and \ref{cor:unique-non-labeled-cluster-2}, the set
\[\{v\in \mathrm{CM}\mathbb T(k_1,k_2,k_3,\sigma)\mid \sigma\in \mathfrak S_3\}\cup \{v\in \mathrm{CM}\mathbb T^\dag(k_1,k_2,k_3,\sigma)\mid \sigma\in \mathfrak S_3\}\]
coincides with the set of all non-labeled clusters (and its parities) in $\Acal(k_1,k_2,k_3)$.
\subsection{Fraction Labeling of Cluster Variables}\label{subsec:LabelClVar}
In this subsection, we recall the Farey tree, and we label cluster variables with irreducible fractions.

\begin{definition}
For $\frac{a}{b}$ and $\frac{c}{d}$, we denote $ad-bc$ by $\det\left(\frac{a}{b},\frac{c}{d}\right)$. A triple $\left(\frac{a}{b},\frac{c}{d},\frac{e}{f}\right)$ is called a \emph{Farey triple} if
\begin{itemize}\setlength{\leftskip}{-15pt}
    \item [(1)] $\frac{a}{b},\frac{c}{d}$ and $\frac{e}{f}$ are irreducible fractions and
    \item [(2)] $\left|\det\middle(\frac{a}{b},\frac{c}{d}\middle)\middle|=\middle|\det\middle(\frac{c}{d},\frac{e}{f}\middle)\middle|=\middle|\det\middle(\frac{e}{f},\frac{a}{b}\middle)\right|=1$.
\end{itemize}
\end{definition}
We define the \emph{Farey tree} $\mathrm{F}\mathbb T$ as follows.
\begin{itemize}\setlength{\leftskip}{-15pt}
\item [(1)] The root vertex is $\left(\frac{0}{1},\frac{1}{1},\frac{1}{0}\right)$.
\item[(2)]Every vertex $\left(\frac{a}{b},\frac{c}{d},\frac{e}{f}\right)$ has the following two children.
\[\begin{xy}(0,0)*+{\left(\dfrac{a}{b},\dfrac{c}{d},\dfrac{e}{f}\right)}="1",(-30,-15)*+{\left(\dfrac{a}{b},\dfrac{a+c}{b+d},\dfrac{c}{d}\right)}="2",(30,-15)*+{\left(\dfrac{c}{d},\dfrac{c+e}{d+f},\dfrac{e}{f}\right).}="3", \ar@{-}"1";"2"\ar@{-}"1";"3"
\end{xy}\]
\end{itemize}
The first few vertices of $\mathrm{F}\mathbb{T}$ are given by the following.
\relsize{+1}
\begin{align*}
\begin{xy}(0,0)*+{\left(\frac{0}{1},\frac{1}{1},\frac{1}{0}\right)}="1",(20,-14)*+{\left(\frac{0}{1},\frac{1}{2},\frac{1}{1}\right)}="2",(20,14)*+{\left(\frac{1}{1},\frac{2}{1},\frac{1}{0}\right)}="3",
(50,-24)*+{\left(\frac{0}{1},\frac{1}{3},\frac{1}{2}\right)}="4",(50,-8)*+{\left(\frac{1}{2},\frac{2}{3},\frac{1}{1}\right)}="5",(50,8)*+{\left(\frac{1}{1},\frac{3}{2},\frac{2}{1}\right)}="6",(50,24)*+{\left(\frac{2}{1},\frac{3}{1},\frac{1}{0}\right)}="7",(85,-28)*+{\left(\frac{0}{1},\frac{1}{4},\frac{1}{3}\right)\cdots}="8",(85,-20)*+{\left(\frac{1}{3},\frac{2}{5},\frac{1}{2}\right)\cdots}="9",(85,-12)*+{\left(\frac{1}{2},\frac{3}{5},\frac{2}{3}\right)\cdots}="10",(85,-4)*+{\left(\frac{2}{3},\frac{3}{4},\frac{1}{1}\right)\cdots}="11",(85,4)*+{\left(\frac{1}{1},\frac{4}{3},\frac{3}{2}\right)\cdots}="12",(85,12)*+{\left(\frac{3}{2},\frac{5}{3},\frac{2}{1}\right)\cdots}="13",(85,20)*+{\left(\frac{2}{1},\frac{5}{2},\frac{3}{1}\right)\cdots}="14",(85,28)*+{\left(\frac{3}{1},\frac{4}{1},\frac{1}{0}\right)\cdots}="15",\ar@{-}"1";"2"\ar@{-}"1";"3"\ar@{-}"2";"4"\ar@{-}"2";"5"\ar@{-}"3";"6"\ar@{-}"3";"7"\ar@{-}"4";"8"\ar@{-}"4";"9"\ar@{-}"5";"10"\ar@{-}"5";"11"\ar@{-}"6";"12"\ar@{-}"6";"13"\ar@{-}"7";"14"\ar@{-}"7";"15"
\end{xy}
\end{align*}
\relsize{-1}
\begin{proposition}[see {\cite[Section 3.2]{aigner2013Markov}}]\label{prop:property-farey}\indent
\begin{itemize}\setlength{\leftskip}{-15pt}
    \item [(1)]If $\left(\frac{a}{b},\frac{c}{d},\frac{e}{f}\right)$ is a Farey triple, then so are $\left(\frac{a}{b},\frac{a+c}{b+d},\frac{c}{d}\right)$ and $\left(\frac{c}{d},\frac{c+e}{d+f},\frac{e}{f}\right)$. In particular, each vertex in $\mathrm{F}\mathbb T$ is a Farey triple.
    \item [(2)] For every irreducible fraction $\frac{a}{b} \in \mathbb Q_{> 0}$, there exists a unique Farey triple $F$ in $\mathrm{F}\mathbb T$ such that $\frac{a}{b}$ is the second entry of $F$.
    \item [(3)] All Farey triples $\left(\frac{a}{b},\frac{c}{d},\frac{e}{f}\right)$ in $\mathrm{F}\mathbb T$ satisfy $\frac{a}{b}<\frac{c}{d}<\frac{e}{f}$.
\end{itemize}
\end{proposition}

The canonical graph isomorphisms from the Farey tree to the cluster generalized Markov trees provide a correspondence from positive rational numbers to pairs consisting of a cluster variable and its parity. First, consider the canonical graph isomorphsim from $\mathrm{F}\mathbb{T}$ to $\mathrm{CM}\mathbb{T}(k_1,k_2,k_3,(1\ 3\ 2))$. If this graph isomorphism maps $(\frac{p}{q}, \frac{p+r}{q+s},\frac{r}{s})$ to $((a,\alpha),(b,\beta),(c,\gamma))$, we denote $b = x_{1,\frac{p+r}{q+s}}$ and $\beta = k_{1,\frac{p+r}{q+s}}$. We extend this labeling at the root and set $x_{1,\frac01} = x_3$ and $k_{1,\frac01} = k_3$.
Using Proposition \ref{prop:isomorphism-CMT-kT}, this correspondence induces the map from irreducible fractions in $[0,\infty)$ to pairs of a cluster variable in $k_1\mathbb{T}(k_1,k_2,k_3)$ and its parity.
This map is called the \emph{fraction labeling to cluster variables}. In the same way, by considering the canonical isomorphism from $\mathrm{F}\mathbb T$ to $\mathrm{CM}\mathbb T(k_1,k_2,k_3,\mathrm{id})$ (resp. $\mathrm{CM}\mathbb T(k_1,k_2,k_3,(1\ 2\ 3))$), we can construct the map from irreducible fractions in $[0,\infty)$ to pairs of cluster variables in $k_2\mathbb{T}(k_1,k_2,k_3)$ (resp. $k_3\mathbb{T}(k_1,k_2,k_3)$) and their parities. We denote by $(x_{2,{\frac{p}{q}}},k_{2,{\frac{p}{q}}})$ (resp. $(x_{3,{\frac{p}{q}}},k_{3,{\frac{p}{q}}})$) the pair of a cluster variable and it parity labeled with a fraction $\frac{p}{q}$. We again extend the correspondences to include $\frac01$ by setting $x_{2,\frac01} = x_1$ and $x_{3,\frac01} = x_2$, with the $k_i$ defined accordingly.  We have
\[\mathcal X=\{x_{i,\frac{p}{q}}\mid i\in\{1,2,3\},\ \frac{p}{q}\in [0,\infty)\cap\mathbb Q\}\]
where $\Xcal$ is the set of cluster variables.

\begin{remark}
 Just as we extended the fractional labeling to $\frac01$, we could also extend to $\frac10$. Explicitly, this would yield \[
x_{1,\frac10} = x_2 \qquad x_{2,\frac10} = x_3 \qquad x_{3,\frac10} = x_1.
 \]
 This labeling can be convenient at times, but many arguments will fail to extend to the $\frac10$ case. For example, see Remark \ref{rem:1/0badg-vector}. For this reason, later in Sections \ref{sec:CombinatorialCohnMatrices} and \ref{sec:CombinatorialMM}, matrices labeled $\frac10$ will be treated separately.
\end{remark}

\begin{example}
We give a few examples of the fraction labeling of cluster variables in $\Acal(0,1,2)$. Considering the correspondence between $\mathrm{F}\mathbb{T}$ and $\mathrm{CM}\mathbb T(0,1,2,\textrm{id})$ (Example \ref{ex:CMT(0,1,2,id)}), we have the following.
\begin{align*}
 x_{2,\frac{1}{1}}&=\frac{x_1^2+x_1x_3+x_3^2}{x_2},\\
 x_{2,\frac{1}{2}}&=\frac{x_1^5 + 2x_1^2x_3^3 + x_1x_3^4 + 2x_1^3x_2 + x_1^2x_2^2 + (3x_1^3 + 2x_1x_2)x_3^2 + 2(x_1^4 + x_1^2x_2)x_3}{x_2^2x_3},\\
 x_{2,\frac{2}{1}}&=\frac{x_1^4 + 2x_1^3x_3 + 2x_1x_3^3 + x_3^4 + (3x_1^2 + x_2^2)x_3^2}{x_1x_2^2}.
\end{align*}
\end{example}
\begin{remark}\label{rem:Permutations}
In this paper, we use the permutation $(1\ 3\ 2), \mathrm{id},(1\ 2\ 3)\in \mathfrak{S}_3$ when we give the fraction labeling to cluster variables (and its parities).
The paper \cite{chavez2011c} also defines the same fraction labeling for cluster variables in $\Acal(0,0,0)$, but in a way that maps fractions to cluster patterns, as opposed to our method of mapping cluster variables to the Farey tree.

We can use $(1\ 2), (1\ 3), (2\ 3)\in \mathfrak{S}_3$ instead of $(1\ 3\ 2), \mathrm{id},(1\ 2\ 3)$, and in this case, the correspondence between irreducible fractions and cluster variables changes. If $(\tilde x_{i,\frac{p}{q}},\tilde{k}_{i,\frac{p}{q}})$ is the fraction labeling given by$(1\ 2), (1\ 3),(2\ 3)$, then these labelings are related by  $(\tilde x_{i,\frac{p}{q}},\tilde{k}_{i,\frac{p}{q}})=(x_{i,\frac{q}{p}},k_{i,\frac{q}{p}})$.
\end{remark}

We denote by $\gg_{i,\frac{p}{q}}$ the $g$-vector of $x_{i,\frac{p}{q}}$.

\begin{theorem}\label{thm:g-vector-description}
For any generalized Markov cluster algebra $\Acal(k_1,k_2,k_3)$ and $\frac{p}{q}\in [0,\infty)\cap\mathbb Q$, we have the following.
\begin{align*}
\gg_{1,\frac{p}{q}}=\begin{bmatrix}q-p-1\\-q+1\\p+1\end{bmatrix},\quad
\gg_{2,\frac{p}{q}}=\begin{bmatrix}p+1\\q-p-1\\-q+1\end{bmatrix},\quad
\gg_{3,\frac{p}{q}}=\begin{bmatrix}-q+1\\p+1\\q-p-1\end{bmatrix}.
\end{align*}
In particular, $\gg_{i,\frac{p}{q}}=\gg_{j,\frac{p'}{q'}}$ implies $i=j$ and $\frac{p}{q}=\frac{p'}{q'}$.
\end{theorem}

\begin{proof}
When $k_1=k_2=k_3=0$, this statement is due to\cite[Theorem 3.1.5]{chavez2011c} . For all other triples, we can use Proposition \ref{prop:NakanishiGvectors}, since for all possible matrices $B$, $DB$ is equal to the matrix associated to the Markov cluster algebra.
\end{proof}

In our setting, this means that the set of $g$-vectors is the same regardless of our choices of $k_i$.

\begin{remark}\label{rem:1/0badg-vector}
Note that this result cannot be applied by $x_{i,\frac10}$. For example, we have $\gg_{2,\frac{1}{0}}=(
0,0,1)^\top$ since $x_{2,\frac10} = x_3$. On the other hand, if we apply Theorem \ref{thm:g-vector-description} to $p=1,q=0$, we have $\gg_{2,\frac{1}{0}}=
(2,-2,1)^\top$.
\end{remark}

 Proposition \ref{thm:ximpliesg} and Theorem \ref{thm:g-vector-description} yield the following corollary.

\begin{corollary}\label{cor:uniqueness-of-fraction}
For any $i,j\in\{1,2,3\}$ and $\frac{p}{q}, \frac{p'}{q'}\in[0,\infty)\cap\mathbb Q$, $x_{i,\frac{p}{q}}=x_{j,\frac{p'}{q'}}$ implies $i=j$ and $\frac{p}{q}=\frac{p'}{q'}$.
\end{corollary}

\begin{proof}[Proof of Theorem \ref{position-variable}]
From the results of Proposition \ref{prop:isomorphism-CMT-kT} and Corollary \ref{cor:uniqueness-of-fraction}, we can see that all the newly given cluster variables in $\mathbb T(k_1,k_2,k_3)$ are non-initial cluster variables, and they are different when mutating in a direction that increases the distance from $t_0$. This means that the cluster variables that disappear when mutating in a direction that increases the distance from $t_0$ do not reappear. From this fact, we can see that the position of each cluster variable in the labeled cluster that contains it is uniquely determined.
\end{proof}

N\'ajera Chavez classified the $c$-vectors of $\mathcal{A}(0,0,0)$,  that is,  the Markov cluster algebra. We update the statement of the result to fit the conventions used here.

Let $F=(\frac{p}{q},\frac{p+r}{q+s},\frac{r}{s})\in \mathrm{F}\mathbb T$ be a Farey triple.
We denote by $\xx_{i,F}$ the non-labeled cluster $\{x_{i,\frac{p}{q}},x_{i,\frac{p+r}{q+s}},x_{i,\frac{r}{s}}\}$. We denote by $\cc_{i,\frac{p}{q},F},\cc_{i,\frac{p+r}{q+s},F},\cc_{i,\frac{r}{s},F}$ the corresponding $c$-vectors.

\begin{proposition}[{\cite[Theorem 3.1.2]{chavez2011c}}]\label{prop:MarkovCVectors}
For a cluster $\xx_{i,F}$ be a non-initial cluster in $\mathcal{A}(0,0,0)$,  where $F=(\frac{p}{q},\frac{p+r}{q+s},\frac{r}{s})\in \mathrm{F}\mathbb T$, we have the following:
\begin{itemize}\setlength{\leftskip}{-15pt}
\item [(1-I)]if $i=1$ and $\frac{r}{s}=\frac{1}{0}$ (that is, $F=(\frac{p}{1},\frac{p+1}{1},\frac{1}{0})$),
\[
\cc_{1,\frac{p}{1},F} = \begin{bmatrix} p+2\\0\\p+1\end{bmatrix}, \quad
\cc_{1,\frac{p+1}{1},F} =  \begin{bmatrix} -p-1\\0\\-p\end{bmatrix},\quad
\cc_{1,\frac{1}{0},F} = \begin{bmatrix}0\\1\\0\end{bmatrix},
\]
\item [(1-II)]if $i=1$ and $\frac{r}{s}\neq\frac{1}{0}$,
\[
\cc_{1,\frac{p}{q},F} = \begin{bmatrix} p+q+1\\q+1\\p+1\end{bmatrix}, \quad
\cc_{1,\frac{p+r}{q+s},F} =  \begin{bmatrix} r+s-p-q-1\\s-q-1\\r-p-1\end{bmatrix},\quad
\cc_{1,\frac{r}{s},F} = \begin{bmatrix}-r-s+1\\-s+1\\-r+1\end{bmatrix},
\]
\item [(2-I)]if $i=2$ and $\frac{r}{s}=\frac{1}{0}$ (that is, $F=(\frac{p}{1},\frac{p+1}{1},\frac{1}{0})$),
\[
\cc_{2,\frac{p}{1},F} = \begin{bmatrix} p+1\\p+2\\0\end{bmatrix}, \quad
\cc_{2,\frac{p+1}{1},F} =  \begin{bmatrix} -p\\-p-1\\0\end{bmatrix},\quad
\cc_{2,\frac{1}{0},F} = \begin{bmatrix}0\\0\\1\end{bmatrix},
\]
\item [(2-II)]if $i=2$ and $\frac{r}{s}\neq\frac{1}{0}$,
\[
\cc_{2,\frac{p}{q},F} = \begin{bmatrix} p+1\\p+q+1\\q+1\end{bmatrix}, \quad
\cc_{2,\frac{p+r}{q+s},F} =  \begin{bmatrix} r-p-1\\r+s-p-q-1\\s-q-1\end{bmatrix},\quad
\cc_{2,\frac{r}{s},F} = \begin{bmatrix}-r+1\\-r-s+1\\-s+1\end{bmatrix},
\]
\item [(3-I)]if $i=3$ and $\frac{r}{s}=\frac{1}{0}$ (that is, $F=(\frac{p}{1},\frac{p+1}{1},\frac{1}{0})$),
\[
\cc_{3,\frac{p}{1},F} = \begin{bmatrix} 0\\p+1\\p+2\end{bmatrix}, \quad
\cc_{3,\frac{p+1}{1},F} =  \begin{bmatrix} 0\\-p\\-p-1\end{bmatrix},\quad
\cc_{3,\frac{1}{0},F} = \begin{bmatrix}1\\0\\0\end{bmatrix},
\]
\item [(3-II)]if $i=3$ and $\frac{r}{s}\neq\frac{1}{0}$,
\[
\cc_{3,\frac{p}{q},F} = \begin{bmatrix} q+1\\p+1\\p+q+1\end{bmatrix}, \quad
\cc_{3,\frac{p+r}{q+s},F} =  \begin{bmatrix} s-q-1\\r-p-1\\r+s-p-q-1\end{bmatrix},\quad
\cc_{3,\frac{r}{s},F} = \begin{bmatrix}-s+1\\-r+1\\-r-s+1\end{bmatrix}.
\]
\end{itemize}
\end{proposition}

We can combine Propositions \ref{prop:NakanishiCvectors} and \ref{prop:MarkovCVectors} to give $c$-vectors for all cluster algebras $\mathcal{A}(k_1,k_2,k_3)$. Note that for any choices of $k_i$, there is a canonical isomorphism between the cluster pattern of $\mathcal{A}(k_1,k_2,k_3)$ and $\mathcal{A}(0,0,0)$.

\begin{lemma}\label{lem:CVectorsGenMarkov}
Let $\xx_t$ be a cluster in $\mathcal{A}(k_1,k_2,k_3)$ with principal coefficients. Let  $C_t^{B;t_0} = (c_{ij;t})_{1 \leq i,j \leq 3}$ denote the $C$-matrix of this seed and let $C_t^{B_M;t_0} =  (c^M_{ij;t})_{1 \leq i,j \leq 3}$ denote the $C$-matrix of the associated seed  in $\mathcal{A}(0,0,0)$. Then,  $c_{ij;t} = \frac{d_i}{d_j} c_{ij;t}^M$.
\end{lemma}

\begin{proof}
Suppose $k_1 \leq k_2 \leq k_3$ without loss of generality.  If $d_1 = d_2 = d_3$,  then $DB$ is equal to the Markov $B$-matrix and the statement follows immediately from Proposition \ref{prop:NakanishiCvectors}.

Suppose that we are in one of the other two cases. In particular, $k_1 = 0$ and $k_3 > 0$, implying $d_1 = 1$ and $d_3 = 2$.  Let $R$ be a matrix such that $RDB$ is skew-symmetric. Let $G_t^{B;t_0}$ denote the $G$-matrix at this seed.  By \cite[Proposition 3.21]{nakanishi2015structure},  we have \[
 (DR)^{-1} (G_t^{B;t_0})^\top(DR) C_t^{B;t_0} = I
\]
where $I$ denotes the $3 \times 3$ identity matrix.  Since $BD$ is equal to the Markov $B$-matrix, from Proposition \ref{prop:NakanishiGvectors} we know $G_t^{B;t_0} = G_t^{B_M;t_0}$, the $G$-matrix of the corresponding seed in the Markov cluster pattern.  Since the Markov cluster algebra is skew-symmetric, from \cite[Equation (3.11)]{nakanishi2012tropical} we know $(G_t^{B_M;t_0})^\top = (C_t^{B_M;t_0})^{-1}$.  Therefore,  we can determine $C_t^{B;t_0}$ from the equation \[
C_t^{B;t_0} = (DR)^{-1} C_t^{B_M;t_0} (DR).
\]

In particular, if $d_1 = d_2 = 1$ and $d_3 = 2$,   we can choose $R = \text{diag}(4,4,1)$ and if $d_1 = 1$ and $d_2 = d_3 = 2$, we can choose $R = \text{diag}(4,1,1)$.  In each case, the claim follows from the effect of conjugating $C_t^{B_M;t_0}$ by the appropriate diagonal matrix.
\end{proof}

Finally, we analyze how the $B$-matrix entries compare with the fraction labeling of cluster variables. In the following, we extend the fractional labeling of cluster variables to a labeling of clusters with Farey triples.

\begin{lemma}\label{lem:BMatrixUsingLabeling}
Consider a non-initial seed $(\mathbf{x}_t,\mathbf{y}_t,B_t,\mathbf Z_t)$ in $\Acal(k_1,k_2,k_3)$ and let $i$ be such that this triple also is in $k_i\mathbb T(k_1,k_2,k_3)$. Let $(\frac{p}{q},\frac{p+r}{q+s},\frac{r}{s})$ be the Farey triple which, up to permutation, is associated to this cluster. If $b_{\frac{p}{q},\frac{p+r}{q+s}}$ denotes the entry in $B_t$ in the row associated to the cluster variable $x_{i,\frac{p}{q}}$ and in the column associated to the cluster variable $x_{i,\frac{p+r}{q+s}}$, then $b_{\frac{p}{q}, \frac{p+r}{q+s}}, b_{\frac{p+r}{q+s}, \frac{r}{s}},$ and  $b_{\frac{r}{s}, \frac{p}{q}}$ are all positive.
\end{lemma}

\begin{proof}
We first set $i=2$. Recall we can index as \[\left(x_1,\frac{x_1^2+k_2x_1x_3+x_3^2}{x_2},x_3\right) = (x_{2,\frac{0}{1}},x_{2,\frac{1}{1}},x_{2,\frac{1}{0}}),\]
associated to the Farey triple $(\frac{0}{1},\frac{1}{1},\frac{1}{0})$. From our choice of initial exchange matrix (see Section 2.1), the corresponding exchange matrix is one of the following.
\[\text{(I) }\begin{bmatrix}
    0&2&-2\\-2&0&2\\2&-2&0
\end{bmatrix},\   \text{(II) }\begin{bmatrix}
    0&2&-1\\-2&0&1\\2&-2&0
\end{bmatrix},\  \text{(III) }\begin{bmatrix}
    0&1&-1\\-2&0&1\\2&-1&0
\end{bmatrix},\ \text{(IV) }\begin{bmatrix}
    0&1&-1\\-1&0&1\\1&-1&0
\end{bmatrix}\]
We see the claim is true in this case.

Assume for a Farey triple $(\frac{p}{q},\frac{p+r}{q+s},\frac{r}{s})$ we have shown our claim.  Recall that $\mu(B) = -B$ for any single mutation $\mu$.  Suppose we mutate $x_{\frac{r}{s}}$. Then,  the row and column associated to $\frac{r}{s}$ will now be associated to $\frac{2p+r}{2q+s}$.  By definition of mutation, we know that $b_{\frac{p+r}{q+s},\frac{2p+r}{2q+s}} < 0$,  implying $b_{\frac{2p+r}{2q+s},\frac{p+r}{q+s}} > 0$. Since $\fpq<\frac{2p+r}{2q+s} < \frac{p+r}{q+s} $,  we see the claim holds for this element.  The sign pattern of the matrices $\pm B$ guarantees the other entries have the desired signs as well. Checking the statement is true for other values of $i$ is similar.
\end{proof}

\section{Cluster Generalized Cohn Matrices}\label{sec:ClusterCohn}
In this and the next section, we will give two matrixizations of the cluster variables of generalized cluster algebra introduced in the previous section which preserve their combinatorial structure. In this section, we will introduce the cluster generalized Cohn matrix.

We recall \[\mathscr{M}=\frac{x_1^2+x_2^2+x_3^2+k_1x_2x_3+k_2x_1x_3+k_3x_1x_2}{x_1x_2x_3}.\]

\begin{definition}\label{def:gen-Cohn-matrix}
For $k_1,k_2,k_3\in \mathbb {Z}_{\geq 0}$, we define a \emph{$(k_1,k_2,k_3)$-cluster generalized Cohn matrix} (abbreviated as \emph{$(k_1,k_2,k_3)$-CGC matrix}) $P=\begin{bmatrix}p_{11}&p_{12}\\p_{21}&p_{22}\end{bmatrix}$  as a matrix satisfying the following conditions:
\begin{itemize}\setlength{\leftskip}{-15pt}
    \item [(1)] $P\in SL(2,\mathbb{Z}[x_1^{\pm1},x_2^{\pm1},x_3^{\pm1}])$,
    \item [(2)] $p_{12}$ is a cluster variable in the $\Acal(k_1,k_2,k_3)$, and
    \item [(3)] $\mathrm{tr}(P)=\mathscr M p_{12}-k_{p_{12}}.$
\end{itemize}
\end{definition}
For $P\in M(2,\mathbb Z[x_1^{\pm1},x_2^{\pm1},x_3^{\pm1}])$, we define \[S_P:=\begin{bmatrix}
    \mathscr M p_{12}-\mathrm{tr}(P)&0\\(\mathscr M p_{12}-\mathrm{tr}(P))\mathscr M  &\mathscr M p_{12}-\mathrm{tr}(P)
    \end{bmatrix}.\] We note that if $P$ is $(k_{1},k_{2},k_{3})$-CGC matrix, then we have $S_P=\begin{bmatrix}
   k_{p_{12}}&0\\k_{p_{12}}\mathscr M  &k_{p_{12}}
    \end{bmatrix}$.
\begin{definition}\label{def:gen-Cohn-triple}
For $(k_1,k_2,k_3)\in \mathbb {Z}_{\geq 0}$, we define a \emph{$(k_1,k_2,k_3)$-cluster generalized Cohn triple} (abbreviated as \emph{$(k_1,k_2,k_3)$-CGC triple}) $(P,Q,R)$ as a triple satisfying the following conditions:
\begin{itemize}\setlength{\leftskip}{-15pt}
    \item [(1)] $P,Q,R$ are $(k_1,k_2,k_3)$-CGC matrices,
    \item[(2)] $Q=PR-S_Q$, and
    \item[(3)] $(p_{12},q_{12},r_{12})$ is a non-labeled cluster in $\Acal(k_1,k_2,k_3)$, where $p_{12},q_{12},r_{12}$ are the $(1,2)$-entries of $P,Q,R$, respectively.
\end{itemize}
\end{definition}

For a non-labeled cluster $(a,b,c)$, CGC-triples satisfying $(p_{12},q_{12},r_{12})=(a,b,c)$  are called those \emph{associated with $(a,b,c)$}.

The definition of the $(k_1,k_2,k_3)$-CGC triple does not refer to the existence of the triple satisfying their conditions. First, we demonstrate the existence of such triples where $(p_{12},q_{12},r_{12})$ is $(x_1,x_2,x_3)$ or one of its permutations. We remark that even if $(P,Q,R)$ is a $(k_1,k_2,k_3)$-CGC triple, a permutation of these matrices is not necessarily a $(k_1,k_2,k_3)$-CGC triple. Therefore, we need to consider all possible permutations of $(x_1,x_2,x_3)$. First, we will consider the case $(p_{12},q_{12},r_{12})=(x_1,x_2,x_3)$. We set
\begin{align*}
    P_{\mathrm{id}}(f)&=\begin{bmatrix}
        f&x_1\\\ast&\frac{-fx_2x_3+k_2x_1x_3+k_3x_1x_2+x_1^2+x_2^2+x_3^2}{x_2x_3}
    \end{bmatrix},\\
    Q_{\mathrm{id}}(f)&=\begin{bmatrix}
        \frac{fx_2+k_1x_2+x_3}{x_1}&x_2\\ *& \frac{-fx_2x_3+k_3x_1x_2+x_1^2+x_2^2}{x_1x_3}
    \end{bmatrix},\\
     R_{\mathrm{id}}(f)&=\begin{bmatrix}
        \frac{fx_2x_3+k_1x_2x_3+k_2x_1x_3+x_1^2+x_3^2}{x_1x_2}&x_3\\ *& \frac{-fx_2x_3+x_2^2}{x_1x_2}
    \end{bmatrix},
\end{align*}
where $f\in \mathbb Z[x_1^{\pm1},x_2^{\pm1},x_3^{\pm1}]$, and we omit $(2,1)$-entries of $P_{\mathrm{id}}(f), Q_{\mathrm{id}}(f),R_{\mathrm{id}}(f)$, which are determined by $\det P_{\mathrm{id}}(f)=\det Q_{\mathrm{id}}(f)=\det R_{\mathrm{id}}(f)=1$.
The following proposition is proved by solving the simultaneous equations derived from the definition of $(k_1,k_2,k_3)$-CGC matrix and $(k_1,k_2,k_3)$-CGC triple:
\begin{proposition}
 The triple $(P_{\mathrm{id}}(f),Q_{\mathrm{id}}(f),R_{\mathrm{id}}(f))$ is a $(k_1,k_2,k_3)$-CGC triple. Conversely, for a $(k_1,k_2,k_3)$-CGC triple $(P,Q,R)$ satisfying $(p_{12},q_{12},r_{12})=(x_1,x_2,x_3)$, there exists $f\in \mathbb Z[x_1^{\pm 1},x_2^{\pm1},x_3^{\pm1}]$ such that $(P,Q,R)=(P_{\mathrm{id}}(f),Q_{\mathrm{id}}(f),R_{\mathrm{id}}(f))$.
\end{proposition}
The other cases satisfying $\{p_{12},q_{12},r_{12}\}=\{x_1,x_2,x_3\}$ are given by a permutation of $(x_1,k_1), (x_2,k_2),(x_3,k_3)$ in $(P_{\mathrm{id}}(f),Q_{\mathrm{id}}(f),R_{\mathrm{id}}(f))$. We set
\begin{align*}
    P_{\sigma}(f)&=\begin{bmatrix}
        f&x_{\sigma(1)}\\\ast&\frac{-fx_{\sigma(2)}x_{\sigma(3)}+k_{\sigma(2)}x_{\sigma(1)}x_{\sigma(3)}+k_{\sigma(3)}x_{\sigma(1)}x_{\sigma(2)}+x_{\sigma(1)}^2+x_{\sigma(2)}^2+x_{\sigma(3)}^2}{x_{\sigma(2)}x_{\sigma(3)}}
    \end{bmatrix},\\
    Q_{\sigma}(f)&=\begin{bmatrix}\frac{fx_{\sigma(2)}+k_{\sigma(1)}x_{\sigma(2)}+x_{\sigma(3)}}{x_{\sigma(1)}}&x_{\sigma(2)}\\ *& \frac{-fx_{\sigma(2)}x_{\sigma(3)}+k_{\sigma(3)}x_{\sigma(1)}x_{\sigma(2)}+x_{\sigma(1)}^2+x_{\sigma(2)}^2}{x_{\sigma(1)}x_{\sigma(3)}}
    \end{bmatrix},\\
     R_{\sigma}(f)&=\begin{bmatrix}
        \frac{fx_{\sigma(2)}x_{\sigma(3)}+k_{\sigma(1)}x_{\sigma(2)}x_{\sigma(3)}+k_{\sigma(2)}x_{\sigma(1)}x_{\sigma(3)}+x_{\sigma(1)}^2+x_{\sigma(3)}^2}{x_{\sigma(1)}x_{\sigma(2)}}&x_{\sigma(3)}\\ *& \frac{-fx_{\sigma(2)}x_{\sigma(3)}+x_{\sigma(2)}^2}{x_{\sigma(1)}x_{\sigma(2)}}
    \end{bmatrix},
\end{align*}
where $\sigma$ is an element of $\mathfrak{S}_3$. Clearly, if $\sigma=\mathrm{id}$, then we have $(P_{\sigma}(f),Q_{\sigma}(f),R_{\sigma}(f))=(P_{\mathrm{id}}(f),Q_{\mathrm{id}}(f),R_{\mathrm{id}}(f))$.

\begin{proposition}\label{cohn-mat-with-111}
The triple $(P_{\sigma}(f),Q_{\sigma}(f),R_{\sigma}(f))$ is a $(k_1,k_2,k_3)$-CGC triple. Conversely, for any $(k_1,k_2,k_3)$-CGC triple $(P,Q,R)$ satisfying $\{p_{12},q_{12},r_{12}\}=\{x_1,x_2,x_3\}$, there exists $f\in \mathbb Z[x_1^{\pm 1},x_2^{\pm1},x_3^{\pm1}]$ and $\sigma\in \mathfrak{S}_3$ such that $(P,Q,R)=(P_\sigma(f),Q_\sigma(f),R_\sigma(f))$.
\end{proposition}

We consider a binary tree, the \emph{$(k_1,k_2,k_3)$-cluster generalized Cohn tree} $\mathrm{CGC}\mathbb T(k_1,k_2,k_3,\sigma,f)$ for $\sigma\in \mathfrak{S}_3$ and $f\in \mathbb Z[x_1^{\pm 1},x_2^{\pm1},x_3^{\pm1}]$:
\begin{itemize}\setlength{\leftskip}{-15pt}
\item [(1)] the root vertex is \[(\widetilde{P_{\sigma}}(f), \widetilde{Q_{\sigma}}(f),\widetilde{R_{\sigma}}(f)):=(P_{\sigma\circ(2\ 3)}(f), P_{\sigma\circ(2\ 3)}(f)Q_{\sigma\circ(2\ 3)}(f)-S_{R_{\sigma\circ(2\ 3)}(f)}, Q_{\sigma\circ(2\ 3)}(f))\]
\item [(2)] each vertex $(P,Q,R)$ has the following two children.
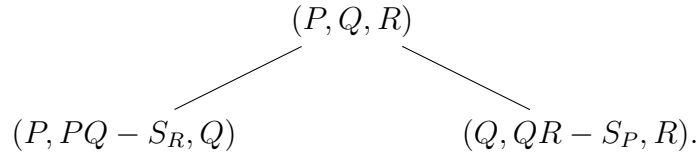
\[\begin{xy}(0,0)*+{(P,Q,R)}="1",(30,-15)*+{(Q,QR-S_P,R).}="2",(-30,-15)*+{(P,PQ-S_R,Q)}="3", \ar@{-}"1";"2"\ar@{-}"1";"3"
\end{xy}\]
\end{itemize}

\begin{remark}\label{rem:right-child-root}
In the definition of $(k_1,k_2,k_3)$-CGC tree, we set the root as $(\widetilde{P_{\sigma}}(f), \widetilde{Q_{\sigma}}(f),\widetilde{R_{\sigma}}(f))$, not $(P_{\sigma}(f), Q_{\sigma}(f), R_\sigma(f))$. This root is associated with $\left( x_1,\frac{x_1^2+k_2x_1x_3+x_3^2}{x_2},x_3\right)$ (if $\sigma=\mathrm{id}$) or a permutations of its indices. The reason for setting the root in this way is to see the duality with the \emph{inverse CGC tree} that we will deal with later. Moreover, we could have instead used $(Q_{\sigma}(f), Q_{\sigma}(f)R_{\sigma}(f)-S_{P_\sigma(f)}, R_{\sigma}(f))$ as the root.  However, we have
\begin{align*}&(Q_{\sigma}(f), Q_{\sigma}(f)R_{\sigma}(f)-S_{P_\sigma(f)}, R_{\sigma}(f))=(\widetilde{P_{\sigma\circ(1\ 2)}}(f'),\widetilde{Q_{\sigma\circ(1\ 2)}}(f'), \widetilde{R_{\sigma\circ(1\ 2)}}(f'),\end{align*}
where $f'=\frac{fx_{\sigma(2)}+k_{\sigma(1)}x_{\sigma(2)}+x_{\sigma(3)}}{x_{\sigma(1)}}$. Therefore, we do not need to consider such trees.
\end{remark}

\begin{theorem}\label{thm:Cohn-well-defined}
Let $(P,Q,R)$ be the $(k_1,k_2,k_3)$-CGC triple associated with $(a,b,c)$.
\begin{itemize}\setlength{\leftskip}{-15pt}
    \item[(1)]
    the $(1,2)$-entry of $PQ-S_R$ (resp. $QR-S_P$) is $\frac{a^2+k_cab+b^2}{c}$ (resp. $\frac{b^2+k_abc+c^2}{a}$),
    \item [(2)] $\mathrm{tr}(PQ-S_R)=\mathscr M m_{12}-k_c$ and $\mathrm{tr}(QR-S_P)=\mathscr M n_{12}-k_a$ hold, where $m_{12}$ (resp. $n_{12}$) is the $(1,2)$-entry of $PQ-S_R$ (resp. $QR-S_P$),
    \item[(3)] $PQ-S_R, QR-S_P\in SL(2,\mathbb Z[x_1^{\pm1},x_2^{\pm1},x_3^{\pm1}])$.
\end{itemize}
In particular, $(P,PQ-S_{R},Q)$ (resp. $(Q,QR-S_{P},R)$) is a $(k_1,k_2,k_3)$-CGC triple.
\end{theorem}

To prove Theorem \ref{thm:Cohn-well-defined}, we introduce four lemmas.

\begin{lemma}[{see \cite[Lemma 4.2]{aigner2013Markov}}]\label{lem:basic-property-trace}
Given $A$ and $B\in SL(2,\mathbb Z[x_1^{\pm1},x_2^{\pm1},x_3^{\pm1}])$, we have
\begin{itemize}\setlength{\leftskip}{-15pt}
    \item [(1)] $\mathrm{tr}(A)=\mathrm{tr}(A^{-1})$,
    \item[(2)] $\mathrm{tr}(AB)= \mathrm{tr}(A)\mathrm{tr}(B)-\mathrm{tr}(AB^{-1})$, and
    \item[(3)] $A^2=\mathrm{tr}(A)A-I$.
\end{itemize}
\end{lemma}
\begin{lemma}\label{lem:MtM}
For a  matrix $P\in SL(2,\mathbb Z[x_1^{\pm1},x_2^{\pm1},x_3^{\pm1}])$ with $\mathrm{tr}(P)=\mathscr M m_{12}-k_{a}$, we have
\begin{align*}
    P\begin{bmatrix}
    0&0\\ \mathscr M & 0
\end{bmatrix}P&=(\mathrm{tr}(P)+k_a)P+\begin{bmatrix}
    0&0\\\mathscr M & 0
\end{bmatrix} \text{ and}\\
 P^{-1}\begin{bmatrix}
    0&0\\ \mathscr M & 0
\end{bmatrix}P^{-1}&=-(\mathrm{tr}(P^{-1})+k_a)P^{-1}+\begin{bmatrix}
    0&0\\\mathscr M & 0
\end{bmatrix}.
\end{align*}
\end{lemma}
\begin{proof}
We only prove the first equality. We set $P=\begin{bmatrix}
    p_{11}&p_{12}\\ p_{21}& p_{22}
\end{bmatrix}$. Then we have
\begin{align*}
    P\begin{bmatrix}
    0&0\\ \mathscr M & 0
\end{bmatrix}P&=P\begin{bmatrix}
    0\\1
\end{bmatrix}
\begin{bmatrix}
    \mathscr M & 0
\end{bmatrix}P=\begin{bmatrix}
    p_{12}\\p_{22}
\end{bmatrix}
\begin{bmatrix}
    \mathscr M p_{11}&\mathscr M p_{12}
\end{bmatrix}\\&=\mathscr M \begin{bmatrix}
    p_{11}p_{12}&p_{12}^2\\p_{11}p_{22}&p_{12}p_{22}
\end{bmatrix}=
\mathscr M p_{12}\begin{bmatrix}
    p_{11}&p_{12}\\p_{21}&p_{22}
\end{bmatrix}+
\begin{bmatrix}
    0&0\\ \mathscr M & 0
\end{bmatrix}\\&=(\mathrm{tr}P+k_a)P+\begin{bmatrix}
    0&0\\\mathscr M & 0\end{bmatrix}.
\end{align*}
Note that $p_{11}p_{22}-p_{21}p_{12}=1$ is used to show that the fourth equal sign is valid.
\end{proof}

\begin{lemma}\label{lem:GSME-transformation}
Given $(P,Q,R)$, $(k_1,k_2,k_3)$-CGC triple associated with $(a,b,c)$, we have
\[\mathrm{tr}(PQ-S_R)=\mathrm{tr}(P)\mathrm{tr}(Q)-\mathrm{tr}(R)-k_ak_b-2k_c\] and \[\mathrm{tr}(QR-S_P)=\mathrm{tr}(Q)\mathrm{tr}(R)-\mathrm{tr}(P)-k_bk_c-2k_a.\]
\end{lemma}

\begin{proof}
We will only prove the first equality.
Since $Q=PR-S_P$, we have
\begin{align*}
    \mathrm{tr}(PQ-S_R)&=\mathrm{tr}(P(PR-S_Q))-2k_c=\mathrm{tr}(P^2R)-\mathrm{tr}(PS_Q)-2k_c\\
    &\overset{\text{Lemma \ref{lem:basic-property-trace} (2)}}{=}\mathrm{tr}(P)\mathrm{tr}(PR)-\mathrm{tr}(PR^{-1}P^{-1})-\mathrm{tr}(PS_Q)-2k_c\\
    &=\mathrm{tr}(P)\mathrm{tr}(PR)-\mathrm{tr}(R)-\mathrm{tr}(PS_Q)-2k_c.
\end{align*}
On the other hand, since
\begin{align*}
\mathrm{tr}(S_QP^{-1})&=\mathrm{tr}\left(\begin{bmatrix}
    p_{22}k_b& \ast\\\ast &p_{11}k_b-p_{12}\mathscr M k_b\end{bmatrix}\right)\\
    &=k_b\cdot\mathrm{tr}(P)-\mathscr M k_bp_{12}=-k_ak_b,
\end{align*}
we have
\begin{align*}
\mathrm{tr}(P)\mathrm{tr}(Q)-\mathrm{tr}(R)-k_ak_b-2k_c&=\mathrm{tr}(P)\mathrm{tr}(PR)-\mathrm{tr}(P)\mathrm{tr}(S_Q)-\mathrm{tr}(R)-k_ak_b-2k_c\\
\overset{\text{Lemma \ref{lem:basic-property-trace} (2)}}{=}&\mathrm{tr}(P)\mathrm{tr}(PR)-\mathrm{tr}(S_QP)-\mathrm{tr}(S_QP^{-1})-\mathrm{tr}(R)-k_ak_b-2k_c\\
&=\mathrm{tr}(P)\mathrm{tr}(PR)-\mathrm{tr}(PS_Q)-\mathrm{tr}(R)-2k_c.
\end{align*}
Therefore, we have the desired equality.
\end{proof}

\begin{lemma}\label{lem:2nd-markov-solution}
If $(a,b,c)$ is a cluster of $\mathcal{A}(k_1,k_2,k_3)$, then $(\mathscr M a-k_1,\mathscr M b-k_2,\mathscr M c-k_3)$ is a solution of the second $(k_1,k_2,k_3)$-GM equation,
\[x^2+y^2+z^2+(2k_1+k_2k_3)x+(2k_2+k_1k_3)y+(2k_3+k_1k_2)z+k_1^2+k_2^2+k_3^2+2k_1k_2k_3=xyz.\]
\end{lemma}

\begin{proof}
We can directly check the statement when $a=x_1,b=x_2,c=x_3$, and then we can show that this quantity is preserved under mutation.
\end{proof}

\begin{proof}[Proof of Theorem \ref{thm:Cohn-well-defined}]
We only prove the statements for $PQ-S_R$. We prove $PQ-S_R$ satisfies (2), i.e., that $\mathrm{tr}(PQ-S_R)=\mathscr M m_{12}-k_c$.
By Lemma \ref{lem:GSME-transformation}, we have
\begin{align*}
\mathrm{tr}(PQ-S_R)&=\mathrm{tr}(P)\mathrm{tr}(Q)-\mathrm{tr}(R)-k_ak_b-2k_c\\
&=\left(\begin{bmatrix}
    \mathscr M & 0
\end{bmatrix}P\begin{bmatrix}
    0 \\ 1
\end{bmatrix}-k_a\right)\left(\begin{bmatrix}
    \mathscr M & 0
\end{bmatrix}Q\begin{bmatrix}
    0 \\ 1
\end{bmatrix}-k_b\right)\\
&-\left(\begin{bmatrix}
    \mathscr M & 0
\end{bmatrix}R\begin{bmatrix}
    0 \\ 1
\end{bmatrix}-k_c\right)-k_ak_b-2k_c\\
&=\begin{bmatrix}
    \mathscr M & 0
\end{bmatrix}P\begin{bmatrix}
    0&0\\\mathscr M & 0
\end{bmatrix}PR\begin{bmatrix}
    0 \\ 1
\end{bmatrix}-k_a\begin{bmatrix}
    \mathscr M & 0
\end{bmatrix}PR\begin{bmatrix}
    0 \\ 1
\end{bmatrix}\\
&-k_b\begin{bmatrix}
    \mathscr M & 0
\end{bmatrix}P\begin{bmatrix}
    0 \\ 1
\end{bmatrix}-\begin{bmatrix}
    \mathscr M & 0
\end{bmatrix}R\begin{bmatrix}
    0 \\ 1
\end{bmatrix}-k_c\\
\overset{\text{Lemma \ref{lem:MtM}}}{=}&\begin{bmatrix}
    \mathscr M & 0
\end{bmatrix}\left((\mathrm{tr}P+k_a)P+\begin{bmatrix}
    0&0\\\mathscr M & 0\end{bmatrix}\right)R\begin{bmatrix}
    0 \\ 1
\end{bmatrix}-k_a\begin{bmatrix}
    \mathscr M & 0
\end{bmatrix}PR\begin{bmatrix}
    0 \\ 1
\end{bmatrix}\\
&-k_b\begin{bmatrix}
    \mathscr M & 0
\end{bmatrix}P\begin{bmatrix}
    0 \\ 1
\end{bmatrix}-\begin{bmatrix}
    \mathscr M & 0
\end{bmatrix}R\begin{bmatrix}
    0 \\ 1
\end{bmatrix}-k_c\\
\overset{\text{Lemma \ref{lem:basic-property-trace} (3)}}{=}&\begin{bmatrix}
    \mathscr M & 0
\end{bmatrix}P^2R\begin{bmatrix}
    0 \\ 1
\end{bmatrix}-k_b\begin{bmatrix}
    \mathscr M & 0
\end{bmatrix}P\begin{bmatrix}
    0 \\ 1
\end{bmatrix}-k_c\\
&=\begin{bmatrix}
    \mathscr M & 0
\end{bmatrix}(P^2R-PS_Q-S_R)\begin{bmatrix}
    0 \\ 1
\end{bmatrix}-k_c\\
&=\begin{bmatrix}
    \mathscr M & 0
\end{bmatrix}(PQ-S_R)\begin{bmatrix}
    0 \\ 1
\end{bmatrix}-k_c=\mathscr M m_{12}-k_c,
\end{align*}
as desired. We remark the third and second-to-last equalities follow from the equality $Q = PR - S_Q$.
Next, we will prove $(P,PQ-S_R,Q)$ satisfies (3). It suffices to show $\det(PQ-S_R)=1$. We set $PQ=\begin{bmatrix}
    x_{11} & x_{12}\\ x_{21} &x_{22}
\end{bmatrix}$.
By the above discussion, we have $\mathrm{tr}(PQ-S_R)=\mathscr M x_{12}-k_c$, and thus $\mathrm{tr}(PQ)=\mathscr M x_{12}+k_c$. Therefore, we have
\begin{align*}
    \det(PQ-S_R)&=(x_{11}-k_c)(x_{22}-k_c)-x_{12}(x_{21}-\mathscr M k_c)\\
    &=\det(PQ)-k_c\cdot\mathrm{tr}(PQ)+k_c^2+\mathscr M k_cx_{12}\\
    &=\det(PQ)=1.
\end{align*}
Finally, we will prove that $(P,PQ-S_R,Q)$ satisfies (1). By assumption and Lemma \ref{lem:2nd-markov-solution}, $(\mathrm{tr}(P),\mathrm{tr}(Q),\mathrm{tr}(R))$ is a solution to the second $(k_1,k_2,k_3)$-GM equation. Substituting $(x,y,z)=(\mathrm{tr}(P),\mathrm{tr}(Q),\mathrm{tr}(R))$ into this equation and transforming it, we obtain
\begin{align*}&\mathrm{tr}(P)\mathrm{tr}(Q)-\mathrm{tr}(R)-k_ak_b-2k_c\\&=\frac{(\mathrm{tr}(P)+k_a)^2+k_c(\mathrm{tr}(P)+k_a)(\mathrm{tr}(Q)+k_b)+(\mathrm{tr}(Q)+k_b)^2}{\mathrm{tr}(R)+k_c}-k_c.\end{align*}
By Lemma \ref{lem:GSME-transformation}, we have the following:
\[\mathrm{tr}(PQ-S_R)=\frac{(\mathrm{tr}(P)+k_a)^2+k_c(\mathrm{tr}(P)+k_a)(\mathrm{tr}(Q)+k_b)+(\mathrm{tr}(Q)+k_b)^2}{\mathrm{tr}(R)+k_c}-k_c.\]
Therefore, by the result (2), we have
\[m_{12}=\frac{1}{\mathscr M }\frac{\mathscr M (a^2+k_cab+b^2)}{c}=\frac{a^2+k_cab+b^2}{c},\]
as desired. This finishes the proof.
\end{proof}

The matrices $PQ-S_R$ and $QR-S_P$ in Theorem \ref{thm:Cohn-well-defined} are new matrices appearing in the two children of $(P,Q,R)$ in $\mathrm{CGC}\mathbb{T}(k_1,k_2,k_3,\sigma,f)$. On the other hand, $\frac{a^2+k_cab+b^2}{c}$ and $\frac{b^2+k_abc+c^2}{a}$ are new cluster variables appearing in the two children of $(a,b,c)$ in $\mathrm{CM}\mathbb T(k_1,k_2,k_3,\sigma)$. Therefore, the following corollary holds.

\begin{corollary}\label{cor:CT-MT}
We fix $f\in \mathbb Z[x_1^{\pm1},x_2^{\pm1},x_{3}^{\pm1}]$ and $\sigma\in \mathfrak S_3$. The correspondence
\[(P,Q,R)\mapsto\left(\left (p_{12},\mathscr M p_{12}-\mathrm{tr}(P)\right),\left(q_{12},\mathscr M q_{12}-\mathrm{tr}(Q)\right),\left(r_{12},\mathscr M r_{12}-\mathrm{tr}(R)\right)\right)\]
induces the canonical graph isomorphism between $\mathrm{CGC}\TT(k_1,k_2,k_3,\sigma, f)$ and $\mathrm{CM}\TT(k_1,k_2,k_3,\sigma)$. In particular, for any (non-labeled) cluster $(a,b,c)$ with $|b|> \max\{|a|,|c|\}$, there is a $(k_1,k_2,k_3)$-CGC triple associated with $(a,b,c)$.
\end{corollary}

The same can be said about the inverse transformation of the transformation treated in Theorem \ref{thm:Cohn-well-defined}.

\begin{theorem}\label{thm:Cohn-well-defined2}
If $(P,Q,R)$ is a $(k_1,k_2,k_3)$-CGC triple associated with $(a,b,c)$, then
\begin{itemize}\setlength{\leftskip}{-15pt}
    \item[(1)]
    the $(1,2)$-entries of $P^{-1}(R+S_R)$ and $(P+S_P)R^{-1}$ are both $\frac{a^2+k_bac+c^2}{b}$;
    \item [(2)] $\mathrm{tr}(P^{-1}(R+S_R))=\mathscr M m_{12}-k_b$ and $\mathrm{tr}((P+S_P)R^{-1})=\mathscr M n_{12}-k_b$ hold, where $m_{12}$ (resp. $n_{12}$) is the $(1,2)$-entry of $P^{-1}(R+S_R)$ (resp. $(P+S_P)R^{-1}$); and
    \item[(3)] $P^{-1}(R+S_R), (P+S_P)R^{-1} \in SL(2,\mathbb Z[x_1^{\pm1},x_2^{\pm1},x_3^{\pm1}])$.
\end{itemize}
In particular, $(P,R,P^{-1}(R+S_R))$ (resp. $((P+S_P)R^{-1},P,R)$) is a $(k_1,k_2,k_3)$-CGC triple.
\end{theorem}

The proof method is the same as Theorem \ref{thm:Cohn-well-defined}, so we will omit it. We note that the correspondence $(P,Q,R)\mapsto(P,R,P^{-1}(R+S_R))$ is the inverse correspondence of $(P,Q,R)\mapsto(P,PQ-S_R,Q)$, and $(P,Q,R)\mapsto((P+S_P)R^{-1},P,R)$ is the inverse of $(P,Q,R)\mapsto(Q,QR-S_P,R)$.

In Corollary \ref{cor:CT-MT}, we found that every non-labeled cluster in $\Acal(k_1,k_2,k_3)$ satisfying $|b|>\max\{|a|,|c|\}$ has a corresponding $(k_1,k_2,k_3)$-CGC triple.
In fact, we can also find that the $(k_1,k_2,k_3)$-CGC triples which correspond to the non-labeled clusters in $\Acal(k_1,k_2,k_3)$ satisfying $|b|>\max\{|a|,|c|\}$ are all those in $\mathrm{CGC}\mathbb T(k_1,k_2,k_3,\sigma,f)$.

\begin{theorem}\label{thm:all-cohn-triple}
Let $(P,Q,R)$ be a $(k_1,k_2,k_3)$-CGC triple associated with $(a,b,c)$. We assume that $|b|>\max\{|a|,|c|\}.$ Then, there exist a unique $f\in \mathbb Z[x_1^{\pm1},x_2^{\pm1},x_3^{\pm1}]$, a unique $\sigma\in \mathfrak{S}_3$ and a unique vertex $v$ of $\mathrm{CGC}\mathbb T(k_1,k_2,k_3,\sigma,f)$ such that $v=(P,Q,R)$.
\end{theorem}

\begin{proof}
By Corollary \ref{cor:unique-non-labeled-cluster}, the non-labeled cluster $(a,b,c)$ together with the parities of its components determines a unique permutation $\sigma$ and a unique vertex $w$ of $\mathrm{CM}\mathbb T(k_1,k_2,k_3,\sigma)$. Let
\[
w_0,w_1,\ldots,w_m=w
\]
be the unique path from the root $w_0$ to $w$. Starting with $(P_m,Q_m,R_m):=(P,Q,R)$, we move backwards along this path. If $w_i$ is the left child of $w_{i-1}$, we apply the first inverse transformation in Theorem \ref{thm:Cohn-well-defined2}; if it is the right child, we apply the second inverse transformation. In this way we obtain a uniquely determined sequence of CGC triples $(P_i,Q_i,R_i)$ associated with $w_i$. Thus the sequence of inverse transformations is uniquely determined by the path $w_0,\ldots,w_m$.

The triple $(P_0,Q_0,R_0)$ is associated with the root of $\mathrm{CM}\mathbb T(k_1,k_2,k_3,\sigma)$. Applying the first inverse transformation of Theorem \ref{thm:Cohn-well-defined2} once more gives a CGC triple $(A,B,C)$ associated with the ordered initial cluster
\[
(x_{\tau(1)},x_{\tau(2)},x_{\tau(3)}),\qquad \tau=\sigma\circ(2\ 3).
\]
By Proposition \ref{cohn-mat-with-111}, $(A,B,C)=(P_\tau(f),Q_\tau(f),R_\tau(f))$ for some $f$. The ordered $(1,2)$-entries force the permutation to be $\tau$, and the explicit formula for $P_\tau(f)$ gives
\[
f=(P_\tau(f))_{11}=A_{11}.
\]
Therefore $f$ is uniquely determined. Reversing the uniquely determined sequence of transformations shows that $(P,Q,R)$ occurs in $\mathrm{CGC}\mathbb T(k_1,k_2,k_3,\sigma,f)$.

Finally, suppose the same matrix triple occurred in $\mathrm{CGC}\mathbb T(k_1,k_2,k_3,\sigma',f')$. Corollary \ref{cor:CT-MT} and Corollary \ref{cor:unique-non-labeled-cluster} imply $\sigma'=\sigma$ and that the two occurrences correspond to the same vertex $w$. Following the unique path back to the root gives the same terminal triple $(A,B,C)$, and hence $f'=A_{11}=f$. Since Corollary \ref{cor:CT-MT} is a graph isomorphism, the vertex itself is also unique.
\end{proof}

\begin{example}
We will consider the case where $k_1=k_2=k_3=0$, $\sigma=\textrm{id}$ and $f=0$.
We have
\begin{align*}
    P_{\mathrm{id}}(0)=\begin{bmatrix}
        0&x_1\\-\frac{1}{x_1}&\frac{x_1^2+x_2^2+x_3^2}{x_2x_3}
    \end{bmatrix}, \quad
    Q_{\mathrm{id}}(0)=\begin{bmatrix}
        \frac{x_3}{x_1}&x_2\\ \frac{x_2}{x_1^2}& \frac{x_1^2+x_2^2}{x_1x_3}
    \end{bmatrix},\quad
     R_{\mathrm{id}}(0)=\begin{bmatrix}
        \frac{x_1^2+x_3^2}{x_1x_2}&x_3\\ \frac{x_3}{x_1^2}& \frac{x_2^2}{x_1x_2}
    \end{bmatrix}.
\end{align*}
Therefore, we have
\begin{align*}
    \widetilde{P_{\mathrm{id}}}(0)&=P_{(2\ 3)}(0)=P_{\mathrm{id}}(0), \\
   \widetilde{Q_{\mathrm{id}}}(0)&= P_{(2\ 3)}(0)Q_{(2\ 3)}(0)=\begin{bmatrix}
        \frac{x_3}{x_1}&\frac{x_1^2+x_3^2}{x_2}\vspace{1mm} \\\frac{x_1^2+x_3^2}{x_1^2x_2}&\frac{(x_1^2+x_3^2)^2+x_1^2x_2^2}{x_1x_2^2x_3}
    \end{bmatrix},\\
     \widetilde{R_{\mathrm{id}}}(0)&=Q_{(2\ 3)}(0)=\begin{bmatrix}
        \frac{x_2}{x_1}&x_3\\ \frac{x_3}{x_1^2}& \frac{x_1^2+x_3^2}{x_1x_2}.
    \end{bmatrix}
\end{align*}
We can check that the triple of $(1,2)$-entries in $(\widetilde{P_{\mathrm{id}}}(0),\widetilde{Q_{\mathrm{id}}}(0),\widetilde{R_{\mathrm{id}}}(0))$ coincides with the root of $\mathrm{CM\mathbb T(0,0,0,id)}$.
The middle component of the left child of the root is
\[\widetilde{P_{\mathrm{id}}}(0)\widetilde{Q_{\mathrm{id}}}(0)=\begin{bmatrix}
 \frac{x_1^2+x_3^2}{x_1x_2}&\frac{(x_1^2+x_3^2)^2+x_1^2x_2^2}{x_2^2x_3}\vspace{1mm}\\\frac{(x_1^2+x_3^2)^2+x_1^2x_2^2}{x_1^2x_2^2x_3} & \frac{(x_1^2+x_3^2)\left((x_1^2+x_3^2)^2+x_1^2x_2^2\right)+x_2^2(x_1^4+x_1^2x_2^2+x_1^2x_3^2)}{x_1x_2^3x_3^2}
\end{bmatrix}.\]
We can check that the triple of $(1,2)$-entries \[\left(x_1,\frac{(x_1^2+x_3^2)^2+x_1^2x_2^2}{x_2^2x_3},\frac{x_1^2+x_3^2}{x_2}\right)\] of $(\widetilde{P_{\mathrm{id}}}(0),\widetilde{P_{\mathrm{id}}}(0)\widetilde{Q_{\mathrm{id}}}(0),\widetilde{Q_{\mathrm{id}}}(0))$ is the left child of the root in $\mathrm{CM\mathbb T(0,0,0,id)}$.
\end{example}

So far, we saw the properties of CGC triple that satisfy $|b|>\max\{|a|,|c|\}$. Next, we will see cases where $|b|\leq\max\{|a|,|c|\}$. We consider another binary tree, the \emph{inverse $(k_1,k_2,k_3)$-cluster generalized Cohn tree} $\mathrm{CGC}\mathbb T^\dag(k_1,k_2,k_3,\sigma,f)$ for $\sigma\in \mathfrak{S}_3$ and $f\in \mathbb Z[x_1^{\pm 1},x_2^{\pm1},x_3^{\pm1}]$.
\begin{itemize}\setlength{\leftskip}{-15pt}
\item [(1)] Set the root vertex to be  $(P_{\sigma}(f), Q_{\sigma}(f),R_{\sigma}(f))$.
\item [(2)] The two children of  a vertex $(P,Q,R)$ are as follows.
\[\begin{xy}(0,0)*+{(P,Q,R)}="1",(-30,-15)*+{(P,R,P^{-1}(R+S_R))}="2",(30,-15)*+{((P+S_P)R^{-1},P,R).}="3", \ar@{-}"1";"2"\ar@{-}"1";"3"
\end{xy}\]
\end{itemize}
Using Theorem \ref{thm:Cohn-well-defined2} instead of Theorem \ref{thm:Cohn-well-defined}, we have the following corollary in parallel with Corollary \ref{cor:CT-MT}.
\begin{corollary}\label{cor:CT-MT2}
We fix $f\in \mathbb Z[x_1^{\pm1},x_2^{\pm1},x_{3}^{\pm1}]$ and $\sigma\in \mathfrak S_3$. The correspondence
\[(P,Q,R)\mapsto\left(\left (p_{12},\mathscr M p_{12}-\mathrm{tr}(P)\right),\left(q_{12},\mathscr M q_{12}-\mathrm{tr}(Q)\right),\left(r_{12},\mathscr M r_{12}-\mathrm{tr}(R)\right)\right)\]
induces the canonical graph isomorphism between $\mathrm{CGC}\TT^\dag(k_1,k_2,k_3,\sigma, f)$ and $\mathrm{CM}\TT^\dag(k_1,k_2,k_3,\sigma)$. In particular, for any (non-labeled) cluster $(a,b,c)$ with $|b|\leq \max\{|a|,|c|\}$, there is a $(k_1,k_2,k_3)$-CGC triple associated with $(a,b,c)$.
\end{corollary}

Exchanging the role of Theorem \ref{thm:Cohn-well-defined} and that of Theorem \ref{thm:Cohn-well-defined2} in Theorem \ref{thm:all-cohn-triple}, we have the following theorem.

\begin{theorem}\label{thm:all-cohn-triple2}
Let $(P,Q,R)$ be a $(k_1,k_2,k_3)$-CGC triple associated with $(a,b,c)$. We assume that $|b|\leq\max\{|a|,|c|\}.$ Then, there exist a unique $f\in \mathbb Z[x_1^{\pm1},x_2^{\pm1},x_3^{\pm1}]$, a unique $\sigma\in \mathfrak{S}_3$ and a unique vertex $v$ in $\mathrm{CGC}\mathbb T^\dag(k_1,k_2,k_3,\sigma,f)$ such that $v=(P,Q,R)$.
\end{theorem}

\begin{proof}
The proof is parallel to that of Theorem \ref{thm:all-cohn-triple}, with the ordinary trees replaced by their inverse versions, Corollaries \ref{cor:unique-non-labeled-cluster} and \ref{cor:CT-MT} replaced by Corollaries \ref{cor:unique-non-labeled-cluster-2} and \ref{cor:CT-MT2}, respectively, and the roles of Theorems \ref{thm:Cohn-well-defined} and \ref{thm:Cohn-well-defined2} interchanged. In this case, the unique path back to the root yields directly
\[
(A,B,C)=(P_\sigma(f),Q_\sigma(f),R_\sigma(f)),
\]
and the explicit formula gives $f=A_{11}$. Thus $\sigma$, $f$, and the vertex are unique, while reversing the path proves existence.
\end{proof}

\section{Cluster Markov-Monodoromy Matrices}\label{sec:ClusterMM}
In this section, we will introduce another matrixization of cluster variables in $\Acal(k_1,k_2,k_3)$ called the cluster Markov-monodromy matrices. These will have properties that are analogue to CGC matrices. Almost all theorems in this section are proved in the next section, where explicit connections between the two matrix families are given.

\begin{definition}
For $k_1,k_2,k_3\in \mathbb {Z}_{\geq 0}$, we define a \emph{$(k_1,k_2,k_3)$-cluster Markov-monodromy matrix} (abbreviated as \emph{$(k_1,k_2,k_3)$-CMM matrix}) $X=\begin{bmatrix}x_{11}&x_{12}\\x_{21}&x_{22}\end{bmatrix}$  as a matrix such that
\begin{itemize}\setlength{\leftskip}{-15pt}
    \item [(1)] $X\in SL(2,\mathbb{Z}[x_1^{\pm1},x_2^{\pm1},x_3^{\pm1}])$,
    \item [(2)] $x_{12}$ is a cluster variable in the $(k_1,k_2,k_3)$-generalized Markov cluster algebra, and
    \item [(3)] $\mathrm{tr}(X)=-k_{x_{12}}$.
\end{itemize}
\end{definition}
\begin{definition}\label{def:markov-monodoromy-triple}
For $(k_1,k_2,k_3)\in \mathbb {Z}_{\geq 0}$, we define a \emph{$(k_1,k_2,k_3)$-cluster Markov-monodromy triple} (abbreviated as \emph{$(k_1,k_2,k_3)$-CMM triple}) $(X,Y,Z)$ as a triple such that
\begin{itemize}\setlength{\leftskip}{-15pt}
    \item [(1)] $X,Y,Z$ are $(k_1,k_2,k_3)$-CMM matrices,
    \item[(2)] $XYZ=T$, where $T=\begin{bmatrix}
        -1&0\\\mathscr M  &-1
    \end{bmatrix}$, and
    \item[(3)] $(x_{12},y_{12},z_{12})$ is a (non-labeled) cluster in $(k_1,k_2,k_3)$-generalized Markov cluster algebra, where $x_{12},y_{12},z_{12}$ are the $(1,2)$-entries of $X,Y,Z$, respectively.
\end{itemize}
\end{definition}

For a non-labeled cluster $(a,b,c)$, CMM-triples satisfying $(x_{12},y_{12},z_{12})=(a,b,c)$ are called those \emph{associated with $(a,b,c)$}.

As for the CGC matrices, the definition of $(k_1,k_2,k_3)$-CMM triples does not refer to the their existence. We begin our discussion by demonstrating the existence of such triples where $(x_{12},y_{12},z_{12})=(x_1,x_2,x_3)$ or one of its permutations. We remark that even if $(P,Q,R)$ is a $(k_1,k_2,k_3)$-CMM triple, a permutation of these matrices is not necessarily a $(k_1,k_2,k_3)$-CMM triples. Therefore, we again need to consider all possible permutations of $(x_1,x_2,x_3)$.
First, we will consider the case $(x_{12},y_{12},z_{12})=(x_1,x_2,x_3)$. We set
\begin{align*}
    X_{\mathrm{id}}(f)&=\begin{bmatrix}
        f&x_1\\\ast&-f-k_1
    \end{bmatrix},\\
    Y_{\mathrm{id}}(f)&=\begin{bmatrix}
        \frac{fx_2-k_2x_1-x_3}{x_1}&x_2\\ *& \frac{-fx_2+x_3}{x_1}
    \end{bmatrix},\\
    Z_{\mathrm{id}}(f)&=\begin{bmatrix}
        \frac{fx_2x_3-k_3x_1x_2-k_2x_1x_3-x_1^2-x_3^2}{x_1x_2}&x_3\\ *& \frac{-fx_2x_3+k_2x_1x_3+x_1^2+x_3^2}{x_1x_2}
    \end{bmatrix},
\end{align*}
where $f\in \mathbb Z[x_1^{\pm1},x_2^{\pm1},x_3^{\pm1}]$, and we omit $(2,1)$-entries of $X_{\mathrm{id}}(f), Y_{\mathrm{id}}(f),Z_{\mathrm{id}}(f)$, which are determined by $\det X_{\mathrm{id}}(f)=\det Y_{\mathrm{id}}(f)=\det Z_{\mathrm{id}}(f)=1$.
The following proposition is proved by solving the simultaneous sytstem of equations derived from the definition of $(k_1,k_2,k_3)$-CMM matrices and $(k_1,k_2,k_3)$-CMM triples.
\begin{proposition}
 The triple $(X_{\mathrm{id}}(f),Y_{\mathrm{id}}(f),Z_{\mathrm{id}}(f))$ is a $(k_1,k_2,k_3)$-CMM triple. Conversely, for a $(k_1,k_2,k_3)$-CMM triple $(X,Y,Z)$ satisfying $(x_{12},y_{12},z_{12})=(x_1,x_2,x_3)$, there exists $f\in \mathbb Z[x_1^{\pm 1},x_2^{\pm1},x_3^{\pm1}]$ such that $(X,Y,Z)=(X_{\mathrm{id}}(f),Y_{\mathrm{id}}(f),Z_{\mathrm{id}}(f))$.
\end{proposition}
The other cases satisfying $\{x_{12},y_{12},z_{12}\}=\{x_1,x_2,x_3\}$ are given by a permutation of $(x_1,k_1), (x_2,k_2),(x_3,k_3)$ in $(X_{\mathrm{id}}(f),Y_{\mathrm{id}}(f),Z_{\mathrm{id}}(f))$. We set
\begin{align*}
    X_{\sigma}(f)&=\begin{bmatrix}
        f&x_{\sigma(1)}\\\ast&-f-k_{\sigma(1)}
    \end{bmatrix},\\
    Y_{\sigma}(f)&=\begin{bmatrix}
        \frac{fx_{\sigma(2)}-k_{\sigma(2)}x_{\sigma(1)}-x_{\sigma(3)}}{x_{\sigma(1)}}&x_{\sigma(2)}\\ *& \frac{-fx_{\sigma(2)}+x_{\sigma(3)}}{x_{\sigma(1)}}
    \end{bmatrix},\\
    Z_{\sigma}(f)&=\begin{bmatrix}
        \frac{fx_{\sigma(2)}x_{\sigma(3)}-k_{\sigma(3)}x_{\sigma(1)}x_{\sigma(2)}-k_{\sigma(2)}x_{\sigma(1)}x_{\sigma(3)}-x_{\sigma(1)}^2-x_{\sigma(3)}^2}{x_{\sigma(1)}x_{\sigma(2)}}&x_{\sigma(3)}\\ *& \frac{-fx_{\sigma(2)}x_{\sigma(3)}+k_{\sigma(2)}x_{\sigma(1)}x_{\sigma(3)}+x_{\sigma(1)}^2+x_{\sigma(3)}^2}{x_{\sigma(1)}x_{\sigma(2)}}
    \end{bmatrix},
\end{align*}
where $\sigma$ is an element of $\mathfrak{S}_3$. Clearly, if $\sigma=\mathrm{id}$, then we have $(X_{\sigma}(f),Y_{\sigma}(f),Z_{\sigma}(f))=(X_{\mathrm{id}}(f),Y_{\mathrm{id}}(f),Z_{\mathrm{id}}(f))$.

\begin{proposition}\label{prop:CMM-initial-triples}
The triple $(X_{\sigma}(f),Y_{\sigma}(f),Z_{\sigma}(f))$ is a $(k_1,k_2,k_3)$-CMM triple. Conversely, for a $(k_1,k_2,k_3)$-CMM triple $(X,Y,Z)$ satisfying $\{x_{12},y_{12},z_{12}\}=\{x_1,x_2,x_3\}$, there exist $f\in \mathbb Z[x_1^{\pm 1},x_2^{\pm1},x_3^{\pm1}]$ and $\sigma\in \mathfrak{S}_3$ such that $(X,Y,Z)=(X_\sigma(f),Y_\sigma(f),Z_\sigma(f))$.
\end{proposition}

We define the \emph{$(k_1,k_2,k_3)$-cluster generalized Cohn tree} $\mathrm{CMM}\mathbb T(k_1,k_2,k_3,\sigma,f)$ for $\sigma\in \mathfrak{S}_3$ and $f\in \mathbb Z[x_1^{\pm 1},x_2^{\pm1},x_3^{\pm1}]$ to be the binary tree as follows.
\begin{itemize}\setlength{\leftskip}{-15pt}
\item [(1)] The root vertex is
\[(\widetilde{X_{\sigma}}(f), \widetilde{Y_{\sigma}}(f),\widetilde{Z_{\sigma}}(f)):=(X_{\sigma\circ(2\ 3)}(f),Y_{\sigma\circ(2\ 3)}(f)Z_{\sigma\circ(2\ 3)}(f)Y_{\sigma\circ(2\ 3)}(f)^{-1}, Y_{\sigma\circ(2\ 3)}(f)).\]
\item [(2)] A vertex $(X,Y,Z)$ has the following two children.
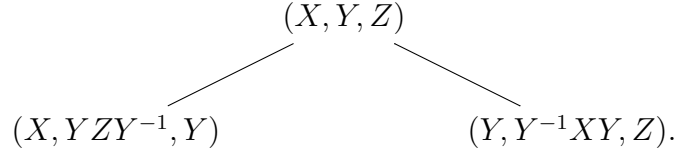
\[\begin{xy}(0,0)*+{(X,Y,Z)}="1",(-30,-15)*+{(X,YZY^{-1},Y)}="2",(30,-15)*+{(Y,Y^{-1}XY,Z).}="3", \ar@{-}"1";"2"\ar@{-}"1";"3"
\end{xy}\]
\end{itemize}

\begin{remark}
As in the case of the CGC-tree, in the definition of the $(k_1,k_2,k_3)$-CMM trees, we set the root as $(\widetilde{X_{\sigma}}(f), \widetilde{Y_{\sigma}}(f),\widetilde{Z_{\sigma}}(f))$, not $(X_{\sigma}(f), Y_{\sigma}(f), Z_\sigma(f))$. Moreover, we could instead of have used $(Y_{\sigma}(f), Y_{\sigma}(f)^{-1}X_{\sigma}(f)Y_{\sigma}(f), Z_{\sigma}(f))$ as the root, but we have
\begin{align*}&(Y_{\sigma}(f), Y_{\sigma}(f)^{-1}X_{\sigma}(f)Y_{\sigma}(f), Z_{\sigma}(f))=(\widetilde{X_{\sigma\circ(1\ 2)}}(f'),\widetilde{Y_{\sigma\circ(1\ 2)}}(f'), \widetilde{Z_{\sigma\circ(1\ 2)}}(f'),\end{align*}
where $f'=\frac{fx_{\sigma(2)}-k_{\sigma(2)}x_{\sigma(1)}-x_{\sigma(3)}}{x_{\sigma(1)}}$. Therefore, we do not need to consider these cases.
\end{remark}

Our next result shows that all triples in each tree $\mathrm{CMM}\mathbb T(k_1,k_2,k_3,\sigma,f)$ are CMM triples.

\begin{theorem}\label{thm:Markov-monodoromy-well-defined}
If $(X,Y,Z)$ is a $(k_1,k_2,k_3)$-CMM triple associated with $(a,b,c)$, then
\begin{itemize}\setlength{\leftskip}{-15pt}
    \item[(1)]
    the $(1,2)$-entry of $YZY^{-1}$ (resp. $Y^{-1}XY$) is $\frac{a^2+k_cab+b^2}{c}$ (resp. $\frac{b^2+k_abc+c^2}{a}$),
    \item [(2)] the trace of $YZY^{-1}$ is $k_c$ and the trace of $Y^{-1}XY$ is $-k_a$, and
    \item[(3)] $YZY^{-1}, Y^{-1}XY\in SL(2,\mathbb Z[x_1^{\pm1},x_2^{\pm1},x_3^{\pm1}])$.
\end{itemize}
In particular, $(X,YZY^{-1},Y)$ and $(Y,Y^{-1}XY,Z)$ are  $(k_1,k_2,k_3)$-CMM triples.
\end{theorem}

It is easy to see that (2) holds by using the fact $\mathrm{tr}(AB)=\mathrm{tr}(BA)$ for general matrices $A$ and $B$, and it is clear that (3) holds. We will prove item (1) in the next section.

The correspondence $(X,Y,Z)\mapsto(X,Z,Z^{-1}YZ)$ is the inverse correspondence of $(X,Y,Z)\mapsto(X,YZY^{-1},Y)$, and $(X,Y,Z)\mapsto(XYX^{-1},X,Z)$ is the inverse of $(X,Y,Z)\mapsto(Y,Y^{-1}XY,Z)$. These inverse transformations also preserve CMM triples.

\begin{theorem}\label{thm:Markov-monodoromy-well-defined2}
If $(X,Y,Z)$ be the $(k_1,k_2,k_3)$-CMM triple associated with $(a,b,c)$, then
\begin{itemize}\setlength{\leftskip}{-15pt}
    \item[(1)]
    the $(1,2)$-entry of $Z^{-1}YZ$ and $XYX^{-1}$ are both $\frac{a^2+k_bac+c^2}{b}$;
    \item [(2)] $\mathrm{tr}(Z^{-1}YZ)=\mathrm{tr}(XYX^{-1})=-k_b$ holds; and
    \item[(3)] $Z^{-1}YZ,XYX^{-1}\in SL(2,\mathbb Z[x_1^{\pm1},x_2^{\pm1},x_3^{\pm1}])$.
\end{itemize}
In particular, $(X,Z,Z^{-1}XZ)$ (resp. $(XYX^{-1},X,Z)$) is a $(k_1,k_2,k_3)$-CMM triple.
\end{theorem}

The following two results, Corollary \ref{cor:MMT-MT} and Theorem \ref{thm:all-Markov-triple}, are in parallel with Corollary \ref{cor:CT-MT} and Theorem \ref{thm:all-cohn-triple}. We prove these in the next section.

\begin{corollary}\label{cor:MMT-MT}
We fix $f\in \mathbb Z[x_1^{\pm1},x_2^{\pm1},x_{3}^{\pm1}]$ and $\sigma\in \mathfrak S_3$. The correspondence
\[(X,Y,Z) \mapsto ((x_{12},-\mathrm{tr}(X)),(y_{12},-\mathrm{tr}(Y)),(z_{12},-\mathrm{tr}(Z)))\]
induces the canonical graph isomorphism between $\mathrm{CMM}\TT(k_1,k_2,k_3,\sigma, f)$ and $\mathrm{CM}\TT(k_1,k_2,k_3,\sigma)$. In particular, for any (non-labeled) cluster $(a,b,c)$ with $b> \max\{|a|,|c|\}$, there is a $(k_1,k_2,k_3)$-CMM matrix associated with $(a,b,c)$.
\end{corollary}

\begin{theorem}\label{thm:all-Markov-triple}
Let $(X,Y,Z)$ be a $(k_1,k_2,k_3)$-CMM triple associated with $(a,b,c)$. We assume that $|b|>\max\{|a|,|c|\}.$ Then, there exist a unique $f\in \mathbb Z[x_1^{\pm1},x_2^{\pm1},x_3^{\pm1}]$, a unique $\sigma\in \mathfrak{S}_3$ and a unique vertex $v$ in $\mathrm{CMM}\mathbb T(k_1,k_2,k_3,\sigma,f)$ such that $v=(X,Y,Z)$.
\end{theorem}

\begin{example}
We will consider the case where $k_1=k_2=k_3=0$, $\sigma=\textrm{id}$ and $f=0$.
We have
\begin{align*}
    X_{\mathrm{id}}(0)=\begin{bmatrix}
        0&x_1\\-\frac{1}{x_1}&0
    \end{bmatrix}, \quad
    Y_{\mathrm{id}}(0)=\begin{bmatrix}
        -\frac{x_3}{x_1}&x_2\\ -\frac{x_1^2+x_3^2}{x_1^2x_2}& \frac{x_3}{x_1}
    \end{bmatrix},\quad
     Z_{\mathrm{id}}(0)=\begin{bmatrix}
        -\frac{x_1^2+x_3^2}{x_1x_2}&x_3\\ -\frac{(x_1^2+x_3^2)^2+x_1^2x_2^2}{x_1^2x_2^2x_3}& \frac{x_1^2+x_3^2}{x_1x_2}
    \end{bmatrix}.
\end{align*}
Therefore, we have
\begin{align*}
    \widetilde{X_{\mathrm{id}}}(0)&=X_{(2\ 3)}(0)=X_{\mathrm{id}}(0), \\
   \widetilde{Y_{\mathrm{id}}}(0)&= Y_{(2\ 3)}(0)Z_{(2\ 3)}(0)Y_{(2\ 3)}(0)^{-1}=\begin{bmatrix}
        -\frac{x_3}{x_1}&\frac{x_1^2+x_3^2}{x_2}\vspace{1mm} \\-\frac{x_2}{x_1^2}&\frac{x_3}{x_1}
    \end{bmatrix},\\
     \widetilde{Z_{\mathrm{id}}}(0)&=Y_{(2\ 3)}(0)=\begin{bmatrix}
        -\frac{x_2}{x_1}&x_3\\ -\frac{x_1^2+x_2^2}{x_1^2x_3}& \frac{x_2}{x_1}.
    \end{bmatrix}
\end{align*}
We can check that the triple of $(1,2)$-entries in $(\widetilde{X_{\mathrm{id}}}(0),\widetilde{Y_{\mathrm{id}}}(0),\widetilde{Z_{\mathrm{id}}}(0))$ coincides with the root of $\mathrm{CM\mathbb T(0,0,0,id)}$.
The middle component of the left child of the root is
\[\widetilde{Y_{\mathrm{id}}}(0)\widetilde{Z_{\mathrm{id}}}(0)\widetilde{Y_{\mathrm{id}}}(0)^{-1}=\begin{bmatrix}
 -\frac{x_1^2+x_3^2}{x_1x_2}&\frac{(x_1^2+x_3^2)^2+x_1^2x_2^2}{x_2^2x_3}\vspace{1mm}\\\frac{x_3}{x_1^2} & \frac{x_1^2+x_3^2}{x_1x_2}
\end{bmatrix}.\]
We can check that the triple of $(1,2)$-entries \[\left(x_1,\frac{(x_1^2+x_3^2)^2+x_1^2x_2^2}{x_2^2x_3},\frac{x_1^2+x_3^2}{x_2}\right)\] of $(\widetilde{X_{\mathrm{id}}}(0),\widetilde{Y_{\mathrm{id}}}(0)\widetilde{Z_{\mathrm{id}}}(0)\widetilde{Y_{\mathrm{id}}}(0)^{-1},\widetilde{Y_{\mathrm{id}}}(0))$ is the left child of the root in $\mathrm{CM\mathbb T(0,0,0,id)}$.
\end{example}
Finally, in order to discuss CMM triples associated to clusters $(a,b,c)$ with $|b| \leq \max\{|a|,|c|\}$, we consider the \emph{inverse $(k_1,k_2,k_3)$-cluster-Markov monodromy tree}, which is denoted by $\mathrm{CMM}\mathbb T^\dag(k_1,k_2,k_3,\sigma,f)$, for $\sigma\in \mathfrak{S}_3$ and $f\in \mathbb Z[x_1^{\pm 1},x_2^{\pm1},x_3^{\pm1}]$. This is the binary tree defined as follows.
\begin{itemize}\setlength{\leftskip}{-15pt}
\item [(1)] The root vertex is $(X_{\sigma}(f), Y_{\sigma}(f),Z_{\sigma}(f))$.
\item [(2)] Each vertex $(X,Y,Z)$ has the following two children of it.
\[\begin{xy}(0,0)*+{(X,Y,Z)}="1",(-30,-15)*+{(X,Z,Z^{-1}YZ)}="2",(30,-15)*+{(XYX^{-1},X,Z).}="3", \ar@{-}"1";"2"\ar@{-}"1";"3"
\end{xy}\]
\end{itemize}

The following corollary and theorem are the inverse-tree counterparts of Corollary \ref{cor:MMT-MT} and Theorem \ref{thm:all-Markov-triple}. They follow by the same argument, using Corollary \ref{cor:CT-MT2} and Theorem \ref{thm:all-cohn-triple2} in place of Corollary \ref{cor:CT-MT} and Theorem \ref{thm:all-cohn-triple}. More precisely, Theorem \ref{thm:BT-CT2-inverse} gives the canonical graph isomorphism
\[
\Psi_0\colon \mathrm{CMM}\mathbb T^\dag(k_1,k_2,k_3,\sigma,f)
\longrightarrow
\mathrm{CGC}\mathbb T^\dag(k_1,k_2,k_3,\sigma,h),
\qquad h=-f-k_{\sigma(1)}.
\]
Combining this with Corollary \ref{cor:CT-MT2} gives the following corollary. Moreover, the uniqueness of $\sigma$, $h$, and the vertex in Theorem \ref{thm:all-cohn-triple2} yields the corresponding uniqueness below, since $f=-h-k_{\sigma(1)}$. These results complete our classification.

\begin{corollary}\label{cor:MMT-MT2}
We fix $f\in \mathbb Z[x_1^{\pm1},x_2^{\pm1},x_{3}^{\pm1}]$ and $\sigma\in \mathfrak S_3$. The correspondence
\[(X,Y,Z) \mapsto ((x_{12},-\mathrm{tr}(X)),(y_{12},-\mathrm{tr}(Y)),(z_{12},-\mathrm{tr}(Z)))\]
induces the canonical graph isomorphism between $\mathrm{CMM}\TT^\dag(k_1,k_2,k_3,\sigma, f)$ and $\mathrm{CM}\TT^\dag(k_1,k_2,k_3,\sigma)$. In particular, for any (non-labeled) cluster $(a,b,c)$ with $|b|\leq \max\{|a|,|c|\}$, there is a $(k_1,k_2,k_3)$-CMM matrix associated with $(a,b,c)$.
\end{corollary}

\begin{theorem}\label{thm:all-markov-triple2}
Let $(X,Y,Z)$ be a $(k_1,k_2,k_3)$-CMM triple associated with $(a,b,c)$. We assume that $|b|\leq\max\{|a|,|c|\}.$ Then, there exist a unique $f\in \mathbb Z[x_1^{\pm1},x_2^{\pm1},x_3^{\pm1}]$, a unique $\sigma\in \mathfrak{S}_3$ and a unique vertex $v$ in $\mathrm{CMM}\mathbb T^\dag(k_1,k_2,k_3,\sigma,f)$ such that $v=(X,Y,Z)$.
\end{theorem}

\section{Relations between CGC Triples and CMM Triples}\label{sec:Relations}
In this section, we introduce two correspondences between CGC triples and CMM triples. We moreover use these correspondences to  prove some of the properties of CMM-triples introduced in the previous section.
\subsection{Isomorphism Maps between CGC Trees and CMM Trees}
For $g\in \mathbb Z[x_1^{\pm1},x_2^{\pm1},x_3^{\pm1}]$, we define the map \[\psi_{g}\colon M(2,\mathbb Z[x_1^{\pm1},x_2^{\pm1},x_3^{\pm1}]) \to M(2,\mathbb Z[x_1^{\pm1},x_2^{\pm1},x_3^{\pm1}])\]
by
\[\begin{bmatrix}
    m_{11}&m_{12}\\m_{21}&m_{22}
\end{bmatrix}\mapsto \begin{bmatrix}1&0\\\mathscr M-g&1\end{bmatrix}\begin{bmatrix}
    m_{22}&m_{12}\\m_{21}& m_{11}
\end{bmatrix}\begin{bmatrix}
    1&0\\g&1
\end{bmatrix}.\]
This map is a bijection. Indeed, we can construct the inverse map $\psi_{{g}}^{-1}$ as
\[\begin{bmatrix}
    m_{11}&m_{12}\\m_{21}&m_{22}
\end{bmatrix}\mapsto \begin{bmatrix}1&0\\-g&1\end{bmatrix}\begin{bmatrix}
    m_{22}&m_{12}\\m_{21}& m_{11}
\end{bmatrix}\begin{bmatrix}
    1&0\\-(\mathscr M-g)&1
\end{bmatrix}.\]
To be precise, note that $\mathscr M=\mathscr M(k_1,k_2,k_3)$ depends on the given a non-negative integer triplet $(k_1,k_2,k_3)$, so that $\psi_g$ is determined by the values of $(k_1,k_2,k_3)$.
We will see the properties of $\psi_g$.

\begin{proposition}\label{pr:Markov-monodromy-cohn-matrix}
Let $g\in \mathbb Z[x_1^{\pm1},x_2^{\pm1},x_3^{\pm1}]$.

\begin{itemize}\setlength{\leftskip}{-15pt}
\item[(1)] If $X$ is a $(k_1,k_2,k_3)$-CMM matrix associated with $a$, then $\psi_g(X)$ is a $(k_1,k_2,k_3)$-CGC matrix associated with $a$.
\item[(2)] If $P$ is a $(k_1,k_2,k_3)$-CGC matrix associated with $a$, then $\psi_g^{-1}(P)$ is a $(k_1,k_2,k_3)$-CMM matrix associated with $a$.
\end{itemize}
In particular, $\psi_g$ induces a bijection between the set of $(k_1,k_2,k_3)$-CMM matrices associated with $a$ and the set of $(k_1,k_2,k_3)$-CGC matrices associated with $a$.
\end{proposition}

\begin{proof}
First, we will prove (1).
Since \[\det\begin{bmatrix}1&0\\g&1\end{bmatrix}=\det\begin{bmatrix}1&0\\\mathscr M-g&1\end{bmatrix}=1\text{ and }\det\begin{bmatrix}
    x_{22}&x_{12}\\x_{21}&x_{11}
\end{bmatrix}=\det X=1,\] we have $\det(\psi_g(X))=1$ and $\psi_{g}(X)\in SL(2,\ZZ[x_1^{\pm1},x_2^{\pm1}.x_3^{\pm1}])$. Next, we evaluate the trace of $\psi_g(X)$. We have
\begin{align*}\mathrm{tr}\left(\begin{bmatrix}1&0\\\mathscr M-g&1\end{bmatrix}\begin{bmatrix}
    x_{22}&x_{12}\\x_{21}& x_{11}
\end{bmatrix}\begin{bmatrix}
    1&0\\g&1
\end{bmatrix}\right)&=\mathrm{tr}\left(\begin{bmatrix}
    x_{22}&x_{12}\\x_{21}& x_{11}
\end{bmatrix}\begin{bmatrix}
    1&0\\g&1
\end{bmatrix}\begin{bmatrix}1&0\\\mathscr M-g&1\end{bmatrix}\right)\\&=\mathrm{tr}\left(\begin{bmatrix}
    x_{22}&x_{12}\\x_{21}& x_{11}
\end{bmatrix}\begin{bmatrix}
    1&0\\\mathscr M&1
\end{bmatrix}\right)\\
&=\mathrm{tr}\left(\begin{bmatrix}
    x_{22}+\mathscr{M}x_{12}&x_{12}\\x_{11}+\mathscr Mx_{11}&x_{11}
\end{bmatrix}\right)\\
&=\mathscr{M}x_{12}+(x_{11}+x_{22})=\mathscr{M}x_{12}-k_a,
\end{align*}
as desired. Therefore, by definition, $\psi_g(X)$ is a $(k_1,k_2,k_3)$-CGC matrix associated with $a$, as desired.

Next, we will prove (2).  Since \[\det\begin{bmatrix}1&0\\-g&1\end{bmatrix}=\det\begin{bmatrix}1&0\\\mathscr -(M-g)&1\end{bmatrix}=1\text{ and }\det\begin{bmatrix}
    p_{22}&p_{12}\\p_{21}&p_{11}
\end{bmatrix}=\det P=1,\] we have $\det(\psi^{-1}_g(P))=1$ and $\psi_{g}^{-1}(P)\in SL(2,\ZZ[x_1^{\pm1},x_2^{\pm1},x_3^{\pm1}])$. Next, we compute the trace of $\psi_{g}^{-1}(P)$. We have
\begin{align*}\mathrm{tr}\left(\begin{bmatrix}1&0\\-g&1\end{bmatrix}\begin{bmatrix}
    p_{22}&p_{12}\\p_{21}& p_{11}
\end{bmatrix}\begin{bmatrix}
    1&0\\-(\mathscr M-g)&1
\end{bmatrix}\right)&=\mathrm{tr}\left(\begin{bmatrix}
    p_{22}&p_{12}\\p_{21}& p_{11}
\end{bmatrix}\begin{bmatrix}
    1&0\\-g&1
\end{bmatrix}\begin{bmatrix}1&0\\-(\mathscr M-g)&1\end{bmatrix}\right)\\&=\mathrm{tr}\left(\begin{bmatrix}
    p_{22}&p_{12}\\p_{21}& p_{11}
\end{bmatrix}\begin{bmatrix}
    1&0\\-\mathscr M&1
\end{bmatrix}\right)\\
&=\mathrm{tr}\left(\begin{bmatrix}
    p_{22}-\mathscr{M}p_{12}&p_{12}\\p_{11}-\mathscr Mp_{11}&p_{11}
\end{bmatrix}\right)\\
&=-\mathscr{M}p_{12}+(p_{11}+p_{22})=-k_a,
\end{align*} as desired. Therefore, by definition, $\psi^{-1}_g(P)$ is a $(k_1,k_2,k_3)$-CMM matrix associated with $a$.
\end{proof}

For a $(k_1,k_2,k_3)$-CMM triple $(X,Y,Z)$ associated with $(a,b,c)$, we set
\[\Psi_g(X,Y,Z):=(\psi_{g}(X),\psi_{g}(Y),\psi_{g}(Z)),\]
and for each $(k_1,k_2,k_3)$-CGC triple $(P,Q,R)$ associated with $(a,b,c)$, we set
\[\Psi^{-1}_g(P,Q,R):=(\psi^{-1}_{g}(P),\psi^{-1}_{g}(Q),\psi^{-1}_{g}(R)).\]

We will show that $\Psi_g$ respects the tree structure of CGC triples and CMM triples.

\begin{proposition}\label{prop:psi-cohn-triples}\label{pr:Markov-monodromy-cohn-triple}
Let $g\in \mathbb Z[x_1^{\pm1},x_2^{\pm1},x_3^{\pm1}]$.
\begin{itemize}\setlength{\leftskip}{-15pt}
\item[(1)] If $(X,Y,Z)$ is a $(k_1,k_2,k_3)$-CMM triple associated with $(a,b,c)$, $\Psi_g(X,Y,Z)$ is a $(k_1,k_2,k_3)$-CGC triple associated with $(a,b,c)$,
\item[(2)] If $(P,Q,R)$ is a $(k_1,k_2,k_3)$-CGC triple associated with $(a,b,c)$, $\Psi_g^{-1}(P,Q,R)$ is a $(k_1,k_2,k_3)$-CMM triple associated with $(a,b,c)$.
\end{itemize}
In particular, $\Psi_g$ induces a bijection between the set of $(k_1,k_2,k_3)$-CMM triples associated with $(a,b,c)$ and the set of $(k_1,k_2,k_3)$-CGC triples associated with $(a,b,c)$.
\end{proposition}
\begin{proof}
We begin with item (1). We assume that $(X,Y,Z)$ is a $(k_1,k_2,k_3)$-CMM triple associated with $(a,b,c)$. By Proposition \ref{pr:Markov-monodromy-cohn-matrix} (1), it suffices to show that $\psi_{g}(Y)=\psi_{g}(X)\psi_{g}(Z)-S_{\psi_{g}(Y)}$.

We set $X=\begin{bmatrix}
    x_{11}&x_{12}\\x_{21}&x_{22}
\end{bmatrix}$ and $Z=\begin{bmatrix}
    z_{11}&z_{12}\\z_{21}&z_{22}
\end{bmatrix}.$
By assumption, we have $\psi_{g}(Y)=\psi_{g}(X^{-1}TZ^{-1})$. Therefore, it is enough to show that $\psi_{g}(X)\psi_{g}(Z)-\psi_{g}(X^{-1}TZ^{-1})=S_{\psi_{g}(Y)}$.
Calculating  $\psi_{g}(X)\psi_{g}(Z)$, we have
\begin{align*}
    &\psi_{g}(X)\psi_{g}(Z)\\&=\begin{bmatrix}
    1 &0\\\mathscr M-g&1
\end{bmatrix}\begin{bmatrix}
    x_{22}&x_{12}\\x_{21}&x_{11}
\end{bmatrix}\begin{bmatrix}
    1 &0\\g&1
\end{bmatrix}\begin{bmatrix}
    1 &0\\\mathscr M-g&1
\end{bmatrix}\begin{bmatrix}z_{22}&z_{12}\\z_{21}&z_{11}\end{bmatrix}\begin{bmatrix}
    1 &0\\g&1
\end{bmatrix}\\&=\begin{bmatrix}
    1 &0\\\mathscr M-g&1
\end{bmatrix}\begin{bmatrix}
    x_{22}&x_{12}\\x_{21}&x_{11}
\end{bmatrix}\begin{bmatrix}
    1 &0\\\mathscr M&1
\end{bmatrix}\begin{bmatrix}z_{22}&z_{12}\\z_{21}&z_{11}\end{bmatrix}\begin{bmatrix}
    1 &0\\g&1
\end{bmatrix}\\
&=\begin{bmatrix}
    1&0\\\mathscr{M}-g&1
\end{bmatrix}\begin{bmatrix}
   (x_{22}+\mathscr M x_{12})z_{22}+x_{12}z_{21} &(x_{22}+\mathscr{M}x_{12})z_{12}+x_{12}z_{11}\\(x_{21}+\mathscr{M}x_{11})z_{22}+x_{11}z_{21}&(x_{21}+\mathscr{M}x_{11})z_{12}+x_{11}z_{11}
\end{bmatrix}\begin{bmatrix}
    1&0\\g&1
\end{bmatrix}.\end{align*}
On the other hand, since \begin{align*}
    X^{-1}TZ^{-1}&=\begin{bmatrix}
    x_{22}&-x_{12}\\-x_{21}&x_{11}
\end{bmatrix}\begin{bmatrix}
    -1 &0\\\mathscr M&-1
\end{bmatrix}\begin{bmatrix}z_{22}&-z_{12}\\-z_{21}&z_{11}\end{bmatrix}\\&=\begin{bmatrix}
    -((x_{22}+\mathscr{M}x_{12})z_{22}+x_{12}z_{21})&(x_{22}+\mathscr{M}x_{12})z_{12}+x_{12}z_{11}\\(x_{21}+\mathscr{M}x_{11})z_{22}+x_{11}z_{21}& -((x_{21}+\mathscr M x_{11})z_{12}+x_{11}z_{11})
\end{bmatrix},\end{align*}
We have
\begin{align*}&\psi_g(X^{-1}TZ^{-1})\\&=\begin{bmatrix}
    1&0\\\mathscr{M}-g&1
\end{bmatrix}\begin{bmatrix}
   -((x_{21}+\mathscr M x_{11})z_{12}+x_{11}z_{11}) &(x_{22}+\mathscr{M}x_{12})z_{12}+x_{12}z_{11}\\(x_{21}+\mathscr{M}x_{11})z_{22}+x_{11}z_{21}&-((x_{22}+\mathscr{M}x_{12})z_{22}+x_{12}z_{21})
\end{bmatrix}\begin{bmatrix}
    1&0\\g&1
\end{bmatrix}.
\end{align*}
Therefore, we have
\begin{align*}
&\psi_{g}(X)\psi_{g}(Z)-\psi_g(X^{-1}TZ^{-1})=\begin{bmatrix}
    1&0\\\mathscr{M}-g&1
\end{bmatrix}\begin{bmatrix}
   \alpha &0\\0&\alpha
\end{bmatrix}\begin{bmatrix}
    1&0\\g&1
\end{bmatrix}=\begin{bmatrix}
\alpha&0\\\alpha\mathscr M& \alpha
\end{bmatrix},
\end{align*}
where $\alpha=(x_{22}+\mathscr{M}x_{12})z_{22}+x_{12}z_{21}+(x_{21}+\mathscr M x_{11})z_{12}+x_{11}z_{11}$.
Since
\[-k_b=\mathrm{tr}(Y)=\mathrm{tr}(X^{-1}TZ^{-1})=-\alpha,\]
we have
\[\psi_{g}(X)\psi_{g}(Z)-\psi_g(X^{-1}TZ^{-1})=\begin{bmatrix}
\alpha&0\\\alpha\mathscr M& \alpha
\end{bmatrix}=\begin{bmatrix}
k_b&0\\k_b\mathscr M& k_b
\end{bmatrix}=S_{\psi_g(Y)},\]
as desired.

Next, we will prove item (2). By Proposition \ref{pr:Markov-monodromy-cohn-matrix} (2), it suffices to show that \[\psi_{g}^{-1}(P)\psi_{g}^{-1}(Q)\psi^{-1}_{g}(R)=T.\] By assumption, it is equivalent to show that \[\psi_{g}^{-1}(PR-S_Q)=(\psi_{g}^{-1}(P))^{-1}T(\psi_{g}^{-1}(R))^{-1}.\]
We set $P=\begin{bmatrix}
    p_{11}& p_{12}\\p_{21}&p_{22}
\end{bmatrix}$ and $R=\begin{bmatrix}
    r_{11}& r_{12}\\r_{21}& r_{22}
\end{bmatrix}$.
By a direct calculation, we have
\begin{align*}
    PR-S_Q=\begin{bmatrix}
        p_{11}r_{11}+p_{12}r_{21}-k_b&p_{11}r_{12}+p_{12}r_{22}\\p_{21}r_{11}+p_{22}r_{21}-k_b\mathscr{M}&p_{21}r_{12}+p_{22}r_{22}-k_b
    \end{bmatrix}.
\end{align*}
On the one hand, by applying $\psi_{g}^{-1}$, we have
\begin{align*}\psi_{g}^{-1}(PR-S_Q)=\begin{bmatrix}
1&0\\-g&1
\end{bmatrix}\begin{bmatrix}
       p_{21}r_{12}+p_{22}r_{22}-k_b&p_{11}r_{12}+p_{12}r_{22}\\p_{21}r_{11}+p_{22}r_{21}-k_b\mathscr{M}& p_{11}r_{11}+p_{12}r_{21}-k_b
    \end{bmatrix}\begin{bmatrix}
1&0\\-(\mathscr M-g)&1
\end{bmatrix}.
\end{align*}
On the other hand, through direct calculation, we have \begin{align*}&(\psi^{-1}(P))^{-1}T(\psi^{-1}(R))^{-1}\\
&=\begin{bmatrix}
    1&0\\\mathscr{M}-g&1
\end{bmatrix}\begin{bmatrix}
    p_{11}&-p_{12}\\-p_{21}&p_{22}
\end{bmatrix}\begin{bmatrix}
    1&0\\g&1
\end{bmatrix}\begin{bmatrix}
    -1&0\\\mathscr{M}&-1
\end{bmatrix}\begin{bmatrix}
    1&0\\\mathscr{M}-g&1
\end{bmatrix}\begin{bmatrix}
    r_{11}&-r_{12}\\-r_{21}&r_{22}
\end{bmatrix}\begin{bmatrix}
    1&0\\g&1
\end{bmatrix}\\
&=-\begin{bmatrix}
    1&0\\\mathscr{M}-g&1
\end{bmatrix}\begin{bmatrix}
    p_{11}&-p_{12}\\-p_{21}&p_{22}
\end{bmatrix}\begin{bmatrix}
    r_{11}&-r_{12}\\-r_{21}&r_{22}
\end{bmatrix}\begin{bmatrix}
    1&0\\g&1
\end{bmatrix}\\
&=\begin{bmatrix}
    1&0\\\mathscr{M}-g&1
\end{bmatrix}\begin{bmatrix}-p_{11}r_{11}-p_{12}r_{21}&p_{11}r_{12}+p_{12}r_{22}\\p_{21}r_{11}+p_{22}r_{21}&-p_{21}r_{12}-p_{22}r_{22}
\end{bmatrix}\begin{bmatrix}
    1&0\\g&1
\end{bmatrix}
\end{align*}
Comparing these two computations yields
\begin{align*}
&\psi_{g}^{-1}(PR-S_Q)-\psi^{-1}(P))^{-1}T(\psi^{-1}(R))^{-1}\\&=(\mathscr Mp_{11}r_{12} + \mathscr Mp_{12}r_{22} - p_{11}r_{11} - p_{21}r_{12} - p_{12}r_{21} - p_{22}r_{22} + k_b)I,\end{align*}
where $I$ is the identity matrix.
Therefore, it suffices to show that
\begin{align}
\mathscr Mp_{11}r_{12} + \mathscr Mp_{12}r_{22} - p_{11}r_{11} - p_{21}r_{12} - p_{12}r_{21} - p_{22}r_{22} + k_b=0.\label{eq:ShowIsZero}
\end{align}

Since $(P,Q,R)$ is a $(k_1,k_2,k_3)$-CGC triple, we have
\[\mathrm{tr}(Q)=\mathscr{M}q_{12}-k_b=\mathscr{M}(p_{11}r_{12}+p_{12}r_{22})-k_b,\]
and by calculation we have \[\mathrm{tr}(PR-S_Q)=(p_{11}r_{11}+p_{12}r_{21}-k_b)+(p_{21}r_{12}+p_{22}r_{22}-k_b).\]

We also know $Q=PR-S_Q$, and thus  $\mathrm{tr}(Q)=\mathrm{tr}(PR-S_Q)$.
Therefore, we conclude
\[\mathscr{M}(p_{11}r_{12}+p_{12}r_{22})-k_b=(p_{11}r_{11}+p_{12}r_{21}-k_b)+(p_{21}r_{12}+p_{22}r_{22}-k_b).\]
This equality is equivalent to the Equation \eqref{eq:ShowIsZero}, thus completing the proof.
\end{proof}

The following lemma implies that $\Psi_g$ is compatible with the generation rules of two trees.

\begin{lemma}\label{lem:psi-cohn}
Let $g\in \mathbb Z[x_1^{\pm1},x_2^{\pm1},x_3^{\pm1}]$.For a $(k_1,k_2,k_3)$-CMM triple $(X,Y,Z)$, we have
\begin{align*}
\Psi_g(X,YZY^{-1},Y)&=(\psi_{g}(X),\psi_{g}(X)\psi_{g}(Y)-S_{\psi_{g}(Z)},\psi_{g}(Y)),\\
\Psi_g(Y,Y^{-1}XY,Z)&=(\psi_{g}(Y),\psi_{g}(Y)\psi_{g}(Z)-S_{\psi_{g}(X)},\psi_{g}(Z)).
\end{align*}
\end{lemma}

\begin{proof}
We assume that $(X,Y,Z)$ is associated with $(a,b,c)$. Our goal is to show \begin{align}\psi_{g}(YZY^{-1})=\psi_{g}(X)\psi_{g}(Y)-S_{\psi_{g}(Z)}\label{eq:Lem6.3}\end{align}
and
\begin{align*} \psi_{g}(Y^{-1}XY)=\psi_{g}(Y)\psi_{g}(Z)-S_{\psi_{g}(X)}.\end{align*}
Due to the symmetry of the two statements, we will only prove Equation \eqref{eq:Lem6.3}. By Proposition \ref{pr:Markov-monodromy-cohn-triple}, we have $\psi_{g}(Y)=\psi_{g}(X)\psi_{g}(Z)-S_{\psi_{g}(Y)}$. By using this to rewrite $\psi_g(X)$ in Equation \eqref{eq:Lem6.3}, we have
\[\psi_{g}(X)\psi_{g}(Y)-S_{\psi_{g}(Z)}=\psi_{g}(Y)\psi_{g}(Z)^{-1}\psi_{g}(Y)+S_{\psi_{g}(Y)}\psi_{g}(Z)^{-1}\psi_{g}(Y)-S_{\psi_{g}(Z)}.\]
Therefore, we will show that
\begin{equation}\label{eq:BT-CT2}
\psi_{g}(YZY^{-1})=\psi_{g}(Y)\psi_{g}(Z)^{-1}\psi_{g}(Y)+S_{\psi_{g}(Y)}\psi_{g}(Z)^{-1}\psi_{g}(Y)-S_{\psi_{g}(Z)}.
\end{equation}
We set $Y=\begin{bmatrix}
    y_{11}&y_{12}\\y_{21}&y_{22}
\end{bmatrix}$ and $Z=\begin{bmatrix}
    z_{11}&z_{12}\\z_{21}&z_{22}
\end{bmatrix}$. By setting the left hand side of \eqref{eq:BT-CT2} as $\begin{bmatrix}
    \alpha_{11}& \alpha_{12}\\\alpha_{21} &\alpha_{22}
\end{bmatrix}$ and the right hand side of \eqref{eq:BT-CT2} as  $\begin{bmatrix}
    \beta_{11}& \beta_{12}\\\beta_{21} &\beta_{22}
\end{bmatrix}$.
Then, by a direct calculation, we have
\begin{align*}
    \alpha_{11}&=-gy_{11}y_{12}z_{11} + gy_{11}^2z_{12} - gy_{12}^2z_{21} + gy_{11}y_{12}z_{22} - y_{12}y_{21}z_{11} + y_{11}y_{21}z_{12} - y_{12}y_{22}z_{21} \\&+ y_{11}y_{22}z_{22} \\
    \alpha_{12}&=-y_{11}y_{12}z_{11} + y_{11}^2z_{12} - y_{12}^2z_{21} + y_{11}y_{12}z_{22},\\
    \alpha_{22}&=-\mathscr My_{11}y_{12}z_{11} + gy_{11}y_{12}z_{11} + \mathscr My_{11}^2z_{12} - gy_{11}^2z_{12} - \mathscr My_{12}^2z_{21} + gy_{12}^2z_{21} + \mathscr My_{11}y_{12}z_{22} \\&- gy_{11}y_{12}z_{22} + y_{11}y_{22}z_{11} - y_{11}y_{21}z_{12} + y_{12}y_{22}z_{21} - y_{12}y_{21}z_{22},\\
    \beta_{11}&=gk_by_{12}z_{11} + gy_{12}y_{22}z_{11} - gk_by_{11}z_{12} - gy_{11}y_{22}z_{12} - gy_{12}^2z_{21} + gy_{11}y_{12}z_{22} + k_by_{22}z_{11} \\&+ y_{22}^2z_{11} - k_by_{21}z_{12} - y_{21}y_{22}z_{12} - y_{12}y_{22}z_{21} + y_{12}y_{21}z_{22} - k_c\\
    \beta_{12}&=k_by_{12}z_{11} + y_{12}y_{22}z_{11} - k_by_{11}z_{12} - y_{11}y_{22}z_{12} - y_{12}^2z_{21} + y_{11}y_{12}z_{22}\\
    \beta_{22}&=\mathscr Mk_by_{12}z_{11} - gk_by_{12}z_{11} + \mathscr{M}y_{12}y_{22}z_{11} - gy_{12}y_{22}z_{11} - \mathscr Mk_by_{11}z_{12} + gk_by_{11}z_{12} \\&- \mathscr My_{11}y_{22}z_{12} + gy_{11}y_{22}z_{12} - \mathscr My_{12}^2z_{21} + gy_{12}^2z_{21} + \mathscr My_{11}y_{12}z_{22} - gy_{11}y_{12}z_{22} + y_{12}y_{21}z_{11} \\&- y_{11}y_{21}z_{12} - k_by_{12}z_{21} - y_{11}y_{12}z_{21} + k_by_{11}z_{22} + y_{11}^2z_{22} - k_c.
\end{align*}
Using $y_{11}+y_{22}=-k_b$, we can check $\alpha_{12}=\beta_{12}$. Using $y_{11}+y_{22}=-k_b$, $z_{11}+z_{22}=-k_c$, and $y_{11}y_{22}-y_{12}y_{21}=1$, we can check $\alpha_{11}=\beta_{11}$ and $\alpha_{22}=\beta_{22}$.
Moreover, $\det(\psi_{g}(YZY^{-1}))$ is $1$ since  $\det(YZY^{-1})=1$ and $\psi_{g}$ preserves determinant. On the other hand, we have
\[\det(\psi_{g}(Y)\psi_{g}(Z)^{-1}\psi_{g}(Y)+S_{\psi_{g}(Y)}\psi_{g}(Z)^{-1}\psi_{g}(Y)-S_{\psi_{g}(Z)})=1\] because $\psi_{g}(X)\psi_{g}(Y)-S_{\psi_{g}(Z)}$ is a $(k_1,k_2,k_3)$-CGC triple.
Therefore, we have verified Equation \eqref{eq:BT-CT2} and this finishes the proof.
\end{proof}
We now conclude that $\Psi_g$ extends to an isomorphisms of the trees studied here.
\begin{theorem}\label{thm:BT-CT2}
The map $\Psi_g$ induces the canonical graph isomorphism\[\mathrm{CMM}\mathbb T(k_1,k_2,k_3,\sigma,f)\to \mathrm{CGC}\mathbb T(k_1,k_2,k_3,\sigma, -f+gx_{\sigma(1)}-k_{\sigma(1)}).\]
\end{theorem}
\begin{proof}
We can check the correspondence between the roots by a direct calculation. We can prove other vertices by Lemma \ref{lem:psi-cohn} and induction.
\end{proof}

Theorem \ref{thm:BT-CT2} is powerful as it allows us to quickly deduce properties of CMM triples by using the analogous properties of CGC triples.

\begin{proof}[Proof of Theorem \ref{thm:Markov-monodoromy-well-defined} (1)]
This follows from Theorem \ref{thm:Cohn-well-defined} (1) and Theorem \ref{thm:BT-CT2}.
\end{proof}

\begin{proof}[Proof of Corollary \ref{cor:MMT-MT}]
This follows from Corollary \ref{cor:CT-MT} and Theorem \ref{thm:BT-CT2}.
\end{proof}

\begin{proof}[Proof of Theorem \ref{thm:all-Markov-triple}]
Apply the bijection $\Psi_0$ to $(X,Y,Z)$. By Proposition \ref{prop:psi-cohn-triples}, its image is a CGC triple associated with the same non-labeled cluster $(a,b,c)$; since $\Psi_0$ preserves the $(1,2)$-entries, the parities of its components are also the same. Theorem \ref{thm:all-cohn-triple} gives a unique permutation $\sigma$, a unique parameter $h$, and a unique vertex in $\mathrm{CGC}\mathbb T(k_1,k_2,k_3,\sigma,h)$ representing this image. By Theorem \ref{thm:BT-CT2},
\[
\Psi_0\colon \mathrm{CMM}\mathbb T(k_1,k_2,k_3,\sigma,f)
\longrightarrow
\mathrm{CGC}\mathbb T(k_1,k_2,k_3,\sigma,-f-k_{\sigma(1)})
\]
is a canonical graph isomorphism. Therefore
\[
f=-h-k_{\sigma(1)}
\]
is uniquely determined, and the inverse graph isomorphism gives a unique vertex representing $(X,Y,Z)$. This proves both existence and all parts of the uniqueness assertion.
\end{proof}

The map $\Psi_g$ is also compatible with generation rules of the inverse trees. These statements can be shown in parallel methods to Lemma \ref{lem:psi-cohn} and Theorem \ref{thm:BT-CT2}.

\begin{lemma}\label{lem:psi-cohn-inverse}
Let $g\in SL(2,\mathbb Z[x_1^{\pm1},x_2^{\pm1},x_3^{\pm1}])$. For a $(k_1,k_2,k_3)$-CMM triple $(X,Y,Z)$ associated with $(a,b,c)$, we have
\begin{align*}
\Psi_g(X,Z,Z^{-1}YZ)&=(\psi_{g}(X),\psi_{g}(Z),\psi_{g}(X)^{-1}(\psi_{g}(Z)+S_{\psi_{g}(Y)})),\\
\Psi_g(XYX^{-1},X,Z)&=((\psi_{g}(X)+S_{\psi_{g}(Y)})\psi_{g}(Z)^{-1},\psi_{g}(X),\psi_{g}(Z)).
\end{align*}
\end{lemma}

\begin{theorem}\label{thm:BT-CT2-inverse}
The map $\Psi_g$ induces the canonical graph isomorphism\[\mathrm{CMM}\mathbb T^\dag(k_1,k_2,k_3,\sigma,f)\to \mathrm{CGC}\mathbb T^\dag(k_1,k_2,k_3,\sigma, -f+gx_{\sigma(1)}-k_{\sigma(1)}).\]
\end{theorem}

\begin{remark}
A similar map is introduced in \cite[Section 5]{gyoda2024sl}, working with similar matrices but under the assumption $k_1 = k_2 = k_3$ and setting all initial cluster variables to 1.  Working with these assumptions, the map provided therein can be seen to be $\psi_\frac{k_\sigma(1)}{x_\sigma(1)}$.
\end{remark}
\subsection{Markov-Monodromy Decompositions of CGC Matrices}
We introduce another connection between $(k_1,k_2,k_3)$-CGC triples and $(k_1,k_2,k_3)$-CMM triples.
\begin{definition}\label{def:Markov-monodromy-dec}
 We fix $k_1,k_2,k_3\in \mathbb Z_{\geq 0}$. For a $(k_1,k_2,k_3)$-CGC triple $(P,Q,R)$ associated with $(a,b,c)$, we consider a triple $(X,Y,Z)$ satisfying the following conditions:
\begin{itemize}\setlength{\leftskip}{-15pt}
\item [(1)] $(X,Y,Z)$ is a $(k_1,k_2,k_3)$-CMM triple,
\item [(2)] $P=-(YZ)^{-1},\quad Q=-(XZ)^{-1},\quad R=-(XY)^{-1}$.
\item[(3)] $\mathrm{tr}(X)=-k_a$, $\mathrm{tr}(Y)=-k_b$, $\mathrm{tr}(Z)=-k_c$.
\end{itemize}
This triple $(X,Y,Z)$ of $(P,Q,R)$ is called a \emph{Markov-monodromy decomposition of} $(P,Q,R)$.
\end{definition}
In this paper, we abbreviate the Markov-monodromy decomposition as the \emph{MM decomposition}.

\begin{remark}
MM decompositions of generalized Cohn triples are also discussed in \cite{gyoda2024sl}. This paper works under the assumption  $k_1=k_2=k_3$. In this special case, it was enough to assume $\mathrm{tr}(X)=\mathrm{tr}(Y)=\mathrm{tr}(Z)$; in particular, it was not necessary to assume that this value is $-k_1$.
However, in the current case, that is, $k_1,k_2,k_3$ are possibly distinct, we can not apply the same arguments as in \cite{gyoda2024sl}.
\end{remark}

\begin{lemma}\label{thm:Markov-monodromy-uniqueness-initial}
The triple $(\widetilde{X_{\sigma}}(f), \widetilde{Y_{\sigma}}(f),\widetilde{Z_{\sigma}}(f))$ is a MM decomposition of $(P_\sigma(f),Q_\sigma(f),R_\sigma(f))$. Moreover, there are no other MM decompositions of $(P_\sigma(f),Q_\sigma(f),R_\sigma(f))$.
\end{lemma}

\begin{proof}
Set
\[
(P,Q,R)=(P_\sigma(f),Q_\sigma(f),R_\sigma(f))
\quad\text{and}\quad
(\widetilde X,\widetilde Y,\widetilde Z)=(\widetilde{X_\sigma}(f),\widetilde{Y_\sigma}(f),\widetilde{Z_\sigma}(f)).
\]
By Proposition \ref{prop:CMM-initial-triples} and Theorem \ref{thm:Markov-monodoromy-well-defined}, $(\widetilde X,\widetilde Y,\widetilde Z)$ is a CMM triple. Its three matrices have traces $-k_{\sigma(1)}$, $-k_{\sigma(2)}$, and $-k_{\sigma(3)}$, respectively, so conditions (1) and (3) in Definition \ref{def:Markov-monodromy-dec} are satisfied. A direct calculation from the explicit formulas gives
\[
\widetilde X=-TP,\qquad \widetilde Z=-RT.
\]
Since $\widetilde X\widetilde Y\widetilde Z=T$, these identities imply
\[
\widetilde Y=P^{-1}T^{-1}R^{-1}.
\]
Consequently,
\[
\widetilde Y\widetilde Z=-P^{-1},\qquad
\widetilde X\widetilde Y=-R^{-1}.
\]
To verify the remaining identity, write $Q=(q_{ij})$ and set $k_Q=\mathscr M q_{12}-\mathrm{tr}(Q)$. Since
\[
S_Q=\begin{bmatrix}k_Q&0\\k_Q\mathscr M&k_Q\end{bmatrix},
\]
a direct multiplication, using $\det Q=1$, gives
\[
T(Q+S_Q)T=
\begin{bmatrix}-q_{22}&q_{12}\\q_{21}&-q_{11}\end{bmatrix}=-Q^{-1}.
\]
The CGC relation $Q=PR-S_Q$ therefore yields
\[
\widetilde X\widetilde Z=TPRT=T(Q+S_Q)T=-Q^{-1}.
\]
Thus all conditions in Definition \ref{def:Markov-monodromy-dec} hold, proving existence.

The same calculation gives uniqueness for any MM decomposition, not only for the initial triples. Indeed, if $(X,Y,Z)$ is an MM decomposition of a CGC triple $(P,Q,R)$, then
\[
YZ=-P^{-1},\qquad XY=-R^{-1},\qquad XYZ=T.
\]
Hence
\[
X=-TP,\qquad Z=-RT,
\]
and then
\[
Y=P^{-1}T^{-1}R^{-1}.
\]
Thus $(X,Y,Z)$ is uniquely determined by $(P,Q,R)$ and must equal $(\widetilde X,\widetilde Y,\widetilde Z)$.
\end{proof}

\setcounter{theorem}{12}
\setcounter{equation}{4}
Our next goal is to characterize MM decompositions.

\begin{theorem}\label{thm:Markov-monodromy-uniqueness}
If $(P,Q,R)$ is a $(k_1,k_2,k_3)$-CGC triple associated with $(a,b,c)$ such that $|b|> \max\{|a|,|c|\}$, then there exists a unique MM decomposition $(X,Y,Z)$ of $(P,Q,R)$.
\end{theorem}

Before proving Theorem \ref{thm:Markov-monodromy-uniqueness}, we provide some preparatory results. The first demonstrates how these decompositions interact with the tree structure on each set of matrices.

\begin{proposition}\label{thm:BT-CT}
Let $(P,Q,R)$ be a $(k_1,k_2,k_3)$-CGC triple. If $(X,Y,Z)$ is an MM decomposition of $(P,Q,R)$, then $(X,Z,Z^{-1}YZ)$ (resp. $(XYX^{-1},X,Z)$) is an MM decomposition of $(P,PQ-S_R,Q)$ (resp. $(Q,QR-S_P,R)$).
\end{proposition}

\begin{proof}
 We assume $(X,Y,Z)$ is associated with $(a,b,c)$. We denote by $(X',Y',Z')$ (resp. $(X'',Y'',Z'')$) the left (resp. right) child of $(X,Y,Z)$ in $\mathrm{CMM}\TT^\dag(k,\ell)$. We can easily see that  $(X',Y',Z')$ and $(X'',Y'',Z'')$ satisfy (1) and (3) in Definition \ref{def:Markov-monodromy-dec}.  We will prove that they satisfy (2). This entails showing
\begin{align*}
    -(Y'Z')^{-1}&=-(YZ)^{-1},\\
    -(X'Z')^{-1}&=(YZ)^{-1}(XZ)^{-1}-S_R,\\
    -(X'Y')^{-1}&=-(XZ)^{-1},\\
    -(Y''Z'')^{-1}&=-(XZ)^{-1},\\
    -(X''Z'')^{-1}&=(XZ)^{-1}(XY)^{-1}-S_P, \text{ and}\\
    -(X''Y'')^{-1}&=-(XY)^{-1}.
\end{align*}
All but the second and fifth equality are clear, and due to the symmetry involved, we will only prove the second. This is equivalent to showing
\begin{align}\label{eq:GC-MM-decomposition-relation}
-Z^{-1}Y^{-1}ZX^{-1}=Z^{-1}Y^{-1}Z^{-1}X^{-1}-S_R.\end{align}
By assumption, $XYZ=T$, and one can also calculate $TS_R=-k_cI$. These two equalities give us
\begin{align*}
    YZS_RX=X^{-1}XYZS_RX=X^{-1}TS_RX=-k_cX^{-1}IX=-k_cI.
\end{align*}
Since $\det(Z) = 1$, we know $Z + Z^{-1} = \mathrm{tr}(Z)I$. Using this and the previous computation, we have
\[Z+Z^{-1}=-k_cI=YZS_RX.\]
 By multiplying $Z^{-1}Y^{-1}$ to the left and $X^{-1}$ to the right by the above equality, we have \eqref{eq:GC-MM-decomposition-relation}, as desired.
\end{proof}

In parallel with Proposition \ref{thm:BT-CT}, we have the following lemma.
\begin{lemma}\label{lem:inverseBT-CT}
If $(X,Y,Z)$ is an MM decomposition of $(P,Q,R)$, then $(X,YZY^{-1},Y)$ (resp. $(Y,Y^{-1}XY,Z)$) is an MM decomposition of $(P,R,P^{-1}(R+S_R))$ (resp. $((P+S_P)R^{-1},P,R)$).
\end{lemma}

\begin{proof}[Proof of Theorem \ref{thm:Markov-monodromy-uniqueness}]
By Theorem \ref{thm:all-cohn-triple}, there exist $\sigma\in \mathfrak{S}_3$ and $f\in \ZZ[x_1^{\pm1},x_2^{\pm1},x_3^{\pm1}]$ such that $(P,Q,R)\in \mathrm{CGC}\mathbb T(k_1,k_2,k_3,\sigma,f)$. The existence of a MM decomposition follows from Lemma \ref{thm:Markov-monodromy-uniqueness-initial} and Proposition \ref{thm:BT-CT}. Uniqueness follows from the general calculation in the proof of Lemma \ref{thm:Markov-monodromy-uniqueness-initial}: every MM decomposition of $(P,Q,R)$ is necessarily
\[
(X,Y,Z)=(-TP,\,P^{-1}T^{-1}R^{-1},\,-RT).
\]
\end{proof}

\begin{lemma}\label{thm:Markov-monodromy-uniqueness-initial2}
 The triple $(X_{\sigma}(f), Y_{\sigma}(f),Z_{\sigma}(f))$ is a MM decomposition of $(\widetilde{P_\sigma}(f),\widetilde{Q_\sigma}(f),\widetilde{R_\sigma
 }(f))$. Moreover, there are no other MM decompositions of $(\widetilde{P_\sigma}(f),\widetilde{Q_\sigma}(f),\widetilde{R_\sigma
 }(f))$.
\end{lemma}

\begin{proof}
It follows from Lemma \ref{thm:Markov-monodromy-uniqueness-initial}, Proposition \ref{thm:BT-CT}, and Theorem \ref{thm:Markov-monodromy-uniqueness}.
\end{proof}

Finally, we introduce an inverse operation of the MM decomposition.

\begin{corollary}\label{cor:tree-iso-BT-CTdag}
We define $\Phi\colon SL(2,\mathbb Z[x_{1}^{\pm1},x_{2}^{\pm1},x_{3}^{\pm1}])^3\to SL(2,\mathbb Z[x_{1}^{\pm1},x_{2}^{\pm1},x_{3}^{\pm1}])^3$ by \[\Phi(X,Y,Z)= (-(YZ)^{-1},-(XZ)^{-1},-(XY)^{-1}).\] The map $\Phi$ induces the canonical graph isomorphism from $\mathrm{CMM}\mathbb T^{\dag}(k_1,k_2,k_3,\sigma,f)$ to $\mathrm{CGC}\mathbb T(k_1,k_2,k_3,\sigma,f)$.
\end{corollary}
\begin{proof}
By Lemma \ref{thm:Markov-monodromy-uniqueness-initial}, $\Phi$ maps from the root of $\mathrm{CMM}\mathbb T^{\dag}(k_1,k_2,k_3,\sigma,f)$ to the root of $\mathrm{CGC}\mathbb T(k_1,k_2,k_3,\sigma,f)$. By using Proposition \ref{thm:BT-CT}, the correspondence between the remaining vertices follows inductively.
\end{proof}

In the previous discussions, we assume that $|b|> \max\{|a|,|c|\}$. We can do the same under the assumption $|b|\leq \max\{|a|,|c|\}$ by considering $\mathrm{CGC}\mathbb T^\dag(k_1,k_2,k_3,\sigma,f)$ and $\mathrm{CMM}\mathbb T(k_1,k_2,k_3,\sigma,f)$ instead of $\mathrm{CGC}\mathbb T(k_1,k_2,k_3,\sigma,f)$ and $\mathrm{CMM}\mathbb T^{\dag}(k_1,k_2,k_3,\sigma,f)$:

\begin{corollary}\label{cor:tree-iso-BT-CT^dag2}
If $(P,Q,R)$ is a $(k_1,k_2,k_3)$-CGC triple associated with $(a,b,c)$ such that $|b|\leq \max\{|a|,|c|\}$, then there is a unique MM decomposition of $(P,Q,R)$.
\end{corollary}

\begin{proof}
Existence is proved as in Theorem \ref{thm:Markov-monodromy-uniqueness}, with the roles of Proposition \ref{thm:BT-CT} and Lemma \ref{lem:inverseBT-CT} exchanged and using Theorem \ref{thm:all-cohn-triple2}. Uniqueness again follows from the formula $(X,Y,Z)=(-TP,P^{-1}T^{-1}R^{-1},-RT)$ established in the proof of Lemma \ref{thm:Markov-monodromy-uniqueness-initial}.
\end{proof}

\begin{corollary}\label{cor:tree-iso-BT-CT^dag3}
The map $\Phi$ induces the canonical graph isomorphism from $\mathrm{CMM}\mathbb T(k_1,k_2,k_3,\sigma,f)$ to $\mathrm{CGC}\mathbb T^\dag(k_1,k_2,k_3,\sigma,f)$.
\end{corollary}
\begin{proof}
 By using Lemma \ref{thm:Markov-monodromy-uniqueness-initial2} in place of Lemma \ref{thm:Markov-monodromy-uniqueness-initial} and Lemma \ref{lem:inverseBT-CT} in place of Proposition \ref{thm:BT-CT} in the proof of Corollary \ref{cor:tree-iso-BT-CTdag}, we have the statement.
\end{proof}

From Corollaries \ref{cor:tree-iso-BT-CTdag} and \ref{cor:tree-iso-BT-CT^dag3}, we have the following:

\begin{corollary}\label{cor:bijection-triples}
The map $\Phi$ induces a bijection from the set of $(k_1,k_2,k_3)$-CMM triples to the set of $(k_1,k_2,k_3)$-CGC triples.
\end{corollary}

To compute the MM decomposition (i.e., the inverse map of $\Phi$), we use the  isomorphism, $\Psi_g$. Since $\Psi_g$ induces the isomorphism from $\mathrm{CMM}\mathbb T^\dag(k_1,k_2,k_3,\sigma,f)$ to $\mathrm{CGC}\mathbb T^\dag(k_1,k_2,k_3,\sigma, -f+gx_{\sigma(1)}-k_{\sigma(1)})$, by setting $g=\frac{k_{\sigma(1)}}{x_{\sigma(1)}}$, the map $\Psi_{\frac{k_{\sigma(1)}}{x_{\sigma(1)}}}^{-1}\circ\Phi$ induces the isomorphism
from $\mathrm{CMM}\mathbb T(k_1,k_2,k_3,\sigma,f)$ to $\mathrm{CMM}\mathbb T^\dag(k_1,k_2,k_3,\sigma,-f)$. Since this is involution, we have the following:

\begin{corollary}\label{cor:Markov-monodromy-decom-algorithm}
    We have $\left(\Psi_{\frac{k_{\sigma(1)}}{x_{\sigma(1)}}}^{-1}\circ\Phi\right)^2=\mathrm{id}$ and $\left(\Phi\circ\Psi^{-1}_{\frac{k_{\sigma(1)}}{x_{\sigma(1)}}}\right)^2=\mathrm{id}$. In particular, the MM decomposition $\Phi^{-1}$ is given by $\Psi_{\frac{k_{\sigma(1)}}{x_{\sigma(1)}}}^{-1}\circ \Phi\circ \Psi_{\frac{k_{\sigma(1)}}{x_{\sigma(1)}}}^{-1}$.
\end{corollary}

\section{Poset Expansion Formulas}\label{sec:Posets}

At this point, we will turn to our second goal of giving an explicit description of all CGC and CMM matrices in one specific tree. In order to do this, we will use and add to a set of combinatorial tools used to understand cluster algebras from surfaces and orbifolds.

These tools are most often used  in the context of principal coefficients, and we will at times also work in this context. We will write $\Acal^\prin(k_1,k_2,k_3)$ for the algebras resulting from using  $\PP=\trop(y_1,y_2,y_3)$ instead of $\PP = \{1\}$ and $\yy=(y_1,y_2,y_3)$ instead of $\yy = (1,1,1)$ in the definition $\mathcal A(k_1,k_2,k_3)$  in Section \ref{subsec:Setting}. In short, this is the principal coefficient version of $\Acal(k_1,k_2,k_3)$.

\subsection{Construction}
Given a triple $(k_1,k_2,k_3) \in \mathbb{Z}_{\geq 0}$ and a rational number $\frac{p}{q} \in \mathbb{Q} \cap [0,\infty)$, we build a poset $\calP_{\frac{p}{q}}(k_1,k_2,k_3)$. We stress that the construction depends on $k_1,k_2,$ and $k_3$, but for brevity we will simply write $\calP_{\fpq}$, as we did earlier with $\mathscr M$.
The construction will be similar to that in \cite{banaian2024skein,ezgieminecluster2024,pilaud2023posets}.  Before we explicitly give the construction, we will describe the features of these posets in general.

The poset $\calP_{\frac{p}{q}}$ will be a \emph{fence poset}, which means it will be a finite poset whose Hasse diagram is isomorphic, as an undirected graph, to a path graph. Every element of $\calP_{\frac{p}{q}}$ will come with both a \emph{label} and a \emph{weight}. The \emph{label} will be one of the initial cluster variables $\{x_1,x_2,x_3\}$. The weight of an element will be a rational multiple of a Laurent monomial in the initial cluster. To this end, we define \[
\widehat{x_1} = \frac{x_2}{x_3}, \quad \widehat{x_2} = \frac{x_3}{x_1}, \quad \widehat{x_3} = \frac{x_1}{x_2}.
\]

We will occasionally work instead with principal coefficients. In this case, define $\widehat{x_i}^{\prin}$ by $\widehat{x_i}y_i$ if $k_i > 0$ and $\widehat{x_i}\sqrt{y_i}$ if $k_i = 0$.

The fact that the Hasse diagram of $\calP_\fpq$ is a path graph invites a natural indexing of the elements reading from left to right. Following \cite{banaian2024skein}, we call this a \emph{chronological ordering}. To be precise, this is a bijection from $\{1,2,\ldots,\vert \calP_\fpq\vert\}$ to $\calP_\fpq$ where $i$ is sent to $\calP_\fpq(i)$ in such a way that $\calP_\fpq(i)$ only has cover relations with $\calP_\fpq(i-1)$ and $\calP_\fpq(i+1)$, when these indices exist. We remark that if $\vert \calP_\fpq \vert > 1$, then there are two possible chronological orderings to place on $\calP_\fpq$, given by reading the Hasse diagram from left to right and from right to left. Note that, unless the poset is a chain, neither of the two chronological orderings is a linear extension.
Given a poset $\calP$ on $n$ elements with chronological ordering, we set $L_\calP:= \calP(1)$ and $R_\calP:= \calP(n)$. If the poset is understood, we may drop the subscript.

Following notation from \cite{ouguz2025oriented,ezgieminecluster2024}, given two posets with chronological ordering, $\calP$ and $\calP'$, we let $\calP \nearrow \calP'$ be the poset on $\calP \sqcup \calP'$ which is the transitive closure of all of the relations in $\calP$ and $\calP'$ as well as $R_\calP < L_{\calP'}$, and define $\calP \searrow \calP'$ similarly.

One advantage to placing a chronological ordering on a poset $\calP$ is that we can now easily reference certain subposets of $\calP_\fpq$. In particular, for $i \leq j$, let $\calP[i,j]$ denote the subposet of $\calP$ resulting from restricting the partial order to $\{\calP(i),\calP(i+1),\ldots,\calP(j)\}$.

We now describe the poset $\calP_{\frac{p}{q}}$.
Throughout this section, let $\mathcal{L}$ be the lattice of all line segments of slope 0, $\infty$, and $-1$ which go through points of $\mathbb{Z}^2$. We label every line segment of slope $-1$ with $x_1$, every line segment of slope $\infty$ with $x_2$, and every line segment of slope $0$ with $x_3$. In this way, we have described the universal cover of the initial triangulation of the once-punctured torus, with lifts of arcs labeled by the Farey triple $(\frac{-1}{1},\frac{0}{1},\frac11)$. Unlike the triples in Section \ref{subsec:LabelClVar}, this triple contains a negative number. Recall that we used positive Farey triples to label cluster variables in $k_i\mathbb{T}(k_1,k_2,k_3)$, and none of these trees contain the initial cluster. This triple can be seen as the parent of $(\frac01, \frac11, \frac10)$ in the Farey tree and hence it is a reasonable way to label the initial cluster. This indexing is standard in the literature as well, as in \cite{chavez2011c,musiker2025super}.

Let $\gamma_\fpq$ be the line segment between $(0,0)$ and $(q,p)$;  this could be viewed as a lift of an arc with rational slope on the once-punctured torus. We consider $\gamma_\fpq$ to be oriented from $(0,0)$ to $(q,p)$. This gives a linear ordering on the intersections of $\gamma_\fpq$ with the line segments in $\mathcal{L}$.

To construct $\calP_\fpq$, we do the following steps in order based on this order of the intersections. We will refer to each line segment's left and right endpoint with respect to the direction of travel of $\gamma_\fpq$. This poset will have one function labeling the elements and another weighting the elements. In the following, we will define a function $wt^\star$ in terms of $\hat{x}^\star$. This will specialize to either $wt$ or $wt^\prin$ based on whether we specialize $\hat{x}^\star$ to $\hat{x}$ or $\hat{x}^\prin$.

\begin{definition}[Construction Algorithm]\label{def:ConstructPpq}
The poset $\calP_{\frac01}$ is $\emptyset$. Now, let $\fpq \in (0,\infty) \cap \mathbb{Q}$. To initialize, we consider the first intersection of $\gamma_\fpq$ with a line segment $\tau_a$ in $\mathcal{L}$. Let $x_a$ be the label of $\tau_a$. If $k_a = 0$, then we introduce one element $\calP_\fpq$ which will be $\calP_\fpq(1)$ with the chronological ordering. We label $\calP_\fpq(1) $ with $x_a$ and set $wt^\star(\calP_\fpq(1)) = (\xhat_a^\star)^2$. If $k_a \neq 0$, then we introduce two elements $\calP_\fpq(1)$ and $\calP_\fpq(2)$. If this intersection occurs strictly closer to the right endpoint of $\tau$ than the left (which is the case when $\fpq < 1$), then we set $\calP_\fpq(1) > \calP_\fpq(2)$; otherwise, we set $\calP_\fpq(1) < \calP_\fpq(2)$. In each case, we label both elements $x_a$, set the weight the smaller of the two elements with $k_a \xhat_a^\star$ and the larger of the two elements with $\frac{1}{k_a} \xhat_a^\star$.

The remainder of the poset is constructed with the following algorithm.

\begin{enumerate}\setlength{\leftskip}{-15pt}
    \item Let $\tau_a$ be the last line segment in $\mathcal{L}$ whose intersection with $\gamma_\fpq$ we have accounted for. If $\tau_a$ is the final line segment intersected by $\gamma_\fpq$, we terminate and return $\calP_\fpq$. Otherwise, let $\tau_b$ be the next line segment intersected by $\gamma_\fpq$. Let $x_a$ and $x_b$ be the labels of $\tau_a$ and $\tau_b$ respectively
and $\mathcal P_\frac{p}{q}(j)$ the last element labeled with $x_a$.
    \item Introduce an element $\calP_\fpq(j+1)$ which is labeled $x_b$.
    \item  Necessarily, $\tau_a$ and $\tau_b$ share an endpoint. If the endpoint shared by $\tau_a$ and $\tau_b$ lies to the right of $\gamma_\fpq$, then we set $\calP_\fpq(j) > \calP_\fpq(j+1)$; otherwise, we set $\calP_\fpq(j) < \calP_\fpq(j+1)$.
    \item If $k_b = 0$, we set the weight of $\calP_\fpq(j+1)$ to be $(\xhat_b^\star)^2$. Return to step 1.
    \item If $k_b \neq 0$, we introduce an element $\calP_\fpq(j+2)$ which is labeled $x_b$. If the intersection between $\gamma_\fpq$ and $\tau_b$ is strictly closer to the right endpoint of $\tau_b$, then we set $\calP_\fpq(j+1) > \calP_\fpq(j+2)$; otherwise, we set $\calP_\fpq(j+1) < \calP_\fpq(j+2)$. In each case, we weight the smaller of the two elements with $k_b \xhat_b^\star$ and the larger of the two elements with $\frac{1}{k_b} \xhat_b^\star$. Return to step 1.
\end{enumerate}
\end{definition}

We will introduce three other posets which we will frequently use.

\begin{definition}\label{def:H}
Define $H$ to be the poset resulting from applying the algorithm in Definition \ref{def:ConstructPpq} to the half circle which travels from $(-\epsilon,-\epsilon)$ to $(\epsilon,\epsilon)$ in the clockwise direction for small $\epsilon>0$.

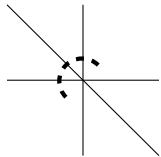
\begin{center}
\begin{tikzpicture}
 \draw[] (1,-1) -- (1,1);
 \draw[] (0,0) -- (2,0);
 \draw[] (0,1) -- (2,-1);
\draw[dashed,ultra thick] (1.2,0.2) arc(45:225:.3);
\end{tikzpicture}
\end{center}

We define $\Ptilde_{\frac01}$ to be the poset resulting from applying the algorithm from Definition \ref{def:ConstructPpq} to the line segment $(\epsilon,\epsilon) - (1+\epsilon,\epsilon)$.

Given $\fpq \in (0,\infty) \cap \mathbb{Q}$, define $\Ptilde_\fpq$ to be the poset resulting from applying the algorithm from Definition \ref{def:ConstructPpq} to the curve which travels along $\gamma_\fpq$ from $(0,0)$ to $(q-\epsilon,p-\epsilon)$, and then travels in a half circle in the clockwise direction from $(q-\epsilon,p-\epsilon)$ to $(q+\epsilon,p+\epsilon)$.

For all $\fpq \in [0,\infty) \cap \mathbb{Q}$, define $\calP^\circ_\fpq$ to be the poset with the same underlying set as $\Ptilde_\fpq$ and whose relations are the transitive closure of the relations in $\Ptilde_\fpq$ as well as $L > R$.
\end{definition}

Since we are at times closely following \cite{banaian2024skein}, we want to point out that the symbol $\Ptilde$ was used in a different way in this source.
We remark that $\Ptilde_\fpq$ could also be described as $\calP_\fpq \nearrow H$. Moreover, note that $\calP^\circ_\fpq$ is not a fence poset; its Hasse diagram is a cycle. Such posets are sometimes referred to as ``circular fence posets''.

It will be useful to have shorthand for the number of elements of $\calP_\fpq$ and $\Ptilde_\fpq$. We set $h_\fpq:= \vert \calP_\fpq \vert$ and $\widetilde{h}_\fpq := \vert \Ptilde_\fpq \vert = \vert \calP^\circ_\fpq \vert$. These quantities are functions of $p,q,k_1,k_2,$ and $k_3$, and we give explicit expressions in Proposition \ref{prop:Computedfpq}.

\begin{example}\label{ex:Construction}
Let $k_1 = 0, k_2 = 2, k_3 = 3$. In the following, we will draw the Hasse diagram of each poset with elements written as $(\ell,w)$ where $\ell$ is the label and $w$ is one of the two possible weights.

Below, on the left we draw $H$ and on the right we draw $\Ptilde_{\frac01}$, each with principal coefficient weights. The coefficient-free weights can be obtained simply by setting all $y_i$ to 1.

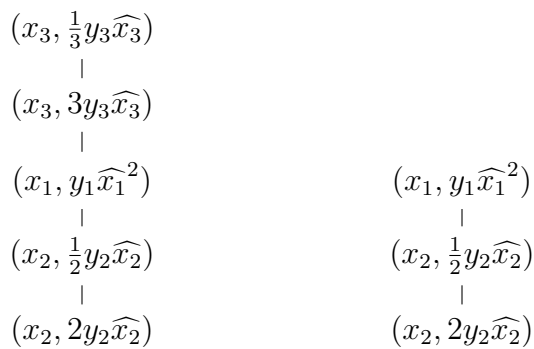
\begin{center}
\begin{tikzpicture}
\node(0) at (0,0){$(x_2, 2 y_2\widehat{x_2})$};
\node(1) at (0,1){$(x_2, \frac12 y_2 \widehat{x_2})$};
\node(2) at (0,2){$(x_1, y_1\widehat{x_1}^2)$};
\node(3) at (0,3){$(x_3, 3 y_3\widehat{x_3})$};
\node(4) at (0,4){$(x_3, \frac13 y_3\widehat{x_3})$};
\draw(0) -- (1);
\draw(1) -- (2);
\draw(2)--(3);
\draw(3) -- (4);
\node(5) at (5,0){$(x_2, 2 y_2\widehat{x_2})$};
\node(6) at (5,1){$(x_2, \frac12 y_2\widehat{x_2})$};
\node(7) at (05,2){$(x_1, y_1\widehat{x_1}^2)$};
\draw(5) -- (6);
\draw(6) -- (7);
\end{tikzpicture}
\end{center}

As a larger example of the algorithm, we construct $\calP_{\frac23}$. Consider $\gamma_{\frac23}$ as drawn below.

\begin{center}
\begin{tikzpicture}
\draw(0,0) -- (3,0) -- (3,2) -- (0,2) -- (0,0);
\draw(1,0) -- (0,1) -- (3,1) --(2,2);
\draw(0,2) -- (2,0) -- (2,2);
\draw(3,0) -- (1,2) -- (1,0);
\draw[red,ultra thick] (0,0) -- (3,2);
\end{tikzpicture}
\end{center}

Then, $\calP_{\frac23}$ is the following labeled, weighted poset. Here, we only include coefficient-free weights for sake of neatness.

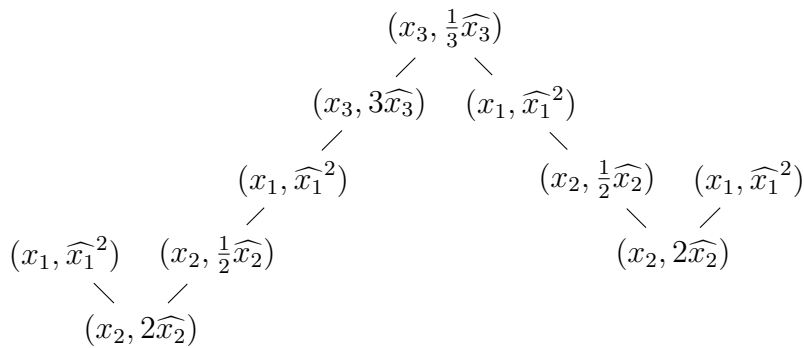
\begin{center}
\begin{tikzpicture}
\node(0) at (0,0){$(x_1,\widehat{x_1}^2)$};
\node(1) at (1,-1){$(x_2,2\widehat{x_2})$};
\node(2) at (2,0){$(x_2,\frac12\widehat{x_2})$};
\node(3) at (3,1){$(x_1,\widehat{x_1}^2)$};
\node(4) at (4,2){$(x_3,3\widehat{x_3})$};
\node(5) at (5,3){$(x_3,\frac13\widehat{x_3})$};
\node(6) at (6,2){$(x_1,\widehat{x_1}^2)$};
\node(7) at (7,1){$(x_2,\frac12\widehat{x_2})$};
\node(8) at (8,0){$(x_2,2\widehat{x_2})$};
\node(9) at (9,1){$(x_1,\widehat{x_1}^2)$};
\draw(0) -- (1);
\draw(1) -- (2);
\draw(2) -- (3);
\draw(3) -- (4);
\draw(4) -- (5);
\draw(5) -- (6);
\draw(6) -- (7);
\draw(7)--(8);
\draw(8) -- (9);
\end{tikzpicture}
\end{center}

Using the alternate definition of $\Ptilde_\fpq$ as $\calP_\fpq \nearrow H$, we can immediately also present $\Ptilde_{\frac23}$.

\begin{center}
\begin{tikzpicture}
\node(0) at (0,0){$\boxed{(x_1,\widehat{x_1}^2)}$};
\node(1) at (1,-1){$(x_2,2\widehat{x_2})$};
\node(2) at (2,0){$(x_2,\frac12\widehat{x_2})$};
\node(3) at (3,1){$(x_1,\widehat{x_1}^2)$};
\node(4) at (4,2){$(x_3,3\widehat{x_3})$};
\node(5) at (5,3){$(x_3,\frac13\widehat{x_3})$};
\node(6) at (6,2){$(x_1,\widehat{x_1}^2)$};
\node(7) at (7,1){$(x_2,2\widehat{x_2})$};
\node(8) at (8,0){$(x_2,\frac12\widehat{x_2})$};
\node(9) at (9,1){$(x_1,\widehat{x_1}^2)$};
\draw(0) -- (1);
\draw(1) -- (2);
\draw(2) -- (3);
\draw(3) -- (4);
\draw(4) -- (5);
\draw(5) -- (6);
\draw(6) -- (7);
\draw(7)--(8);
\draw(8) -- (9);
\node(10) at (11,-1.5){$\boxed{(x_2, 2 \widehat{x_2})}$};
\node(11) at (11,-0.5){$(x_2, \frac12 \widehat{x_2})$};
\node(12) at (11,0.5){$(x_1, \widehat{x_1}^2)$};
\node(13) at (11,1.5){$(x_3, 3 \widehat{x_3})$};
\node(14) at (11,2.5){$(x_3, \frac13 \widehat{x_3})$};
\draw(10) -- (11);
\draw(11) -- (12);
\draw(12)--(13);
\draw(13) -- (14);
\draw(9) --(14);
\end{tikzpicture}
\end{center}

We boxed $L_{\Ptilde_{\frac23}}$ and $R_{\Ptilde_{\frac23}}$. These are, respectively, the boxed element labeled $x_1$ and the boxed element labeled $x_2$. Then, the poset $\calP^\circ_\fpq$ has the same underlying set as $\Ptilde_{\frac23}$ and has all the relations from $\Ptilde_{\frac23}$  as well as $L_{\Ptilde_{\frac23}}> R_{\Ptilde_{\frac23}}$.
\end{example}

\begin{remark}\label{rem:Christoffel}
In the case where $k_1=k_2=k_3=0$ and $\frac{p}{q}\in (0,1]$, we can recover the \emph{Christoffel word} $\mathrm{ch}_{\fpq}$ in the sense of \cite{aigner2013Markov} from $\Ptilde_{\frac{p}{q}}$. Choose the chronological ordering of $\Ptilde_\fpq$ which corresponds with reading the intersections of $\gamma_\fpq$ from bottom-left to top-right. It is straightforward to see that all maximal chains in $\Ptilde_\fpq$ are length 1 or 2. Consider all chains which are decreasing with respect to the chronological ordering. Consider the word in letters $\{A,B\}$ given by writing $A$ for each length 1 maximal decreasing chain and $B$ for each length 2 maximal decreasing chain, using the chronological ordering to order the letters. This word is the Christoffel word $\mathrm{ch}_{\fpq}$. For instance, by constructing $\Ptilde_{\frac23}$ with $k_1 = k_2 = k_3 = 0$, we can calculate $\mathrm{ch}_{\frac23} = ABB$. This construction is given with more detail in \cite{ouguz2025oriented,rabideau2020continued}.
\end{remark}

\subsection{Expansion Formula from a Poset}

In this section, we will describe a way to associate a Laurent polynomial to each poset $\calP_\fpq$.

We recall the construction of a vector $\mathbf{g}_\calP$ to any labeled, weighted poset $\calP$ as in \cite{banaian2024skein}. Let $a_i$ be the number of  elements labeled $x_i$ which are minimal in $\calP$, and let $\mathbf{a}_\calP = (a_1,\ldots,a_n)^\top$. Define an element to be \emph{strictly maximal} element if it is maximal and covers at least two elements. Let $b_i$ be the number of elements labeled $x_i$ which are strictly maximal in $\calP$, and let $\mathbf{b}_\calP = (b_1,\ldots,b_n)^\top$.

In \cite{banaian2024skein}, there is a definition of a third vector $\mathbf{r}$. We give a simplified definition here since we work in a specific setting.  We will refer to an element of $\calP$ as a \emph{leaf} if it is a leaf in the Hasse diagram. Equivalently, an element is a leaf if it only has a cover relation with one other element of the poset.
Given a poset $\calP$, let $r_1$ be the number of leaves labeled $x_2$, $r_2$ be the number of leaves labeled $x_3$ and $r_3$ the number of leaves labeled $x_1$. Let $\mathbf{r}_\calP = (r_1,r_2,r_3)^\top$.

We define $\gb_\calP = -\mathbf{a}_\calP + \mathbf{b}_\calP + \mathbf{r}_\calP$. For illustrations, see Example \ref{ex:GVecFPoly}.  Next, we show that when $\calP = \calP_\fpq$, then this vector $\gb_{\calP_\fpq}$ is equal to the $g$-vector of the corresponding cluster variable. We prove this by analyzing how maximal and minimal elements can appear in $\calP_\fpq$.

\begin{lemma}\label{lem:CompareGVector}
For any $\fpq \in (0,\infty)$, we have $\gb_{\calP_\fpq} = \gb_{1,\fpq}$.
\end{lemma}

\begin{proof}
 From our conventions, every element labeled $x_2$ is either minimal or neither minimal nor maximal. We have exactly one such minimal element for each crossing of $\gamma_\fpq$ with a line of slope $\infty$. Therefore, we have $q-1$ minimal elements, each labeled $x_2$.
 Similarly, every element labeled $x_3$ is either a strict maximal element or neither minimal nor maximal, and overall there are $p-1$ strict maximal elements labeled $x_3$.

 The number of strict maximal elements in a fence poset is one less than the number of minimal elements. The missing minimal or strict maximal elements will be labeled $x_1$. We have two subcases, based on whether $p \leq q$, i.e. $\fpq \leq 1$, or if $p > q$. In the former case, there are $(q-1) - (p-1) - 1 = q-p-1$ strict maximal elements labeled $x_1$ and no minimal elements labeled $x_1$. The latter is true because $\gamma_{\fpq}$ with $\fpq \leq 1$ crosses a line segment of slope $\infty$ before or after each crossing with a line segment of slope $-1$, guaranteeing that an element labeled $x_1$ will never be minimal.  Here, $\mathbf{a}_{\calP_\fpq} = (0,q-1,0)^\top$ and $\mathbf{b}_{\calP_\fpq} = (q-p-1,0,p-1)^\top$.

 If $p > q$, then there are $p-q+1$ minimal elements labeled $x_1$ and no maximal elements labeled $x_1$. In this case, $\mathbf{a}_{\calP_\fpq} = (p-q+1,q-1,0)^\top$ and $\mathbf{b}_{\calP_\fpq} = (0,0,p-1)^\top$.

 Finally, for any $\fpq \in [0,\infty)$, the first and last line segment in $\mathcal{L}$ crossed by $\gamma_\fpq$ is labeled $x_1$, which implies $\mathbf{r}_{\calP_\fpq} = (0,0,2)$. Now the result follows from comparing $\gb_{\calP_\fpq}$ with the expression for $\gb_{1,\fpq}$ in Theorem \ref{thm:g-vector-description}.
\end{proof}

We exclude $\fpq = \frac01$ from the previous since $\mathcal P_{\frac01} = \emptyset$. So, as a convention, we set $\gg_{\mathcal P_{\frac01}} = (0,0,1)$, which again aligns with $\gb_{1,\frac01}$.

The analysis in the previous proof naturally lends itself to enumerating the elements in $\calP_\fpq$. Recall that when $k_i = 0$, this means we have $d_i = 1$ in our associated cluster algebra and when $k_i > 0$, we have $d_i = 2$. As was the case for Lemma \ref{lem:CompareGVector}, this proof relies on a careful analysis of $\calP_\fpq$.

\begin{proposition}\label{prop:Computedfpq}
Given $\fpq \in (0,\infty) \cap \mathbb{Q}$, the number of elements of $\calP_\fpq$ labeled $x_1$ is $d_1(q+p-1)$, the number of elements labeled $x_2$ is $d_2(q-1)$, and the number of labeled $x_3$ is $d_3(p-1)$. In particular, \[
h_\fpq = d_1(q+p-1) + d_2(q-1) + d_3(p-1)
\]
and
\[
\tilde{h}_\fpq = d_1(q+p) + d_2q + d_3p.
\]
\end{proposition}

\begin{proof}
The line segment $\gamma_\fpq$ crosses $q-1$ lines of slope $\infty$ and $p-1$ lines of slope $0$. Each intersection contributes $d_2$ and $d_3$ elements to $\calP_\fpq$ respectively. Among all intersections of  $\gamma_\fpq$ with $\mathcal{L}$, every other intersection is with a line of slope $-1$. These are also the first and last intersections. Therefore, $\gamma_\fpq$ crosses $(p-1) + (q-1) + 1 = p+q-1$ lines of slope $\frac{-1}{1}$, and each intersection contributes $d_1$ elements to $\calP_\fpq$. This verifies our formula for $h_\fpq$.

Since $\Ptilde_\fpq = \calP_\fpq \nearrow H$, the cardinalities are related by $\tilde{h}_\fpq = h_\fpq + \vert H \vert$. We can calculate $\vert H \vert = d_1 + d_2 + d_3$, which implies our expression for $\tilde{h}_\fpq$ is correct.
\end{proof}

For example, in Example \ref{ex:Construction} we have $d_1 = 1, d_2 = d_3 = 2$. We found that $\calP_{\frac23}$ has 10 elements, where $10 = (2+3-1) + 2(1) + 2(2)$. Moreover, $H$ has 5 elements with these values of $k_i$, so that $\Ptilde_{\frac23}$ has $15 = (2+3) + 2(2) + 2(3)$ elements.

 In what follows, we will refer to $\gb_{1,\fpq}$ simply as $\gb_\fpq$ as we will focus on this branch of $\mathbb{T}(k_1,k_2,k_3)$. See Remark \ref{rem:OtherBranches} for dicussion of the other branches. We will also abbreviate $\gb_{\calP^\circ_{\fpq}}$ as $\gb^\circ_\fpq$. The following identity will be useful in later computations.

\begin{lemma}\label{lem:CompareGVecAndGCirc}
The vectors $\gb_\fpq$ and $\gb^\circ_\fpq$ are related by \[
\gb_\fpq^\circ = \gb_\fpq + \begin{bmatrix}1\\-1\\-1\end{bmatrix}.
\]
\end{lemma}

\begin{proof}
Recall the poset $\Ptilde_\fpq$ is related to $\calP_\fpq$ by adding the chain $H$, and then $\calP^\circ_\fpq$ is the result of imposing a relation to the leaves of $\Ptilde_\fpq$. By adding a chain, we include a new minimal element labeled $x_2$ and a new maximal element labeled $x_3$. Moreover, there is an element labeled $x_1$ which is a leaf in $\calP_\fpq$ but is a strictly maximal element in $\calP^\circ_\fpq$. Since all maximal elements in $\calP^\circ_\fpq$ are strictly maximal, the conclusion follows.
\end{proof}

\begin{remark}
 While we set $\gg_{\frac01} = (0,0,1)^\top$ as a convention, we can follow our definition of $g$-vectors and compute $\gg^\circ_{\frac01} = (1,-1,0)^\top$ which indeed is equal to $\gg_{\frac01} + (1,-1,-1)^\top$.
\end{remark}

Notice in particular that the coordinates of all vectors $\gb_\fpq$ sum to 1 and the coordinates of all vectors $\gb_\fpq^\circ$ sum to 0. For an example, see Example \ref{ex:GVecFPoly}.

Recall that an \emph{order ideal} of a poset $(\calP,\leq)$ is $I \subseteq \calP$ such that for all $y \in I$ and $x \in \calP$, if $x \leq y$, then $x \in I$. We let $J(\calP)$ denote the set of all order ideals of a poset $\calP$.

\begin{definition}
 Given $I \subseteq \calP$, let $wt(I) = \prod_{x \in I} wt(x)$, and set $\mathcal{W}(\calP) = \sum_{I \in J(\calP)} wt(I)$.

Similarly, let $wt^\prin(I) = \prod_{x \in I} wt^\prin(x)$ and let $\mathcal{W}^\prin(\calP) = \sum_{I \in J(\calP)} wt^\prin(I)$.

 We also set $\mathcal{W}(\calP;\neg L) = \sum_{\substack{I \in J(\calP)\\ L \notin I}} wt(I)$, $\mathcal{W}(\calP;L) = \sum_{\substack{I \in J(\calP)\\ L \in I}} wt(I)$. We define other symbols, such as $\mathcal{W}(\calP;\neg R)$ and $\mathcal{W}^\prin(\calP; L,R)$, similarly.
\end{definition}

\begin{definition}
Given a labeled, weighted poset, let $X_{\calP}:= \xb^{\gb_\calP} \mathcal{W}(\calP)$ and $X^\prin_{\calP} = \xb^{\gb_\calP} \mathcal{W}^\prin(\calP)$. For shorthand, set $X_\fpq := X_{\calP_\fpq}$, $X^\circ_\fpq := X_{\calP^\circ_\fpq}$, and similarly for principal coefficients.
\end{definition}

In Theorem \ref{thm:CorrectnessOfPosetFormula}, we will show that $X^\prin_{\calP_\fpq}$ is in fact a cluster variable in $\mathcal{A}^{\prin}(k_1,k_2,k_3)$, implying the same is true in the coefficient-free case. This proof requires several preparatory results, which will be presented in the next few sections.

\begin{example}\label{ex:GVecFPoly}
Recall we computed $\calP_{\frac23},\Ptilde_{\frac23}$ and $\calP^\circ_{\frac23}$ in Example \ref{ex:Construction}. We compute the $g$-vector of each.
\[\mathbf{a}_{\calP_{\frac23}} = \begin{bmatrix}0\\2\\0\end{bmatrix},\ \mathbf{b}_{\calP_{\frac23}} = \begin{bmatrix}0\\0\\1\end{bmatrix},\ \mathbf{r}_{\calP_{\frac23}} = \begin{bmatrix}0\\0\\2\end{bmatrix}\implies \mathbf{g}_{\calP_{\frac23}} = \begin{bmatrix}0\\-2\\3\end{bmatrix}\]\vspace{1pt}
\[\mathbf{a}_{\Ptilde_{\frac23}} = \begin{bmatrix}0\\3\\0\end{bmatrix},\
\mathbf{b}_{\Ptilde_{\frac23}} = \begin{bmatrix}0\\0\\2\end{bmatrix},\
\mathbf{r}_{\Ptilde_{\frac23}} = \begin{bmatrix}1\\0\\1\end{bmatrix} \implies \mathbf{g}_{\Ptilde_{\frac23}} = \begin{bmatrix}1\\-3\\3\end{bmatrix}\]\vspace{1pt}
\[\mathbf{a}_{\calP^\circ_{\frac23}} = \begin{bmatrix}0\\3\\0\end{bmatrix},\ \mathbf{b}_{\calP^\circ_{\frac23}} =\begin{bmatrix}1\\0\\2\end{bmatrix},\  \mathbf{r}_{\calP^\circ_{\frac23}} = \begin{bmatrix}0\\0\\0\end{bmatrix} \implies \mathbf{g}_{\calP^\circ_{\frac23}} = \begin{bmatrix}1\\-3\\2\end{bmatrix}\]

The Laurent polynomial $X_{\frac23}^\prin$ is given by \[
X^\prin_\frac23 = \frac{x_3^3}{x_2^2}\big(1 + (2 + 2)y_2\widehat{x_2} + (1+1+4)y_2^2\widehat{x_2}^2 + (2 + 2)y_1y_2\widehat{x_1}^2\widehat{x_2}   + \cdots \big)
\]
and $X_{\frac23}$ is given by setting $y_i = 1$.
\end{example}

\section{Markov Posets and Cluster Variables}\label{sec:MarkovPosets}

From this section, we will sometimes refer to $x_{1,\fpq}$ as simply $x_{\fpq}$.
Here, we will provide several results on the structure of the poset $\calP_\fpq$ and $\Ptilde_\fpq$ and relations on the Laurent polynomials $\mathcal{W}(\calP_\fpq)$ and $X_\fpq$. In particular, we will be able to identify each $X_\fpq$ with a cluster variable $x_\fpq$.  These results will be useful in exhibiting a specific class of cluster Cohn matrices.
\subsection{Skein Relations}

The following definitions are inspired by \cite{canakci2013snake}. In this section, we consider general labeled, weighted fence posets.
We say two posets, $\calP$ and $\calP'$ have an \emph{overlap} if there exists $1 \leq s \leq t \leq \vert \calP \vert$ and $1 \leq s' \leq t' \leq \vert \calP' \vert$ such that  $\calP[s,t]$ and $\calP'[s',t']$ are isomorphic as labeled posets with isomorphism $\Omega$ sending $\calP(s+i) \mapsto \calP'(s'+i)$. We emphasize that the subposet $\calP[s,t]$ (or equivalently $\calP'[s',t']$) will be referred to as the overlap.

Say a set $S \subseteq \calP$ is \emph{on top} of $\calP$ if there does not exist $x \in S, y \in P\backslash S$ such that $x < y$, and define a set being \emph{on bottom} of $\calP$ similarly.

\begin{definition}\label{def:Overlap}
Let $\calP_1$ and $\calP_2$ be a pair of labeled, weighted posets. An overlap $R:=\calP_1[s,t] \cong \calP_2[s',t]'$ is \emph{crossing} if
\begin{enumerate}\setlength{\leftskip}{-15pt}
 \item $R$ is on top of $\calP_1$ and on bottom of $\calP_2$;
 \item we do not have $s = s' = 1$;
 \item we do not have both $t = \vert \calP_1 \vert$ and $t' = \vert \calP_2 \vert$; and
 \item $\calP_1[s,t]$ and $\calP_2[s',t']$ are isomorphic as \emph{labeled, weighted} posets.
 \end{enumerate}
\end{definition}

Given an arc $\gamma$ on a surface with triangulation, there is a general method to construct a labeled, weighted poset $\calP_\gamma$ which encodes the Laurent expansion of the cluster variable associated to $\gamma$ in the seed labeled by $T$ \cite{ezgieminecluster2024, pilaud2023posets}. The inspiration for the definition of a crossing overlap is that, if posets $\calP_\gamma$ and $\calP_{\gamma'}$ associated to two arcs on a surface have a crossing overlap, then the arcs $\gamma$ and $\gamma'$ intersect.

\begin{definition}\label{Def:resolveOverlap}
Let $\calP_1$ and $\calP_2$ be two labeled, weighted posets which have a crossing overlap $\calP_1[s,t] \cong \calP_2[s',t']$. Let $h_1 = \vert \calP_1 \vert$ and $h_2 = \vert \calP_2 \vert$. Given a poset $\calP$ with chronological ordering, we let $\overline{\calP}$ be the same poset with the opposite chronological ordering. We define the \emph{resolution} of a pair of posets with a crossing overlap to be the pair of sets $\{\calP_3 \cup \calP_4\}, \{\calP_5 \cup \calP_6\}$ where the posets are defined as follows.
\begin{itemize}\setlength{\leftskip}{-15pt}
    \item Let $\calP_3 = \calP_1[1,t] \nearrow \calP_2[t'+1,h_2]$.
    \item Let $\calP_4 = \calP_2[1,t'] \searrow \calP_1[t+1,h_1]$.
    \item We have three subcases for $\calP_5$.
    \begin{itemize}\setlength{\leftskip}{-15pt}
        \item If $s > 1$ and $s' > 1$, define $\calP_5$ to be $\calP_1[1,s-1] \searrow \overline{\calP_2[s'-1,1]}$.
        \item If $s = 1$, implying $s' > 1$, let $u'<s'-1$ be the largest integer such that $\calP_2(u') \nless \calP_2(s'-1)$, if it exists,  and otherwise let $u' = 0$.  Let $\calP_5 = \calP_2[1,u']$.
        \item If $s' = 1$, implying $s > 1$, let $u < s-1$ be the largest integer such that $\calP_1(u) \ngtr \calP_1(s-1)$, if it exists,  and otherwise $u = 0$.  Let $\calP_5 = \calP_1[1,u]$.
    \end{itemize}
    \item We have three subcases for $\calP_6$.
    \begin{itemize}\setlength{\leftskip}{-15pt}
    \item If $t < h_1$ and $t' < h_2$,  define $\calP_6$ to be $\overline{\calP_1[t+1,h_1]} \searrow \calP_2[t'+1,h_2]$.
    \item If $t = h_1$, implying $t'<h_2$,  let $v'>t'+1$ be the smallest integer such that $\calP_2(v')\nless \calP_2(t+1)$, if it exists,  and otherwise set $v' = h_2+1$. Let  $\calP_6 = \calP_2[v',h_2]$.
    \item If $t'=h_2$,  implying $t < h_1$,  let $v>t+1$ be the smallest integer such that $\calP_1(v) \ngtr \calP(t+1)$, if it exists,  and otherwise set $v = h_1+1$. Let $\calP_6 = \calP_1[v,h_1]$.
   	\end{itemize}
\end{itemize}
\end{definition}

To be precise, we define $\calP_i[1,0]=\calP_i[h_i+1,h_i]=\emptyset$, $\gb_\emptyset = (0,0,0)^\top$, and  $\mathcal{W}(\emptyset) = 1$. We will want to apply the following results to both weighting functions we use on $\calP_\fpq$. To emphasize these hold for either weight, as in Definition \ref{def:ConstructPpq}, we will let $wt^\star$ denote an arbitrary weighting function and $\mathcal{W}^\star(\calP)$ to denote the Laurent polynomial of the poset $\calP$ with respect to $wt^\star$.

\begin{example}\label{ex:CrossingOverlap}
Consider $\calP_1$ and $\calP_2$ drawn below on the left and right respectively. We denote each element by its label.

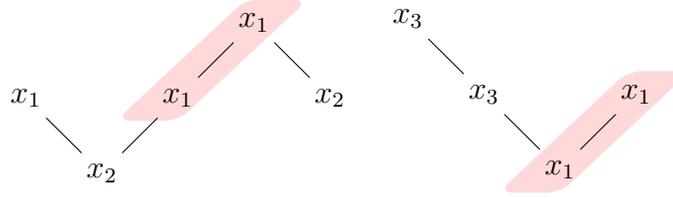
\begin{center}
\begin{tabular}{cc}
\begin{tikzpicture}
\filldraw [rounded corners=5pt, red!15]
(1.2,-.3)
-- (2.9,1.3)
-- (3.7,1.3)
-- (2,-.3)
-- cycle;
\node(0) at (0,0){$x_1$};
\node(1) at (1,-1){$x_2$};
\node(2) at (2,0){$x_1$};
\node(3) at (3,1){$x_1$};
\node(4) at (4,0){$x_2$};
\draw(0) -- (1);
\draw(1) -- (2);
\draw(2) -- (3);
\draw(3) -- (4);
\end{tikzpicture}&
\begin{tikzpicture}
\filldraw [rounded corners=5pt, red!15]
(1.2,-2.3)
-- (2.9,-.7)
-- (3.7,-.7)
-- (2,-2.3)
-- cycle;
\node(0) at (0,0){$x_3$};
\node(1) at (1,-1){$x_3$};
\node(2) at (2,-2){$x_1$};
\node(3) at (3,-1){$x_1$};
\draw(0) -- (1);
\draw(1) -- (2);
\draw(2) -- (3);
\end{tikzpicture}
\end{tabular}
\end{center}

For each poset, we use the chronological ordering induced by reading the Hasse diagram left to right.  According to Definition \ref{def:Overlap}, if $wt^*(\calP_1(3)) = wt^*(\calP_2(3))$ and $wt^*(\calP_1(4)) = wt^*(\calP_2(4))$, then $\calP_1$ and $\calP_2$ have a crossing overlap in $\calP_1[3,4] \cong \calP_2[3,4]$. We can use Definition \ref{Def:resolveOverlap} to construct a resolution. First, we draw $\calP_3$ and $\calP_4$ on the left and right respectively.

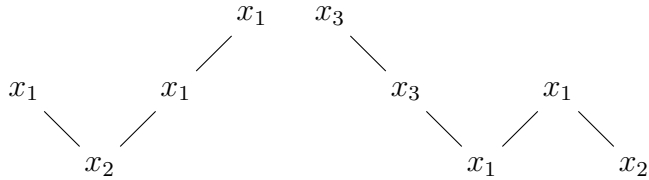
\begin{center}
\begin{tabular}{cc}
\begin{tikzpicture}
\node(0) at (0,0){$x_1$};
\node(1) at (1,-1){$x_2$};
\node(2) at (2,0){$x_1$};
\node(3) at (3,1){$x_1$};
\draw(0) -- (1);
\draw(1) -- (2);
\draw(2) -- (3);
\end{tikzpicture}&
\begin{tikzpicture}
\node(0) at (0,0){$x_3$};
\node(1) at (1,-1){$x_3$};
\node(2) at (2,-2){$x_1$};
\node(3) at (3,-1){$x_1$};
\node(4) at (4,-2){$x_2$};
\draw(0) -- (1);
\draw(1) -- (2);
\draw(2) -- (3);
\draw(3) -- (4);
\end{tikzpicture}
\end{tabular}
\end{center}

Since $s = s' = 3 > 1$, there is a nonempty poset $\calP_5$, which we draw below. However, $t' = 4 = h_2$ and there is no $v > t+1 = 5$ with $\calP_1(v) \ngtr \calP_1(5)$, so we have $\calP_6 = \emptyset$.

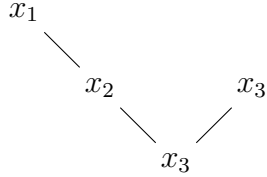
\begin{center}
\begin{tikzpicture}
\node(0) at (0,0){$x_1$};
\node(1) at (1,-1){$x_2$};
\node(2) at (2,-2){$x_3$};
\node(3) at (3,-1){$x_3$};
\draw(0) -- (1);
\draw(1) -- (2);
\draw(2) -- (3);
\end{tikzpicture}
\end{center}
\end{example}

Our first result comes directly from \cite{banaian2024skein,canakci2013snake}.

\begin{proposition}\label{prop:CrossingOverlap}
Given $\calP_1$ and $\calP_2$ with a crossing overlap in $\calP_1[s,t] \cong \calP_2[s',t']$,  define $Z^\star_R,Z^\star_s$ and $Z^\star_t$ by \[
Z^\star_R = \prod_{i = s}^t wt^\star(\calP_1(i)),
\]
\[
Z^\star_s = \begin{cases} \prod_{i=u'+1}^{s'-1} wt^\star(\calP_2(i)) & s = 1\\
1 & s>1\end{cases},
\qquad
Z^\star_t= \begin{cases} \prod_{i=t'+1}^{v'-1} wt^\star(\calP_2(i)) & t = d\\
1 & t<d\end{cases}.
\]
We have the following relation amongst the weight polynomials of $\calP_1,\calP_2$ and their resolution,\[
\mathcal{W}^\star(\calP_1)\mathcal{W}^\star(\calP_2) = \mathcal{W}^\star(\calP_3)\mathcal{W}^\star(\calP_4) +
Z^\star_RZ^\star_sZ^\star_t\mathcal{W}^\star(\calP_5)\mathcal{W}^\star(\calP_6).
\]
\end{proposition}

\begin{proof}
The proof of \cite[Proposition 6]{banaian2024skein} exhibited a partition $A \sqcup B$ of $J(\calP_1) \times J(\calP_2)$ and bijections $A \cong J(\calP_3) \times J(\calP_4)$ and $B \cong J(\calP_5) \times J(\calP_6)$. The first bijection was weight-preserving while the second bijection changed the weight by exactly the monomial $Z^\star_RZ^\star_sZ^\star_t$.
\end{proof}

\begin{example}
Applying Proposition \ref{prop:CrossingOverlap} to Example \ref{ex:CrossingOverlap} gives the following relation,\[
\mathcal{W}^\star(\calP_1)\mathcal{W}^\star(\calP_2) = \mathcal{W}^\star(\calP_3)\mathcal{W}^\star(\calP_4) + wt^\star(\calP_1(3))wt^\star(\calP_1(4))\mathcal{W}^\star(\calP_5).
\]
Using the bijection between order ideals and perfect matchings of a snake graph (see \cite{ezgieminecluster2024, musiker2013bases}) and the theory of snake graphs and continued fractions (see \cite{canakci2020snake}), one can compute the number of order ideals of each poset,
\begin{align*}
\vert J(\calP_1) \vert = 12, \quad  \vert J(\calP_2) \vert = 7,\quad  \vert J(\calP_3) \vert = 7,\quad \vert J(\calP_4)\vert = 11,\quad \vert J(\calP_5) \vert = 7,
\end{align*}
which supports the validity of this relation.
\end{example}

It is possible that two arcs $\gamma_1,\gamma_2$ intersect but the corresponding posets $\calP_{\gamma_1}$ and $\calP_{\gamma_2}$ do not have a crossing overlap.  We recall one type of operation which would be the result of resolving such an intersection. This result can be seen as coming from snake graph calculus \cite{canakci2013snake}. We include a proof for completeness.

\begin{proposition}\label{prop:SkeinType2}
Given $\calP_1,\calP_2$ with $h_i:= \vert \calP_i \vert$,  define $\calP_{34} = \calP_1 \nearrow \calP_2$.  Let $u$ be the largest integer such that $\calP_1(u) \ngtr \calP_1(h_1)$
if it exists, and otherwise let $u = 0$.  Similarly,  let $v$ be the smallest integer such that $\calP_2(v) \nless \calP_2(1)$ if it exists and otherwise let $v = h_2+1$.  Let $\calP_{5} = \calP_1[1,u]$ and let $\calP_6 = \calP_2[v,h_2]$.
Let $Z^\star = \prod_{i=1}^{v-1} wt^\star(\calP_2(i))$.  Then,  we have the following relation amongst the weight polynomials of these posets, \[
\mathcal{W}^\star(\calP_1) \mathcal{W}^\star(\calP_2) = \mathcal{W}^\star(\calP_{34}) + Z^\star \mathcal{W}^\star(\calP_5)\mathcal{W}^\star(\calP_6).
\]
\end{proposition}

\begin{proof}
We will follow the proof method which was described in \cite[Section 3.4]{banaian2024skein}. We partition $J(\calP_1) \times J(\calP_2)$ into $A$ and $B$ where $A$ consists of the subset of all pairs $(I_1,I_2)$ such that $\calP_2(1) \in I_2$ only if $\calP_1(h_1) \in I_1$ and $B$ consists of the complement, namely,  the set of all $(I_1,I_2)$ such that $\calP_2(1) \in I_2$ and $\calP_1(h_1) \notin I_1$.

Any pair $(I_1,I_2) \in A$ can be naturally mapped to an order ideal of $\calP_{34}$ in a manner such that the weight of the new order ideal is equal to the product of weights of $I_1$ and $I_2$. This establishes a weight-preserving bijection between $A$ and $J(\calP_{34})$.

Given $(I_1,I_2) \in B$,  we know $I_1 \in J(\calP_5)$ and $I_2 \backslash \langle \calP_2(1) \rangle \in J(\calP_6)$.  This gives a bijection between $B$ and $J(\calP_5) \times J(\calP_6)$ which changes the weight of the image exactly by $\prod_{x \in \langle \calP_2(1)\rangle} wt^\star(x) = \prod_{i=1}^{v-1} wt^\star(\calP_2(i))$.
\end{proof}

\begin{example}\label{ex:SkeinType2}
Recall the posets provided in Example \ref{ex:Construction} in the setting $k_1 = 0, k_2 = 2, k_3 = 3$. If we choose the chronological ordering on $H$ which reads the poset from top to bottom as drawn, then the poset $\Ptilde_{\frac23}$ is equal to $\calP_{\frac23} \nearrow H$. Following the statement of Proposition \ref{prop:SkeinType2}, $u = 9$, $v = 6$, $\calP_6 = \emptyset$, and $\calP_5$ is the result of removing the rightmost element from $\calP_{\frac23}$. Setting $\calP_1 = \calP_{\frac23}$ and $\calP_2 = H$ in Proposition \ref{prop:SkeinType2} yields  \[
\mathcal{W}^\star(\calP_{\frac23})\mathcal{W}^\star(H) = \mathcal{W}^\star(\Ptilde_{\frac23}) + Z^\star\mathcal{W}^\star(\calP_5),
\]
where if $wt^\star = wt$, the monomial $Z^\star$ is $1$ and if $wt^\star = wt^\prin$,  $Z^\star = y_1y_2^2y_3^3$.
\begin{center}

\end{center}
\end{example}

We remark that Propositions \ref{prop:CrossingOverlap} and \ref{prop:SkeinType2} can also be determined by reinterpreting the work in \cite{canakci2013snake} in terms of fence posets.
We will need one more type of relation which does not appear when considering posets from a pair of intersecting arcs on a surface.  This will involve a \emph{self-overlap},  i.e. an overlap where both subposets are contained in the same poset.

\begin{definition}\label{def:ReverseOverlap}
Let $\calP$ be a labeled, weighted poset such that there exist values $1 \leq r \leq s \leq t \leq h:= \vert \calP \vert$ such that $R:=\calP[r,s]$ and  $R':=\overline{\calP[s+1,t]}$ are overlaps. We say this overlap is a \emph{reverse, kissing self-overlap} if
\begin{enumerate}\setlength{\leftskip}{-15pt}
\item $R$ is on top of $\calP$ and $R'$ is on bottom of $\calP$;
\item we do not have both $r = 1$ and $t = h$;
\item $R$ and $R'$ are isomorphic as labeled posets; and
\item the isomorphism between $R$ and $R'$ is weight-preserving when restricted to $\calP[r,s-1]$ and $\calP[s+2,t]$.
\end{enumerate}
\end{definition}

\begin{definition}\label{def:ResolutionReverseSelfOverlap}
Let $\calP$ be a labeled, weighted poset on $h$ elements with a reverse, kissing self-overlap $\calP[r,s] \cong \overline{\calP[s+1,t]}$.  We define the \emph{resolution} of $\calP$ with respect to this overlap to be $\{\calP_{34}\}$ and $\{\calP_{56}\}$ where the posets are defined as follows.

\begin{itemize}\setlength{\leftskip}{-15pt}
    \item Let $\calP_{34} = \calP[1,r-1] \nearrow \overline{\calP[t,r]} \nearrow \calP[t+1,h]$.
    \item We have three subcases for $\calP_{56}$.
    \begin{itemize}\setlength{\leftskip}{-15pt}
        \item If $r > 1$ and $t < h$, define $\calP_{56}$ to be  $\calP[1,s-1] \searrow \calP[t+1,h]$.
        \item If $r = 1$, implying $t < d$, let $v>t+1$ be the smallest integer such that $\calP(v) \nless \calP(t+1)$, if it exists, and otherwise let $v = h+1$. Let $\calP_{56} = \calP[v,h]$.
        \item If $t = h$, implying $r > 1$, let $u < s-1$ be the largest integer such that $\calP(u) \ngtr \calP(s-1)$, if it exists, and otherwise let $u = 0$. Let $\calP_{56} = \calP[1,u]$.
    \end{itemize}
\end{itemize}
\end{definition}

To prove our relation on posets with reverse, kissing overlaps, we will  use the idea of a \emph{switching-position},  which was introduced in \cite{canakci2013snake} in terms of perfect matchings and was translated to the poset setting in \cite{banaian2024skein}.  We recall the definition and existence of switching positions here.

\begin{definition}
Let $\calP_1,\calP_2$ be two weighted posets with overlaps $\calP_1[s,t] \cong \calP_2[s',t']$,  and let $I_1 \in J(\calP_1)$ and $I_2 \in J(\calP_2)$. Suppose there exists a value $0 \leq i \leq t-s$ such that $\calP_1(s+i) \in I_1$ if and only if $\calP_2(s'+i) \in I_2$.  Then,  we say the \emph{switching position} of $I_1$ and $I_2$ is the smallest such $i$.  If no such value $i$ exists, then we say there is no switching position between $I_1$ and $I_2$ with respect to this overlap.
\end{definition}

\begin{lemma}[{\cite[Lemma 3]{banaian2024skein}}]\label{lem:SwitchingPosition}
Let $\calP_1,\calP_2$ be two labeled, weighted posets with overlaps $R_1:=\calP_1[s,t] \cong \calP_2[s',t']=:R_2$,  and let $I_1 \in J(\calP_1)$ and $I_2 \in J(\calP_2)$. Then,  a switching position between $I_1$ and $I_2$ exists unless $R_1 \cap I_1 = \emptyset $ and $R_2 \subseteq I_2$ or $R_1 \subseteq I_1$ and $R_2 \cap I_2 = \emptyset$.
\end{lemma}

We now exhibit a relation concerning the weight polynomial of a poset with a reverse self-kissing overlap. This new result can be seen as an addition to the poset skein relations discussed in \cite{banaian2024skein}, which built on the original snake graph calculus \cite{canakci2013snake}.

\begin{proposition}\label{prop:KissingReverseOverlap}
Given a poset $\calP$ with a reverse, kissing self-overlap in $\calP[r,s] \cong \overline{\calP[s+1,t]}$, define $Z^\star_R$ and $Z^\star_t$ by \[
Z^\star_R = \prod_{i = s+1}^t wt^\star(\calP(x))
\quad
Z^\star_t= \begin{cases} \prod_{i=t+1}^{v-1} wt^\star(\calP(i)) & r = 1\\
1 & r>1\end{cases}.
\]
We can rewrite $\mathcal{W}(\calP)$ in the following way using the resolution of $\calP$,\[
\mathcal{W}^\star(\calP) = \mathcal{W}^\star(\calP_{34}) + Z^\star_RZ^\star_t \mathcal{W}^\star(\calP_{56}).
\]
\end{proposition}

\begin{proof}
 For shorthand, let $R$ denote the subposet $\calP[r,s]$ and let $R'$ denote the subposet $\calP[s+1,t]$.  We will follow the proof method which was described in \cite[Section 3.4]{banaian2024skein}.

 We first assume $r > 1$ and $t < d$.  Let $A$ denote the subset of $J(\calP)$ consisting of all order ideals $I$ such that if $R' \subseteq I$ and $R \cap I = \emptyset$, then $\calP(r-1) \in I$ and $\calP(t+1) \notin I$.  Let $B$ denote the complement of $A$ in $J(\calP)$; equivalently, $B$ is the set of $I \in J(\calP)$ such that $R' \subseteq I$, $R \cap I = \emptyset$, and $\calP(r-1) \in I$ only if $\calP(t+1) \in I$.

 We begin with a weight-preserving bijection $\Omega$ between $A$ and $J(\calP_{34})$. Let $A = A_1 \sqcup A_2$ where $A_1$ consists of $I \in J(\calP)$ such that we do not have  $R' \subseteq I$ and $R \cap I = \emptyset$ and $A_2$ is the complement, namely, the set of $I \in J(\calP)$ such that $R' \subseteq I$, $R \cap I = \emptyset$,  $\calP(r-1) \in I$ and $\calP(t+1) \notin I$.

 Suppose we have $I \in A_1$. From Lemma \ref{lem:SwitchingPosition}, we know there exists a switching position,  $0 \leq i \leq s-r$, between $I$ and itself with respect to the reverse, kissing self-overlap. We define $\Omega(I)$ to be the order ideal of $\calP_{34}$ which is determined by the following conditions.
 \begin{enumerate}\setlength{\leftskip}{-15pt}
 \item The subset of $\Omega(I)$ supported on  $\calP_{34}[1,r+i-1] \cup \calP_{34}[t-i+1,h]$ is equivalent to the subset of $I$ supported on $\calP[1,r+i-1] \cup \calP[t-i+1,h]$.
 \item The subset of $\Omega(I)$ supported on $\calP_{34}[r+i,s]$ is equivalent to the subset of $I$ supported on $\calP[s+1,t-i]$.
\item The subset of $\Omega(I)$ supported on $\calP_{34}[s+1,t-i]$ is equivalent to the subset of $I$ supported on $\calP[r+i,s]$.
 \end{enumerate}

From part 3 of Definition \ref{def:ResolutionReverseSelfOverlap} and the fact that  $\calP(s) \in I$, is sent to $\calP_{34}(s+1) \in \Omega(I)$ and vice versa,  we conclude  that if $I \in A_1$, then $wt^\star(I) = wt^\star(\Omega(I))$.

 Now, suppose we have $I \in A_2$. We define $\Omega(I)$ to be the element of $J(\calP_{34})$ which contains $\calP_{34}[r,s]$ and whose restriction to $\calP[1,r-1] \cup \calP[t+1,h]$ is the same as that of $I$.  Part 3 of Definition \ref{def:ResolutionReverseSelfOverlap} again confirms that if $I \in A_2$,  then $wt^\star(I) = wt^\star(\Omega(I))$.

Notice that the image $\Omega(A_1)$ consists of all $I \in J(\calP_{34})$ such that we do not have $\calP_{34}[r,s] \subseteq I$ and $\calP_{34}[s+1,t] \cap I = \emptyset$ while $\Omega(A_2)$ exactly contains the complement.  Each map is reversible,  yielding a bijection between $A$ and $J(\calP_{34})$.

 Given $I \in B$, we can naturally construct an order ideal of $\calP_{56}$ from $I \backslash \calP[s,t]$, and this yields a bijection between $B$ and $J(\calP_{56})$.  Since we delete the elements in $\calP[r,s]$,  we see that $wt^\star(I) = Z^\star_R wt^\star(\Omega(I))$ for all $I \in B$.  Since $J(\calP) = A \sqcup B$,  our formula follows in this case.

 Next,  suppose $r = 1$,  which implies $t < h$.  If we define $A = A_1 \sqcup A_2$ where $A_1$ is identical to the previous case and $A_2$ consists of all $I \in J(\calP)$ such that $R' \subseteq I$, $R \cap I = \emptyset$ and $\calP(t+1) \notin I$,  then we can use a similar weight-preserving bijection $\Omega: A \to J(\calP_{34})$.  Let $B$ again be the complement of $A$; that is,  $B$ consists of all $I \in J(\calP)$ such that $R' \cup \langle \calP(t+1) \rangle \subseteq I$ and $R \cap I = \emptyset$.  Given $I \in B$,  $I \backslash (R' \cup \langle \calP(t+1) \rangle)$ is an element of $J(\calP_{56})$,  and this defines a bijection between the two sets.  Since for all $I \in B$,  $wt^\star(I) = Z^\star_R Z^\star_t wt^\star(\Omega(I))$,  our claim again follows.

 Finally,  suppose $t = h$,  implying $r > 1$.  If we define $A = A_1 \sqcup A_2$ where $A_1$ is identical to the previous case and $A_2$ consists of all $I \in J(\calP)$ such that $R' \subseteq I$, $R \cap I = \emptyset$ and $\calP(r-1) \in I$,  then we can use a similar weight-preserving bijection $\Omega: A \to J(\calP_{34})$. Let $B$ again be the complement of $A$; that is,  $B$ consists of all $I \in J(\calP)$ such that $R'  \subseteq I$ and $(R \cup \{\calP(r-1)\}) \cap I = \emptyset$.  Given $I \in B$,  $I \backslash R'$ is an element of $J(\calP_{56})$,  and this defines a bijection between the two sets.  Since for all $I \in B$,  $wt^\star(I) = Z^\star_R  wt^\star(\Omega(I))$,  our claim again follows.
\end{proof}

\begin{example}\label{ex:ReverseKissingSelfOverlap}
Consider the following poset $\calP$, where elements are expressed by a pair which first denotes their label and then denotes their weight. Use the chronological ordering induced by reading the Hasse diagram left to right.

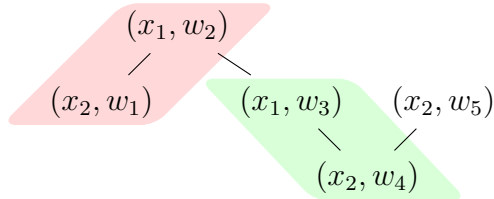
\begin{center}
\begin{tikzpicture}
\filldraw [rounded corners=5pt, red!15]
(-1.3,-.3)
-- (0.3,1.3)
-- (2.3,1.3)
-- (0.7,-.3)
-- cycle;
\filldraw [rounded corners=5pt, green!15]
(1.3,.3)
-- (3.3,.3)
-- (4.8,-1.3)
-- (2.8,-1.3)
-- cycle;
\node(0) at (0,0){$(x_2,w_1)$};
\node(1) at (1,1){$(x_1,w_2)$};
\node(2) at (2.5,0){$(x_1,w_3)$};
\node(3) at (3.5,-1){$(x_2,w_4)$};
\node(4) at (4.5,0){$(x_2,w_5)$};
\draw(0) -- (1);
\draw(1) -- (2);
\draw(2) -- (3);
\draw(3) -- (4);
\end{tikzpicture}
\end{center}

If $w_1 = w_4$,  then according to Definition \ref{def:ReverseOverlap}, this poset has a reverse kissing self-overlap in $\calP[1,2] \cong \overline{\calP[3,4]}$. Using Definition \ref{def:ResolutionReverseSelfOverlap}, we construct the poset $\calP_{34}$ below. The same Definition tells us that $\calP_{56} = \emptyset$.

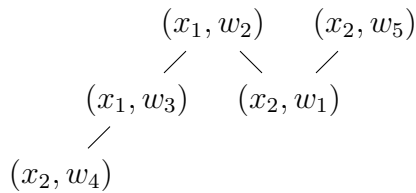
\begin{center}
\begin{tikzpicture}
\node(0) at (0,0){$(x_2,w_4)$};
\node(1) at (1,1){$(x_1,w_3)$};
\node(2) at (2,2){$(x_1,w_2)$};
\node(3) at (3,1){$(x_2,w_1)$};
\node(4) at (4,2){$(x_2,w_5)$};
\draw(0) -- (1);
\draw(1) -- (2);
\draw(2) -- (3);
\draw(3) -- (4);
\end{tikzpicture}
\end{center}

Applying Proposition \ref{prop:KissingReverseOverlap}, we have \[
\mathcal{W}^\star(\calP) = \mathcal{W}^\star(\calP_{34}) + w_3w_4w_5 .
\]
\end{example}

\begin{remark}
The resolution in Definition \ref{def:ResolutionReverseSelfOverlap} and Proposition \ref{prop:KissingReverseOverlap} may seem surprising since this corresponds to a resolution of a self-intersection, and usually such resolutions include a closed curve.  In our setting, a reverse, kissing self-overlap arises from a self-intersection near an orbifold point.  In \cite{banaian2024generalization}, the first author and Sen discussed $\mathcal{A}(1,1,1)$ as the generalized cluster algebra associated to a sphere with one puncture and three orbifold points of order 3, so it is natural that such relations appear in our combinatorial model. An example of the sort of intersection which gives rise to a reverse, kissing self-overlap can be seen in the third column in \cite[Table 1]{banaian2023snake}. In this setting, the closed curve in the resolution would be contractible onto the orbifold point and would only contribute a scalar.
If one were to define such poset combinatorics abstractly,  outside of the scope of cluster algebras from surfaces and orbifolds,  then possibly  Definition \ref{def:ResolutionReverseSelfOverlap} and Proposition \ref{prop:KissingReverseOverlap} would need to be broadened.
\end{remark}

\subsection{A Relation on Laurent Polynomials from Posets}\label{subsec:ExtendReverseKissingToXs}

This subsection is independent from the remainder of the paper. Here we will rewrite Proposition \ref{prop:KissingReverseOverlap} in terms of the variables $X_\calP$. To do this in a general setting, we need to introduce a relationship between the labels and weights of the poset elements compare. A broad setting where this will apply and be of interest is for a poset coming from a \emph{gentle pair}.

Recall a quiver is a directed graph.  Denote the start and tail of each arrow by $\mathrm{s}(\alpha) \xrightarrow{\alpha} \mathrm{t}(\alpha)$.  Let $Q$ be a finite quiver on $n$ vertices, labeled $\{1,\ldots,n\}$, such that every vertex has at most two incoming arrows and two outgoing arrows. Let a 2-path in $Q$ be two consecutive (not necessarily distinct) arrows, i.e. $\alpha_1\alpha_2$ such that $\mathrm{t}(\alpha_1) = \mathrm{s}(\alpha_2)$.  Let $Z$ be a set of length 2-paths from $Q$. We say $(Q,Z)$ is a \emph{gentle pair} if for every vertex $v$ of $Q$, the following is true. \begin{itemize}\setlength{\leftskip}{-15pt}
    \item For each arrow $\alpha$ which points into $v$, there is at most one arrow $\beta$ which points out of $v$ such that $\alpha \beta \in Z$ and there is at most one arrow $\gamma$ which points out of $v$ such that $\alpha \gamma \notin Z$.
    \item For each arrow $\beta$ which points out of $v$, there is at most one arrow $\alpha$ ($\delta$) which points into $v$ such that $\alpha \beta \in Z$ and there is at most one arrow $\delta$ which points into $v$ such that $\delta \beta \notin Z$.
\end{itemize}

We moreover say that a gentle pair is \emph{orbifold type} if all elements of $Z$  either are of the form $\mu^2$ for a loop $\mu$ or come in a triple $\{\alpha \beta,\beta \gamma,\gamma\alpha\}$ for a triangle $\alpha \beta\gamma$ in $Q$.  The rationale for this term comes from the first author's work with Valdivieso \cite{banaian2023snake}. Examples of gentle pairs of orbifold type are given in Table \ref{table:GentlePairs}.

For each arrow $\alpha$ in $Q$ let $\alpha^{-1}$ denote a formal inverse,  i.e. $\mathrm{s}(\alpha^{-1}) = \mathrm{t}(\alpha)$ and $\mathrm{t}(\alpha^{-1}) = \mathrm{s}(\alpha)$. We say a \emph{walk} $w$ in a quiver $Q$ is a sequence of arrows and their formal inverses $w = \alpha_1^{\epsilon_1} \alpha_2^{\epsilon_2} \cdots \alpha_\ell^{\epsilon_\ell}$ for $\epsilon_i \in \{+,-\}$ such that $\mathrm{t}(\alpha_i^{\epsilon_i}) = \mathrm{s}(\alpha_{i+1}^{\epsilon_{i+1}})$ for all $1 \leq i \leq \ell-1$.  Given the data of $Z$,  we say a walk is \emph{admissible} there is no subsequence of the form $\alpha\beta$ or $\beta^{-1}\alpha^{-1}$ for an element $\alpha \beta \in Z$ nor any subsequence $\alpha^\pm\alpha^\mp$. Every admissible walk of length $\ell$ naturally gives rise to a poset $\calP_w$ on $\ell + 1$ elements with a chronological ordering  which is the transitive closure of the relations \[
\calP_w(i) > \calP_w(i+1) \text{ for all } \epsilon_i = + \qquad \text{and} \qquad   \calP_w(i) < \calP_w(i+1) \text{ for all } \epsilon_i = - .
\]

Let $x_1,\ldots,x_n$ be $n$ indeterminants, one for each vertex of $Q$. Each element of $\calP_w$ corresponds to a vertex of $Q$. If $\calP_w(i)$ corresponds to vertex $v$, we label this element with $x_v$ and set its weight to be $\lambda_i \frac{\prod_{v \to u} x_u}{\prod_{z \to v} x_z}$ for some choice of $\lambda_i \in \mathbb{C}$.

All of our posets $\calP_\fpq$ and $\Ptilde_\fpq$ come from admissible walks in a gentle pair of orbifold type. These are displayed in Table \ref{table:GentlePairs}. Note that when some but not all $k_i$ are nonzero,  some indices will correspond to multiple vertices of the quiver. This is a consequence of the fact that in these cases, the $B$-matrix of the generalized cluster algebra $\mathcal{A}(k_1,k_2,k_3)$ is not skew-symmetric.

\begin{remark}
The quivers in Table \ref{table:GentlePairs} bear strong resemblance to the ``exceptional $s$-blocks'' listed in \cite[Table 3.2]{FeliksonShapiroTumarkin}. The gentle pair for all $k_i > 0$ appears in \cite[Example 2.3]{ford2024homological}.
The ordinary case ($k_i = 0$) has shown up in many places in the context of gentle algebras and representation theory as well. For example, as explained in \cite[Section 9]{RepType}, this quiver is a rare example of a quiver of finite mutation type which has multiple \emph{non-degenerate potentials}. A potential is a linear combination of cycles in a quiver which gives rise to a set of relations, and one of the non-degenerate potentials gives rise to the relations listed here.
\end{remark}

\begin{table}[ht]
\begin{center}
\begin{tabular}{cc}
$k_1 = k_2 = k_3 = 0$ & $k_1 = k_2 = 0, k_3 > 0$\\
\begin{tikzcd}[arrow style=tikz,>=stealth,row sep=3em]
& 2 \arrow[dr,shift left=1.25ex,"\beta_1"] \arrow[dr,swap,"\beta_2"] & \\
1 \arrow[ur, shift left = 1.25ex, "\alpha_1"]  \arrow[ur, swap,"\alpha_2"] & & 3 \arrow[ll,swap, "\gamma_2"]  \arrow[ll,shift left = 1.25ex, "\gamma_1"] \\
\end{tikzcd}\vspace{-.5cm}  &
\begin{tikzcd}[arrow style=tikz,>=stealth,row sep=3em]
& 2 \arrow[dl,"\beta_1"]\arrow[dr,swap,"\beta_2"] & \\
3\arrow[dr,"\gamma_1"] \arrow[out=180,in=90,loop, swap,"\nu_1"]& & 3\arrow[dl,swap,"\gamma_2"]\arrow[out=0,in=90,loop, "\nu_2"]  \\
& 1 \arrow[uu,swap, shift right = .6ex,"\alpha_2"] \arrow[uu,shift left = .6ex,"\alpha_1"] & \\
\end{tikzcd}\vspace{-.5cm} \\
$Z = \{\alpha_1\beta_1,\beta_1\gamma_1,\gamma_1\alpha_1,\alpha_2\beta_2,\beta_2\gamma_2,\gamma_2\alpha_2\}$ & $Z = \{\alpha_1\beta_1,\beta_1\gamma_1,\gamma_1\alpha_1,\alpha_2\beta_2,$\\
& $\beta_2\gamma_2,\gamma_2\alpha_2,\nu_1^2,\nu_2^2\}$ \\ \hline
&\\
$k_1 = 0,  k_2 >0, k_3 > 0$ & $k_1 > 0, k_2 > 0, k_3 > 0$\\
\begin{tikzcd}[arrow style=tikz,>=stealth,row sep=3em]
2 \arrow[dd,"\beta_1"] \arrow[ out = 180, in =90,loop, swap,"\rho_1"] & & 3  \arrow[ out = 0, in =90,loop, "\nu_1"]\arrow[dl,"\gamma_2"]\\
& 1\arrow[ul,"\alpha_1"] \arrow[dr,"\alpha_2"] & \\
3 \arrow[ur,"\gamma_1"]  \arrow[out = 180, in =-90,loop,"\nu_2"] & & 2\arrow[uu,"\beta_2"]\arrow[out = 0, in =-90,loop,swap,"\rho_2"]\\
\end{tikzcd}  &
\begin{tikzcd}[arrow style=tikz,>=stealth,row sep=3em]
& 2 \arrow[dr,"\beta"]\arrow[out = 135, in = 45, loop,swap,"\rho"]  & \\
1 \arrow[ur, "\alpha"]\arrow[out = 180, in =90,loop,swap,"\mu"]  & & 3 \arrow[ll,"\gamma"] \arrow[out = 0, in =90,loop, "\nu"]\\
\end{tikzcd} \\
$Z = \{\alpha_1\beta_1,\beta_1\gamma_1,\gamma_1\alpha_1,\alpha_2\beta_2,$ & $Z = \{\alpha\beta,\beta\gamma,\gamma\alpha,\mu^2,\rho^2,\nu^2\}$\\
$\beta_2\gamma_2,\gamma_2\alpha_2,\rho_1^2,\rho_2^2,\nu_1^2,\nu_2^2\}$ & \\
\end{tabular}
\end{center}
\caption{Gentle pairs associated to the four types of generalized cluster algebras we consider here. }\label{table:GentlePairs}
\end{table}

Now that we are considering more general posets, we must include a more general definition of the vector $\mathbf{r}_{\calP}$.  Let $\calP_w$ come from an admissible walk $w = \alpha_1^{\epsilon_1}\cdots\alpha_\ell^{\epsilon_\ell}$.  If there exist arrows $\beta$ and $\gamma$ such that $\beta \alpha_1^{\epsilon_1}\cdots\alpha_\ell^{\epsilon_\ell}\gamma^-$ is an admissible walk,  then $\mathbf{r}_{\calP_w} = \ee_{s(\beta)} + \ee_{s(\gamma)}$.  If one or both of these arrows do not exist, then we ignore the corresponding term.  From the definition of a gentle pair,  if such a $\beta$ or $\gamma$ exists, it is unique.

For the following lemma, note that given a poset $\calP_w$ associated to an admissible walk,  every subposet of the form $\calP_w[u,v]$ is equivalent to a poset associated to an admissible subwalk of $w$.  The case $u = v$ is the case of a trivial walk of length 0.

\begin{lemma}\label{lem:ProductOfWeightsGeneral}
Let $(Q,Z)$ be a gentle pair of orbifold type and let $w = \alpha_1^{\epsilon_1} \cdots \alpha_\ell^{\epsilon_\ell}$ be an admissible walk.  Let $\calP_w$ be the corresponding poset.   Given $1 \leq j \leq k \leq \ell + 1$,  let $w' = \alpha_j^{\epsilon_j} \cdots \alpha_{k-1}^{\epsilon_{k-1}}$ and define $\gb'$ to be the $g$-vector associated to the poset $\calP_{w'}$. Let $\beta,\gamma,\delta,\rho$ be arrows such that $\beta w, \gamma^{-1}w,w\delta,w\rho^{-1}$ are admissible walks,  if they exist.  Let $\xb_{\mathrm{out}} = x_{\mathrm{t}(\gamma)}x_{\mathrm{t}(\delta)}$ and $\xb_{\mathrm{in}} = x_{\mathrm{s}(\beta)}x_{\mathrm{s}(\rho)}$
The product of weights of all elements of $\calP_w[j,k]$ is given by \[
\prod_{i=j}^k wt(\calP_w(i)) = \lambda \xx^{-2\gb'} \xb_{\mathrm{out}} \xb_{\mathrm{in}} x_{\mathrm{s}(\alpha_j^{\epsilon_j})}^{-1}x_{\mathrm{t}(\alpha_{k-1}^{\epsilon_{k-1}})}^{-1}.
\]for a scalar $\lambda\in \mathbb C$.
\end{lemma}
\begin{proof}
We will fix an admissible walk $w$ and induct on $k-j$.  Let $\calP:=\calP_w$. First, suppose $k-j = 0$ so that $\calP_{w'}$ is a poset with one element $\{\calP(j)\}$.  Let $v$ be the vertex in $Q$ associated to $\calP(j)$. Then,  by definition $wt(\calP(j)) = \lambda \frac{\prod_{v \to u} x_u}{\prod_{z \to v} x_z} = \lambda \frac{\xb_{\mathrm{out}}}{\xb_{\mathrm{in}}}$ for some scalar $\lambda$.  If $\gb'$ is the $g$-vector associated to this singleton poset,  by definition we have $\xx^{-2\gb'} = \frac{x_v^2}{\xb_{\mathrm{in}}^2}$. Since in this case we interpret $s(\alpha_j^{\epsilon_j}) = t(\alpha_{k-1}^{\epsilon_{k-1}}) = v$,  the claim follows.

Now, suppose we have shown the claim for all subwalks of $w$ of length $m-1$, where $m-1 < \ell$, and consider a subwalk $w' = \alpha_j^{\epsilon_j} \cdots \alpha_{k-1}^{\epsilon_{k-1}}$ where $k-j = m$.  Let $w'' =  \alpha_j^{\epsilon_j} \cdots \alpha_{k-2}^{\epsilon_{k-2}}$; if $k-j = 1$,  then $w''$ is a trivial walk as in the base case.  The inductive hypothesis applies to $w''$,  so we have
\begin{align}
\prod_{i=j}^k wt(\calP(i)) = (\lambda'' \xx^{-2\gb''} \xb_{\mathrm{out}}'' \xb_{\mathrm{in}}'' x_{\mathrm{s}(\alpha_j^{\epsilon_j})}^{-1}x_{\mathrm{t}(\alpha_{k-2}^{\epsilon_{k-2}})}^{-1}) \lambda_k \frac{\prod_{\calP(k) \to u}x_u}{\prod_{z \to \calP(k)}x_z}, \label{eq:InductiveStep}
\end{align}

where $\gb''$ is the $g$-vector associated to $\calP[j,k-1]$ and in the products we conflate the poset elements and their associated vertices. We will in fact continue to conflate these. Note that, using this abuse of notation, $\mathrm{t}(\alpha_{k-2}^{\epsilon_{k-2}})$ is the vertex $\calP(k-1)$.

Suppose first that $\epsilon_{k-1} = +$.  Then,  in $Q$,  the most complicated neighborhood around the vertex $\calP(k)$ is the following.  Either of these triangles could in fact be a loop, in which case all vertices and all arrows coincide.

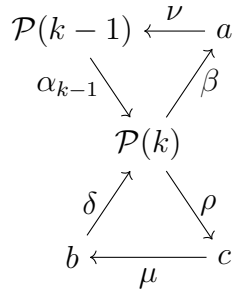
\begin{center}
\begin{tikzpicture}
\node(t1) at (-1,1.5){$\calP(k-1)$};
\node(t) at (0,0){$\calP(k)$};
\node(a) at (1,1.5){$a$};
\node(b) at (-1,-1.5){$b$};
\node(c) at (1,-1.5){$c$};
\draw[->] (t1) to node[left]{$\alpha_{k-1}$}(t);
\draw[->] (t) to node[right]{$\beta$} (a);
\draw[->] (a) to node[above]{$\nu$} (t1);
\draw[->] (t) to node[right]{$\rho$}(c);
\draw[->] (c) to node[below]{$\mu$} (b);
\draw[->] (b) to node[left]{$\delta$} (t);
\end{tikzpicture}
\end{center}

If the neighborhood of $\calP(k)$ is as above, then we have \[
\{\alpha_{k-1}\beta,\beta\nu,\nu\alpha_{k-1},\delta\rho,\rho\mu,\mu\delta\} \subseteq Z.
\]

The existence of $\nu$ and $\mu$ and the form of $Z$ comes from the fact that $(Q,Z)$ is of orbifold-type.
In all cases, the neighborhood of $\calP(k)$ is a subdiagram of the following, resulting from deleting some vertices and incident arrows. Necessarily the corresponding terms would also be removed from $Z$.  We will assume we are in this most complicated setting. In the other cases, one can take the expressions here and specialize variables associated to non-existent vertices to 1.

From this diagram, we have $wt(\calP(k)) = \frac{x_ax_c}{x_{\calP(k-1)}x_{b}}$.  Now,  we revisit Equation \eqref{eq:InductiveStep} with our more detailed information regarding the neighborhood of $\calP(k)$,\[
\prod_{i=j}^k wt(\calP(i)) = \big(\lambda'' Z \frac{ x_{\calP(k)}}{x_ax_{\calP(k-1)}}\big)\lambda_k  \frac{x_ax_c}{x_{\calP(k-1)}x_{b}} = \lambda Z \frac{x_{\calP(k)} x_c}{x_{\calP(k-1)}^2 x_b}
\]
where $\lambda = \lambda'' \lambda_k$ and $Z$ represents the remainder of the monomial $\prod_{i=j}^{k-1} wt(\calP(i))$ after factoring out $\lambda'' \frac{ x_{\calP(k)}}{x_ax_{\calP(k-1)}}$.

Now,  we verify this monomial has the desired form, that is, that \[
\lambda Z \frac{x_{\calP(k)} x_c}{x_{\calP(k-1)}^2 x_b} = \lambda \xx^{-2\gb'} \xb_{\mathrm{out}} \xb_{\mathrm{in}} x_{\mathrm{s}(\alpha_j^{\epsilon_j})}^{-1}x_{\mathrm{t}(\alpha_{k-1}^{\epsilon_{k-1}})}^{-1}.
\]

Since $\calP(k)$ is minimal in $\calP_{w'}$,  we know that $x_{\calP(k)}$ appears with multiplicity 2 in the numerator of $\xx^{-2\gb'}$. Since $\mathrm{t}(\alpha_{k-1}^{\epsilon_{k-1}} = \calP(k)$ as well, the $x_{\calP(k)}$ matches on each side. Next,  we see that arrows $\rho$ and $\delta$ in the figure above match the arrows $\rho$ and $\delta$ from the statement of the lemma, implying $x_c$ is a factor of $\mathbf{x}_{\textrm{out}}$ and $x_b$ is a factor of $\mathbf{x}_{\textrm{in}}$.  Since $x_b$ appears with multiplicity 2 in the denominator of $\xx^{-2\gb'}$, the exponents of each also match.

Finally,  when analyzing the exponent of $x_{\calP(k-1)}$ on each side, there are two subcases based on whether $\calP(k-1)$ was minimal or maximal in $\calP[j,k-1]$.  If it was minimal,  then $x_{\calP(k-1)}$ appeared with multiplicity 2 in the numerator of $\xx^{-2\gb''}$, i.e. in $Z$.  Since $\calP(k-1)$ would not be minimal or maximal in $\calP[j,k]$, it does not appear on the righthand side, which agrees with the lefthand side. Otherwise,  if $\calP(k-1)$ was maximal in $\calP[j,k-1]$ then it is still maximal in $\calP[j,k]$.  It is only strictly maximal in the latter though,  so while there is no factor $x_{\calP(k-1)}$ in $\xx^{-2\gb''}$, i.e. in $Z$, it appears with multiplicity 2 in the denominator of $\xx^{-2\gb'}$.  From this analysis,  we conclude that $ \lambda Z \frac{x_{\calP(j)} x_c}{x_{\calP(k-1)}^2 x_b}$ is exactly the monomial from the statement of the lemma.

If $\epsilon_{k-1} = -$,  one can perform similar checks to verify the lemma.
\end{proof}

When a gentle pair of orbifold type does not have any loops,  Lemma \ref{lem:ProductOfWeightsGeneral} is implied by \cite[Lemma 2]{banaian2024skein}.  This is true because all such gentle pairs (i.e. those of orbifold type without loops) arise from a triangulation of a surface.
With an understanding of how products of weights compare with $g$-vectors of subposets, we can extend Proposition \ref{prop:KissingReverseOverlap} to the values $X_\calP$.

When the poset $\calP_{56}$ from Definition \ref{def:ResolutionReverseSelfOverlap} is empty, we must take more care in defining $X_{\calP_{56}}$. Let $\calP$ be a poset from an admissible walk $w$ in a gentle pair of orbifold type $(Q,Z)$ with a reverse, kissing self-overlap. Let $h = \vert \calP \vert$. Recall the parameters $r$ and $t$ from Definition \ref{def:ReverseOverlap}.

\begin{itemize}
    \item Suppose $r = 1$ and $\calP_{56} = \emptyset$. If there exists an arrow $\beta$ such that $w\beta$ is an admissible walk, let $X_{\calP_{56}} = x_{t(\beta)}$. Otherwise, let $X_{\calP_{56}} = 1$.
    \item Suppose $t = h$ and $\calP_{56} = \emptyset$. If there exists an arrow $\alpha$ such that $\alpha^{-1}w$ is an admissible walk, let $X_{\calP_{56}} = x_{s(\alpha^{-1})}$. Otherwise, let $X_{\calP_{56}} = 1$.
\end{itemize}

These special cases for $X_{\calP_{56}}$ correspond to a geometric resolution of a self-intersection where one arc in the resolution lies in the initial triangulation that determines $(Q,Z)$.

The following builds on Proposition \ref{prop:KissingReverseOverlap}. The comparison of $g$-vectors in the proof is akin to calculations done in \cite{banaian2024skein} although the setting is new.

\begin{theorem}\label{thm:DecomposeVariables}
Let $\calP$ be a poset which comes from an admissible walk in a gentle pair of orbifold type $(Q,Z)$. Suppose that $\calP$ has a reverse, kissing self-overlap in $\calP[r,s] \cong \overline{\calP[s+1,t]}$.  We can rewrite $X_\calP$ using the resolution of $\calP$,\[
X_{\calP} = X_{\calP_{34}} + \lambda X_{\calP_{56}},
\]
for a scalar $\lambda\in \mathbb C$.
\end{theorem}
\begin{proof}
Recall in the proof of Proposition \ref{prop:KissingReverseOverlap}, we partitioned $J(\calP)$ as $A \sqcup B$ and provided maps $\Omega_A: A \to J(\calP_{34})$ and $\Omega_B: B \to J(\calP_{56})$. These maps were such that, for all $I \in A$, $wt(I) = wt(\Omega_A(I))$ and for all $I \in B$, $wt(I) = Z_RZ_t wt(\Omega_B(I))$. Therefore, we can rewrite $X_\calP$ as \begin{align*}
X_\calP &= \xx^{\gg_\calP} \sum_{I \in J(\calP)} wt(I)\\
&= \xx^{\gg_\calP} \big(\sum_{I \in A} wt(I) + \sum_{I \in B} wt(I)\big)\\
&= \xx^{\gg_\calP} \big(\sum_{I \in J(\calP_{34})} wt(I) + \sum_{I \in J(\calP_{56})} Z_RZ_twt(I)\big).
\end{align*}
Therefore, we will be done if we can show $\xx^{\gg_\calP} = \xx^{\gg_{\calP_{34}}}$ and $\xx^{\gg_\calP}Z_RZ_t = \xx^{\gg_{\calP_{56}}}$. The former equality is immediate from the relationship between $\calP$ and $\calP_{34}$.

We turn to comparing $\gg_\calP$ and $\gg_{\calP_{56}}$. We will assume that  all labels are distinct in $\calP$. One can always begin with this assumption and specialize at the end of the process. Assume first $r > 1$. Given a statement $G$, let $\delta_G = 1$ if $G$ is true and $\delta_G = 0$ if $G$ is false.  As before, assume $Q$ has vertices $\{1,\ldots,n\}$ so that the following vectors are in $\mathbb{R}^n$.  By analyzing the behavior of $\calP$ close to the overlap, we calculate \begin{align*}
\gg_\calP &= \gg' + 2\gg_R - (1 - \delta_{\calP(r-1) > \calP(r-2)})\ee_{\calP(r-1)} + \delta_{\calP(t+1) > \calP(t+2)} \ee_{\calP(t+1)}\\
&+  (-1)^{\delta_{\calP(t) < \calP(t-1)}} \ee_{\calP(t)}  + (-1)^{\delta_{\calP(s+1)<\calP(s+2)}}\ee_{\calP(s+1)}
\end{align*}
where we say that a statement is false if the elements in the inequality do not all exist, $\gg'$ consists of all summands $\ee_x$ with $x \in \calP[1,r-2] \cup \calP[t+2,\vert \calP \vert]$, and $\gg_R$ consists of all summands $\ee_x$ with $x \in \calP[r+1,s-1]$.  The third and fourth term come from the fact that $\calP(r-1)$ contributes if it is a minimal element and $\calP(t+1)$ contributes if it is a strict maximal element. We know that either $\calP(r)$ is a strict maximal element or $\calP(t)$ is a minimal element, and from the definition of a reverse, kissing self-overlap, we know these are labeled identically, explaining  the fifth term. Similarly, either $\calP(s)$ is a strict maximal element or $\calP(s+1)$ is a minimal element, explaining the sixth term.

Given a monomial $\xx^{\vv}$, let $deg(\xx^\vv) = \vv$.  We claim that from Lemma \ref{lem:ProductOfWeightsGeneral},  we have \begin{align*}
\deg(\prod_{i=s}^t wt(\calP(x))) &= -2\gg_R + (-1)^{\delta_{\calP(s+1) > \calP(s+2)}} \ee_{\calP(s+1)} + (-1)^{\delta_{\calP(t) > \calP(t-1)}} \ee_{\calP(t)} \\
& -\ee_{\calP(t+1)} + \ee_{\calP(r-1)}.
\end{align*}

 In the language of this lemma, $\deg(\xb_{\mathrm{out}}) = \ee_{\calP(s)} + \ee_{\calP(r-1)}$ and $\deg(\xb_{\mathrm{in}} = \ee_{\calP(s)} + \ee_{\calP(t+1)}$. If $\gb':= \gb_{\calP[s,t]}$, then following the definitions we have $\gb' = \gg_R - (1 - \delta_{\calP(t) > \calP(t-1)}) \ee_{\calP(t)} - (1 - \delta_{\calP(s+1) > \calP(s+2)}) \ee_{\calP(s+1)} + \ee_{\calP(s)} + \ee_{\calP(t+1)}$. Comparing these values shows us the expression above is valid. By adding in $\gb_\calP$, we have
 \[
\gg_\calP + \deg(\prod_{i=s}^t wt(\calP(x))) = \gg' + (-1 + \delta_{\calP(t+1) > \calP(t+2)}) \ee_{\calP(t+1)} + \delta_{\calP(r-1) > \calP(r-2)}\ee_{\calP(r-1)},
\]
which can be shown to be identical to $\gg_{\calP_{56}}$ following the definition.

Now, suppose $r = 1$, which implies $t < \vert \calP \vert$. We update our initial calculation, \begin{align*}
\gg_\calP &= \gg' + 2\gg_R +  2\ee_{\calP(t+1)}
-\delta_{\calP(t) < \calP(t-1)} \ee_{\calP(t)}-\delta_{\calP(1) < \calP(2)} \ee_{\calP(1)}  \\
&+ (-1)^{\delta_{\calP(s+1)<\calP(s+2)}}\ee_{\calP(s+1)} - \ee_{\calP(v-1)} + \delta_{\calP(v) > \calP(v+1)} \ee_{\calP(v)}
\end{align*}
where $v$ is as in Definition \ref{def:ResolutionReverseSelfOverlap} and now $\gg'$ consists of all summands $\ee_x$ of $\gg_{\calP}$ with $x \in \calP[v+1,\vert \calP \vert]$. Now that $r = 1$, recall we have $\calP(1)$ and $\calP(t)$ are labeled in the same way. Recall from Definition \ref{def:ResolutionReverseSelfOverlap} we defined $v> t+1$ to be the smallest integer such that $\calP(v) \nless \calP(t+1)$ if it exists, and if such an integer does not exist $v = h+1$.  Similar to before,  from Lemma \ref{lem:ProductOfWeightsGeneral} we have
\begin{align*}
\deg(\prod_{i=s}^{\min(v,h)} wt(\calP(x))) &= -2\gg_R + (-1)^{\delta_{\calP(s+1) > \calP(s+2)}}\ee_{\calP(s+1)} + 2 \delta_{\calP(t) > \calP(t+1)} \ee_{\calP(t)}\\
&- 2 \ee_{\calP(t+1)} + \ee_{\calP(v-1)} -\ee_{\calP(v)} + \ee_c
\end{align*}
where if $w$ is the admissible walk associated to $\calP$,  $\ee_c$ is the degree vector for factor from $\xb_{\mathrm{out}}$ corresponding to the beginning of $w$.. One can again now verify that $\gg_\calP +  \deg(\prod_{i=s}^v wt(\calP(x)))  = \gg_{\calP_{56}}$.
\end{proof}

\begin{example}
The poset in Example \ref{ex:ReverseKissingSelfOverlap} can be seen as coming from the string $\alpha^{-1}\mu\alpha\rho^{-1}$ from the gentle pair given at the bottom right of Table \ref{table:GentlePairs}. We must use the more general definition of the  vector $\mathbf{r}_{\calP}$ since $\calP$ is not of the form $\calP_\fpq$. Using this, we compute \[
\mathbf{g}_{\calP} = -\mathbf{a}_\calP + \mathbf{b}_\calP + \mathbf{r}_\calP = -\begin{bmatrix}0\\2\\0\end{bmatrix} + \begin{bmatrix}1\\0\\0\end{bmatrix} + \begin{bmatrix}1\\1\\0\end{bmatrix} =\begin{bmatrix}2\\-1\\0\end{bmatrix}
\]
and we can similarly compute $\mathbf{g}_{\calP_{34}} = (2,-1,0)^\top$. In this setting, the weights are of the form \[
w_1 = \lambda_1 \frac{x_2}{x_3} \quad w_2 = \lambda_2 \frac{x_1}{x_2} \quad w_2 = \lambda_3 \frac{x_1}{x_2} \quad w_4 = \lambda_4 \frac{x_2}{x_3} \quad \text{ and} \quad w_5 = \lambda_5 \frac{x_2}{x_3}.
\]

Applying the proof method of Proposition \ref{prop:KissingReverseOverlap}, we have a weight-preserving bijection between $J(\calP) \backslash \{\calP[3,5]\}$ and $J(\calP_{34})$. Since $\mathbf{g}_\calP = \mathbf{g}_{\calP_{34}}$, we have \[
X_{\calP} - X_{\calP_{34}} = \mathbf{x}^\gg wt(\calP[3,5]) = \lambda_3\lambda_4\lambda_5x_3.
\]

Here, $r = 1$ and $\calP_{56} = \emptyset$. Since the string $\alpha^{-1}\mu\alpha\rho^{-1}$ can be extended on the right by $\beta$, we have $X_{\calP_{56}} = x_{t(\beta)} = x_3$, verifying the statement of Theorem \ref{thm:DecomposeVariables}.
\end{example}

\begin{remark}\label{rem:CC}
 Let $kQ$ denote the path algebra associated to a quiver $Q$. Given a gentle pair $(Q,Z)$, the quotient of $kQ$  by the ideal generated by $Z$ is called a \emph{gentle algebra}.  In this language,  our posets are \emph{strings}, combinatorial objects which index some of the indecomposable modules in \emph{gentle algebras} \cite{butler1987auslander}.  Moreover, the formula for $X_\calP$ is similar to applying the Caldero-Chapoton functions associated to the corresponding string module \cite{cerulli2015caldero}.  In this correspondence, a fence poset having a reverse, kissing self-overlap implies the associated module is not rigid. In some cases, this results in an algebraic dependence involving the Caldero-Chapoton functions. In particular, we remark that the resolution given in Definition \ref{def:ResolutionReverseSelfOverlap} and the identity given in Theorem \ref{thm:DecomposeVariables} bear strong similarities to the third, fifth, and sixth relations on Caldero-Chapoton functions given in the proof of \cite[ Proposition 9.4]{cerulli2015caldero}; see also \cite[Remark 11]{labardini2019family}.
\end{remark}

\subsection{Structure of Markov Posets}
In the case $k_1=k_2=k_3=0$, i.e. the ordinary Markov case, the shape of the poset $\calP_\fpq$ can be deduced from the continued fraction associated to the Markov number associated to $\fpq$. These continued fractions are, for example, computed in \cite{rabideau2020continued} (see also Remark \ref{rem:Christoffel}), and they are always palindromic. However, as noted in \cite{banaian2024generalization,gyoda2024sl}, once some of the $k_i$ are nonzero, the continued fractions are not guaranteed to be palindromic.
This asymmetry will be reflected in our poset formulas. Before we explain this fully, we will provide more details about the posets $\calP_\fpq$.

We introduce a simple rule to determine the parities $k_{1; \fpq}$. For $\frac{p}{q}\in [0,\infty]\cap\mathbb Q$, we set
\[k_{\fpq}:=\begin{cases}k_1\quad \text{if $p$ and $q$ are odd,}\\k_2\quad \text{if $p$ is odd and  $q$ is even,}\\k_3 \quad \text{if $p$ is even and $q$ is odd.}\end{cases}\]

This correspondence is determined based on the fact that $x_1$, $x_2$, and $x_3$ correspond to segments with slopes of $\frac{-1}{1}$, $\frac{1}{0}$, and $\frac{1}{1}$, respectively.

We have the following lemma:
\begin{lemma}\label{lem:parity}
The value $k_{\fpq}$ coincides with the parity $k_{1,\fpq}$ of $x_{1,\fpq}.$
\end{lemma}

\begin{proof}
We can prove this by inducting on the shortest path from the root of $\mathrm{F} \mathbb{T}$ to the first Farey triple containing $\fpq$. First, we check directly all three fractions in the root.
By the definition of the fraction labeling, we have $x_{\frac01}=x_3$, $x_{\frac11}=\frac{x_2^2+k_1x_2x_3+x_3^2}{x_1}$, and $x_{\frac{1}{0}}=x_2$. Since the parities of $x_3$, $\frac{x_2^2+k_1x_2x_3+x_3^2}{x_1}$ and $x_2$ are $k_3$, $k_1$ and $k_2$, respectively, we have the claim for $\frac{p}{q}\in \{\frac01,\frac11,\frac10\}$.

One can prove inductively that for every Farey triple,  $(\frac{a}{b}, \frac{c}{d}, \frac{e}{f})$, the parity types are all distinct. In particular, the parity type of $\frac{e}{f}$ is the same as that of $\frac{a + c}{b + d}$, and the parity type of $\frac{a}{b}$ is the same as that of $\frac{c + e}{d + f}$. Therefore, since the rotation rules of the components in triples when tracing down the tree are the same in both the $\mathrm F\mathbb T$ and $\mathrm{CM}\mathbb T(k_1,k_2,k_3, (1\ 3\ 2))$, we have the claim inductively for any $\fpq\in[0,\infty]\cap \mathbb Q$.
\end{proof}

For the remainder of this article, we will use $k_\fpq$ in place of $k_{1,\fpq}$.

\begin{lemma}\label{lem:WhenIsPfpqPalindromic}
The number of elements of  $\vert \calP_\fpq \vert$ is odd if and only if $k_\fpq = 0$. Moreover, the elements $\calP_\fpq(i)$ and $\calP_\fpq(h_\fpq + 1 - i)$ have the same label and weight for any $1 \leq i <  \frac12 h_\fpq$.
\end{lemma}

\begin{proof}
The sequence of labels of intersections between $\gamma_\fpq$ and $\mathcal{L}$ is palindromic when we consider all lines in $\mathcal{L}$ of the same slope as being equivalent.   This sequence is also odd,  where the middle crossing occurs with the line of the same ``parity'' as $\fpq$. For example, if $p$ and $q$ are both odd, then this middle crossing is a line of slope $\frac{-1}{1}$.  From the Construction Algorithm (Definition \ref{def:ConstructPpq}),  we see that the intersections strictly before this middle intersection contribute the same number of elements of $\calP_\fpq$ as the intersections strictly after the middle intersection. Therefore,  the total number of elements of $\calP_\fpq$ is even exactly when the middle intersection contributes two elements, that is, when $k_\fpq > 0$.  The palindromic nature of the intersections also explains the second statement.
\end{proof}

We define an important $y$-monomial which will show up frequently in the principal coefficient case.

\begin{definition}
Let $\fpq \in \mathbb{Q} \cap (0,\infty)$. Define $Y_\fpq = y_1^{d_1(q+p-1)}y_2^{d_2(q-1)}y_3^{d_3(p-1)}$.
\end{definition}

\begin{lemma}\label{lem:ProductOfAllWeights}
The product of the weights of all elements in $\calP_\fpq$ is given by
\[
\prod_{x \in \calP_\fpq} wt(x) = \xb^{-2 \gb_\fpq^\circ} \qquad \prod_{x \in \calP_\fpq} wt^\prin(x) = Y_\fpq\xb^{-2 \gb_\fpq^\circ}
\]
and if $k_\fpq > 0$, then \[
\prod_{x \in \calP_\fpq[1,\frac12 \cdot h_\fpq]} wt(x) = k_\fpq \xb^{- \gb_\fpq^\circ} \qquad \prod_{x \in \calP_\fpq[1,\frac12 \cdot h_\fpq]} wt^\prin (x) = k_\fpq \sqrt{Y_\fpq} \xb^{- \gb_\fpq^\circ}
\]
\end{lemma}

\begin{proof}
In Proposition \ref{prop:Computedfpq}, we count the number of elements in $\calP_\fpq$ with each label. With similar reasoning, we see that $\gamma_\fpq$ crosses $p+q-1$ line segments of slope $-1$, $q-1$ with slope $\infty$ and $p-1$ with slope 0. When we take the product of all elements of $\calP_\fpq$,  we will cancel all contributions of $k_i$ for each $i$, and each crossing contributes the same $x$-monomial. Therefore, we compute  \[
\prod_{x \in \calP_\fpq} wt^\prin(x)  = \bigg(y_1^{d_1} \frac{x_2^2}{x_3^2}\bigg)^{p+q-1}\bigg(y_2^{d_2} \frac{x_3^2}{x_1^2}\bigg)^{q-1}\bigg( y_3^{d_3}\frac{x_1^2}{x_2^2}\bigg)^{p-1} = Y_\fpq x_1^{2(p-q)}x_2^{2q} x_3^{-2p}
\]
From Theorem \ref{thm:g-vector-description} and Lemma \ref{lem:CompareGVecAndGCirc},  we can conclude that this monomial is $\mathbf{x}^{-2\gb^\circ_\fpq}$.  The second statement follows from the symmetry described in Lemma \ref{lem:WhenIsPfpqPalindromic}.
\end{proof}

\begin{remark}
The symmetry described in Lemma \ref{lem:WhenIsPfpqPalindromic} also guarantees that, if $k_\fpq > 0$, then $\sqrt{Y_\fpq}$ has only integer powers.
\end{remark}

The following lemma describes the asymmetry present in the poset $\calP_\fpq$ when $k_\fpq$ is nonzero. This asymmetry was previously explored for the $k_1 = k_2 = k_3 = 1$ with all $x$-variables set to 1 in \cite{banaian2024generalization}, but here we add more detail. We will only use this lemma in the coefficient-free case, so we only list this version.

\begin{lemma}\label{lem:Asymmetry}
The two polynomials $\mathcal{W}(\calP_\fpq;\neg R) $ and $ \mathcal{W}(\calP_\fpq;\neg L)$ are related by
\[
\mathcal{W}(\calP_\fpq;\neg R) = \mathcal{W}(\calP_\fpq;\neg L) + k_\fpq \xb^{-\gb^\circ_\fpq}.
\]
Similarly, we have \[
\mathcal{W}(\calP_\fpq;L) = \mathcal{W}(\calP_\fpq;R) + k_\fpq \xb^{-\gb^\circ_\fpq}
\]
and
\[
\mathcal{W}(\calP_\fpq;L, \neg R) = \mathcal{W}(\calP_\fpq;R, \neg L) + k_\fpq \xb^{-\gb^\circ_\fpq}
\]
\end{lemma}

\begin{proof}
Each statement is trivial if $\fpq = 1$. Now, for $\fpq \neq 1$, there are six cases, coming from showing each of the two statement in the cases $\fpq < 1$ and $\fpq > 1$.

We will only prove the first statement in the case $\fpq < 1$ as all other proofs follow similarly. Since $\fpq < 1$, we know $\calP_\fpq(1) > \calP_\fpq(2)$ and $\calP_\fpq(h_\fpq) > \calP_\fpq(h_\fpq - 1)$.  Therefore, we can rewrite $\mathcal{W}(\calP_\fpq;\neg R) = \mathcal{W}(\calP_\fpq[1,h_\fpq-1])$ and $\mathcal{W}(\calP_\fpq;\neg L) = \mathcal{W}(\calP_\fpq[2,h_\fpq])$.

Suppose first that $k_\fpq = 0$.  Then,  it follows from Lemma \ref{lem:WhenIsPfpqPalindromic} that $\calP_\fpq[1,h_\fpq-1]$ and $\calP_\fpq[2,h_\fpq]$ are isomorphic as weighted labeled posets,  and the claim follows.

Next, suppose that $k_\fpq > 0$.  In this case,  the poset $\calP_\fpq[1,h_\fpq-1]$ has a reverse, kissing self-overlap in $\calP[2,\frac12 \cdot h_\fpq]$ and  $\calP[\frac12\cdot h_\fpq+1,h_\fpq-1]$.  The fact that this is a self-overlap and that $\frac{h_\fpq}{2}$ is an integer again follows from Lemma \ref{lem:WhenIsPfpqPalindromic}.  The poset $\calP_{34}$ from Proposition \ref{prop:KissingReverseOverlap} is identical to $\calP_\fpq[2,h_\fpq]$.  In this case,  since the largest index of the overlap is equal to $h_\fpq$,  we have that $\calP_{56} = \emptyset$.  Therefore,  the result follows from Proposition \ref{prop:KissingReverseOverlap}  and Lemma \ref{lem:ProductOfAllWeights}.
\end{proof}

\begin{example}\label{ex:Asymmetry}
Let $k_1,k_2$ be positive integers and consider $\calP_{\frac12}$:
\begin{center}
\begin{tikzpicture}[scale=1.25]
\node(a) at (0,2){$x_1$};
\node(b) at (0.5,1.5){$x_1$};
\node(c) at (1,1){$x_2$};
\node(d) at (1.5,1.5){$x_2$};
\node(e) at (2,2){$x_1$};
\node(f) at (2.5,2.5){$x_1$};
\draw(a) -- (b);
\draw(b)--(c);
\draw(c)--(d);
\draw(d)--(e);
\draw(e)--(f);
\end{tikzpicture}
\end{center}
We compute the two polynomials,\[
\mathcal{W}(\calP_{\frac12};\neg R) = 1 + k_2\widehat{x_2}(1 + k_1\widehat{x_1} + \widehat{x_1}^2)(1 + \frac{1}{k_2}\widehat{x_2} + \frac{k_1}{k_2}\widehat{x_1}\widehat{x_2})
\]
\[
\mathcal{W}(\calP_{\frac12};\neg L) = 1 + k_2\widehat{x_2}(1 + k_1\widehat{x_1})(1 + \frac{1}{k_2}\widehat{x_2} + \frac{k_1}{k_2}\widehat{x_1}\widehat{x_2}+ \frac{1}{k_2}\widehat{x_1}^2\widehat{x_2})
\]
and we can see \[
\mathcal{W}(\calP_{\frac12};\neg R) - \mathcal{W}(\calP_{\frac12};\neg L)  = k_2 \widehat{x_1}^2\widehat{x_2} = k_2 \frac{x_2^2}{x_3x_1} = k_2 \mathbf{x}^{(-1,2,-1)^\top}.
\]
\end{example}

It will be useful when showing $X_\fpq$ is a cluster variable to understand how the posets associated to a Farey triple are related. In order to do this precisely, we define $\overline{\calP_\fpq}$ to be the poset $\calP_\fpq$ with the reverse chronological ordering.

\begin{lemma}\label{lem:ShapeOfPosetInFareyTriple}
Let $(\fpq,\frac{p+r}{q+s},\frs)$ be a Farey triple written in increasing order.
\begin{enumerate}
    \item If $\fpq > \frac01$ and $\frs < \frac10$, then the poset $\calP_{\frac{p+r}{q+s}}$ is of the form $\overline{\calP}_{\frs} \searrow \overline{H} \searrow \overline{\mathcal{P}}_{\fpq}$.
    Equivalently, $\check {\calP}_{\frac{p+r}{q+s}}$ is of the form $\Ptilde_{\fpq} \nearrow \calP_{\frs}$.
    \item If $\fpq = \frac01$, then $\calP_{\frac{p+r}{q+s}}$ is of the form $\Ptilde_{\frac01} \nearrow \overline{\calP}_{\frs}$.
\end{enumerate}
\end{lemma}

\begin{proof}
To prove item (1), note that if $(\fpq, \frac{p+r}{q+s},\frs)$ is a Farey triple, then there are no interior lattice points in the triangle with vertices $(0,0), (r,s),$ and $(q+s,p+r)$.
The proof comes from the fact that $\gamma_{\frac{p+r}{q+s}}$ will follow first $\gamma_{\frs}$ and then $\gamma_{\fpq}$, where the added chain corresponds to when $\gamma_{\frac{p+r}{q+s}}$ crosses arcs which are incident to $(s,r)$. The reversal of the posets comes from the fact that $\gamma_{\frac{p+r}{q+s}}$ lies to the right of $\gamma_{\fpq}$ and to the right of the line segment $(s,r) - (q+s,p+r)$, which is equivalent to $\gamma_{\fpq}$, so the middle crossings from the smaller posets are now to the right.
Item (2) can be seen directly from the pattern in the posets $\calP_{\frac{1}{s+1}}$.
\end{proof}

The relationship amongst posets $\Ptilde_\fpq$ in Lemma \ref{lem:ShapeOfPosetInFareyTriple} can be seen as a poset-theoretic version of the relationship between Christoffel words associated to a Markov triple; see Remark \ref{rem:Christoffel}. The symmetry in the posets $\calP_{\frac{1}{s+1}}$ was explored explicitly in \cite[Section 4.2]{banaian2024generalization} using the language of continued fractions.
The $g$-vectors associated to a cluster also have a nice relationship, implied by Lemmas \ref{lem:CompareGVecAndGCirc} and \ref{lem:ShapeOfPosetInFareyTriple}.

\begin{corollary}\label{cor:g_vecs_in_cluster}
Let $(\fpq,\frac{p+r}{q+s},\frs)$ be a Farey triple.  The $g$-vectors of the associated cluster are related by \[
\gg_{\frac{p+r}{q+s}}= \gg_{\fpq} + \gg_\frs + \begin{bmatrix}
    1\\-1\\-1
\end{bmatrix}= \gg^\circ_{\fpq} + \gg_\frs = \gg_{\fpq} + \gg_\frs^\circ,
\]
and in particular, \[
\gg^\circ_{\frac{p+r}{q+s}} = \gg^\circ_{\frac{p}{q}} + \gg_\frs^\circ.
\]
\end{corollary}

\begin{example}\label{ex:PosetsInFareyTriple}
Suppose $k_1 = 0$ and $k_2$ and $k_3$ are nonzero. Then, $\calP_\frac11$ is a single element, labeled $x_1$. In Example \ref{ex:Construction}, we drew $\calP_{\frac12}$ assuming $k_1$ was also nonzero. If $k_1 = 0$, this poset is instead the following.
\begin{center}
\begin{tikzpicture}
\node(1) at (0,0){$x_1$};
\node(2) at (1,-1){$x_2$};
\node(3) at (2,0){$x_2$};
\node(4) at (3,1){$x_1$};
\draw(1) -- (2);
\draw(2) -- (3);
\draw(3) -- (4);
\end{tikzpicture}
\end{center}

The numbers $\frac12$ and $\frac11$ sit together in a Farey triple, $(\frac12, \frac23, \frac11)$. We calculated $\calP_{\frac23}$ in Example \ref{ex:Construction}, and we can observe that it is of the form $\overline{\calP}_{\frac11} \searrow \overline{H} \searrow \overline{\calP}_{\frac12}$.

\begin{center}
\begin{tikzpicture}
\filldraw [rounded corners=10pt, pink!20]
(-1.4,.1)
-- (.4,.1)
-- (.4,.9)
-- (-1.4,.9)
-- cycle;
\filldraw [rounded corners=10pt, green!15]
(-.3,-1.4)
-- (1.7,-1.4)
-- (6.2,3.4)
-- (4.4,3.4)
-- cycle;
\filldraw [rounded corners=10pt, yellow!20]
(5.6,1.8)
-- (5.6,1)
-- (7.6,-.9)
-- (9.5,-.9)
-- (11,.9)
-- (8.5,.9)
-- (7.5,2)
-- cycle;
\node(0) at (-.5,.5){$(x_1,\widehat{x_1}^2)$};
\node(1) at (1,-1){$(x_2,2\widehat{x_2})$};
\node(2) at (2,0){$(x_2,\frac12\widehat{x_2})$};
\node(3) at (3,1){$(x_1,\widehat{x_1}^2)$};
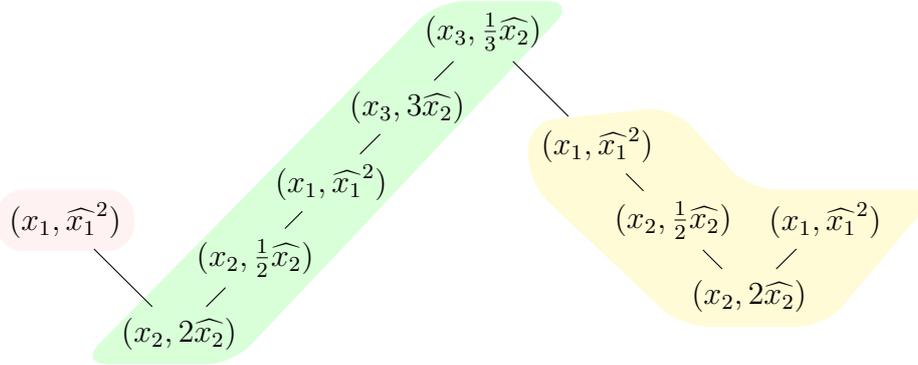
\node(4) at (4,2){$(x_3,3\widehat{x_2})$};
\node(5) at (5,3){$(x_3,\frac13\widehat{x_2})$};
\node(6) at (6.5,1.5){$(x_1,\widehat{x_1}^2)$};
\node(7) at (7.5,0.5){$(x_2,\frac12\widehat{x_2})$};
\node(8) at (8.5,-0.5){$(x_2,2\widehat{x_2})$};
\node(9) at (9.5,0.5){$(x_1,\widehat{x_1}^2)$};
\draw(0) -- (1);
\draw(1) -- (2);
\draw(2) -- (3);
\draw(3) -- (4);
\draw(4) -- (5);
\draw(5) -- (6);
\draw(6) -- (7);
\draw(7)--(8);
\draw(8) -- (9);
\end{tikzpicture}
\end{center}

We can also compare the $g$-vectors.  These are $\gg_{\frac11} = (-1,0,2)^\top,  \gg_{\frac12} = (0,-1,2)^\top$,  and $\gg_{\frac23} = (0,-2,3)^\top$,  and we indeed see $\gg_{\frac23} = \gg_{\frac11} + \gg_{\frac12} + (1,-1,-1)^\top$.

\end{example}

\subsection{Cluster Variables}
Here, we will show that each $X_\fpq$ is in fact the cluster variable $x_\fpq$ in the (possibly generalized) cluster algebra $\mathcal{A}^\prin(k_1,k_2,k_3)$ labeled by $\fpq$.
This is known when $k_1 = k_2 = k_3 = 0$, i.e. when we have an ordinary cluster algebra, by \cite{pilaud2023posets}.
By setting $y_1=y_2=y_3 = 1$, we also will get an interpretation of cluster variables in $\mathcal{A}(k_1,k_2,k_3)$. Note that we can extend our fractional labeling from $\mathcal{A}(k_1,k_2,k_3)$ to $\mathcal{A}^\prin(k_1,k_2,k_3)$ using the canonical isomorphism between their cluster patterns.

Our proof of this result will rely on the separation formula in Theorem \ref{thm:SeparationFormulaNakanishi}. As we have already discussed $g$-vectors in Lemma \ref{lem:CompareGVector}, we will here focus on $F$-polynomials and the polynomials $\mathcal{W}(\calP)$.

The separation formula has been used in the past when discussing these combinatorial models; see for instance \cite[Lemma 11.4]{Schemes} and \cite[Theorem 10]{banaian2023snake}. Here, our main technique is to show that the combinatorial polynomials $\mathcal{W}(\calP)$ satisfy mutation-like formulas.  Such an idea was used by Schiffler to verify the $T$-path formulas for cluster variables in Type A \cite{SchifflerTypeATPath}.

We begin by specifying how $F$-polynomial mutation behaves in our setting. Throughout this section, let $\hat{y}_i = y_i x_1^{b_{1i}} x_2^{b_{2i}} x_3^{b_{3i}}$ and let $F_\fpq = F_{1,\fpq}(\hat{y}_1,\hat{y}_2,\hat{y}_3)$. Notice that $\hat{y}_i$ is equal to the weight we associate to an element of a poset labeled by $x_i$. Let $d_\fpq$ be the degree of the mutation polynomial associated to $x_{1,\fpq}$. Given a Farey triple $(\fpq,\frac{p+r}{q+s},\frs)$ written in increasing order, define  \[
\cc_{\frac{p}{q}} =\begin{cases} \bigg(\frac{d_1 (p+q+1)}{d_{\frac{p}{q}}} ,\frac{d_2 (q+1)}{d_{\frac{p}{q}}},  \frac{d_3 (p+1)}{d_{\frac{p}{q}}}\bigg)^\top & \frs \neq \frac10\\
\bigg(\frac{d_1 (p+2)}{d_{\frac{p}{q}}} ,0,  \frac{d_3 (p+1)}{d_{\frac{p}{q}}}\bigg)^\top & \frs = \frac10.\end{cases}
\]

\begin{lemma}\label{lem:MutateFPoly}
Let $(\fpq,\frac{p+r}{q+s},\frs)$ be a Farey triple written in increasing order.
If $d_\fpq = 1$, then \begin{align*}
F_{\frac{p+2r}{q+2s}}F_{\frac{p}{q}} = F_{\frac{p+r}{q+s}}^2  +  \hat{\yy}^{\cc_{\frac{p}{q}}}  F_{\frac{r}{s}}^2
\end{align*}
and if $d_{\frac{p}{q}} = 2$,   then \begin{align*}
F_{\frac{p+2r}{q+2s}}F_{\frac{p}{q}} = F_{\frac{p+r}{q+s}}^2 + k_\frac{p}{q} \hat{\yy}^{\cc_{\frac{p}{q}}} F_{\frac{p+r}{q+s}} F_{\frac{r}{s}} +  \hat{\yy}^{2\cc_{\frac{p}{q}}}  F_{\frac{r}{s}}^2.
\end{align*}
\end{lemma}

\begin{proof}
The mutation of $\mathbf{x}_{1,\fpq}$ in the cluster $(\mathbf{x}_{1,\fpq},\mathbf{x}_{1,\frac{p+r}{q+s}},\mathbf{x}_{1,\frs})$ produces $\mathbf{x}_{1,\frac{p+2r}{q+2s}}$. We aim to show that in each case, the righthand side is exactly the formula from the alternate definition of $F$-polynomials in Definition \ref{def:FPoly}. Therefore, we must study the $c$-vector and  the column of the $B$-matrix  associated to $\mathbf{x}_{1,\fpq}$  at the seed associated to this Farey triple.

Recall Proposition \ref{prop:MarkovCVectors} classified $c$-vectors for the ordinary Markov cluster algebra and Lemma \ref{lem:CVectorsGenMarkov}  extended the classification to the generalized setting. Combining these results, we see that in each case $\cc_{\frac{p}{q}}$ is the $c$-vector $c_{\fpq,1,\mathcal{F}}$ where $\mathcal{F}$ here denotes the Farey triple $(\fpq,\frac{p+r}{q+s},\frs)$. Next, recall that from Lemma \ref{lem:BMatrixUsingLabeling}, we know $b_{\frs,\fpq} > 0$ and $b_{\frac{p+r}{q+s},\fpq} < 0$. With these observations, we see that the provided formulas are coincide with from  Definition \ref{def:FPoly}.

\end{proof}

We now compare the weight polynomials associated to a pair of adjacent Farey triples. We will see that, when we are replacing a fraction $\fpq$ with $d_\fpq = 2$, the process requires two steps. This is a common phenomenon in the theory of generalized cluster algebras, which reflects the fact that generalized cluster mutation can be seen as a sequence of ordinary mutations  \cite{chekhov2014teichmuller}. This second step can be understood as resolving a self-intersection, see \cite[Remark 7.2]{BanaianKelley}.

\begin{lemma}\label{lem:MutateWPoly}
Let $(\fpq,\frac{p+r}{q+s},\frs)$ be a Farey triple written in increasing order.
If $d_\fpq = 1$, then \begin{align*}
\mathcal{W}(\calP_{\frac{p+2r}{q+2s}})\mathcal{W}(\calP_{\frac{p}{q}}) = \mathcal{W}(\calP_{\frac{p+r}{q+s}})^2  +  \hat{\yy}^{\cc_{\frac{p}{q}}}  \mathcal{W}(\calP_{\frac{r}{s}})^2
\end{align*}
and if $d_{\frac{p}{q}} = 2$,   then \begin{align*}
\mathcal{W}(\calP_{\frac{p+2r}{q+2s}})\mathcal{W}(\calP_{\frac{p}{q}}) = \mathcal{W}(\calP_{\frac{p+r}{q+s}})^2 + k_\frac{p}{q} \hat{\yy}^{\cc_{\frac{p}{q}}} \mathcal{W}(\calP_{\frac{p+r}{q+s}}) \mathcal{W}(\calP_{\frac{r}{s}}) +  \hat{\yy}^{2\cc_{\frac{p}{q}}}  \mathcal{W}(\calP_{\frac{r}{s}})^2.
\end{align*}
\end{lemma}

\begin{proof}
Assume first that $\fpq \neq \frac01$ and $\frs \neq \frac10$. From Lemma \ref{lem:ShapeOfPosetInFareyTriple}, we know $\calP_{\frac{p+2r}{q+2s}}$ is of the form $\overline{\mathcal P}_{\frac{r}{s}} \searrow \overline{H} \searrow \overline{\mathcal P}_{\frac{p+r}{q+s}}$. Using this Lemma again on $\overline{\mathcal P}_{\frac{p+q}{r+s}}$ shows us that $\calP_{\frac{p+2r}{q+2s}}$ is of the following form where the chains $x_2 \cdots x_3$ represent, from left to right,  $\overline{H}$ and $H$.

\begin{center}
\begin{tikzpicture}
\node(pq1) at (0,0){$\overline{\calP}_{\frac{r}{s}}$};
\node(2) at (1,-1){$x_2$};
\node(dots1) at (2,0){$\iddots$};
\node(3) at (3,1){$x_3$};
\node(rs) at (4,0){$\calP_{\frac{p}{q}}$};
\node(33) at (5,1){$x_3$};
\node(dots2) at (6,0){$\ddots$};
\node(22) at (7,-1){$x_2$};
\node(pq2) at (8,0){$\calP_{\frac{r}{s}}$};
\draw(pq1) -- (2);
\draw(2) -- (dots1);
\draw(dots1) -- (3);
\draw(3) -- (rs);
\draw(rs) -- (33);
\draw(33) -- (dots2);
\draw(dots2) -- (22);
\draw(22) -- (pq2);
\end{tikzpicture}
\end{center}

The depiction shows us  $\calP_{\frac{p+2r}{q+2s}}$ and $\calP_{\frac{p}{q}}$ have a crossing overlap in $\calP_{\frac{p}{q}}$.  This means we can invoke Proposition \ref{prop:CrossingOverlap}, which implies
\begin{equation}\label{eq:FirstResolutionCase1}
\mathcal{W}^\prin(\calP_{\frac{p+2r}{q+2s}})\mathcal{W}^\prin(\calP_{\frac{p}{q}}) = \mathcal{W}^\prin(Q_1)\mathcal{W}^\prin(Q_2)  +
Z\mathcal{W}^\prin(\calP_{\frac{r}{s}})^2
\end{equation}
where $Z$ is the product of the weights of the elements on $\overline{H},\calP_{\fpq},$ and $H$, and the posets $Q_1$ and $Q_2$ are below on the left and right respectively.

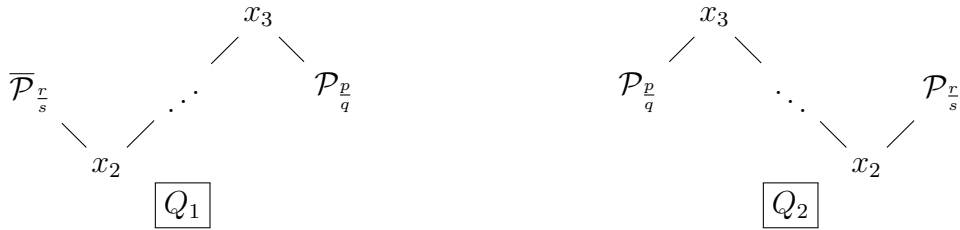
\begin{center}
\begin{tikzpicture}
\node(pq1) at (0,0){$\overline{\calP}_{\frac{r}{s}}$};
\node(2) at (1,-1){$x_2$};
\node(dots1) at (2,0){$\iddots$};
\node(3) at (3,1){$x_3$};
\node(rs) at (4,0){$\calP_{\frac{p}{q}}$};
\node(rs1) at (8,0){$\calP_{\frac{p}{q}}$};
\node(33) at (9,1){$x_3$};
\node(dots2) at (10,0){$\ddots$};
\node(22) at (11,-1){$x_2$};
\node(pq2) at (12,0){$\calP_{\frac{r}{s}}$};
\draw(pq1) -- (2);
\draw(2) -- (dots1);
\draw(dots1) -- (3);
\draw(3) -- (rs);
\draw(rs1) -- (33);
\draw(33) -- (dots2);
\draw(dots2) -- (22);
\draw(22) -- (pq2);
\node[] at (2,-1.5){$\boxed{Q_1}$};
\node[] at (10,-1.5){$\boxed{Q_2}$};
\end{tikzpicture}
\end{center}

By Lemma \ref{lem:ProductOfAllWeights}, the product of weights of elements along $H$ (equivalently, $\overline{H}$) is $\hat{y}_1^{d_1} \hat{y}_2^{d_2} \hat{y}_3^{d_3}$. Let  $\hat{Y}_\fpq:=\prod_{x \in \calP_{\fpq}} wt^\prin(x)$.  By Proposition \ref{prop:Computedfpq}, $\hat{Y}_\fpq = \hat{y}_1^{d_1(q+p-1)}\hat{y}_2^{d_2(q-1)}\hat{y}_3^{d_3(p-1)}$, so $Z = \hat{y}_1^{d_1(q+p+1)}\hat{y}_2^{d_2(q+1)}\hat{y}_3^{d_3(p+1)}$.

From Lemma \ref{lem:ShapeOfPosetInFareyTriple}, we recognize $Q_2$ as $\overline{\calP}_{\frac{p+r}{q+s}}$, so that $\mathcal{W}^\prin(Q_2) = \mathcal{W}^\prin(\calP_{\frac{p+r}{q+s}})$.
If $d_{\frac{p}{q}} = 1$, then $\calP_{\frac{p}{q}} = \overline{\calP}_{\frac{p}{q}}$, so $Q_1$ is identical to $\calP_{\frac{p+r}{q+s}}$ as well. In this case, we can see that the degree vector of $Z$ in the variables $\hat{y}_i$ is exactly $\cc_\fpq$, and overall we have reached the desired formula.

 Next, suppose that $d_\fpq = 2$, so that $\calP_\fpq$ has an even number of elements. Then, we can refine our depiction of $Q_1$, which illustrates that this poset has a reverse-kissing overlap in $\calP_\fpq[1,\frac12 \cdot h_\fpq]$ and $\calP_\fpq[\frac12 \cdot h_\fpq + 1,h_\fpq]$.

\begin{center}
\begin{tikzpicture}[scale=0.8]
\node(pq1) at (-0.5,0){$\overline{\calP}_{\frac{r}{s}}$};
\node(2) at (0.5,-1){$x_2$};
\node(dots1) at (1.5,0){$\iddots$};
\node(3) at (3,1.5){$x_3$};
\node(rs) at (4,0){$\calP_{\frac{p}{q}}(1)$};
\node(dots2) at (5,0){$\cdots$};
\node(hpq) at (6.5,0){$\calP_\fpq(\frac12 \cdot h_\fpq)$};
\node(hpq1) at (8.5,1.5){$\calP_\fpq(\frac12 \cdot h_\fpq+1)$};
\node(dots3) at (10.5,1.5){$\cdots$};
\node(end) at (11.8,1.5){$\calP_\fpq(h_\fpq)$};
\draw(pq1) -- (2);
\draw(2) -- (dots1);
\draw(dots1) -- (3);
\draw(3) -- (rs);
\draw(hpq) -- (hpq1);
\end{tikzpicture}
\end{center}

Applying Proposition \ref{prop:KissingReverseOverlap} to $\mathcal{W}^\prin(Q_1)$ yields
\[
\mathcal{W}^\prin(Q_1) = \mathcal{W}^\prin(Q_3)  + \hat{y}_1^{d_1}\hat{y}_2^{d_2}\hat{y}_3^{d_3}\big(\prod_{x \in \calP_\fpq[1,\frac12 \cdot h_\fpq]} wt^\prin(x) \big)\mathcal{W}^\prin(\calP_\frs)
\]

where $Q_3$ is the poset on the same set as $Q_1$ but with the relation between  $\calP_\fpq(\frac12 \cdot h_\fpq)$ and  $\calP_\fpq(\frac12 \cdot h_\fpq+1)$ reversed, i.e. $\calP_\fpq(\frac12 \cdot h_\fpq) > \calP_\fpq(\frac12 \cdot h_\fpq+1)$. Using Proposition \ref{prop:Computedfpq} again,  we recognize the degree vector of $\hat{y}_1^{d_1}\hat{y}_2^{d_2}\hat{y}_3^{d_3}\big(\prod_{x \in \calP_\fpq[1,\frac12 \cdot h_\fpq]} wt^\prin(x) \big)$ to match $\cc_\fpq$. From Lemma \ref{lem:ProductOfAllWeights}, we see that this monomial has a factor of $k_\fpq$.

 The poset $Q_3$ is the same as $\overline{\calP_\frs} \searrow \overline{H} \searrow \overline{\calP_\fpq}$, i.e., $\calP_{\frac{p+r}{q+s}}$, so $\mathcal{W}^\prin(Q_3) = \mathcal{W}^\prin(\calP_{\frac{p+r}{q+s}})$. Therefore, substituting this expression for $\mathcal{W}^\prin(Q_1)$ into Equation \eqref{eq:FirstResolutionCase1} gives us the desired formula.

Now, we assume $\frs = \frac10$. Here, we can assume $\fpq > \frac01$ since we can check the Farey triple $(\frac01, \frac11, \frac10)$ directly. Therefore, we are considering a initial Farey triple of the form $(\frac{p}{1},\frac{p+1}{1},\frac10)$ for $p \geq 1$, and replacing $\frac{p}{1}$ gives us $(\frac{p+1}{1},\frac{p+2}{1},\frac10)$. The poset $\calP_{\frac{p+2}{1}}$ is as below. We can see from this depiction that $\calP_{\frac{p+2}{1}}$ and $\calP_{\frac{p}{1}}$ have a crossing overlap in $\calP_{\frac{p}{1}}$. The proof in this case then follows the general case.

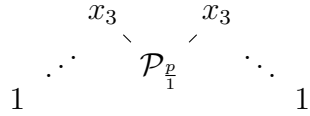
\begin{center}
\begin{tikzpicture}[scale = 0.75]
\node(1) at (0.5,0.5){$x_1$};
\node(dots1) at (1.25,1.25){$\iddots$};
\node(2) at (2,2){$x_3$};
\node(3) at (3,1){$\calP_\frac{p}{1}$};
\node(d) at (4,2){$x_3$};
\node(dots2) at (4.75,1.25){$\ddots$};
\node(e) at (5.5,0.5){$x_1$};
\draw(2) -- (3);
\draw(3) -- (d);
\end{tikzpicture}
\end{center}

Finally, we assume $\fpq = \frac01$. As mentioned before, we now can assume $\frs \neq \frac10$. Then, our original Farey triple is of the form $(\frac01, \frac{1}{r+1}, \frac1r)$ with $r \geq 1$ and replacing $\frac01$ gives us the triple $(\frac{1}{r+1},\frac{2}{2r + 1},\frac1r)$. In this special case, $\cc_\frac01 = (2\frac{d_1}{d_3},2\frac{d_2}{d_3},1)^\top$. We draw $\calP_{\frac{2}{2r+1}}$ in the case where $d_1 = d_2 = 1$ and $d_3 = 2$. In general, there are $d_3$ elements labeled $x_3$ and to either side of this element/these elements there are $d_1$ and $d_2$ elements labeled $x_1$ and $x_2$ respectivey.

\begin{center}
\begin{tikzpicture}[scale = 0.75]
\node(1) at (0,0){$\calP_\frac{1}{r+1}$};
\node(a) at (1,-1){$x_2$};
\node(c) at (2,0){$x_1$};
\node(d) at (3,1){$x_3$};
\node(e) at (4,2){$x_3$};
\node(f) at (5,1){$x_1$};
\node(g) at (6,0){$x_2$};
\node(2) at (7,1){$\calP_\frac{1}{r+1}$};
\draw(1) -- (a);
\draw(a) -- (c);
\draw(c) -- (d);
\draw(d) -- (e);
\draw(e) -- (f);
\draw(f) -- (g);
\draw(g) -- (2);
\end{tikzpicture}
\end{center}

If $d_{\frac01} = 2$, since $\calP_\frac01 = \emptyset$, we want to show \[
\mathcal{W}^\prin(\calP_{\frac{2}{2r + 1}}) = \mathcal{W}^\prin(\calP_{\frac{1}{r + 1}})^2 + k_3\hat{y}_1^{d_1}\hat{y}_2^{d_2}\hat{y}_3\mathcal{W}^\prin(\calP_{\frac{1}{r + 1}})\mathcal{W}^\prin(\calP_{\frac{1}{r}})^2 + \hat{y}_1^{2d_1}\hat{y}_2^{2d_2}\hat{y}_3^2\mathcal{W}^\prin(\calP_{\frac{1}{r}})^2.
\]

From left to right, we can realize these three terms as the weight generating functions of order ideals of $\calP_{\frac{2}{2r + 1}}$ which (1) do not contain either element labeled $x_3$, (2) contain the smaller but not the larger element labeled $x_3$, and (3)  contain both elements labeled $x_3$. When $d_{\frac01} = 1$, one can do a similar analysis.
\end{proof}

\begin{theorem}\label{thm:CorrectnessOfPosetFormula}
Let $\fpq \in \mathbb{Q} \cap [0,\infty)$. The cluster variable $x^{\prin}_{1,\frac{p}{q}}$ is given by \[
x^{\prin}_{1,\frac{p}{q}} = \mathbf{x}^{\gb_{\frac{p}{q}}} \mathcal{W}^{\prin}(\calP_{\frac{p}{q}}).
\]
\end{theorem}

\begin{proof}
By Theorem \ref{thm:SeparationFormulaNakanishi}, $x^{\prin}_{1,\frac{p}{q}}$ can be expressed in terms of its $g$-vector and $F$-polynomial.  Lemma \ref{lem:CompareGVector} has already shown that the $g$-vector associated to the cluster variable $x_{1,\fpq}$ matches the $g$-vector associated to the poset $\calP_\fpq$. Proposition 3.19 and  Theorem 3.20 in \cite{nakanishi2015structure}, show that, when using principal coefficients, $F$-polynomials are polynomials with constant term 1; therefore, we need not consider the denominator in Theorem \ref{thm:SeparationFormulaNakanishi}.

Therefore, we will be done once we show $F_\fpq = \mathcal{W}^\prin(\calP_\fpq)$.
We do this inductively, checking the base cases directly.

 Let $(\frac{p}{q},\frac{p+r}{q+s},\frs)$ be a Farey triple such that we already know $F_{\frac{p}{q}} = \mathcal{W}^\prin(\calP_{\frac{p}{q}})$,  $F_{\frac{p+r}{q+s}} = \mathcal{W}^\prin(\calP_{\frac{p+r}{q+s}})$, and $F_{\frac{r}{s}} = \mathcal{W}^\prin(\calP_{\frac{r}{s}})$.
 If we mutate at $x^{\prin}_{1,\frac{p+r}{q+s}}$, we will reach a cluster closer to the initial cluster, and by induction we also know the desired statement for the resulting cluster variable. Therefore, we need only consider mutating at $x^{\prin}_{1,\frac{p}{q}}$ and at $x^{\prin}_{1,\frac{r}{s}}$. Lemmas \ref{lem:MutateFPoly} and \ref{lem:MutateWPoly} together show that, if $d_\fpq = 1$, then \[
 F_{\frac{p+2r}{q+2s}} = \frac{F_{\frac{p+r}{q+s}}^2  +  \hat{\yy}^{\cc_{\frac{p}{q}}}  F_{\frac{r}{s}}^2}{F_\fpq} = \frac{\mathcal{W}(\calP_{\frac{p+r}{q+s}})^2  +  \hat{\yy}^{\cc_{\frac{p}{q}}}  \mathcal{W}(\calP_{\frac{r}{s}})^2}{\mathcal{W}(\calP_\fpq)} = \mathcal{W}(\calP_\frac{p+2r}{q+2s}) ,
 \]
 and similarly if $d_\fpq = 2$. A parallel argument would show the statement for $F_{\frac{2p+r}{2q+d}}$.
\end{proof}

By substituting $y_1=y_2=y_3=1$, we have the following corollary of the trivial coefficient case.

\begin{corollary}\label{cor:CorrectnessOfPosetFormula-trivial}
Let $\fpq \in \mathbb{Q} \cap [0,\infty)$. Let $x_{1,\frac{p}{q}}$ be the cluster variable in $k_1 \mathbb{T}(k_1,k_2,k_3)$ labeled by $\frac{p}{q}$ in $\mathcal{A}(k_1,k_2,k_3)$. Then,\[
x_{1,\frac{p}{q}} = \mathbf{x}^{\gb_{\frac{p}{q}}} \mathcal{W}(\calP_{\frac{p}{q}})
\]
\end{corollary}

This form will be used frequently in the discussion that follows.

\begin{remark}\label{rem:OtherBranches}
 Theorem \ref{thm:CorrectnessOfPosetFormula} and Corollary \ref{cor:CorrectnessOfPosetFormula-trivial} only addressed cluster variables of the form $x_{1,\fpq}$. Recall these are the cluster variables reachable from the result of mutating the initial cluster at position 1. One can easily extend this to reach all cluster variables by applying permutations to the indices.
\end{remark}

It is particularly nice to work in the case $k:=k_1 = k_2 = k_3$. In this setting, the coefficient in the formulas in Lemmas \ref{lem:MutateFPoly} and \ref{lem:MutateWPoly} does not depend on $\fpq$. If we furthermore set $x_1 = x_2 = x_3 = 1$, we recover ``$k$-Markov numbers''. These have been conjectured to satisfy a $k$ version of Frobenius's Uniqueness conjecture \cite{gyoda2023uniqueness}.  It is hopeful that Theorem \ref{thm:CorrectnessOfPosetFormula} could aid in making progress towards this conjecture.

\section{Combinatorial Cluster Generalized Cohn Matrices}\label{sec:CombinatorialCohnMatrices}

To each positive rational number, we associate a $2 \times 2$ matrix whose entries are in $\mathbb{Z}[x_1^\pm,x_2^\pm,x_3^\pm]$. We discuss the absence of $y$-variables in Remark \ref{rmk:YVar}. This definition was inspired by the $q$-deformation of Cohn matrices given in \cite{ouguz2025oriented}.

\begin{definition}
Given $\fpq \in \mathbb{Q} \cap [0,\infty)$, define $C_\fpq$ to be the matrix given by
\[
C_{\frac{p}{q}} = \begin{bmatrix} \xx^{\gg^\circ_{\frac{p}{q}}} \mathcal{W}(\Ptilde_{\frac{p}{q}}; R) & \xx^{\gg_{\frac{p}{q}}} \mathcal{W}(\Ptilde_{\frac{p}{q}}; \neg R)\\
\xx^{2\gg^\circ_\fpq - \gg_\fpq} \mathcal{W}(\Ptilde_{\frac{p}{q}}; R, \neg L)& \xx^{\gg^\circ_\fpq} \mathcal{W}(\Ptilde_{\frac{p}{q}}; \neg R, \neg L)\end{bmatrix}.
\]

We also define $C_{\frac10}$ as follows. Let $\calP_{\frac10}$ be the poset associated to the line segment $(0,1) - (2,0)$ and let $\Ptilde_{\frac10}$ be $\calP_{\frac10} \searrow H$. The poset $\Ptilde_{\frac10}$ when all $k_i > 0$ is drawn in Figure \ref{fig:Ptilde10}. We define $C_{\frac10}$ in terms of $\Ptilde_{\frac10}$.

\[
C_{\frac10} = \begin{bmatrix} \frac{x_1}{x_3} \mathcal{W}(\Ptilde_{\frac10}; R, \neg L) & x_2 \mathcal{W}(\Ptilde_{\frac10}; \neg L, \neg R) \\
 -\frac{x_1^2}{x_2x_3^2} \mathcal{W}(\Ptilde_{\frac10}; L,R) &  -\frac{x_1}{x_3}  \mathcal{W}(\Ptilde_{\frac10}; L, \neg R)
\end{bmatrix}
\]
\end{definition}

Note that $\mathcal{W}(\Ptilde_{\frac10}; \neg L, \neg R) = wt(\emptyset) = 1$ so that the $(1,2)$-entry of $C_{\frac10}$ is simply $x_2$. Moreover, the $(2,2)$ entry is simply $-k_2$. We write the matrix in the language of posets simply to demonstrate how its structure compares to the other combinatorial cluster generalized Cohn matrices.

While the matrix $C_{\frac10}$ has a different structure than the other combinatorial cluster generalized Cohn matrices, we can quickly verify some of its properties from the following lemma, which can be observed directly.

\begin{lemma}\label{lem:C10}
The matrix $C_{\frac10}$ is equal to $Q_{(1,3)}(f)$ where $f = \frac{x_1}{x_3} \bigg(1 +  k_2 \widehat{x_2} + \widehat{x_2}^2 + k_1\widehat{x_1}\widehat{x_2}^2 + \widehat{x_1}^2\widehat{x_2}^2 + k_3 \widehat{x_1}^2\widehat{x_2}^2\widehat{x_3}\bigg)$.
\end{lemma}
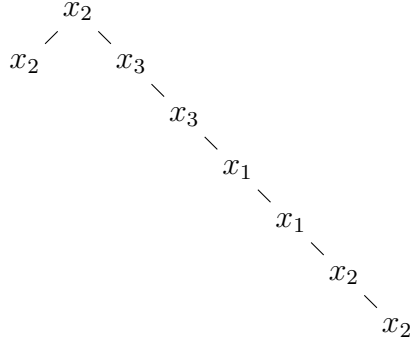
\begin{figure}
    \centering
\begin{tikzpicture}[scale = 0.7]
\node(1) at (0,0){$x_2$};
\node(2) at (1,1){$x_2$};
\node(3) at (2,0){$x_3$};
\node(4) at (3,-1){$x_3$};
\node(5) at (4,-2){$x_1$};
\node(6) at (5,-3){$x_1$};
\node(7) at (6,-4){$x_2$};
\node(8) at (7,-5){$x_2$};
\draw(1) -- (2);
\draw(2) -- (3);
\draw(3) -- (4);
\draw(4) -- (5);
\draw(5) -- (6);
\draw(6) -- (7);
\draw(7) -- (8);
\end{tikzpicture}
    \caption{The poset $\Ptilde_{\frac10}$ when all $k_i > 0$. Note that if $k_2 = 0$, then this is a chain.}
    \label{fig:Ptilde10}
\end{figure}

The goal of this section is to prove the following. The key to our proofs will be to use the machinery developed in \cite{ouguz2025oriented} for similar types of matrices.

\begin{theorem}\label{thm:Combinatorial_Cohn_Matrix_and_Triple}
For each $\fpq \in \mathbb{Q} \cap [0,\infty)$ or for $\fpq = \frac10$,  $C_\fpq$ is a $(k_1,k_2,k_3)$-CGC matrix, and for each Farey triple $(\frac{p}{q}, \frac{p+r}{q+s},\frs)$ in $\mathrm{F}\mathbb T$,  the triple $(C_{\frac{p}{q}}, C_{\frac{p+r}{q+s}},C_{\frs})$ is a $(k_1,k_2,k_3)$-CGC triple.
\end{theorem}

To this end, we refer to the family of matrices $C_\fpq$ as \emph{combinatorial cluster generalized Cohn matrices} and triples $(C_{\frac{p}{q}}, C_{\frac{p+r}{q+s}},C_{\frs})$ associated to a Farey triple as \emph{combinatorial cluster generalized Cohn triples}.
We divide the verification that combinatorial cluster generalized Cohn matrices and triples satisfy  Definitions \ref{def:gen-Cohn-matrix} and \ref{def:gen-Cohn-triple} into a series of intermediate results. The first Proposition is proven using poset skein relations.

\begin{proposition}\label{prop:Combinatorial_Cohn_determinant_1}
For each $\fpq\in \mathbb{Q} \cap [0,\infty)$,  $\det(C_\fpq) = 1$.
\end{proposition}

\begin{proof}
Notice the determinant of $C_\fpq$ is given by \[
\xx^{2\gg^\circ_\fpq}\big(\mathcal{W}(\Ptilde_{\frac{p}{q}}; R)\mathcal{W}(\Ptilde_{\frac{p}{q}}; \neg L, \neg R) - \mathcal{W}(\Ptilde_{\frac{p}{q}}; \neg R)\mathcal{W}(\Ptilde_{\frac{p}{q}}; R, \neg L)\big)
\]
so it suffices to show that $\mathcal{W}(\Ptilde_{\frac{p}{q}}; R)\mathcal{W}(\Ptilde_{\frac{p}{q}}; \neg L, \neg R) - \mathcal{W}(\Ptilde_{\frac{p}{q}}; \neg R)\mathcal{W}(\Ptilde_{\frac{p}{q}}; R, \neg L)$ is equal to $\xx^{-2\gg^\circ_\fpq}$.

We begin by considering the exceptional case when $\Ptilde_{\frac{p}{q}}$ is a chain, which occurs only when $\frac{p}{q} = \frac01$. Recall we drew an example in Example \ref{ex:Construction}. In this case, the chain is decreasing with respect to the chronological ordering. This implies $L$ is the maximum and $R$ is the minimum. Since only the empty order ideal does not contain $R$, we have $\mathcal{W}(\Ptilde_{\frac{0}{1}}; \neg R) = \mathcal{W}(\Ptilde_{\frac{p}{q}}; \neg L, \neg R) =1$. Moreover, one can quickly see $\mathcal{W}(\Ptilde_{\frac{p}{q}}; R) -\mathcal{W}(\Ptilde_{\frac{p}{q}}; R, \neg L)$ is equal to the weight of the maximal order ideal, i.e. $I = \Ptilde_{\frac{0}{1}}$. This weight is $\widehat{x_1}^2\widehat{x_2}^2 = \frac{x_2^2}{x_1^2} = \mathbf{x}^{(2,-2,0)^\top}$. Since $(2,-2,0) = -2 \gg^\circ_{\frac{0}{1}}$, we are done.

Now, for the rest of the proof, let $\fpq > 0$ and let $u>1$ be the smallest number such that $\Ptilde_\fpq(u) \nless \Ptilde_\fpq(1)$.

Since $\Ptilde_\fpq(\widetilde{h}_\fpq) < \Ptilde_\fpq(\widetilde{h}_\fpq - 1)$,  we can rewrite our polynomial as \[
\widehat{x_2}\mathcal{W}(\Ptilde_{\frac{p}{q}}[1,\widetilde{h}_\fpq - 1])\mathcal{W}(\Ptilde_{\frac{p}{q}}[u,h_\fpq]) -\widehat{x_2}\mathcal{W}(\Ptilde_{\frac{p}{q}}[1,h_\fpq])\mathcal{W}(\Ptilde_{\frac{p}{q}}[u,\widetilde{h}_\fpq-1]).
\]
The posets $\Ptilde_{\frac{p}{q}}[u,h_\fpq]$ and $\Ptilde_{\frac{p}{q}}[1,\widetilde{h}_\fpq - 1]$ have a crossing overlap in $\Ptilde_{\frac{p}{q}}[2,h_\fpq]$. This is because this subposet is trivially on top in $\Ptilde_{\frac{p}{q}}[u,h_\fpq]$ and since $\widetilde{P_\fpq}(h_\fpq) < \widetilde{P_\fpq}(h_\fpq+1)$,  it is on bottom in $\Ptilde_{\frac{p}{q}}[1,\widetilde{h}_\fpq - 1]$.  Using Proposition \ref{prop:CrossingOverlap}, we have \[
\widehat{x_2}\mathcal{W}(\Ptilde_{\frac{p}{q}}[1,\widetilde{h}_\fpq - 1])\mathcal{W}(\Ptilde_{\frac{p}{q}}[u,h_\fpq]) - \widehat{x_2} \mathcal{W}(\Ptilde_{\frac{p}{q}}[1,h_\fpq])\mathcal{W}(\Ptilde_{\frac{p}{q}}[u,\widetilde{h}_\fpq-1]) = \widehat{x_2} \prod_{i=2}^{\widetilde{h}_\fpq} wt(\Ptilde_\fpq)(i),
\]
and by Lemma \ref{lem:ProductOfAllWeights} and the fact that $\widehat{x}_1\widehat{x}_2\widehat{x}_3 = 1$, the righthand side is $\xx^{-2\gg^\circ_\fpq}$,  as desired.

\end{proof}
In several proofs,  we will find it convenient to express $\mathscr M $ in terms of a weight-generating function of a weighted poset.  We record this fact for reference.

\begin{lemma}\label{lem:RewriteMinv}
The expression $\mathscr M $ is equal to both $\frac{x_1}{x_2x_3} \mathcal{W}(H; \neg L)$ and $\frac{x_1}{x_2x_3} \mathcal{W}(H; R)$. Explicitly, $\mathscr M $  can be rewritten as \begin{align*}
\mathscr M  &= \frac{x_1}{x_2x_3}(1 + k_2 \widehat{x_2} + \widehat{x_2}^2 + k_1\widehat{x_1}\widehat{x_2}^2 + \widehat{x_1}^2\widehat{x_2}^2 + k_3 \widehat{x_1}^2\widehat{x_2}^2\widehat{x_3})\\
&=  \frac{x_1}{x_2x_3}(k_2 \widehat{x_2} + \widehat{x_2}^2 + k_1\widehat{x_1}\widehat{x_2}^2 + \widehat{x_1}^2\widehat{x_2}^2 + k_3 \widehat{x_1}^2\widehat{x_2}^2\widehat{x_3} + \widehat{x_1}^2\widehat{x_2}^2\widehat{x_3}^2 )
\end{align*}
\end{lemma}

\begin{proof}
Verifying each of these is a straightforward computation. In particular, the last equality is immediate from the fact that $\widehat{x_1}\widehat{x_2}\widehat{x_3} = 1$.
\end{proof}

\begin{definition}\label{def:alphaShorthand}
For shorthand in the ensuing proofs, let $\alpha:= \mathcal{W}(H;R, \neg L)$. Explicitly, $\alpha = k_2 \widehat{x_2} + \widehat{x_2}^2 + k_1\widehat{x_1}\widehat{x_2}^2 + \widehat{x_1}^2\widehat{x_2}^2 + k_3 \widehat{x_1}^2\widehat{x_2}^2\widehat{x_3}$.
\end{definition}

We highlight some ways in which $\alpha$ will simplify expressions.

\begin{lemma}\label{lem:alphaSimplifying}
Given $\alpha$ as in Definition \ref{def:alphaShorthand}, we have \[
\mathscr M  = \frac{x_1}{x_2x_3}(1 + \alpha ),
\]
\[
\mathcal{W}(\Ptilde_{\fpq}) = (1+\alpha)\mathcal{W}(\calP_{\fpq}) + \mathcal{W}(\calP_{\fpq};R),
\]
and
\[
\mathcal{W}(\Ptilde_{\fpq};R) = \alpha\mathcal{W}(\calP_{\fpq}) + \mathcal{W}(\calP_{\fpq}).
\]
\end{lemma}

The following Proposition, verifying that the matrix $C_\fpq$ satisfies (2) and (3) from Definition \ref{def:gen-Cohn-matrix}, can be proven fairly directly.

\begin{proposition}\label{prop:Combinatorial_Cohn_Trace_Formula}
For each $\fpq\in \mathbb{Q} \cap [0,\infty)$,  let $p_{12}$ denote the $(1,2)$-entry of $C_\fpq$. We have $\text{tr}(C_\fpq) = \mathscr M  p_{12} - k_\fpq$.
\end{proposition}

\begin{proof}
We first consider the $(1,1)$ entry of $C_\fpq$.  To simplify notation, we assume $k_2 > 0$.  This entry can be written as \[
\xx^{\gg^\circ_\fpq} \mathcal{W}(\Ptilde_\fpq;R_{\Ptilde_\fpq}) = k_2\widehat{x_2}\xx^{\gg^\circ_\fpq} \mathcal{W}(\Ptilde_\fpq[1,\widetilde{h}_\fpq-1]).
\]
We partition the order ideals of $\Ptilde_\fpq[1,\widetilde{h}_\fpq-1]$ into two sets, based on whether they contain $\Ptilde_\fpq(h_\fpq+1)$.  Note that this element is maximal.  The contribution from order ideals which do not contain $\calP_\fpq(h_\fpq+1)$ is \[
k_2\widehat{x_2}\xx^{\gg^\circ_\fpq}\sum_{\substack{I \in J(\calP_\fpq[1,\widetilde{h}_\fpq-1])\\ \calP_\fpq(h_\fpq+1) \notin I}} wt(I)= \xx^{\gg^\circ_\fpq} \alpha \mathcal{W}(\calP_\fpq)
\]
and the contribution from order ideals which contain $\calP_\fpq(h_\fpq+1)$ is \[
k_2\widehat{x_2}\xx^{\gg^\circ_\fpq}\sum_{\substack{I \in J(\calP_\fpq[1,\widetilde{h}_\fpq-1])\\ \calP_\fpq(h_\fpq+1) \in I}} wt(I)= \xx^{\gg^\circ_\fpq}  \widehat{x_1}^2\widehat{x_2}^2\widehat{x_3}^2 \mathcal{W}(\calP_\fpq;R_{\calP_\fpq})=  \xx^{\gg^\circ_\fpq}  \mathcal{W}(\calP_\fpq;R_{\calP_\fpq})
\]
where we use the fact that  $ \widehat{x_1}\widehat{x_2}\widehat{x_3} = 1$ to simplify the right hand side. We remark that $R_{\calP_\fpq} = \calP_\fpq(h_\fpq)$.

Next, we consider the (2,2) entry, which is $ \xx^{\gg^\circ_\fpq}\mathcal{W}(\Ptilde_{\frac{p}{q}}; \neg R_{\Ptilde_{\frac{p}{q}}}, \neg L_{\Ptilde_{\frac{p}{q}}}) $. By the definition of $\Ptilde_\fpq$, this Laurent polynomial is equal to $ \xx^{\gg^\circ_\fpq}\mathcal{W}(\calP_{\frac{p}{q}}; \neg L_{\calP_{\frac{p}{q}}}) $,  and from Lemma \ref{lem:Asymmetry} this is equal to $ \xx^{\gg^\circ_\fpq}(\mathcal{W}(\calP_{\frac{p}{q}}; \neg R_{\calP_{\frac{p}{q}}}) )- k_\fpq $.

Now,  we have \begin{align*}
\text{tr}(C_\fpq) &= \xx^{\gg^\circ_\fpq}\big(  \mathcal{W}(\Ptilde_\fpq;R_{\Ptilde_\fpq}) + \mathcal{W}(\Ptilde_{\frac{p}{q}}; \neg R_{\Ptilde_\fpq}, \neg L_{\Ptilde_\fpq})\big)\\
&=  \xx^{\gg^\circ_\fpq}\big(\alpha \mathcal{W}(\calP_\fpq)+  \mathcal{W}(\calP_\fpq;R_{\calP_\fpq}) + \mathcal{W}(\calP_\fpq;\neg R_{\calP_\fpq})\big) - k_\fpq\\
&=\frac{x_1}{x_2x_3} \xx^{\gg_\fpq}(1+ \alpha)\mathcal{W}(\calP_\fpq) - k_\fpq\\
&= \mathscr M  \xx^{\gg_\fpq} \mathcal{W}(\calP_\fpq) - k_\fpq
\end{align*}

where in the second-to-last equality,  we use  Lemma \ref{lem:CompareGVecAndGCirc} and in the last equality we use Lemma \ref{lem:RewriteMinv}.  If $k_2 = 0$,  then we follow a similar proof but use the fact that the weight of the element $R$ is not $\widehat{x_2}^2$.
\end{proof}

As a corollary of Proposition \ref{prop:Combinatorial_Cohn_Trace_Formula},  we show how the cluster variable $X_\fpq$ relates to the corresponding the variable $X^\circ_\fpq$. These latter Laurent polynomials are related to \emph{band graphs} as in \cite{musiker2013bases}, which associates Laurent polynomials to closed curves in a surface. In the ordinary case ($k_1=k_2=k_3=0$),  with the initial cluster variables set to 1, this recovers the well-known fact that the length of the closed curve of slope $\fpq$ on the once-punctured torus is three times the length of the associated arc (for example, see \cite[Theorem 4.9]{canakci2020snake}).  The case $k_1=k_2=k_3=1$ with the initial cluster set to 1 was given in \cite[Theorem 8]{banaian2024generalization} using a combinatorial proof involving snake graphs.

\begin{corollary}\label{cor:Band_Graph_Formula_Conserved}
For each $\fpq \in \mathbb{Q} \cap [0,\infty)$,  we have
\[
X_\fpq^\circ = \mathscr M  x_\fpq - k_\fpq.
\]
\end{corollary}

\begin{proof}
Recall we form $\calP^\circ_\fpq$ by adding the relation $L_{\Ptilde_\fpq} > R_{\Ptilde_\fpq}$ to $\Ptilde_\fpq$.  Therefore,  order ideals of $\calP^\circ_\fpq$ correspond to the subset of order ideals of $\Ptilde_\fpq$ which contain $L_{\Ptilde_\fpq}$ only if they contain $R_{\Ptilde_\fpq}$.  This shows that the trace of $C_\fpq$ is equal to  $\xx^{\gg^\circ_\fpq} \mathcal{W}(\calP^\circ_\fpq)$.  Then,  the equality follows from Proposition \ref{prop:Combinatorial_Cohn_Trace_Formula}.
\end{proof}

The previous results concerned individual positive rational numbers, i.e. $\fpq \in \mathbb{Q} \cap [0,\infty)$. The following concerns Farey triples, and thus we will also consider $\fpq = \frac10$. We will rely on previous results in \cite{ouguz2025oriented} concerning how to multiply the matrices $C_\fpq$.

\begin{proposition}\label{prop:Mutation_Of_Combinatorial_Cohn}
Let $(\frac{p}{q},\frac{p+r}{q+s},\frs)$ be a Farey triple in $\mathrm{F}\mathbb T$ with $\frac{p}{q} < \frac{p+r}{q+s} < \frs$.  The combinatorial cluster generalized Cohn matrices satisfy \[
C_{ \frac{p+r}{q+s}} = C_{\frac{p}{q}}C_\frs - S_{ \frac{p+r}{q+s} }
\]
where \[
S_{ \frac{p+r}{q+s}}:= S_{C_{ \frac{p+r}{q+s}} } = \begin{bmatrix} k_{\frac{p+r}{q+s}} & 0 \\ k_{\frac{p+r}{q+s}} \mathscr M  &k_{\frac{p+r}{q+s}}\end{bmatrix}.
\]
\end{proposition}

\begin{proof}
\textbf{Case 1: $\frs \in \mathbb{Q} \cap [0,\infty)$}

First, we note that from Corollary \ref{cor:g_vecs_in_cluster},  the leading terms in each entry of $C_{\frac{p}{q}}C_\frs - S_{ \frac{p+r}{q+s}}$ agree with those in $C_{ \frac{p+r}{q+s}} $; by ``leading term'', we mean the terms associated with empty order ideals in each poset.
The matrices $C_\fpq$ are closely related to the matrices used in \cite{ouguz2025oriented,ezgieminecluster2024}. In particular, if we ignore the leading monomial in each term, we recover the \emph{dual rank matrix} in \cite{ouguz2025oriented}. As a consequence, from \cite[Proposition 4.2]{ouguz2025oriented} as well as an analysis of the leading monomials $\xx^\gg$, we have \[
C_{\frac{p}{q}}C_\frs = \begin{bmatrix} \xx^{\gg^\circ_\frac{p+r}{q+s}} \mathcal{W}(\Ptilde_{\frac{p}{q}} \nearrow \Ptilde_\frs;R_{\Ptilde_{\frac{r}{s}}})
&  \xx^{\gg_\frac{p+r}{q+s}} \mathcal{W}(\Ptilde_{\frac{p}{q}} \nearrow \Ptilde_\frs;\neg R_{\Ptilde_{\frac{r}{s}}})\\
 \xx^{2\gg^\circ_\frac{p+r}{q+s} - \gg_\frac{p+r}{q+s} } \mathcal{W}(\Ptilde_{\frac{p}{q}} \nearrow \Ptilde_\frs;R_{\Ptilde_{\frac{r}{s}}}, \neg L_{\Ptilde_{\frac{p}{q}}})
 &  \xx^{\gg^\circ_\frac{p+r}{q+s}} \mathcal{W}(\Ptilde_{\frac{p}{q}} \nearrow \Ptilde_\frs;\neg R_{\Ptilde_{\frac{r}{s}}}, \neg L_{\Ptilde_{\frac{p}{q}}})\end{bmatrix},\]
and now we will compare entries of $N:=C_{\frac{p}{q}}C_\frs - S_{ \frac{p+r}{q+s}}$ with those of $C_{ \frac{p+r}{q+s}} $.

 In order to compare these two matrices, which involve different posets, we define $Q_{11}$ to be the poset $\overline{H} \searrow \overline{P}_{\frac{p+r}{q+s}} \nearrow H$.  Using Lemma \ref{lem:ShapeOfPosetInFareyTriple}, we draw a sketch of this poset below, where each diagonal $\cdots$ denotes a copy of $H$. Throughout, we assume $\fpq > \frac01$, so that $\calP_\fpq$ is nonempty. However, one can use item (2) in Lemma \ref{lem:ShapeOfPosetInFareyTriple} when $\fpq = \frac01$ to get the desired result in this case.

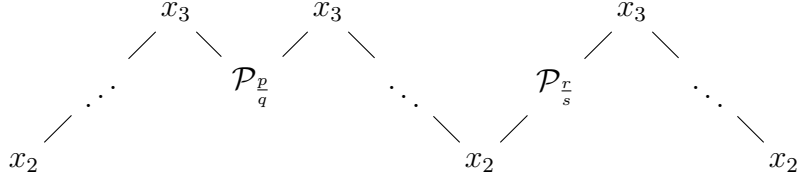
\begin{center}
\begin{tikzpicture}
\node(21) at (0,0){$x_2$};
\node(dots1) at (1,1){$\iddots$};
\node(31) at (2,2){$x_3$};
\node(rs) at (3,1){$\calP_{\frac{p}{q}}$};
\node(32) at (4,2){$x_3$};
\node(dots2) at (5,1){$\ddots$};
\node(22) at (6,0){$x_2$};
\node(pq2) at (7,1){$\calP_{\frac{r}{s}}$};
\node(33) at (8,2){$x_3$};
\node(dots3) at (9,1){$\ddots$};
\node(23) at (10,0){$x_2$};
\draw(21) -- (dots1);
\draw(dots1) -- (31);
\draw(31) -- (rs);
\draw(rs) -- (32);
\draw(32) -- (dots2);
\draw(dots2) -- (22);
\draw(22) -- (pq2);
\draw(pq2) -- (33);
\draw(33) -- (dots3);
\draw(dots3)-- (23);
\end{tikzpicture}
\end{center}

 Removing the left chain yields  $\Ptilde_{\frac{p}{q}} \nearrow \Ptilde_\frs$, and since the element $L_{Q_{11}}$ is minimal, we can write $\mathcal{W}(\Ptilde_{\frac{p}{q}} \nearrow \Ptilde_\frs;R_{\Ptilde_{\frac{r}{s}}}) = \mathcal{W}(Q_{11};R_{Q_{11}},\neg L_{Q_{11}})$.   From Lemma \ref{lem:ShapeOfPosetInFareyTriple}, removing the chains on the left and right yields $\overline{\calP}_{\frac{p+r}{q+s}}$, and so by reversing chronological orderings we have $\mathcal{W}(\Ptilde_{\frac{p+r}{q+s}};R_{\Ptilde_{\frac{p+r}{q+s}}}) = \mathcal{W}(Q_{11};L_{Q_{11}},\neg R_{Q_{11}})$. We are now prepared to compare the $(1,1)$ entries of $N$ and $C_{\frac{p+r}{q+s}}$, \begin{align*}
\mathcal{W}(\Ptilde_\fpq \nearrow \Ptilde_\frs;R_{\Ptilde_{\frac{r}{s}}}) - k_{\frac{p+r}{q+s}} \xx^{-\gb_{\frac{p+r}{q+s}}^\circ} &=  \mathcal{W}(Q_{11};R_{Q_{11}},\neg L_{Q_{11}})- k_{\frac{p+r}{q+s}} \xx^{-\gb_{\frac{p+r}{q+s}}}\\
&=  \alpha \mathcal{W}(\calP_{\frac{p+r}{q+s}}) +\mathcal{W}(\calP_{\frac{p+r}{q+s}}; L_{\calP_{\frac{p+r}{q+s}}}) - k_{\frac{p+r}{q+s}} \xx^{-\gb_{\frac{p+r}{q+s}}^\circ}\\
&=  \alpha \mathcal{W}(\calP_{\frac{p+r}{q+s}}) +\mathcal{W}(\calP_{\frac{p+r}{q+s}}; R_{\calP_{\frac{p+r}{q+s}}})\\
&= \mathcal{W}(\Ptilde_{\frac{p+r}{q+s}};R_{\Ptilde_{\frac{p+r}{q+s}}})
 \end{align*}
 where we use logic similar to Lemma \ref{lem:alphaSimplifying} in the second equality and use this Lemma again in the last equality while in the third equality we use Lemma \ref{lem:Asymmetry}.
Once we multiply through by $\xx^{\gb_{\frac{p+r}{q+s}}}$, we can conclude that the $(1,1)$ entries of $N$ and $C_{\frac{p+r}{q+s}}$ agree.

We again use Lemmas \ref{lem:Asymmetry} and \ref{lem:alphaSimplifying} when comparing the $(2,1)$ entries of these matrices. For concision, we begin by comparing the $(2,1)$ entries of $C_\fpq C_\frs$ and $C_{\frac{p+r}{q+s}}$,    \begin{align*}
 \mathcal{W}(\Ptilde_{\frac{r}{s}} \nearrow \Ptilde_\frs; R_{\Ptilde_{\frs}},\neg L_{\Ptilde_{\frac{p}{q}}}) &=
 \alpha \mathcal{W}(\calP_{\frac{p+r}{q+s}}; \neg R_{\calP_{\frac{p+r}{q+s}}}) + \mathcal{W}(\calP_{\frac{p+r}{q+s}}; L_{\calP_{\frac{p+r}{q+s}}}, \neg R_{\calP_{\frac{p+r}{q+s}}})\\
 &= \alpha \big(\mathcal{W}(\calP_{\frac{p+r}{q+s}}; \neg L_{\calP_{\frac{p+r}{q+s}}}) + k_{\frac{p+r}{q+s}} \xb^{-\gg^\circ_{\frac{p+r}{q+s}}}\big) \\&+ \mathcal{W}(\calP_{\frac{p+r}{q+s}}; R_{\calP_{\frac{p+r}{q+s}}}, \neg L_{\calP_{\frac{p+r}{q+s}}}) + k_{\frac{p+r}{q+s}} \xb^{-\gg^\circ_{\frac{p+r}{q+s}}}\\
 &= \mathcal{W}(\Ptilde_{\frac{p+r}{q+s}}; R_{\Ptilde_{\frac{p+r}{q+s}}}, \neg L_{\Ptilde_{\frac{p+r}{q+s}}}) + (1+\alpha)  k_{\frac{p+r}{q+s}} \xb^{-\gg^\circ_{\frac{p+r}{q+s}}}
 \end{align*}

By multiplying through by $\xx^{2\gb^\circ_{\frac{p+r}{q+s}} - \gb_{\frac{p+r}{q+s}}} = \xx^{\gb^\circ_{\frac{p+r}{q+s}} + (1,-1,-1)}$ we see the difference between the $(2,1)$ entries in $N$ and $C_{\frac{p+r}{q+s}}$ is exactly $k_{\frac{p+r}{q+s}} \frac{x_1}{x_2x_3}(1+\alpha)  = k_{\frac{p+r}{q+s}} \mathscr M $.

For the $(1,2)$-entry,  we have that \[
\mathcal{W}(\Ptilde_{\frac{p}{q}} \nearrow \Ptilde_\frs;\neg R) =\mathcal{W}(\Ptilde_{\frac{p}{q}} \nearrow \calP_\frs)\]
and from Lemma   \ref{lem:ShapeOfPosetInFareyTriple}, this latter poset is $\calP_{\frac{p+r}{q+s}}$ with reverse chronological orientation.  Since $\mathcal{W}(\calP) = \mathcal{W}(\overline{\calP})$ for all posets $\calP$,  we see the (1,2) entries of our two matrices agree.

Finally,  we consider the $(2,2)$ entry.  With similar reasoning as for the $(1,2)$-entry,  we have that the $(2,2)$ entry of $N$ is given by \begin{align*}
\xx^{\gg_\frac{p+r}{q+s}} \mathcal{W}(\Ptilde_{\frac{p}{q}} \nearrow \Ptilde_\frs;\neg R_{\Ptilde_\frs}, \neg L_{\Ptilde_{\frac{p}{q}}}) - k_\frs &= \xx^{\gg_\frac{p+r}{q+s}} \mathcal{W}(\calP_{\frac{p+r}{q+s}};\neg R_{\calP_{\frac{p+r}{q+s}}}) - k_\frs\\
& = \xx^{\gg_\frac{p+r}{q+s}} \mathcal{W}(\calP_{\frac{p+r}{q+s}};\neg L_{\calP_{\frac{p+r}{q+s}}})\\
&=\xx^{\gg_\frac{p+r}{q+s}} \mathcal{W}(\Ptilde_{\frac{p+r}{q+s}};\neg R_{\Ptilde_{\frac{p+r}{q+s}}}, \neg L_{\Ptilde_{\frac{p+r}{q+s}}})\\
\end{align*}
where in the second equality we use Lemma \ref{lem:Asymmetry} and in the last we use the fact that order ideal of $\Ptilde_{\frac{p+r}{q+s}}$ which do not include the rightmost element are in bijection with order ideals of $\calP_{\frac{p+r}{q+s}}$.

\textbf{Case 2: $\frs =\frac10$}

Recall the poset $\Ptilde_{\frac10}$ is drawn in Figure \ref{fig:Ptilde10}. In this case, our Farey triple is of the form $(\frac{p}{1},\frac{p+1}{1},\frac10)$. We can check the case $p = 0$ directly, so we will assume $p > 0$ for convenience. Let $Q_{\frac10}$ denote the poset associated to the line segment $(\epsilon,\epsilon) - (\epsilon,1+\epsilon)$ for small $\epsilon$. The posets $\calP_{\frac{p}{1}}$ are related by \begin{equation}
\calP_{\frac{p+1}{1}} = Q_{\frac10} \searrow \overline{\calP}_{\frac{p}{1}}.\label{eq:Pr1Relationship}
\end{equation}

As in the other case, we will check each component separately.
First, consider the $(1,1)$ entry of $C_{\frac{p}{1}}C_{\frac10}$. This is given by \begin{equation}
\xx^{\gg^\circ_{\frac{p}{1}} + (1,0,-1)^\top} \bigg(\mathcal{W}(\Ptilde_{\frac{p}{1}};R)\mathcal{W}(\Ptilde_{\frac{1}{0}};R,\neg L) - \mathcal{W}(\calP_{\frac{p}{1}})\mathcal{W}(\Ptilde_{\frac{1}{0}};L,R)\bigg)\label{eq:(1,1)In10Case}
\end{equation}

We can see that $\mathcal{W}(\Ptilde_{\frac{1}{0}};R,\neg L) = \mathcal{W}(H;R) = 1 + \alpha$. Similarly, $\mathcal{W}(\Ptilde_{\frac{1}{0}};L,R) = k_2\widehat{x_2}(1+\alpha) + \widehat{x_2}^2$. Moreover, from comparing the structure of $\Ptilde_{\frac{p}{1}}$ and $\calP_{\frac{p+1}{1}}$, from Equation \eqref{eq:Pr1Relationship}, we have \begin{align*}
\mathcal{W}(\Ptilde_{\frac{p}{1}};R) &= k_2\widehat{x_2}\mathcal{W}(\calP_{\frac{p}{1}}) + \widehat{x_2}^2 \mathcal{W}(\overline{\calP}_{\frac{p+1}{1}})\\
&= k_2\widehat{x_2}\mathcal{W}(\calP_{\frac{p}{1}}) + \widehat{x_2}^2 \mathcal{W}(\calP_{\frac{p+1}{1}}).
\end{align*}

Substituting these into Equation \eqref{eq:(1,1)In10Case}, we have that the $(1,1)$ entry of $C_{\frac{p}{1}}C_{\frac10}$ is \begin{align}
&\xx^{\gg^\circ_{\frac{p}{1}} + (1,0,-1)^\top} \bigg((k_2\widehat{x_2}\mathcal{W}(\calP_{\frac{p}{1}}) + \widehat{x_2}^2 \mathcal{W}(\overline{\calP}_{\frac{p+1}{1}}))(1+\alpha) - \mathcal{W}(\calP_{\frac{p}{1}})(k_2\widehat{x_2}(1+\alpha) + \widehat{x_2}^2)\bigg)\nonumber\\
&= \xx^{\gg^\circ_{\frac{p}{1}} + (-1,0,1)^\top}\bigg(\mathcal{W}(\calP_{\frac{p+1}{1}})\alpha + \mathcal{W}(\calP_{\frac{p+1}{1}})-\mathcal{W}(\calP_{\frac{p}{1}}))\bigg)\label{eq:Halfwaythrough}\\
&= \xx^{\gg^\circ_{\frac{p+1}{1}}}\bigg(\mathcal{W}(\calP_{\frac{p+1}{1}})\alpha + \mathcal{W}(\calP_{\frac{p+1}{1}};L)\bigg)\nonumber\\
&= \xx^{\gg^\circ_{\frac{p+1}{1}}}\bigg(\mathcal{W}(\calP_{\frac{p+1}{1}})\alpha + \mathcal{W}(\calP_{\frac{p+1}{1}};R) + k_{\frac{p+1}{1}}\xx^{-\gg^\circ_{\frac{p+1}{1}}}\bigg)\nonumber\\
&= \xx^{\gg^\circ_{\frac{p+1}{1}}}\mathcal{W}(\Ptilde_{\frac{p+1}{1}};R)\alpha + k_{\frac{p+1}{1}}\nonumber
\end{align}
where the first and last equality come from purely algebraic manipulations, the second follows from the relationship in Equation \eqref{eq:Pr1Relationship} and the third follows from Lemma \ref{lem:Asymmetry}. One can also see that $\gg^\circ_{\frac{p}{1}} + (-1,0,1)^\top = \gg^\circ_{\frac{p+1}{1}}$ from Equation \eqref{eq:Pr1Relationship} and our method of computing $\gb_\calP$.

Next, consider the $(1,2)$-entry of $C_{\frac{p}{1}}C_{\frac10}$. Using our previous analysis of the $(1,1)$ entry, we see this is given by \[
\xx^{\gg^\circ_{\frac{p}{1}} + (0,1,0)^\top}\bigg(k_2\widehat{x_2}\mathcal{W}(\calP_{\frac{p}{1}}) + \widehat{x_2}^2 \mathcal{W}(\overline{\calP}_{\frac{p+1}{1}}) - \mathcal{W}(\calP_{\frac{p}{1}})k_2\widehat{x_2}\bigg) = \xx^{\gg^\circ_{\frac{p}{1} + (-2,1,2)^\top}}\mathcal{W}(\calP_{\frac{p+1}{1}}).
\]
Now, from Lemma \ref{lem:CompareGVecAndGCirc}, we have $\gg^\circ_{\frac{p}{1}} + (-2,1,2)^\top = \gg_{\frac{p}{1}} + (-1,0,1)^\top$, and from Equation \eqref{eq:Pr1Relationship} we see this is $\gg_{\frac{p+1}{1}}$.

Analyzing the $(2,1)$ entry is very similar to analyzing the $(1,1)$, except in the analogue to Equation \eqref{eq:(1,1)In10Case}, both terms involving $\calP_{\frac{p}{1}}$ or $\Ptilde_{\frac{p}{1}}$ only consider order ideals without the corresponding element $L$. The leading monomial is also slightly different but follows the same pattern. As a result, following similar reasoning up to Equation \eqref{eq:Halfwaythrough}, we reach the conclusion that the $(2,1)$ entry of $C_{\frac{p}{1}}C_{\frac10}$ is \[
\xx^{2\gg^\circ_{\frac{p+1}{1}}-\gg_{\frac{p+1}{1}}}\bigg(\mathcal{W}(\calP_{\frac{p+1}{1}};\neg R)\alpha + \mathcal{W}(\calP_{\frac{p+1}{1}};\neg R)-\mathcal{W}(\calP_{\frac{p}{1}};\neg L))\bigg)
\]
where the condition $\neg R$ comes from the condition $\neg L$ in the $(2,1)$ and $(2,2)$ entries of $C_{\frac{p}{1}}$ after reversing the chronological ordering. Now, again by Lemma \ref{lem:Asymmetry}, this is equal to \begin{align*}
&\xx^{2\gg^\circ_{\frac{p+1}{1}}-\gg_{\frac{p+1}{1}}}\bigg((1+\alpha)k_{\frac{p+1}{1}}\xx^{-\gg^\circ_{\frac{p+1}{1}}}  + \mathcal{W}(\calP_{\frac{p+1}{1}};\neg L)(1+\alpha) - \mathcal{W}(\calP_{\frac{p}{1}};\neg L)\bigg)\\\\
&=\frac{x_1}{x_2x_3}k_{\frac{p+1}{1}} \alpha + \xx^{2\gg^\circ_{\frac{p+1}{1}}-\gg_{\frac{p+1}{1}}}\bigg(\mathcal{W}(\calP_{\frac{p+1}{1}};\neg L)\alpha +\mathcal{W}(\calP_{\frac{p+1}{1}};R,\neg L)\bigg)\\
&= k_{\frac{p+1}{1}}\mathscr M  + \xx^{2\gg^\circ_{\frac{p+1}{1}}-\gg_{\frac{p+1}{1}}}\mathcal{W}(\Ptilde_{\frac{p+1}{1}};R,\neg L).
\end{align*}
 The second equality follows from $\mathcal{W}(\calP_{\frac{p+1}{1}};\neg L) - \mathcal{W}(\calP_{\frac{p}{1}};\neg L) = \mathcal{W}(\calP_{\frac{p+1}{1}};R,\neg L)$, which can be seen as another consequence of Equation \eqref{eq:Pr1Relationship}, and the last follows from Lemma \ref{lem:alphaSimplifying}.

 Verifying the final entry, i.e., $(2,2)$ is the same in $C_{\frac{p}{1}}C_{\frac10} - S_{\frac{p+1}{1}}$ and $C_{\frac{p+1}{1}}$ uses very similar reasoning to the previous cases, so we omit the details.
\end{proof}

\begin{proof}[Proof of Theorem \ref{thm:Combinatorial_Cohn_Matrix_and_Triple}]
First,  we check that $C_\fpq$ satisfies the three requirements of Definition \ref{def:gen-Cohn-matrix}. If $\fpq = \frac10$, then this follows from Lemma \ref{lem:C10}. Now, suppose $\fpq \in \mathbb{Q} \cap [0,\infty)$.  We note that $C_\fpq \in SL(2,\mathbb{Z}[x_1^\pm, x_2^\pm,x_3^\pm])$ by Proposition \ref{prop:Combinatorial_Cohn_determinant_1}.  The $(1,2)$-entry of $C_\fpq$ is equal to $\xx^{\gg_\fpq}\mathcal{W}(\calP_\fpq)$,  which is the cluster variable $x_\fpq$ in the cluster algebra $\mathcal{A}(k_1,k_2,k_3)$ by Corollary \ref{cor:CorrectnessOfPosetFormula-trivial}.  Finally,  Proposition \ref{prop:Combinatorial_Cohn_Trace_Formula} and Corollary \ref{cor:Band_Graph_Formula_Conserved} show that $C_\fpq$ satisfies the third condition.

Now,  let $(\frac{p}{q},\frac{p+r}{q+s},\frs)$ be a Farey triple.  The triple $(C_\frac{p}{q},C_{\frac{p+r}{q+s}},C_\frs)$ clearly satisfies part (1) of Definition \ref{def:gen-Cohn-triple}. Proposition \ref{prop:Mutation_Of_Combinatorial_Cohn} verifies part (2). The indexing of cluster variables by rational numbers (as in Section \ref{subsec:LabelClVar}) and Corollary \ref{cor:CorrectnessOfPosetFormula-trivial} shows that the triple also satisfies part (3).
\end{proof}

\begin{remark}\label{rmk:YVar}
As our poset expansion formula is valid in the context of principal coefficients, one might expect to see our Cohn matrices also in this context. While such a family of matrices is interesting (see Section \ref{sec:RelationshipsToOthers}), the obvious extension of the definition to principal coefficients loses some of the nice properties enjoyed in the coefficient-free case. A key issue is our reliance in several proof on the identity  $\widehat{x_1}\widehat{x_2}\widehat{x_3} = 1$. When we work with principal coefficients, we no longer have such cancellation. In particular, with principal coefficients, there is no clear candidate for $\mathscr M $ since $\mathcal{W}^{\prin}(H;R) \neq \mathcal{W}^{\prin}(H;\neg L)$.
\end{remark}

\begin{example}\label{ex:CohnTriple}
Here, we give several examples of combinatorial cluster generalized Cohn matrices. First, consider $\fpq = \frac{0}{1}$. In this case, for general $k_1,k_2,k_3$, we have \[
C_{\frac01} = \begin{bmatrix} \frac{x_1}{x_2}\bigg(k_2 \frac{x_3}{x_1} + \frac{x_3^2}{x_1^2} + k_1 \frac{x_2x_3}{x_1^2} + \frac{x_2^2}{x_1^2}\bigg) & x_3\\
\frac{x_1^2}{x_2^2x_3}\bigg(k_2 \frac{x_3}{x_1} + \frac{x_3^2}{x_1^2} + k_1 \frac{x_2x_3}{x_1^2}\bigg) & \frac{x_1}{x_2}\end{bmatrix}.\]

Notice that this coincides with $P_{(1,3,2)}(f)$ when $f = \frac{x_1}{x_2}\bigg(k_2 \frac{x_3}{x_1} + \frac{x_3^2}{x_1^2} + k_1 \frac{x_2x_3}{x_1^2} + \frac{x_2^2}{x_1^2}\bigg)$.

Next, consider $C_{\frac11}$ when all $k_i = 0$. Recall the $g$-vector in this case is $(-1,0,2)^\top$, and the poset $\Ptilde_{\frac11}$ is $x_1 < x_3 > x_1 > x_2$. The matrix is then \[
C_{\frac11} = \begin{bmatrix} \frac{x_3}{x_2} \bigg(\big(1 + \frac{x_2^2}{x_3^2}\big)\big(\frac{x_3^2}{x_1^2}+ \frac{x_2^2}{x_1^2}\big) + \frac{x_2^2}{x_3^2}\bigg) &
\frac{x_3^2}{x_1} \bigg( 1 + \frac{x_2^2}{x_3^2}\bigg) \\
\frac{x_1}{x_2^2} \bigg( \frac{x_3^2}{x_1^2}+ \frac{x_2^2}{x_1^2}\bigg) & \frac{x_3}{x_2}  \end{bmatrix}.
\]

Suppose that instead $k_1 = 0, k_2 = 1$ and $k_3 = 2$. Then, the poset $\Ptilde_{\frac11}$ is $x_1 < x_3 > x_3 > x_1 > x_2 > x_2$, and the combinatorial cluster generalized Cohn matrix is more complicated.

\[
C_{\frac11} = \begin{bmatrix} \frac{x_3}{x_2} \bigg(\big(1 + \frac{x_2^2}{x_3^2}\big)\big(\frac{x_3}{x_1} + \frac{x_3^2}{x_1^2}+ 2\frac{x_3^2}{x_1^2}+  2\frac{x_3^2}{x_1x_2}\big) + \frac{x_2^2}{x_3^2}\bigg) &
\frac{x_3^2}{x_1} \bigg( 1 + \frac{x_2^2}{x_3^2}\bigg) \\
\frac{x_1}{x_2^2} \bigg(\frac{x_3}{x_1} + \frac{x_3^2}{x_1^2}+ 2\frac{x_3^2}{x_1^2}+  2\frac{x_3^2}{x_1x_2} \bigg) & \frac{x_3}{x_2}  \end{bmatrix}.
\]

Finally, we construct $C_{\frac12}$ in the case where all $k_i = 0$, i.e., in the Markov cluster algebra. The poset $\Ptilde_{\frac12}$ is as below, where each $x_i$ is weighted $\widehat{x_i}^2$.

\begin{center}
\begin{tikzpicture}[scale=0.8]
\node(1) at (0,0){$x_1$};
\node(2) at (1,-1){$x_2$};
\node(3) at (2,0){$x_1$};
\node(4) at (3,1){$x_3$};
\node(5) at (4,0){$x_1$};
\node(6) at (5,-1){$x_2$};
\draw(1) -- (2);
\draw(2) -- (3);
\draw(3) -- (4);
\draw(4) -- (5);
\draw(5) -- (6);
\end{tikzpicture}
\end{center}

Recall $\gb_{\frac12} = (0,-1,2)^\top$, so that $\gb^\circ_{\frac12} = (1,-2,1)^\top$. From Corollary \ref{cor:CorrectnessOfPosetFormula-trivial}, we can calculate, \[
x_{\frac12} = \frac{x_3^2}{x_2}\bigg(1 + \frac{x_3^2}{x_1^2} + 2\frac{x_2^2}{x_1^2} + \frac{x_2^4}{x_1^2x_3^2}\bigg) = \frac{x_1^2x_3^2 + x_3^4 + 2x_2^2x_3^2 + x_2^4}{x_1^2x_2}.
\]

Let $z_{\frac12} = \frac{x_2}{x_3^2} x_{\frac12}$. Define also $w_{\frac12} = \mathcal{W}(\calP_{\frac12};\neg L) =  1 + \frac{x_3^2}{x_1^2} + \frac{x_2^2}{x_1^2}$. Then, the combinatorial cluster generalized Cohn matrix associated to $\frac12$ is \[
C_\frac12 = \begin{bmatrix} \frac{x_1x_3}{x_2^2} \bigg( z_{\frac12} \big(\frac{x_3^2}{x_1^2} + \frac{x_2^2}{x_1^2}\big) + \frac{x_2^2}{x_1^2} + \frac{x_2^4}{x_1^2x_3^2}\bigg) & \frac{x_3^2}{x_2}z_{\frac12} \\
\frac{x_1^2}{x_2^3}\bigg(w_{\frac12}\big(\frac{x_3^2}{x_1^2} + \frac{x_2^2}{x_1^2}\big) +  \frac{x_2^2}{x_1^2}\bigg) & \frac{x_1x_3}{x_2^2} w_{\frac12}
\end{bmatrix}.
\]

We check here that the trace is as expected, \begin{align*}
\mathrm{tr}(C_{\frac12}) &= \frac{x_1x_3}{x_2^2} \bigg( z_{\frac12} \big(\frac{x_3^2}{x_1^2} + \frac{x_2^2}{x_1^2}\big) +  \frac{x_2^2}{x_1^2} + \frac{x_2^4}{x_1^2x_3^2} + w_{\frac12} \bigg)\\
&= \frac{x_1x_3}{x_2^2} \bigg(1 + \frac{x_3^2}{x_1^2} + \frac{x_2^2}{x_1^2}\bigg) z_{\frac12}\\
&= \frac{x_1x_3}{x_2^2} \bigg(1 + \frac{x_3^2}{x_1^2} + \frac{x_2^2}{x_1^2}\bigg) \frac{x_2}{x_3^2}x_{\frac12} \\
&= \mathscr M x_{\frac12}.
\end{align*}

In the case $k_1 = k_2 = k_3 = 0$, we have now collected all combinatorial cluster generalized Cohn matrices for the Farey triple $(\frac01, \frac12, \frac11)$. One can check that $C_{\frac12} = C_{\frac01} C_{\frac11}$; the matrix $S_{\frac12}$ is trivial here since $k_2 = 0$.

\end{example}

Finally, we show how the tree of combinatorial cluster generalized Cohn triples fits into the general story of $(k_1,k_2,k_3)$-cluster generalized Cohn trees. Let the \emph{combinatorial cluster generalized Cohn tree} be the tree with root given by $(C_{\frac{0}{1}},C_{\frac{1}{1}},C_{\frac10})$ and propagation rule given by \[\begin{xy}(0,0)*+{(C_{\frac{p}{q}},C_{\frac{p+r}{q+s}},C_{\frs})}="1",(30,-15)*+{(C_{\frac{p+r}{q+s}},C_{\frac{p+r}{q+s}}C_{\frs} - S_{\frac{2p+r}{2q+s}},C_{\frs}).}="2",(-30,-15)*+{(C_{\frac{p}{q}},C_{\frac{p}{q}}C_{\frac{p+r}{q+s}}-S_{\frac{p+2r}{q+2s}},C_{\frac{p+r}{q+s}})}="3", \ar@{-}"1";"2"\ar@{-}"1";"3"
\end{xy}\]

\begin{proposition}\label{prop:WhichTreeIsCombCohn}
The combinatorial cluster generalized Cohn tree is \\$\mathrm{CGC}\mathbb T(k_1,k_2,k_3,(1,3,2),f)$ where $f = \frac{x_1}{x_2}\bigg(k_2 \frac{x_3}{x_1} + \frac{x_3^2}{x_1^2} + k_1 \frac{x_2x_3}{x_1^2} + \frac{x_2^2}{x_1^2}\bigg)$.
\end{proposition}

\begin{proof}
This follows from checking that the roots of the two trees agree, which is straightforward.
\end{proof}

\section{Combinatorial Cluster Markov-Monodromy Matrices}\label{sec:CombinatorialMM}

The goal of this section is to provide an explicit set of Markov-monodromy matrices and triples. This will be achieved by applying $\psi^{-1}_{g}$ to each combinatorial cluster generalized Cohn matrix for a certain value of $g$. Accordingly, we will refer to the matrices constructed as \emph{combinatorial cluster Markov-monodromy matrices}. These matrices appear to be unique in the literature.

\begin{theorem}
Define \[
M_{\frac01} = \begin{bmatrix} \frac{x_1}{x_2} & x_3 \\ -\frac{x_1^2}{x_2^2x_3}(1 + k_3 \frac{x_2}{x_1} + \frac{x_2^2}{x_1^2}) & -\frac{x_1}{x_2} - k_3 \end{bmatrix} \qquad M_{\frac10} = \begin{bmatrix} -k_2 & x_2 \\ -\frac{1}{x_2} & 0 \end{bmatrix} ,
\]
and for each $\fpq \in \mathbb{Q} \cap (0,\infty)$, let $M_\fpq$ be defined by \[
M_\fpq = \begin{bmatrix} \xx^{\gg^\circ_\fpq} \mathcal{W}(\calP_\fpq; \neg L) & \xx^{\gg_\fpq} \mathcal{W}(\calP_\fpq)\\
-\xx^{2\gg^\circ_\fpq - \gg_\fpq} \mathcal{W}(\calP_\fpq; \neg L, \neg R) & -\xx^{\gg^\circ_\fpq} \mathcal{W}(\calP_\fpq;\neg R)\end{bmatrix}.
\]
Each matrix $M_\fpq$ is a $(k_1,k_2,k_3)$-CMM matrix. Moreover, if $(\frac{p}{q},\frac{p+r}{q+s},\frs)$ is a Farey triple, then $(M_{\frac{p}{q}},M_{\frac{p+r}{q+s}},M_\frs)$ is a $(k_1,k_2,k_3)$-CMM triple.
\end{theorem}

\begin{proof}
 Fix $(k_1,k_2,k_3) \in \mathbb{Z}_{\geq 0}^3$.  We will show that  \[
M_\fpq = \psi_{\mathscr M}^{-1}(C_\fpq),
\]
which by Theorem \ref{thm:BT-CT2} will guarantee that $M_\fpq$ is a Markov monodromy matrix and that given a Farey triple $(\frac{p}{q},\frac{p+r}{q+s},\frs)$,  the triple $(M_{\frac{p}{q}},M_{\frac{p+r}{q+s}},M_\frs)$ is a cluster monodromy triple.

Note that the setting $g = \mathscr M$ in the definition of $\psi_{g}^{-1}$ simplifies the expression, \begin{align*}
\psi_{\mathscr M}^{-1}:  \begin{bmatrix} m_{11} & m_{12} \\ m_{21} & m_{22} \end{bmatrix} \mapsto \begin{bmatrix}1&0\\-\mathscr M&1\end{bmatrix}\begin{bmatrix}
    m_{22}&m_{12}\\m_{21}& m_{11}
\end{bmatrix}\begin{bmatrix}
    1&0\\0&1
\end{bmatrix}
=\begin{bmatrix} m_{22} & m_{12}\\
m_{21} - \mathscr M m_{22} & m_{11} - \mathscr M m_{12} \end{bmatrix}.
\end{align*}

We analyze each entry of $\psi_{\mathscr M}^{-1}(C_\fpq)$ individually. We can check $\fpq = \frac01$ and $\fpq = \frac10$ directly, so we assume throughout that $0 < \fpq < \infty$.
Since $\mathcal{W}(\Ptilde_{\fpq}; \neg L, \neg R) = \mathcal{W}(\calP_{\fpq}; \neg L)$, it is clear that the $(1,1)$ entry of $\psi_{\mathscr M}^{-1}(C_\fpq)$ matches that of $M_\fpq$. Moreover, by definition the $(1,2)$-entries agree.

Now, we compute the $(2,1)$ entry of $\psi_{\mathscr M}^{-1}(C_\fpq)$, \begin{align*}
&\xx^{2 \gg^\circ_\fpq - \gg_\fpq}\mathcal{W}(\Ptilde_\fpq; R_{\Ptilde_\fpq},  \neg L_{\Ptilde_\fpq}) - \mathscr M\xx^{\gg^\circ_\fpq} \mathcal{W}(\Ptilde_{\fpq};\neg L_{\calP_{\fpq}},\neg R_{\calP_{\fpq}})\\
&=\xx^{2 \gg^\circ_\fpq - \gg_\fpq}\big(\alpha \mathcal{W}(\calP_\fpq;  \neg L_{\calP_\fpq}) + \mathcal{W}(\calP_\fpq; R_{\calP_\fpq},  \neg L_{\calP_\fpq}) \big) \\
&- \frac{x_1}{x_2x_3}\xx^{\gg^\circ_\fpq}(1+\alpha)\mathcal{W}(\calP_{\fpq};\neg L_{\calP_\fpq})\\
&=\xx^{2 \gg^\circ_\fpq - \gg_\fpq} (\mathcal{W}(\calP_\fpq; R_{\calP_\fpq},  \neg L_{\calP_\fpq}) - \mathcal{W}(\calP_\fpq;  \neg L_{\calP_\fpq}) \big),\\
&= -\xx^{2 \gg^\circ_\fpq - \gg_\fpq} \mathcal{W}(\calP_\fpq  ;\neg L_{\calP_\fpq}, \neg R_{\calP_\fpq},)
\end{align*}

where we use Lemmas \ref{lem:CompareGVecAndGCirc}, \ref{lem:RewriteMinv}, and \ref{lem:alphaSimplifying}. We see this agrees with the $(2,1)$ entry of $M_{\fpq}$.

The proof to show that the $(2,2)$ entries agree is identical to the proof for the $(2,1)$ entries after removing the condition $\neg L$ from each polynomial $\mathcal{W}$. We have now shown  $\psi_{\mathscr M}^{-1}(C_\fpq) = M_\fpq$, and the statement follows.
\end{proof}

\begin{example}\label{ex:ClusterMonodromy}
Recall in Example \ref{ex:CohnTriple}, we gave several examples of Cohn matrices. Here, we give the corresponding combinatorical Markov-monodromy. When $k_1 = k_2 = k_3 = 0$, $\mathscr M  = \frac{x_1^2 + x_2^2 + x_3^3}{x_1x_2x_3}$, and in this case \[
M_{\frac11} = \psi_{\mathscr M}^{-1}(C_{\frac11}) = \begin{bmatrix} \frac{x_3}{x_2} & \frac{x_3^2}{x_1}\bigg(1 + \frac{x_2^2}{x_3^2}\bigg) \\ -\frac{x_1}{x_2^2} & -\frac{x_3}{x_2} \end{bmatrix}.
\]

whereas, with all $k_i$ still 0, $M_{\frac12}$ is given by
\[
M_{\frac12} = \psi_{\mathscr M}^{-1}(C_{\frac12}) = \begin{bmatrix}
\frac{x_1x_3}{x_2^2}\bigg(1 + \frac{x_3^2}{x_1^2} + \frac{x_2^2}{x_1^2}\bigg)
&  \frac{x_3^2}{x_2}\bigg(1 + \frac{x_3^2}{x_1^2} + 2\frac{x_2^2}{x_1^2} + \frac{x_2^4}{x_1^2x_3^2}\bigg) \\ -\frac{x_1^2}{x_2^3}\bigg(1 + \frac{x_3^2}{x_1^2}\bigg) & -\frac{x_1x_3}{x_2^2}\bigg(1 + \frac{x_3^2}{x_1^2} + \frac{x_2^2}{x_1^2}\bigg) \end{bmatrix}.
\]

In each of these examples, the entries on the main diagonal were the same since the corresponding $k_i$ value was 0. We could also consider $M_{\frac12}$ when $k_1 = k_2 = k_3 = 1$. The associated poset $\calP_{\frac12}$ as well as the polynomials $\mathcal{W}(\calP_{\frac12};\neg L)$ and $\mathcal{W}(\calP_{\frac12};\neg R)$ were given in Example \ref{ex:Asymmetry}. We write the matrix $M_{\frac12}$ with $k_1 = k_2 = k_3 = 1$ below.

\[
\scalemath{0.9}{\begin{bmatrix}
\frac{x_1x_3}{x_2^2}\bigg(1 + \frac{x_3}{x_1}\big(1 + \frac{x_2}{x_3}\big)\big(1 + \frac{x_3}{x_1} + \frac{x_2}{x_1} + \frac{x_2^2}{x_1x_3}\big)\bigg)
&  \frac{x_3^2}{x_2}\bigg(1 + \frac{x_3}{x_1}\big(1 + \frac{x_2}{x_3} + \frac{x_2^2}{x_3^2}\big)\big(1 + \frac{x_3}{x_1} + \frac{x_2}{x_1} + \frac{x_2^2}{x_1x_3}\big)\bigg)  \\
-\frac{x_1^2}{x_2^3}\bigg(1 + \frac{x_3}{x_1}\big(1 + \frac{x_2}{x_3}\big)\big(1 + \frac{x_3}{x_1} + \frac{x_2}{x_1} \big)\bigg) &
-\frac{x_1x_3}{x_2^2}\bigg(1 + \frac{x_3}{x_1}\big(1 + \frac{x_2}{x_3} + \frac{x_2^2}{x_3^2}\big)\big(1 + \frac{x_3}{x_1} + \frac{x_2}{x_1}\big)\bigg) \end{bmatrix}}
\]
\end{example}

Now, we show how the tree of combinatorial cluster Markov-monodromy triples fits into the general story of $(k_1,k_2,k_3)$-cluster Markov-monodromy trees. Let the \emph{combinatorial cluster Markov-monodromy (MM) tree} be the tree with root given by $(M_{\frac{0}{1}},M_{\frac{1}{1}},M_{\frac10})$ and propagation rule given by \[\begin{xy}(0,0)*+{(M_{\frac{p}{q}},M_{\frac{p+r}{q+s}},M_{\frs})}="1",(30,-15)*+{(M_{\frac{p+r}{q+s}},M_{\frac{p+r}{q+s}}^{-1}M_{\frac{p}{q}}M_{\frac{p+r}{q+s}} ,M_{\frs}).}="2",(-30,-15)*+{(M_{\frac{p}{q}},M_{\frac{p+r}{q+s}}^{-1}M_{\frac{r}{s}}M_{\frac{p+r}{q+s}},M_{\frac{p+r}{q+s}})}="3", \ar@{-}"1";"2"\ar@{-}"1";"3"
\end{xy}\]

\begin{proposition}\label{prop:WhichTreeIsCombMM}
The combinatorial MM tree is $\mathrm{CMM}\mathbb T(k_1,k_2,k_3,(1,3,2),\frac{x_1}{x_2})$.
\end{proposition}

\begin{proof}
The combinatorial cluster Markov-monodromy tree is the image of the combinatorial cluster generalized Cohn tree under $\Psi_\mathscr M^{-1}$. By Proposition \ref{prop:WhichTreeIsCombCohn} this latter tree is exactly $\mathrm{CGC}\mathbb T(k_1,k_2,k_3,(1,3,2),f)$ where $f = \frac{x_1}{x_2} \bigg(k_2 \frac{x_3}{x_1} + \frac{x_3^2}{x_1^2} + k_1 \frac{x_2x_3}{x_1^2} + \frac{x_2^2}{x_1^2}\bigg)$. Therefore, by applying Theorem \ref{thm:BT-CT2} to $\Psi_M^{-1}$, the   combinatorial cluster Markov-monodromy tree is exactly $\mathrm{CMM}\mathbb T(k_1,k_2,k_3,(1,3,2),-f + x_3M - k_2)$. The statement follows after computing $-f + x_3M - k_2$.
\end{proof}

\section{Relationship to Other Matrix Formulas}\label{sec:RelationshipsToOthers}

In this concluding section, we explain a few ways in which our matrices compare with previously known families of matrices related to cluster algebras or to Markov numbers.

\subsection{Matrices with integer entries}

When we specialize $x_1=x_2=x_3 = 1$, our Cohn matrices and cluster-Markov monodromy matrices will be in $SL(2,\mathbb{Z})$. If $k_1 = k_2 = k_3 = 0$, the specialized Cohn matrices will be exactly recover the ordinary Cohn matrices, as in \cite{Cohn1,Cohn2,remak1924indefinite}.

Both families of matrices, with $x_i = 1$ (and all $y_i = 1$) and all $k_1 = k_2 = k_3 = k$ for a general nonnegative integer $k$, were studied by the second author, Maruyama and Sato in \cite{gyoda2024sl}. The matrices which appear there have entries which are in terms of continued fractions. The numerator of a continued fraction enumerates the number of order ideals in a fence poset. However, even when specializing all $x_i = 1$, if $k > 1$, our posets have some elements which have a weight other than 1. In this setting, one can replace each pair of elements with weights $k,\frac1k$ with a chain of $k+1$ elements. This will recover the fence posets whose order ideals are enumerated by the methods in \cite{gyoda2024sl}. Precisely, the matrix $C_{\frs}\vert_{x_i = 1}$ recovers their matrix $C_\frs(k,2k+2)$ and the matrix $M_{\frs}\vert_{x_i = 1}$ recovers their matrix $M_\frs(k,1)$. See Theorems 7.34 and 7.32 respectively in \cite{gyoda2024sl}.

By making different choices for our combinatorial cluster generalized Cohn matrices, we can construct other families of matrices from \cite{gyoda2024sl}. We detail one different choice we could make and how it would affect our matrices. Every line segment $\gamma_\frs$ crosses one line segment from our lattice $\mathcal{L}$ at its midpoint. If the $k$-value for this line segment is nonzero, this intersection yields two elements of $\calP_\frs$. We made a choice in the Construction Algorithm to always set the smaller indexed element less than the greater indexed element. However, we could have chosen the opposite convention, which would yield an equivalent poset. However, with this choice, we need to amend $\Ptilde$ to instead correspond to adding a \emph{counterclockwise circle} at the end of $\gamma_\frs$ (recall Definition \ref{def:H}). Let the associated poset be denoted $\Ptilde^\vee_\frs$.
We will use $\Ptilde^\vee$ to define the \emph{dual combinatorial cluster generalized Cohn matrices}. We give explicit formulas for the matrices at $\frac01$ and $\frac10$, which are similar $C_{\frac10}$ and $C_{\frac01}$ respectively.  \[
C^\vee_{\frac01} = \begin{bmatrix}
-\frac{x_1}{x_2}\bigg(k_3 \frac{x_2}{x_1}\bigg) & x_3 \\ -\frac{x_1^2}{x_2x_3} \bigg(1 + k_3 \frac{x_2}{x_1}(1 + \alpha)\bigg) & \frac{x_1}{x_2}\bigg(1+\alpha\bigg) \end{bmatrix}\]
\[
C^\vee_{\frac10} = \begin{bmatrix} \frac{x_1}{x_3} & x_2\\
\frac{x_1^2}{x_2x_3^2}\bigg(k_3 \frac{x_2}{x_1} + \frac{x_2^2}{x_1^2} + k_1 \frac{x_2x_3}{x_1^2}\bigg)&\frac{x_1}{x_3}\bigg(k_3 \frac{x_2}{x_1} + \frac{x_2^2}{x_1^2} + k_1 \frac{x_2x_3}{x_1^2} + \frac{x_3^2}{x_1^2}\bigg)\end{bmatrix}
\]

For $\frs \in \mathbb{Q} \cap (0,\infty)$, we define the dual combinatorial cluster generalized Cohn matrix by \[
C^\vee_\frs = \begin{bmatrix} \xx^{\gg^\circ_{\frac{r}{s}}} \mathcal{W}(\Ptilde^\vee_{\frac{r}{s}}; L,R) & \xx^{\gg_{\frac{r}{s}}} \mathcal{W}(\Ptilde^\vee_{\frac{r}{s}};  R)\\
\xx^{2\gg^\circ_\frs - \gg_\frs} \mathcal{W}(\Ptilde^\vee_{\frac{r}{s}}; L, \neg R)& \xx^{\gg^\circ_\frs} \mathcal{W}(\Ptilde^\vee_{\frac{r}{s}}; \neg R)\end{bmatrix}.
\]
Using the same techniques as in Section \ref{sec:CombinatorialCohnMatrices}, one can show that these are indeed Cohn matrices and that for each Farey triple the dual combinatorial cluster generalized Cohn matrices form a Cohn triple. The matrix $C^\vee_{\frs}\vert_{x_i = 1}$ recovers the matrix $C_\frs(k,-k)$ from \cite{gyoda2024sl}. Moreover, if $M^\vee_\frs = \psi_g^{-1}(C^\vee_\frs)$, then the matrix $M^\vee_\frs\vert_{x_i = 1}$ recovers $M_\frs(k,0)$ from \cite{gyoda2024sl}. With variables included, $M^\vee_\frs$ for $\frs \in \mathbb{Q} \cap (0,\infty)$ can be interpreted as follows. \[
M^\vee_\frs = \begin{bmatrix}
- \xb^{\gb_\frs^\circ} \mathcal{W}(\calP_\frs;L) &  \xb^{\gb_\frs} \mathcal{W}(\calP_\frs) \\
- \xb^{2\gb_\frs^\circ - \gb_\frs} \mathcal{W}(\calP_\frs;L,R) &  \xb^{\gb_\frs^\circ} \mathcal{W}(\calP_\frs;R)
\end{bmatrix}
\]

\subsection{Matrices with Laurent polynomial entries}

In Remark \ref{rmk:YVar}, we explain the difficulty of extending Theorem \ref{thm:Combinatorial_Cohn_Matrix_and_Triple} to the principal coefficient case. Nonetheless, consider a \emph{principal coefficient combinatorial cluster generalized Cohn matrix}, \[
C^\prin_\frs = \begin{bmatrix} \xx^{\gg^\circ_{\frac{r}{s}}} \mathcal{W}^\prin(\Ptilde_{\frac{r}{s}}; R) & \xx^{\gg_{\frac{r}{s}}} \mathcal{W}^\prin(\Ptilde_{\frac{r}{s}}; \neg R)\\
\xx^{2\gg^\circ_\frs - \gg_\frs} \mathcal{W}^\prin(\Ptilde_{\frac{r}{s}}; R, \neg L)& \xx^{\gg^\circ_\frs} \mathcal{W}^\prin(\Ptilde_{\frac{r}{s}}; \neg R, \neg L)\end{bmatrix},
\]
which matrices $C_\frac01^\prin$ and $C_\frac10^\prin$ defined analogously to their coefficient-free analogues.
The determinant of $C^\prin_\frs$ is a positive $y$-monomial and by Corollary  \ref{cor:CorrectnessOfPosetFormula-trivial}, the $(1,2)$-entry is the cluster variable $x_{1,\frs}$.

Several $q$-analogues of Cohn matrices have been introduced. A critical inspiration for our combinatorial cluster generalized Cohn matrix was the family of $q$-Cohn matrices, given by Kanatarc{\i} O\u{g}uz in \cite{ouguz2025oriented} extending work in \cite{Leclereqdeform}. We remark these are the result of taking $C_\frs^\prin$ in the case $k_i = 0$ and specializing all $x_i$ to 1 and all $y_i$ to $q$. We can also multiply each component of this specialized matrix by $q^{-r-s}$ to obtain Kogiso's $q$-deformation of Cohn matrix introduced in \cite{kogiso2020q-deformations}.

Next, we discuss the similarity between our combinatorial cluster Markov-monodromy matrices and certain matrices associated with \emph{snake graphs}. Snake graphs are certain labeled, bipartite graphs which can be used to give Laurent expansions for cluster variables in a surface-type cluster algebra. Musiker, Schiffler, and Williams showed that the set of perfect matchings of a snake graph can be given a partial order, and this in fact turns out to be a distributive lattice \cite{musiker2013bases}. The poset of join irreducibles can be easily determined by the snake graph, and the resulting poset is always a fence poset. In particular, our expansion formula $\xx^{\gb_\calP} \mathcal{W}(\calP)$ could easily be rephrased in terms of snake graphs. The correspondence between snake graphs and fence posets is fully explored in \cite{ezgieminecluster2024}. Below, we draw the snake graph associated to $\calP_{\frac12}$ when $k_1 = 0$ and $k_2 > 0$. We drew the associated poset in Example \ref{ex:PosetsInFareyTriple}. Explicitly, the weighted sum of all perfect matchings of this graph, divided by $x_1^2x_2^2$, is equal to $\xx^{\gb_{\frac12}}\mathcal{W}(\calP_{\frac12})$.

\begin{center}
\begin{tikzpicture}[scale = 1.7]
\draw(0,0) to node[left]{$x_2$} (0,1) to node[left]{$x_1$} (0,2) to node[above]{$x_2$} (1,2) to node[above]{$x_1$} (2,2) to  node[above] {$x_2$} (3,2) to node[right]{$x_3$} (3,1) to node[below]{$x_2$} (2,1) to node[below]{$x_2$} (1,1) to node[right]{$x_2$} (1,0) to node[below]{$x_3$} (0,0);
\draw(0,1) to node[above]{$x_3$} (1,1) to node[right, yshift =-10pt]{$k_2 x_2$} (1,2);
\draw(2,1) to node[right,yshift =-10pt]{$x_3$} (2,2);
\draw[dashed, gray](0,1) to node[right, text=black]{$x_1$} (1,0);
\draw[dashed, gray] (0,2) to node[right, text=black]{$x_2$} (1,1);
\draw[dashed, gray] (1,2) to node[right, text=black]{$x_2$} (2,1);
\draw[dashed, gray] (2,2) to node[right, text=black]{$x_1$} (3,1);
\end{tikzpicture}
\end{center}

When studying multiplication formulas for cluster variables coming from skein relations, Musiker and Williams built a family of $2 \times 2$ matrices to associate to snake graphs; see  \cite[Section 5.1]{MusikerMatrix}. The formulas rely on cluster variables denoted $x_a,x_b,x_w,x_z, x_{i_1},$ and $x_{i_d}$, which correspond to the labels of the extreme tiles of a snake graph. When considering posets $\calP_\frs$ for $\frs \in (0,\infty) \cap \QQ$, the first and last element of the poset (or, tile in the snake graph) is labeled by $x_1$. Therefore, we have $x_a=x_w=x_3, x_b=x_z=x_2,$ and $x_{i_1}=x_{i_d}=x_1$. Moreover, the matrix product in \cite[Proposition 5.5]{MusikerMatrix} partitions the set of perfect matchings into four sets. These can be naturally rewritten in the language of order ideals. For instance, $S_A$ corresponds to order ideals of $\calP$ which do not contain $R$ or $L$ while $S_B$ corresponds to order ideals which contain $L$ but not $R$. Remarkably, when we compute the matrix in \cite[Corollary 5.6]{MusikerMatrix} in the language of posets and in the case of the cluster algebra $\mathcal{A}(k_1,k_2,k_3)$, we get a matrix very similar to the combinatorial cluster Markov-monodromy matrix. We refer to this as the \emph{Musiker-Williams-Markov (MWM) matrix}:

\[
MW_{\frs} = \begin{bmatrix} -\xx^{\gg_\frs^\circ} \mathcal{W}(\calP_\frs; L) & \xx^{\gg_\frs} \mathcal{W}(\calP_{\frs})\\
\xx^{2\gg_\frs^\circ-\gg_\frs} \mathcal{W}(\calP_\frs; L,\neg R) & -\xx^{\gg_\frs^\circ} \mathcal{W}(\calP_\frs; \neg R)
\end{bmatrix} .
\]

From Lemma \ref{lem:Asymmetry}, we can see that $\mathrm{tr}(MW_{\frs}) = \frac{x_1}{x_2x_3} x_\frs + k_\frs$. However, it is not clear how the MWM-matrices associated to a Farey triple relate. In other words, it is not clear how to describe a propagation rule in a tree of MWM-matrices.

\begin{question}
Given a Farey triple $(\frac{p}{q},\frac{p+r}{q+s},\frac{r}{s})$, how do the MWM-matrices $M_{\frac{p}{q}}, M_{\frac{p+r}{q+s}},$ and $M_{\frac{r}{s}}$ relate?
\end{question}
\bibliography{bibliography}
\end{document}